\documentstyle[12pt]{article}
\newtheorem{thm}{Theorem}[section]
\newtheorem{prop}[thm]{Proposition}

\newtheorem{conj}[thm]{Conjecture}
\newtheorem{thm-defi}[thm]{Theorem/Definition}
\newtheorem{example}[thm]{Example}
\newtheorem{cor}[thm]{Corollary}
\newtheorem{question}[thm]{Question}
\newtheorem{new-lemma}[thm]{Lemma}
\newtheorem{defi}[thm]{Definition}
\newtheorem{rem}[thm]{Remark}

\newtheorem{condition}[thm]{Condition}

\newcommand{\A}{{\cal A}}
\newcommand{\B}{{\cal B}}

\newcommand{\D}{{\cal D}}
\newcommand{\E}{{\cal E}}
\newcommand{\F}{{\cal F}}
\newcommand{\G}{{\cal G}}

\renewcommand{\H}{{\cal H}}

\newcommand{\LB}{{\cal L}}
\newcommand{\R}{{\cal R}}
\renewcommand{\S}{{\cal S}}

\newcommand{\M}{{\cal M}}

\newcommand{\U}{{\cal U}}

\newcommand{\W}{{\cal W}}
\newcommand{\X}{{\cal X}}

\newcommand{\Z}{{\cal Z}}
\renewcommand{\P}{{\cal P}}

\newcommand{\PP}{{\Bbb P}}

\newcommand{\RealNumbers}{{\Bbb R}}
\newcommand{\Integers}{{\Bbb Z}}
\newcommand{\ComplexNumbers}{{\Bbb C}}
\newcommand{\RationalNumbers}{{\Bbb Q}}
\newcommand{\LieAlg}[1]{{\frak #1}}
\newcommand{\linsys}[1]{{\mid}#1{\mid}}
\newcommand{\Mon}{{\rm Mon}}
\newcommand{\monrep}{{mon}}

\newcommand{\IsomRightArrow}{\stackrel{\cong}{\rightarrow}}
\newcommand{\LongIsomRightArrow}{\stackrel{\cong}{\longrightarrow}}
\newcommand{\RightArrowOf}[1]{\stackrel{#1}{\rightarrow}}

\newcommand{\LongRightArrowOf}[1]{\stackrel{#1}{\longrightarrow}}

\newcommand{\StructureSheaf}[1]{{\cal O}_{#1}}
\newcommand{\EndProof}{\hfill  $\Box$}
\newcommand{\restricted}[2]{#1_{\mid_{#2}}}

\newcommand{\rank}{{\rm rank}}

\newcommand{\Pic}{{\rm Pic}}
\newcommand{\Sym}{{\rm Sym}}
\newcommand{\Ext}{{\rm Ext}}

\newcommand{\Hom}{{\rm Hom}}
\newcommand{\Aut}{{\rm Aut}}
\newcommand{\End}{{\rm End}}

\newcommand{\Abs}[1]{|\!#1\!|}

\newcommand{\RelExt}{{\cal E}xt}

\newcommand{\Wedge}[1]{\stackrel{#1}{\wedge}}

\input amssymb.sty

\begin{document}
%\rightline{\today}
\begin{center}
\begin{Large}
{\bf 
\noindent
On the monodromy of moduli spaces of sheaves on K3 surfaces 
}
\end{Large}
\\
Eyal Markman
\footnote{Partially supported by NSF grant number DMS-9802532}
\end{center}
%\author{Eyal Markman}
%\maketitle
\begin{abstract}
Let $S$ be a K3 surface and $\Aut D(S)$ 
the group of auto-equivalences of the derived category of $S$. 
We construct a natural representation of $\Aut D(S)$ 
on the cohomology of all moduli spaces 
of stable sheaves (with primitive Mukai vectors) on $S$. 
The main result of this paper is the precise 
relation of this action with the monodromy of 
the Hilbert schemes $S^{[n]}$ of points on the surface.  
A formula is provided for the monodromy representation, in terms of 
the Chern character of the universal sheaf. 
Isometries of the second cohomology of $S^{[n]}$ are lifted, via this formula,
to monodromy operators of the whole cohomology ring of $S^{[n]}$. 
\end{abstract}

{\scriptsize 
\tableofcontents 
} 
%***************************************************************************
% Introduction 
%***************************************************************************
\section{Introduction}
%\label{sec-introduction}

A variety $X$, with an ample canonical or anti-canonical line bundle, 
is completely determined by the bounded derived category $D(X)$ 
of coherent sheaves on  $X$. Moreover, $D(X)$ admits only the obvious 
auto-equivalences \cite{bondal-orlov}. In contrast, a 
projective variety $X$ with trivial canonical bundle, often admits 
interesting auto-equivalences. These symmetries have been used
by Mukai \cite{mukai-fourier} and others 
to study moduli spaces of sheaves on abelian varieties and 
Calabi-Yau varieties. 
Equivalences of derived categories of K3 or abelian surfaces 
have been used to relate moduli spaces of sheaves on 
the surface $S$, to its Hilbert scheme. 
Any smooth and compact moduli space $\M$, of stable sheaves on $S$,
is deformation equivalent to the product $\Pic^0(S)\times S^{[n]}$
of the identity component of the Picard group, with the 
Hilbert scheme $S^{[n]}$, 
of length $n$-subschemes
%, where $2n=\dim(M)$ 
\cite{yoshioka-irreducibility,yoshioka-abelian-surface}. 
In particular, $\M$ and $\Pic^0(S)\times S^{[n]}$ have isomorphic 
cohomology rings.

We concentrate in this paper on the case of a projective $K3$ surface $S$. 
Let $\Aut D(S)$ be the group of auto-equivalences of $D(S)$.
One may ask, to what extent $\Aut D(S)$ acts on the collection of 
moduli spaces of stable sheaves on $S$? 
The most naive expectation fails; $\Aut D(S)$ does not act via isomorphisms. 
We will show that $\Aut D(S)$ does act on the level of 
cohomology of the moduli spaces. Moreover, the action is related to 
the monodromy representation of a fixed moduli space. 
%The relationship, with the monodromy representation, is best explained in 
%terms of the Mukai lattice. 
We will consider a more general action (\ref{eq-composite-functor}) 
of a groupoid, whose morphisms are 
equivalences between derived categories of two, possibly non-isomorphic, 
$K3$ surfaces.
These ideas are illustrated in section
\ref{sec-elliptic-curve} in the simple case of 
moduli spaces of stable sheaves on an elliptic curve.

%**********************************************************
% Statements of the results
%**********************************************************
\subsection{Statements of the results}

%**********************************************************
% Symmetries of moduli spaces of stable sheaves
%**********************************************************
\subsubsection{Symmetries of moduli spaces of stable sheaves}
%\subsubsection{The Mukai lattice}
{\bf The Mukai Lattice:}
Let $K_{top}(S)$ be  the Grothendieck $K$-ring of topological complex
vector bundles on $S$ \cite{atiyah-book}. The
Euler characteristic $\chi: K_{top}(S) \rightarrow \Integers$ is given by
\[
\chi(v) = 2\rank(v) + \int_S\left[c_1(v)^2/2-c_2(v)\right].
\]
%\begin{eqnarray*}
%\chi \ : \ K_{top}(S) & \rightarrow & \Integers
%\\
%\chi(v) & = & 2\rank(v) + \int_S\left[c_1(v)^2/2-c_2(v)\right]
%\end{eqnarray*}
%the Euler characteristic. 
$K_{top}(S)$ is a free abelian group
of rank $24$. Any class $x\in K_{top}(S)$ is the difference 
$[E]-[F]$ of classes of complex vector bundles, and
we denote by $x^\vee$ the class $[E^\vee]-[F^\vee]$ obtained from  the
dual vector bundles.
The bilinear form on $K_{top}(S)$, given by 
\begin{equation}
\label{eq-mukai-pairing-on-K-top}
(x,y) \ \ \ := \ \ \ -\chi(x^\vee\otimes y),
\end{equation}
is called the {\em Mukai pairing}. The pairing is symmetric, unimodular,
of signature $(4,20)$, and the resulting lattice is called the 
{\em Mukai lattice}.
We define a polarized weight $2$ 
Hodge structure on $K_{top}(S)_\ComplexNumbers:=
K_{top}(S)\otimes_\Integers\ComplexNumbers$,
by setting $K_{top}^{1,1}(S)\subset K_{top}(S)_\ComplexNumbers$
to be the subspace of classes $v$ with $c_1(v)$ of type $(1,1)$,
and $K_{top}^{2,0}(S)$ the subspace of classes $v$,
orthogonal to $K_{top}^{1,1}(S)$, with $c_1(v)$ of type $(2,0)$.
Note that the Chern character 
\begin{eqnarray}
\label{eq-chern-character}
ch \ : \ K_{top}(S) & \rightarrow & H^*(S,\Integers)
\\
\nonumber
v & \mapsto & \rank(v)+c_1(v)-c_2(v)+c_1(v)^2/2
\end{eqnarray}
is an isomorphism of free abelian groups preserving the ring structure 
%($ch(v)=\rank(v)+c_1(v)-c_2(v)+c_1(v)^2/2$ 
($ch(v)$ is integral, 
since the lattice $H^2(S,\Integers)$ is even, injectivity of $ch$ 
follows from Theorem 3.25 in \cite{karoubi}, 
and surjectivity is easily verified). 
{\em Caution:} From section \ref{sec-weight-2-hodge-structure}
onward, we will use Mukai's normalization
(\ref{eq-isometry-between-K-top-and-Mukai-lattice}) 
of the isomorphism (\ref{eq-chern-character}).

\begin{defi}
\label{def-effective-class}
{\rm
A non-zero class $v\in K_{top}^{1,1}(S)$ will be called {\em effective},
if $(v,v)\geq -2$, 
$\rank(v)\geq 0$, and the following conditions hold.
If $\rank(v)=0$, then $c_1(v)$ is the
class of an effective (or trivial) divisor on $S$.
If both $\rank(v)$ and $c_1(v)$ vanish, then $\chi(v)>0$.
}
\end{defi}
Note that if $v\in K_{top}^{1,1}(S)$, $v\neq 0$, and
$(v,v)\geq -2$, then either $v$ or $-v$ is effective
(\cite{BHPV} chapter VIII Proposition 3.7). 
%We will say in this case that $\pm v$ is effective.
The class $v$ is {\em primitive}, if it is not a multiple of a class
in $K_{top}(S)$ by an integer larger than $1$.
If $v$ is a primitive and effective class, 
%and one of $\rank(v)$ or $\chi(v)$ is non-zero, 
then there almost always exists an ample line bundle $H$ on $S$, called
{\em $v$-suitable} in Definition \ref{def-v-suitable}, such that
the moduli space $\M_H(v)$, of $H$-stable coherent sheaves on $S$ with
class $v$, is non-empty,
smooth, connected, projective, and holomorphic symplectic
of dimension $2+(v,v)$
(\cite{mukai-symplectic-structure,ogrady-hodge-str,yoshioka-irreducibility,
yoshioka-abelian-surface}
and Theorem \ref{thm-irreducibility} below). 
The only exception, where a $v$-suitable ample line bundle $H$ 
need not exist, is the 
case $\rank(v)=\chi(v)=0$, and $\rank(\Pic(S))>1$
(see Condition \ref{condition-minimality} in that case).
%provided 
%$\Pic(S)$ is cyclic or if $c_1(v)$ satisfies a certain minimality condition 
%with respect to some ample line-bundle 
%(Condition \ref{condition-minimality}). 
We use the $H$-stability of sheaves due to Gieseker, Maruyama, and Simpson
\cite{huybrechts-lehn-book}.
If a non-zero class $v$ is not effective, then every sheaf 
with class $v$ has an endomorphism algebra of dimension $\geq 1$
\cite{mukai-symplectic-structure}.
In particular, there does not exist any $H$-stable sheaf with class $v$, 
for any polarization $H$. 
{\em Caution:} There may be an unstable such sheaf,
so the term stably-effective is more accurate. We will use the shorter term.

%\subsubsection{Universal classes}
{\bf Universal classes:}
A universal sheaf $\E_v$ exists over $\M_H(v)\times S$, 
if there exists a class $w\in K_{top}^{1,1}(S)$ satisfying $(v,w)=1$
\cite{mukai-hodge}.
In that case $\E_v$ is determined uniquely up to tensorization by a 
line bundle on $\M_H(v)$. We denote by
\[
e_v \ \ \ \in \ \ \ K_{top}(\M_H(v)\times S)
\]
the class of $\E_v$ in the topological $K$-ring. 
In general, a universal sheaf need not exist. 
Nevertheless, a universal class $e_v$ is constructed in
a sequel to this paper \cite{markman-integral-generators}, where
the integral structure of the cohomology
of moduli spaces is studied
(see also Definition \ref{def-e-v} below). 
The class $e_v$ is canonical, up to tensorization
with the pullback of the class of a topological line-bundle on
$\M_H(v)$. 
In this 
introduction we state the results of this paper using 
the existence of $e_v$, 
as it provides a simpler formulation of the results. 
%We will then indicate 
%the actual results proven 
The proofs we provide will use only the existence of the 
class $ch(e_v)$, which was introduced earlier in 
\cite{markman-diagonal} without relying on the existence of 
the integral class $e_v$ 
(see equations (\ref{eq-chern-character-of-rational-universal-class}) and 
(\ref{eq-comparizon-between-universal-class-and-twisted-universal-sheaf}) 
and Lemma \ref{lemma-two-definitions-of-gamma-g-v} below).

%There is always a holomorphic 
%projective bundle $\PP_v$ on $\M_H(v)$ and a sheaf 
%$\widetilde{\E}_v$ on $\PP_v\times S$, flat over $\PP_v$,
%with the following property. Let $\PP_{(t,s)}$ be the fiber of 
%$\PP_v\times S$ over $(t,s)\in \M_H(v)\times S$.
%The restriction of $\widetilde{\E}_v$ to 
%$\PP_{(t,s)}$ is isomorphic to 
%$E_t\boxtimes \StructureSheaf{\PP_{(t,s)}}(-1)$, where
%$E_t$ is the $H$-stable sheaf on $S$ corresponding to $t$. 
%We prove in \cite{markman-integral-generators} the existence of a
%topological complex line bundle $L$ on $\PP_v$,
%such that the class $[L]\otimes [\widetilde{\E}_v]$
%in $K_{top}(\PP_v\times S)$ is the pullback of a unique class
%$e_v$ in $K_{top}(\M_H(v)\times S)$, called a {\em universal class}. 
%The class $e_v$ is independent of the choice of
%$\PP_v$, and is unique up to a product with the class of a 
%topological line bundle on $\M_H(v)$.

%\subsubsection{Cohomological Fourier-Mukai transformations}
{\bf Cohomological Fourier-Mukai transformations:}
Let $S_1$ and $S_2$ be projective $K3$ surfaces,
$g:K_{top}(S_1)\rightarrow K_{top}(S_2)$ an isometry of Mukai lattices,
and $v_1\in K_{top}(S_1)$ a primitive and effective class. 
Set $v_2:=g(v_1)$ and assume that $v_2$ is effective. 
We do not assume that $g$ preserves the Hodge structures.
Choose $v_i$-suitable polarizations $H_i$, $i=1,2$, and
universal classes $e_{v_i}$ in $K_{top}(\M_{H_i}(v_i)\times S_i)$.
The exterior product induces an 
isomorphism
\[
K_{top}(\M_{H_i}(v_i))\otimes K_{top}(S_i)
\ \ \ \cong \ \ \ 
K_{top}(\M_{H_i}(v_i)\times S_i)
\]
by the K\"{u}nneth Theorem 
(\cite{atiyah-book} Corollary 2.7.15).
We can thus regard $(1\otimes g)$ as a homomorphism from 
$K_{top}(\M_{H_1}(v_1)\times S_1)$ to $K_{top}(\M_{H_1}(v_1)\times S_2)$.
Set $m:=2+(v_1,v_1)$. 
Let $\pi_{ij}$ be the projection from 
$\M_{H_1}(v_1)\times S_2 \times \M_{H_2}(v_2)$ to
the product of the $i$-th and $j$-th factors. 
We define a class in the middle cohomology 
$H^{2m}(\M_{H_1}(v_1)\times\M_{H_2}(v_2),\Integers)$
of the product: 
\begin{equation}
\label{eq-intro-class-gamma-g-v}
\gamma(g,v_1) \ \ \ := \ \ \
c_m\left(-\pi_{{13}_!}\left\{\pi_{12}^![(1\otimes g)(e_{v_1})]^\vee
\otimes \pi_{23}^!(e_{v_2})
\right\}
\right).
\end{equation}
Above, $\pi_{ij}^!$ is the pull-back homomorphism of $K$-rings 
and $\pi_{ij_!}$ is the Gysin homomorphisms
(\cite{bfm}, \cite{karoubi} Proposition IV.5.24). 
When the classes $e_{v_i}$ are defined in the algebraic 
$K$-group $K_{alg}(\M_{H_i}(v_i)\times S_i)$, 
the above formula yields the same class when 
$\pi_{ij}^!$ and $\pi_{ij_!}$ are the usual pull-back and push-forward
homomorphisms of algebraic $K$-groups (see \cite{fulton} or
\cite{hartshorne} for their definition).
Denote by
\begin{equation}
\label{eq-gamma-g-v}
\gamma_{g,v_1} \ : \
H^*(\M_{H_1}(v_1),\Integers)_{free} \ \ \ \longrightarrow \ \ \
H^*(\M_{H_2}(v_2),\Integers)_{free}
\end{equation}
the graded homomorphism, between the cohomology groups modulo their 
torsion subgroups, induced by the class $\gamma(g,v_1)$ using
the K\"{u}nneth and Poincare-Duality theorems.
Isometries between two Mukai lattices come in two flavors;
{\em orientation-preserving} and {\em orientation-reversing}
(Remark \ref{rem-natural-orientation}).

\begin{thm}
\label{thm-action}
Let $g$, $S_i$, $v_i$, $H_i$, $e_{v_i}$, $i=1,2$, be as above.
%Let $S_1$, $S_2$ be projective $K3$ surfaces,
%$g:K_{top}(S_1)\rightarrow K_{top}(S_2)$ an isometry of Mukai lattices,
%and $v_1\in K_{top}(S_1)$ a primitive and effective class. 
%Set $v_2:=g(v_1)$ and assume that $v_2$ is effective. 
%Choose $v_i$-suitable polarizations $H_i$, $i=1,2$, and
\begin{enumerate}
\item
\label{thm-item-gamma-g-v-is-an-isomorphism}
$\gamma_{g,v_1}$ is an isomorphism of cohomology rings,
which is independent of the choice of the universal classes $e_{v_i}$.
\item
\label{thm-item-gamma-g-in-K-alg}
If $g$ is an isomorphism of Hodge structures,  
then so is $\gamma_{g,v_1}$. If $g$ is also orientation preserving,
then $\gamma(g,v_1)$ is represented by a class in 
$K_{alg}[\M_{H_1}(v_1)\times \M_{H_2}(v_2)]\otimes_\Integers\RationalNumbers$,
in general, and in $K_{alg}[\M_{H_1}(v_1)\times \M_{H_2}(v_2)]$,
if universal sheaves exist over both moduli spaces.
\item
\label{thm-item-multiplicative-properties}
The equality $\gamma_{g^{-1},v_2}=(\gamma_{g,v_1})^{-1}$ holds.
Assume $S_3$ is a $K3$ surface, 
$f:K_{top}(S_2)\rightarrow K_{top}(S_3)$ an isometry,
$v_3:=f(v_2)$ an effective class, and $H_3$ a $v_3$-suitable polarization.
Then
\[
\gamma_{f,v_2}\circ \gamma_{g,v_1} \ \ \ = \ \ \ \gamma_{fg,v_1}.
\]
\item
\label{thm-item-invariance-of-the-universal-class}
There exists a topological line bundle $\ell$
on $\M_{H_2}(v_2)$, depending uniquely (up to isomorphism) on the triple 
$(g,c_1(e_{v_1}),c_1(e_{v_2}))$, such that the following equality holds
\begin{equation}
\label{eq-gamma-g-takes-universal-class-to-universal-class}
(\gamma_{g,v_1}\otimes \bar{g})(ch(e_{v_1})) \ \ \ = \ \ \ 
ch(e_{v_2})ch(\ell),
\end{equation}
where $\bar{g}:H^*(S_1,\Integers)\rightarrow H^*(S_2,\Integers)$
is the conjugate $ch\circ g\circ (ch)^{-1}$ of $g$ via the Chern character 
isomorphism (\ref{eq-chern-character}).
\item 
\label{thm-item-characterization-of-gamma-g}
Let $f  : 
H^*(\M_{H_1}(v_1),\Integers)_{free} \rightarrow 
H^*(\M_{H_2}(v_2),\Integers)_{free}$ be an isomorphism of cohomology rings 
satisfying the analogue
\begin{equation}
\label{eq-characterization-of-gamma-g}
(f\otimes \bar{g})(ch(e_{v_1})) \ \ \ = \ \ \ 
ch(e_{v_2})ch(\ell)
\end{equation}
of equation (\ref{eq-gamma-g-takes-universal-class-to-universal-class}),
for some topological line-bundle $\ell$ on $\M_{H_2}(v_2)$. 
Then $f=\gamma_{g,v_1}$. 
\item
\label{thm-item-invariance-of-chern-classes}
Let $c_k(\M_{H_i}(v_i))$ be the image of $c_k(T\M_{H_i}(v_i))$
in $H^{2k}(\M_{H_i}(v_i),\Integers)_{free}$. Then 
the following equality holds:
\[
\gamma_{g,v_1}(c_k(\M_{H_1}(v_1))) \ \ \ = \ \ \ c_k(\M_{H_2}(v_2)). 
\]
\end{enumerate}
\end{thm} 

%Note: Using the definition of Chern classes in $K$-theory, 
%we can interpret that class $\gamma(g,v)$ accordingly,
%view $\gamma_{g,v}$ as a homomorphism of $K$-groups.

%\subsubsection{Relation with equivalences of derived categories}
{\bf Relation with equivalences of derived categories:}
An equivalence $\Phi:D(S_1)\rightarrow D(S_2)$,
of the bounded derived categories of coherent sheaves,
induces an isometry $\phi:K_{top}(S_1)\rightarrow K_{top}(S_2)$
(section \ref{seq-Fourier-Mukai-transformations}).
In some cases, the image $\Phi(F)$, of every 
$H_1$-stable sheaf $F$ on $S_1$ with class $v$, 
can be represented by an $H_2$-stable sheaf on $S_2$.
In such cases, the auto-equivalence induces an isomorphism
\[
f \ : \ \M_{H_1}(v) \ \  \longrightarrow \ \  \M_{H_2}(\phi(v)).
%\\
%\hspace{4ex} F & \mapsto & \Phi(F).
\]
Furthermore, $(f_!\otimes \phi)(e_v)$ is a universal class on
$\M_{H_2}(\phi(v))\times S_2$.
Part \ref{thm-item-characterization-of-gamma-g} of Theorem \ref{thm-action}
implies the equality $f_*=\gamma_{\phi,v}$ of isomorphisms
of cohomology rings 
(section \ref{sec-conceptual-interpretation-of-the-formula-gamma-g}). 

Consider instead the contravariant functor from $D(S_1)^{op}$ to $D(S_2)$
taking an object $F$ to the derived dual $\Phi(F)^\vee$ of $\Phi(F)$.
Assume that $\Phi(F)^\vee$ 
can be represented by an $H_2$-stable sheaf on $S_2$, for every 
$H_1$-stable sheaf $F$ on $S_1$ with class $v$.
Then we get an isomorphism 
$f:\M_{H_1}(v)\rightarrow \M_{H_2}(\phi(v)^\vee)$. The induced
isomorphism of $K$-rings $f_!$ satisfies the equation 
$(f_!\otimes \phi)(e_v)=(e_{\phi(v)^\vee})^{\vee}$, 
for a suitable choice of a universal class $e_{\phi(v)^\vee}$. 
The pushforward $f_*$ on the level of cohomology is identified as follows.
Let 
\begin{equation}
\label{eq-duality-involution}
D \ : \ K_{top}(S) \ \ \ \longrightarrow \ \ \ K_{top}(S)
\end{equation}
be the {\em duality involution} sending $w$ to 
$w^\vee$.
The analogue for the cohomology of a moduli space $\M_H(v)$ 
\begin{equation}
\label{eq-duality-M}
D_\M \ : \ H^*(\M_H(v),\Integers)_{free} \ \longrightarrow \ 
H^*(\M_H(v),\Integers)_{free}
\end{equation}
acts by $(-1)^i$ on $H^{2i}(\M_H(v),\Integers)$.
$D_\M$ is a ring isomorphism since the odd Betti numbers of $\M_H(v)$ 
vanish (Theorem \ref{thm-irreducibility}). 
Part \ref{thm-item-characterization-of-gamma-g} of Theorem \ref{thm-action}
implies the equality 
\begin{equation}
\label{eq-isomorphism-of-cohomology-for-contravariant-functor}
f_* \ \ \ = \ \ \ D_\M\circ \gamma_{D\circ\phi,v}.
\end{equation} 
Symmetries of the cohomology of moduli spaces of sheaves, 
arising from such contravariant functors, 
are analogous to symmetries of Grassmannians 
arising from outer-automorphisms of $GL(n)$ 
(see Example 
\ref{example-two-components-of-automorphisms-of-grassmannian}).

In most cases an equivalence sends some stable sheaves to 
complexes, and does not give rise to an isomorphism of moduli spaces. 
In many cases $\Phi$ determines a birational isomorphism $f$, 
but in that case $f$ does not determine the isomorphism
$\gamma_{\phi,v}$, as the following example illustrates. 
Let $S$ be a $K3$ surface with a smooth rational curve $\Sigma \subset S$.
Choose $v$ to be the class of a sky-scraper sheaf of a point, so that
$\M(v)=S$. An auto-equivalence $\Phi$ of $D(S)$ is constructed in
Example \ref{example-relative-fourir-mukai-transform-wrt-E-0}
from any rigid sheaf with a trivial endomorphisms group. Such is 
the push-forward to $S$ of $\StructureSheaf{\Sigma}(-1)$. The resulting 
$\Phi$ maps the skyscraper sheaf of a point in $S\setminus \Sigma$
to itself. The class $\gamma(\phi,v)$ is {\em not}
Poincare-dual to the diagonal $\Delta$ in $S\times S$ but is rather
Poincare-dual to $[\Delta]+[\Sigma\times\Sigma]$. 
Theorems 
\ref{thm-class-of-correspondence-in-stratified-elementary-trans} and 
\ref{thm-reflection-sigma-satisfies-main-conj} and Example 
\ref{example-fiber-square-of-Hilbert-scheme-over-symmetric-product} 
concern higher dimensional analogues.
In the proofs of 
Theorems \ref{thm-class-of-correspondence-in-stratified-elementary-trans} 
and \ref{thm-reflection-sigma-satisfies-main-conj},
$\gamma(\phi,v)$ is shown to be Poincare-dual to
the reducible lagrangian Steinberg correspondence in 
$\M_\LB(v)\times \M_\LB(v)$, when $\phi$ is a suitable reflection. 

Whether an equivalence $\Phi$ induces an isomorphism
$f: \M_{H_1}(v)\rightarrow \M_{H_2}(\phi(v))$ depends, of course,
also on the choice of polarizations. 
When $\Phi$ is the identity, $H_i$, $i=1,2$, 
are two $v$-suitable polarization, 
and $c_1(v)$ is a multiple of an ample class, the class
$\gamma(id,v)$ is Poincare-dual to a correspondence in 
$\M_{H_1}(v)\times \M_{H_2}(v)$, which deforms to an isomorphism,
under simultaneous small deformations of $\M_{H_i}(v)$, $i=1,2$
(see Proposition \ref{prop-yoshioka-5.1} 
and Lemma \ref{lemma-the-class-of-f-is-gamma-E1-E2}).

%The reflection of a $K3$ surface $S$ 
%by a smooth rational curve $\Sigma \subset S$ provides such an example
%(Example \ref{example-relative-fourir-mukai-transform-wrt-E-0} 
%with the sheaf $E_0$ being the push-forward to $S$ of 
%$\StructureSheaf{\Sigma}(-1)$).
%In that case $v$ is the class of a sky-scraper sheaf of a point,
%$\M(v)=S$, $\phi$ is the reflection with respect to the $-2$ class 
%of $E_0$, $\phi(v)=v$, the the birational isomorphism $f$ is the identity
%on $S\setminus \Sigma$, but the class $\gamma(\phi,v)$ is
%Poincare-dual to $[\Delta]+[\Sigma\times\Sigma]$, where
%$\Delta$ is the diagonal in $S\times S$.

{\bf The moduli space $\M_H(v)$ for $v$ with negative rank:}
The isometry $\phi$, associated to an equivalence $\Phi$,
may send effective classes to non-effective classes.
The shift auto-equivalence, for example, corresponds to the isometry
$-id$. 
If $w\in K^{1,1}_{top}(S)$ is a primitive class and
$(w,w)\geq -2$, but $w$ is not effective, then 
the moduli space $\M_H(w)$ is empty while $\M_H(-w)$ is not, since
$-w$ is effective. Greater symmetry is obtained if we reset
\[
(\M_H(w),e_w) \ \ \ := \ \ \ (\M_H(-w),-e_{-w})
\]
(Definition \ref{def-M-minus-v}).
Then Theorem \ref{thm-action} extends (tautologically) 
relaxing the efficacy condition on $v_i$ by 
the condition that $v_i$ is primitive of type $(1,1)$. 

%*********************************************************************
% A groupoid representation
%*********************************************************************
%\subsubsection{A groupoid representation}
{\bf A groupoid representation:}
Parts \ref{thm-item-gamma-g-v-is-an-isomorphism}
and \ref{thm-item-multiplicative-properties}
of Theorem \ref{thm-action} may be reformulated as the 
construction of 
a representation for a groupoid 
\begin{equation}
\label{eq-groupoid}
\G.
\end{equation} 
Recall, that a groupoid is a category all of which morphisms are isomorphisms.
The objects of $\G$ are triples $(S,v,H)$, where $S$ is a projective $K3$ 
surface, $v\in K_{top}^{1,1}(S)$ is a primitive class, and 
$H$ is $v$-suitable ample line bundle. 
A morphism from $(S_1,v_1,H_1)$ to $(S_2,v_2,H_2)$ is an isometry
$g:K_{top}(S_1)\rightarrow  K_{top}(S_2)$ satisfying $g(v_1)=v_2$.
If $(v_1,v_1)= (v_2,v_2)$, then 
$\Hom_{\G}((S_1,v_1,H_1),(S_2,v_2,H_2))$ is a torsor for
each of the two subgroups $\Gamma_{v_i}\subset OK_{top}(S_i)$ 
consisting of isometries stabilizing $v_i$  
(Lemma \ref{lemma-primitive-embeddings-of-rank-2-lattices}). 
If $(v_1,v_1)\neq (v_2,v_2)$, then $\Hom_{\G}((S_1,v_1,H_1),(S_2,v_2,H_2))$ 
is empty. 
%Let $\G$ be the groupoid consisting of
%septets $(S_1,v_1,H_1,g,S_2,v_2,H_2)$, such that $S_i$ is a projective
%$K3$ surface, $v_i\in K_{top}^{1,1}(S_i)$ is a primitive class, 
%$H_i$ is a $v_i$-suitable polarization, 
%$g:K_{top}(S_1)\rightarrow  K_{top}(S_2)$ is an isometry, and
%$g(v_1)=v_2$. The operation
%\begin{equation}
%\label{eq-groupoid}
%(S_1,v_1,H_1,g,S_2,v_2,H_2)(S'_1,v'_1,H'_1,g',S'_2,v'_2,H'_2) \ \ = \ \ 
%(S'_1,v'_1,H'_1,gg',S_2,v_2,H_2)
%\end{equation}
%is defined, whenever $(S_1,v_1,H_1)=(S'_2,v'_2,H'_2)$.
%Let $\T$ be the set of triples $(S,v,H)$, where $v\in K_{top}^{1,1}(S)$ is a 
%primitive class and $H$ is $v$-suitable. 
%Note, that $\G$ acts transitively on each subset $\T_n$, $n\in \Integers$,
%consisting of triples with $(v,v)=n$ 
%(Lemma \ref{lemma-the-covariant-isometry-group-of-a-K3-acts-transitively}). 

The representation constructed is the functor
\begin{equation}
\label{eq-the-representation-functor-H}
\H:\G\rightarrow \A
\end{equation}
into the category $\A$ of 
commutative $\Integers$-algebras with a unit. 
The object $(S,v,H)$ is sent to the cohomology ring
$H^*(\M_H(v),\Integers)_{free}$ (or to $(0)$ if $(v,v)<-2$). 
A morphism $g\in \Hom_{\G}((S_1,v_1,H_1),(S_2,v_2,H_2))$
is sent to the isomorphism $\gamma_{g,v_1}$ given in
(\ref{eq-gamma-g-v}).

%a cohomology bundle ${\cal H}$
%over $\T$. The fiber of  ${\cal H}$ is the cohomology
%ring $H^*(\M(v)_H,\Integers)_{free}$.
%The latter is set to be $(0)$, if $(v,v)<-2$. 

%The groupoid consists of
%quintuples $(S_1,v_1,g,S_2,v_2)$, such that $S_i$ is a projective
%$K3$ surface, $v_i\in K_{top}^{1,1}(S_i)$ is a primitive class, 
%$g:K_{top}(S_1)\rightarrow  K_{top}(S_2)$ is an isometry, and
%$g(v_1)=v_2$. The representation constructed is a cohomology bundle 
%${\cal H}$ over the space of pairs $(S,v)$, where $v\in K_{top}^{1,1}(S)$ 
%is a primitive class. 
%The fibers of  ${\cal H}$ are the cohomology
%rings $H^*(\M(v),\Integers)_{free}$.
%The isomorphism
%$\gamma_{id,v}$ provides a natural identification of the cohomology rings 
%of $\M_{H_1}(v)$ and $\M_{H_2}(v)$, for any two $v$-suitable polatizations,
%and allows us to omit 
%the polarization in the notation $H^*(\M(v),\Integers)_{free}$. 
%The latter is set to be $(0)$, if $(v,v)<-2$. 

Consider the analogous groupoid $\D$, with the same objects, whose 
morphisms $\Phi\in \Hom_{\D}((S_1,v_1,H_1),(S_2,v_2,H_2))$ are equivalences
%consisting of
%septets $(S_1,v_1,H_1,\Phi,S_2,v_2,H_2)$, where $S_i$, $v_i$, and $H_1$ 
%are as above, and 
$\Phi:D(S_1)\rightarrow D(S_2)$ of derived categories,
sending $v_1$ to $v_2$. There is a natural functor $\D\rightarrow \G$, 
sending $\Phi$ to the associated isometry. 
We get the composite functor 
\begin{equation}
\label{eq-composite-functor}
\D\rightarrow \G\rightarrow \A
\end{equation}
representing the groupoid $\D$.
%This groupoid acts on the 
%cohomology bundle ${\cal H}$, via the groupoid homomorphism to $\G$,
%sending $\Phi$ to the associated isometry.

Let $\Gamma$ be the isometry group of the Mukai lattice of a
$K3$ surface $S$. Fix an effective class $v\in K_{top}(S)$ 
and a $v$-suitable polarization $H$. 
The subgroup $\Gamma_v\subset \Gamma$ stabilizing $v$ 
is also the automorphism group of the object $(S,v,H)$ of $\G$. 
As a corollary of Theorem \ref{thm-action} we get that $\Gamma_v$ 
acts on the cohomology of the corresponding moduli space:

\begin{cor} (Corollary \ref{cor-group-homomorphism}
and Lemma \ref{lemma-conjecture-holds-for-isotropic-mukai-vectors})
The natural map
\begin{eqnarray}
\label{eq-action-of-stabilizer}
\gamma \ : \ \Gamma_v & \longrightarrow & 
\Aut\left[H^*(\M_H(v),\Integers)_{\rm free}\right] 
\\
\nonumber
g & \mapsto & \gamma_{g,v}
\end{eqnarray}
is  group homomorphism. It is injective if $(v,v)\geq 2$.
\end{cor}

Let $\G_n$ be the (full) subgroupoid of $\G$, whose objects 
consist of triples with class $v$ satisfying $2+(v,v)=2n$.
The restriction $\H_n$ of the functor $\H$ to $\G_n$ 
is determined by each of the representations
(\ref{eq-action-of-stabilizer}), since any two objects of 
$\G_n$ are isomorphic in $\G_n$. In this sense, 
the functor $\H_n$ may be regarded as the ``induced representation''
$Ind_{\Gamma_v}^{\G_n}H^*(\M_H(v),\Integers)_{\rm free}$.
%$\H_n$, of the representation $\H$ to $\T_{n}$, 
%is isomorphic to the representation (???)
%$Ind_{\Gamma_v}^{\G}H^*(\M_H(v),\Integers)_{\rm free}$
%of $\G$ induced from the representation
%(\ref{eq-action-of-stabilizer}) of $\Gamma_v$, since 
%$\G$ acts transitively on $\T_n$.

%************************************************************
% The monodromy group of a moduli space
%************************************************************
\subsubsection{The monodromy group of a moduli space}
%{\bf The monodromy group of a moduli space:}
The second main result of this paper, is the relation between the action 
(\ref{eq-action-of-stabilizer}) of the stabilizer $\Gamma_v$, and
the monodromy representation of the moduli space $\M_H(v)$
(Theorem \ref{introduction-thm-Gamma-v-acts-motivicly}).
Roughly, the representation (\ref{eq-action-of-stabilizer}) 
acts via monodromy operators, modulo a sign change. 

\begin{defi}
\label{def-irreducible-symplectic}
{\rm 
An {\em irreducible holomorphic symplectic manifold}
is a simply-connected compact K\"{a}hler manifold $X$, such that 
$H^0(X,\Omega_X^2)$ is generated by an everywhere non-degenerate 
holomorphic two-form.
}
\end{defi}

An irreducible holomorphic symplectic manifold of real 
dimension $4n$ admits a Riemannian metric with holonomy group $Sp(n)$
\cite{beauville-varieties-with-zero-c-1}. 
Such a metric is called {\em hyperk\"{a}hler}.

\begin{defi}
\label{def-monodromy}
{\rm
Let $M$ be an irreducible holomorphic symplectic manifold. 
An automorphism $g$ of the cohomology ring 
$H^*(M,\RationalNumbers)$ is called a {\em monodromy operator}, 
if there exists a 
family $\M \rightarrow B$ (which may depend on $g$) 
of irreducible holomorphic symplectic manifolds, having $M$ as a fiber
over a point $b_0\in B$, 
and such that $g$ belongs to the image of $\pi_1(B,b_0)$ under
the monodromy representation. 
The {\em monodromy group} $\Mon(M)$ of $M$ is the subgroup, 
of the automorphism group of the cohomology ring of
$M$, generated by all the monodromy operators. 
}
\end{defi}

The monodromy group is interesting even when 
the moduli space is the Hilbert scheme $S^{[n]}$ of 
length $n$ subschemes on $S$. The Hilbert scheme is the moduli
space $\M_H(v)$, where $v$ is the class of the ideal sheaf of a 
length $n$ subscheme. When $n\geq 2$,
the Hilbert scheme $S^{[n]}$ admits deformations,
which are not Hilbert schemes of points on any $K3$ surface.
In fact, the moduli space, of irreducible holomorphic symplectic manifolds
deformation equivalent to $S^{[n]}$, 
has dimension one larger than the corresponding moduli space 
of the surface $S$. 
Consequently, the monodromy group of $S^{[n]}$ is larger. 

The isometry group $\Gamma$ of $K_{top}(S)$ has a natural character 
$cov:\Gamma\rightarrow \Integers/2\Integers$
(see (\ref{eq-top-homology-character})). 
This character sends reflections by $-2$ vectors to $0$
and reflections by $+2$ vectors to $1$. Note also that $cov(-1)=0$
and $cov(D)=1$, where $D$ is the duality involution 
(\ref{eq-duality-involution}). 
Isometries in the kernel of $cov$ are called
{\em orientation preserving}. 
Modify the representation $\gamma$ in 
(\ref{eq-action-of-stabilizer}) 
to get the representation 
\begin{eqnarray}
\label{introduction-eq-monodromy-is-gamma-times-cov}
\monrep  \ : \  \Gamma_v & \rightarrow & 
\Aut(H^*(\M_H(v),\Integers)_{\rm free}).
\\
\nonumber
g & \mapsto & (D_{\M_H(v)})^{cov(g)}\circ \gamma_{g,v},
\end{eqnarray}
where the duality involution $D_{\M_H(v)}$ is given in 
equation (\ref{eq-duality-M}).
This representation is faithful, if $(v,v)\geq 4$. If 
$(v,v)=2$, then the kernel of $\monrep$ is 
$\{id,-\sigma_v\}$, where $\sigma_v$ is the reflection in $v$
(see section \ref{sec-reflections-with-respect-to-plus-2-vectors}).

\begin{thm}
\label{introduction-thm-Gamma-v-acts-motivicly}
Let $v$ be an effective and primitive class in $K_{top}(S)$ 
and $H$ a $v$-suitable polarization.
Then the image $\monrep(\Gamma_v)$ 
is a normal subgroup of finite index\footnote{
$\Mon(\M_H(v))$ is shown in \cite{markman-integral-generators}
to be the direct product $K\times \monrep(\Gamma_v)$, where
$K$ is a finite abelian group of exponent $2$, which acts trivially on
$H^2(\M_H(v))$.
}
in the monodromy group
$\Mon(\M_H(v))$.
\end{thm}

%We will further show, that the image of $\monrep$ is
%a normal subgroup of finite index in the monodromy group
%(Lemma \ref{lemma-bounding-K}).
%Theorem \ref{thm-Gamma-v-acts-motivicly} is proven in section
%\ref{sec-hodge-isometries-of-hilbert-schemes}. 

%The universal class $e_v$ is not $\Mon(\M_H(v))$-equivariant, nor is $e_v$
%canonical. 
The $\Gamma_v$-invariance property 
(\ref{eq-gamma-g-takes-universal-class-to-universal-class})
of the class $e_v$ no longer holds, if the representation
$\gamma$ of $\Gamma_v$ is replaces by $\monrep$. 
A diagonal monodromy representation of $\Gamma_v$
decomposes a canonical normalization $u_v$ of $e_v$
as a sum $u_v^+ + u_v^-$,
of an invariant class $u_v^+$ and a class $u_v^-$, acted upon via the 
character $cov$ of $\Gamma_v$ 
(Lemma 
\ref{lemma-normalized-chern-character-is-not-monodromy-invariant}).

Sharper results are obtained about the weight $2$ monodromy representation.
We recall first the fundamental result about 
the second cohomology of a moduli space.
Let $v$ and $H$ be as in 
Theorem \ref{introduction-thm-Gamma-v-acts-motivicly} and
$v^\perp\subset K_{top}(S)$ the sublattice orthogonal to $v$. 
Let $\pi_\M$ and $\pi_S$ be the two projections from $\M_H(v)\times S$.
Mukai introduced the natural homomorphism
\begin{eqnarray}
\label{eq-introduction-theta-v-from-v-perp}
\theta_v \ : \ v^\perp & \longrightarrow & H^2(\M_H(v),\Integers)
\\
\label{eq-introduction-mukai-homomorphism}
\theta_v(x) & = & c_1\left[
\pi_{\M_!}\left(
\pi_S^!(x^\vee)\otimes e_v
\right)
\right].
\end{eqnarray}

The homomorphism $\theta_v$ is $\Gamma_v$-equivariant with respect to the 
representation $\gamma$ in (\ref{eq-action-of-stabilizer}) (see 
diagram (\ref{eq-g-gamma-g-conjugate-theta-to-theta})).
Beauville constructed an integral symmetric bilinear form on
the second cohomology of an irreducible holomorphic symplectic manifold
\cite{beauville-varieties-with-zero-c-1}. 
We will not need the intrinsic formula, but only the fact that it is 
monodromy invariant and its identification in 
the following theorem due to Mukai, Huybrechts, O'Grady, and Yoshioka:

\begin{thm} (\cite{ogrady-hodge-str}, 
\cite{yoshioka-abelian-surface} Theorem 8.1,
and \cite{yoshioka-note-on-fourier-mukai} Corollary 3.15)
\label{thm-irreducibility}
Let $v$ be a primitive and effective class and $H$ a $v$-suitable ample
line-bundle.
\begin{enumerate}
\item 
$\M_H(v)$ is a smooth, non-empty, irreducible symplectic, 
projective variety of dimension $\dim(v)=\langle v,v\rangle+2$. 
It is obtained by  deformations 
from the Hilbert scheme of $\frac{<v,v>}{2}+1$ points on $S$. 
\item
The homomorphism (\ref{eq-introduction-mukai-homomorphism}) 
is an isomorphism of weight 2 
Hodge structures with respect to Beauville's bilinear form on 
$H^2(\M_H(v),\Integers)$ when $\dim(v)\geq 4$.  When $\dim(v)=2$,
(\ref{eq-introduction-mukai-homomorphism}) factors through an isomorphism 
from $v^\perp/\Integers\cdot v$. 
\end{enumerate}
\end{thm} 

{\em Note:} Most of the cases of Theorem \ref{thm-irreducibility} are
proven in \cite{yoshioka-abelian-surface}
Theorem 8.1 under the additional assumption that $\rank(v)>0$
or $c_1(v)$ is ample. Corollary 3.15 in \cite{yoshioka-note-on-fourier-mukai}
proves the remaining cases ($\rank(v)=0$ and $c_1(v)$ not ample),
under the assumption that $\M_H(v)$ is non-empty. The non-emptyness
is proven in \cite{yoshioka-note-on-fourier-mukai} Remark 3.4,
under the assumption that $c_1(v)$ is nef and in complete generality in an 
unpublished note communicated to the author by K. Yoshioka.

Given an element $u$ of $H^2(\M_H(v),\Integers)$ with
$(u,u)$ equal $2$ or $-2$, let $\rho_u$ be the isometry given by
$
\rho_u(w)=\frac{-2}{(u,u)}w+(w,u)u. 
$
Then $\rho_u$ is the
reflection in $u$, when $(u,u)=-2$, and $-\rho_u$ is the reflection in $u$,
when $(u,u)=2$. 
Set 
\begin{equation}
\label{eq-W}
\W \ \ \ := \ \ \ \langle
\rho_u \ : \ u \in H^2(\M_H(v),\Integers) \ \
\mbox{and} \ \ (u,u) = \pm2
\rangle
\end{equation} 
to be the subgroup of $O(H^2(\M_H(v),\Integers))$ generated by the 
elements $\rho_u$. Then $\W$ is a normal subgroup of finite index
in $O(H^2(\M_H(v),\Integers))$.
The image of $\Gamma_v$, via the weight $2$ monodromy representation
of Theorem \ref{introduction-thm-Gamma-v-acts-motivicly}, is equal to
$\W$, by Lemma \ref{lemma-EO-is-maximal}. 

\begin{cor}
\label{cor-W-consists-of-monodromies}
$\W$ is contained in the image $\Mon^2$ in $O(H^2(\M_H(v),\Integers))$
of the whole monodromy group $\Mon(\M_H(v))$. 
Assume $(v,v)\geq 2$. Then 
the representation (\ref{introduction-eq-monodromy-is-gamma-times-cov})
factors as the composition of the surjective homomorphism 
\begin{eqnarray*}
\Gamma_v & \longrightarrow & \W 
\\
g & \mapsto & (-1)^{cov(g)}\cdot \theta_v\circ g \circ \theta_v^{-1}
\end{eqnarray*}
and a $\Mon(\M_H(v))$-equivariant injective homomorphism
\begin{equation}
\label{eq-homomorphism-mu-from-W-to-Mon}
\mu \ : \ \W  \ \ \ \rightarrow \ \ \ \Mon(\M_H(v))\subset 
\Aut[H^*(\M_H(v),\Integers)_{free}].
\end{equation}
\end{cor}

Above $\Mon(\M_H(v))$ acts by conjugation on $\W$ and itself. The 
equivariance of $\mu$ follows from the normality of $\W$ in 
$OH^2(\M_H(v),\Integers)$, the normality of
$\monrep(\Gamma_v)$ in Theorem
\ref{introduction-thm-Gamma-v-acts-motivicly}, and
the equivariance of the inclusion of the second cohomology in
the total cohomology. Being monodromy equivariant, $\mu$ is well defined for
any compact K\"{a}hler manifold $X$ deformation equivalent to $\M_H(v)$.

In a sequel to this paper we prove the opposite inclusion
$\Mon^2\subset \W$ \cite{markman-integral-generators}. 
Consequently, $\W$ is {\em equal} to $\Mon^2$. 
It follows, when $n-1$ is not a prime power, 
that $\Mon^2$ does {\em not} surject onto the 
quotient $O(H^2(S^{[n]},\Integers))/(\pm 1)$, of the whole 
isometry group by its center. This, in turn, 
leads to a counter example to a version of the the Generic-Torelli 
question: 
The weight $2$ Hodge structure of a generic 
irreducible holomorphic symplectic manifold $X$, 
deformation equivalent to $S^{[n]}$, does not determine 
the bimeromorphic class of $X$, when $n-1$ is not a prime power 
\cite{markman-integral-generators}. 
The index of $\Mon^2$ in $O(H^2(S^{[n]},\Integers))/(\pm 1)$
is a lower bound for the number of bimeromorphic classes 
with the same generic weight $2$ Hodge structure.
This index is calculated in Lemma \ref{lemma-the-index-of-Gamma-v}.

%***************************************************************
% Related works
%***************************************************************
\subsection{Related works}
%{\bf Related works by other authors:}
Let $X$ be an irreducible holomorphic symplectic manifold 
deformation equivalent to the Hilbert scheme $S^{[n]}$.
Corollary \ref{cor-W-consists-of-monodromies} describes the 
monodromy representation, of the finite index subgroup $\W$ of
the isometry group $O(H^2(X,\Integers))$, in
the automorphism group of the cohomology ring 
$H^*(X)$.
A related representation of $SO(H^2(X,\Integers))$, 
on the cohomology of a hyperk\"{a}hler variety $X$, 
was studied by Verbitsky via
different techniques (see \cite{verbitsky-mirror-symmetry} and Theorem 
\ref{thm-so-action-determines-hodge-structure} below, for Verbitsky's
result). Lemma \ref{lemma-comparing-two-representations}  
compares the two representations. 

There are many parallels between i) the action 
(\ref{eq-action-of-stabilizer}), of the group $\Gamma_v$, on the
cohomology of $\M_H(v)$,
and ii) the action of a Weyl group $W$, of a semisimple Lie group, 
on the cohomology of the cotangent bundle $T^*\B$ of the flag variety $\B$. 
The latter is of course the same as the cohomology of $\B$. 
$W$ is the reflection group of the root lattice, while 
the group $\Gamma_v$ is equal to the reflection group of the 
sublattice $v^\perp$ (Corollary 
\ref{cor-stabilizer-is-generated-by-reflections}). 
Both $\M_H(v)$ and $T^*\B$ have a holomorphic symplectic structure. 
Both actions are not realized by automorphisms. 
(The general Weyl group does act on the affine bundle $G/T\rightarrow B$,
which is a twisted version of $T^*\B$). 
In both cases, the action is realized geometrically via 
lagrangian correspondences. $W$ acts naturally,
via the Steinberg correspondences, on the cohomology 
of $\B$ (and, more generally, of any Springer fiber \cite{chriss-ginzburg}). 
In the Hilbert scheme case, lagrangian correspondences, 
%(\ref{eq-total-brill-noether-correspondence}) 
%and (\ref{eq-correspondence-inducing-involution-of-hilbert-scheme}), 
analogous to the Steinberg variety, 
%of Definition \ref{def-Z-for-cotangent-bundles}, 
play a major role in the proof of Theorem 
\ref{introduction-thm-Gamma-v-acts-motivicly} via local monodromy operators
(Theorem \ref{thm-image-of-mu-is-generated-by-local-monodromies}).

Similar results were obtained by Nakajima, when the 
K3 surface is replaced 
with the resolution of a simple surface singularity
\cite{nakajima-reflections}. 
A relationship with representations of affine Lie algebras, in the K3 case, 
is exhibited in \cite{nakajima-representations}. 
%*****************************************************************
% Moduli spaces of stable sheaves on an elliptic curve
%*****************************************************************
\subsubsection{Moduli spaces of stable sheaves on an elliptic curve}
\label{sec-elliptic-curve}

We verify the analogue of Theorem \ref{thm-action}, 
about the symmetries of moduli spaces of stable sheaves, 
in the elliptic curve case.
We also explain the analogue of Theorem 
\ref{introduction-thm-Gamma-v-acts-motivicly}, relating 
the above symmetries to the monodromy of a fixed moduli space. 

Let $\Sigma$ be an elliptic curve.  
The analogue $H^*(\Sigma,\Integers)$, 
of the Mukai lattice, decomposes as the direct sum
of $H^1(\Sigma,\Integers)$ with the even cohomology
$H^{{\rm even}}(\Sigma,\Integers):=
H^0(\Sigma,\Integers)\oplus H^2(\Sigma,\Integers)$. 
The odd summand is endowed with the intersection pairing. 
The Mukai vector of a sheaf $F$ is $v(F):=(\rank(F),c_1(F))$
and belongs to $H^{{\rm even}}(\Sigma,\Integers)$.
The Mukai pairing on the even summand is given by
(\ref{eq-Mukai-pairing-in-terms-of-sheaves}). 
%is $(v(E),v(F)):=-\chi(E^\vee\stackrel{L}{\otimes}F)$.
Explicitly, it is the anti-symmetric bilinear form
$((r',d'),(r'',d''))=r''d'-r'd''$. 
The symmetry group $\Gamma_\Sigma$ of $H^*(\Sigma,\Integers)$
is thus $SL(2,\Integers)\times SL(2,\Integers)$. 
This factorization, of the symmetry group, separates the monodromy
action from the action on $H^{{\rm even}}(\Sigma,\Integers)$
of the group $\Aut D(\Sigma)$ of auto-equivalences of the
derived category. $\Aut D(\Sigma)$ surjects onto the even $SL(2,\Integers)$
factor \cite{mukai-fourier}. Any auto-equivalence in 
$\Aut D(\Sigma)$ maps a simple sheaf to a shift of a simple sheaf
(Orlov's determination of $\Aut D(\Sigma)$ \cite{orlov-abelian-varieties}  
reduces the check to the case of Mukai's original Fourier-Mukai transform,
in which it follows from \cite{mukai-fourier} Corollary 2.5). 
Over an elliptic curve, a sheaf is stable if and only if it is simple.
$\Aut D(\Sigma)$ acts on the collection of moduli spaces
of stable sheaves (with a primitive Mukai vector) via isomorphisms, provided
we use the convention in Definition \ref{def-M-minus-v},
relating the moduli spaces $\M_\Sigma(r,d)$ and $\M_\Sigma(-r,-d)$.  

Let $g:H^*(\Sigma_1,\Integers)\rightarrow H^*(\Sigma_2,\Integers)$ 
be an isomorphism, compatible with the bilinear forms. 
Denote by $g=(g_0,g_1)$ the decomposition into even and odd factors.
Define the class $\gamma(g,v)$ using the topological formula 
$\gamma(g,\E_v,\E_{g(v)})$ given below in equation (\ref{eq-gamma-g-E1-E2}).
The class $\gamma(g,v)$ 
depends only on $v$, $g(v)$, and the odd factor $g_1$ of $g$ 
in the following sense. 
Let $\M_{\Sigma_i}(r,d)$ be a moduli space of stable sheaves on $\Sigma_i$, 
with a primitive vector $(r,d)$ with non-negative rank. 
Atiyah proved, that $\M_{\Sigma_i}(r,d)$ 
is naturally isomorphic to $\Pic^d(\Sigma_i)$, via $F\mapsto \det(F)$. 
Thus, the cohomologies of all moduli spaces $\M_{\Sigma_i}(r,d)$, with 
rank and degree coprime, are naturally identified with that of $\Sigma_i$. 
Under the above identifications, 
the class $\gamma(g,v)$ 
induces the unital ring isomorphism determined by 
$g_1:H^1(\Sigma_1,\Integers)\rightarrow H^1(\Sigma_2,\Integers)$. 
The proof is similar to that of Lemma 
\ref{lemma-conjecture-holds-for-isotropic-mukai-vectors}.

The analogue of our main Theorem 
\ref{introduction-thm-Gamma-v-acts-motivicly} 
is thus trivial in the elliptic curve case. Let $\Gamma_{\Sigma,v}$ be 
the subgroup of $\Gamma_\Sigma$ stabilizing a primitive vector $v$.
Then $\Gamma_{\Sigma,v}$ is the product 
$\Gamma_{\Sigma,v}^{\rm even}\times SL(2,\Integers)$, 
where the even factor is infinite cyclic. 
The even factor $\Gamma_{\Sigma,v}^{\rm even}$ 
acts trivially on $H^*(\M_\Sigma(v),\Integers)$,
while the odd $SL(2,\Integers)$ factor acts via the natural action on 
$H^1(\Sigma,\Integers)$. We see that the image of 
$\Gamma_{\Sigma,v}$ in $\Aut H^*(\M_\Sigma(v),\Integers)$ is the
monodromy group. 

Note, that {\em we have broken Mirror Symmetry}, which is supposed to 
interchange the two factors of $\Gamma_\Sigma$. 
The reason for this break of symmetry is, that we are considering only 
graded automorphisms of the cohomology ring of moduli spaces.
This is expressed in formula
(\ref{eq-gamma-delta}), by the fact that the class $\gamma(g,v)$
is the {\em middle dimensional} Chern class of a complex
(see also remark 
\ref{rem-auto-equivalences-of-the-derived-category-of-moduli-spaces}). 
The infinite cyclic even 
factor $\Gamma_{\Sigma,v}^{\rm even}$ is generated by the endomorphism 
$\tau_v(w):= w + (w,v)\cdot v$. The endomorphism 
$\tau_v$ can be lifted to an auto-equivalence of the derived category, via 
transvections in stable sheaves parametrized by $\M_\Sigma(v)$
(these transvections are described in 
Theorem \ref{thm-reflection-by-a-sherical-object} under the name reflections). 
Mirror symmetry would require the even factor $\Gamma_{\Sigma,v}^{\rm even}$ 
to act non-trivially. 
Such transvections are mirror symmetric to Dehn twists
\cite{seidel-thomas} in the odd monodromy factor of $\Gamma_\Sigma$. 

%*************************************************************
% Outline of the proofs
%*************************************************************
\subsection{Outline of the proofs}
We reduce the proofs of Theorems 
\ref{thm-action} and \ref{introduction-thm-Gamma-v-acts-motivicly}
to smaller results, which are proven in the rest of the paper.

\subsubsection{Reduction of the proof of Theorem \ref{thm-action}}
\label{sec-reduction-of-proof-of-thm-action}
%{\bf Outline of the proof of Theorem \ref{thm-action}:}  
The special case of the Theorem, where 
$(S_1,v_1,H_1)=(S_2,v_2,H_2)=(S,v,H)$
and $g$ is the identity, was proven in
\cite{markman-diagonal} (Theorem 
\ref{thm-graph-of-diagonal-in-terms-of-universal-sheaves} below).
The class $\gamma(id,v)$ is Poincare dual to the class of the diagonal in
$\M_{H}(v)\times \M_{H}(v)$. Consequently, $H^*(\M_{H}(v),\RationalNumbers)$
is generated by the K\"{u}nneth factors of the Chern classes of the
universal class $e_v$ (Corollary \ref{cor-kunneth-factors-generate}).
The general characterization of the homomorphism 
$\gamma_{g,v}$, in part \ref{thm-item-characterization-of-gamma-g} 
of Theorems \ref{thm-action}, follows formally from the
case $g=id$ and Corollary \ref{cor-kunneth-factors-generate}
(see Lemma \ref{lemma-recovering-f}).
The multiplicative properties of the isomorphisms $\gamma_{g,v}$, 
listed in part \ref{thm-item-multiplicative-properties}
of Theorems \ref{thm-action}, follow from this characterization 
of $\gamma_{g,v}$ 
and the $\gamma$-invariance of the universal class $e_v$ expressed in 
part \ref{thm-item-invariance-of-the-universal-class} of 
Theorem \ref{thm-action} 
(Lemma \ref{lemma-conjecture-holds-for-compositions}). 

We prove next parts \ref{thm-item-gamma-g-v-is-an-isomorphism} and 
\ref{thm-item-invariance-of-the-universal-class} 
of Theorem \ref{thm-action} (restated in cohomological terms
as Theorem \ref{thm-trancendental-reflections}). 
Let $\G$ be the groupoid (\ref{eq-groupoid}).
%consisting of septets $(S_1,v_1,H_1,g,S_2,v_2,H_2)$ satisfying the 
%hypothesis of Theorem \ref{thm-action}. 
Denote by 
$\G_{OK}\subset \G$ the sub-groupoid, with the same objects, 
whose morphisms  
are those isometries $g\in\Hom_{\G}((S_1,v_1,H_1),(S_2,v_2,H_2))$, 
for which $\gamma_{g,v_1}$
is a ring isomorphism, and 
$\gamma_{g,v_1}\otimes \bar{g}$ maps a universal class to a universal class
(equation (\ref{eq-gamma-g-takes-universal-class-to-universal-class})). 
Lemma \ref{lemma-conjecture-holds-for-compositions}
implies that $\G_{OK}$ is indeed a sub-groupoid. 

First we show that $\G_{OK}$ contains a sub-groupoid $\G_{mon}\subset \G$. 
We say that the two triples $(S_i,v_i,H_i)$, $i=1,2$,
are {\em deformation-equivalent}, if
the pairs $(S_i,v_i)$, $i=1,2$, are related by a 
deformation, satisfying mild conditions 
(Definition \ref{def-deformation-equivalent-triples}).
The morphisms in $\Hom_{\G_{mon}}((S_1,v_1,H_1),(S_2,v_2,H_2))$ 
are isometries 
$g:K_{top}(S_1)\rightarrow K_{top}(S_2)$, which arise as monodromy
operators for such a deformation. 
The inclusion $\G_{mon}\subset \G_{OK}$ is proven in 
Lemma \ref{lemma-the-class-of-f-is-gamma-E1-E2}. 

Next we show that $\G_{OK}$ contains another sub-groupoid $\G_{Mukai}$ 
of $\G$ (Lemma \ref{lemma-G-Mukai-is-in-G-OK}). The morphisms in
$\Hom_{\G_{Mukai}}((S_1,v_1,H_1),(S_2,v_2,H_2))$ 
arise in one of the following two ways (section
\ref{sec-conceptual-interpretation-of-the-formula-gamma-g}):

1) There exists an equivalence functor
$\Phi:D(S_1)\rightarrow D(S_2)$, taking $v_1$ to $v_2$, 
and $\Phi$ induces an isomorphism between the moduli spaces 
$\M_{H_1}(v_1)$ and $\M_{H_2}(v_2)$. 
We further assume that the functor $\Phi$ induces an 
{\em orientation-preserving} 
isometry $\phi$, i.e., that $\phi$ is 
in the kernel of the functor $cov$ defined below in (\ref{eq-functor-cov}). 
(The last assumption is automatically satisfied
by Lemma \ref{lemma-on-the-sign-conjecture}, which depends on
Theorems \ref{thm-action} and 
\ref{introduction-thm-Gamma-v-acts-motivicly}). 
In that case the isometry $\phi$ is in 
$\Hom_{\G_{Mukai}}((S_1,v_1,H_1),(S_2,v_2,H_2))$, by definition. 

2)  There exists an equivalence functor
$\Phi:D(S_1)\rightarrow D(S_2)$, taking $v_1$ to $(v_2)^\vee$, 
and the composite functor $F\mapsto (\Phi(F))^\vee$ induces
an isomorphism between the moduli spaces $\M_{H_1}(v_1)$ and 
$\M_{H_2}(v_2)$. 
We further assume that the functor $\Phi$ induces an 
{\em orientation-preserving} 
isometry $\phi$. (The last assumption is automatically satisfied
by Lemma \ref{lemma-on-the-sign-conjecture}). 
In that case the isometry $g=D\circ \phi$ is in 
$\Hom_{\G_{Mukai}}((S_1,v_1,H_1),(S_2,v_2,H_2))$, by definition. 

Two triples $(S_i,v_i,H_i)$, $i=1,2$, related by a morphism in
$\G_{Mukai}$, are said to be {\em Mukai-equivalent}. 
We conclude that $\G_{OK}$ contains the sub-groupoid of $\G$ generated by
$\G_{mon}$ and $\G_{Mukai}$. By this we mean that any morphism in $\G$, 
which is a composition of a finite set of 
morphisms in $\G_{mon}$ and  $\G_{Mukai}$, is a morphism in $\G_{OK}$. 

Let $Ob(\G)_n$ be the set of objects $(S,v,H)$ of $\G$ satisfying
the equality $2+(v,v)=2n$, $n\geq 0$. 
Yoshioka proved that any two triples $(S_i,v_i,H_i)$, $i=1,2$,
in $Ob(\G)_n$ are related by the equivalence relation generated by 
deformation and Mukai equivalences 
(Theorem \ref{thm-irreducibility}). 
It follows that $\Hom_{\G_{OK}}(x,y)$ is non-empty,
for any two objects $x$, $y$ in $Ob(\G)_n$.
%the groupoid $\G_{OK}$ acts {\em transitively} on the set $B$. 
Identify the automorphisms in $\G$ of $(S,v,H)$ with the subgroup $\Gamma_v$ 
of isometries of $K_{top}(S)$ which stabilize $v$.
It remains to prove that $\Aut_{\G_{OK}}(S,v,H)$ contains  $\Gamma_v$, 
for some object $(S,v,H)\in Ob(\G)_n$, for each $n\geq 1$.
This is done in two steps:

Step 1: 
We choose $v$ to be the class of the ideal sheaf of a length
$n$ subscheme of $S$, so that $\M_H(v)$ is the Hilbert scheme $S^{[n]}$
(regardless of $H$). 
We prove in this step that $\Aut_{\G_{OK}}(S,v,H)$ contains the index $2$
subgroup $\Gamma_v^{cov}\subset \Gamma_v$, consisting of orientation
preserving isometries
(Definition \ref{def-covariant-subgroups}).

Sub-step 1.1: 
The class $v$ belongs to the sublattice 
$U\subset K_{top}(S)$,
consisting of classes $x$ with vanishing $c_1(x)$. 
$U$ is isometric to the rank $2$ hyperbolic lattice,
and  $K_{top}(S)$ is the orthogonal direct sum $H^2(S,\Integers)\oplus U$. 
%Furthermore, the class $v$ belongs to $U$. 
Consequently, the isometry group $O(H^2(S,\Integers))$
of $H^2(S,\Integers)$ is contained in
the stabilizer $\Gamma_v$. 
Let $\Gamma_0^{cov}$ be the index two subgroup
of $O(H^2(S,\Integers))$ consisting of orientation preserving isometries. 
Every isomerty in $\Gamma_0^{cov}$ is a monodromy operator for some 
deformation of the $K3$ surface $S$ 
(Corollary \ref{cor-Mon-S}). 
This follows easily from the Torelli Theorem for K3 surfaces.
Thus, $\Gamma_0^{cov}$ is contained in 
$\Aut_{\G_{mon}}(S,v,H)$  and in particular in $\Aut_{\G_{OK}}(S,v,H)$ 
(Theorem \ref{thm-conj-holds-for-hilbert-schemes-and-isometries-of-H-2}). 

Sub-step 1.2:
$\Gamma_v^{cov}$ is shown to be generated by 
$\Gamma_0^{cov}$ and reflections $\rho_u$ in
$-2$ classes $u\in K_{top}(S)$ of topological line bundles
with $c_1(u)$ a primitive class, and $c_1(u)^2=2n-4$
(Proposition \ref{prop-compare-Gamma-v}). Such a class $u$
belongs to $v^\perp$. 
%$\Gamma_0^{cov}$ acts transitively on the set of such classes $u$
Any two such reflections are conjugate under $\Gamma_0^{cov}$
(Lemma \ref{lemma-primitive-embeddings-of-rank-2-lattices}).
The inclusion $\Gamma_v^{cov}\subset \Aut_{\G_{OK}}(S,v,H)$ would follow,
once we find one such reflection in $\Aut_{\G_{OK}}(S,v,H)$.
Proposition \ref{reflection-of-Hilbert-schemes-with-respect-to-topological-lb}
exhibits such a reflection $\rho_u$
in $\Aut_{\G_{OK}}(S,v,H)$, for each $n$
(Theorem 
\ref{thm-class-of-correspondence-in-stratified-elementary-trans}
exhibits another such reflection). 
%Hence, $\Gamma_v^{cov}$ is contained in $\Aut_{\G_{OK}}(S,v,H)$. 

Step 2:
It remains to prove that for each $n\geq 1$ there exists 
$(S,v,H)\in Ob(\G)_n$ and $g\in \Aut_{\G_{OK}}(S,v,H)$, such that $cov(g)=1$.
This follows from Theorem 
\ref{thm+2-reflection-induces-an-isomorphism}
(see also Theorem \ref{thm-reflection-sigma-satisfies-main-conj}).
This completes the proof of 
parts \ref{thm-item-gamma-g-v-is-an-isomorphism} and 
\ref{thm-item-invariance-of-the-universal-class} 
of Theorem \ref{thm-action} (as well as Theorem 
\ref{thm-trancendental-reflections}).

Proof of part \ref{thm-item-gamma-g-in-K-alg} of Theorem \ref{thm-action}:
Assume that $g\in \Hom_\G((S_1,v_1,H_1),(S_2,v_2,H_2))$ is an
orientation-preserving Hodge isometry.
%Note, that when the isometry $g$ in the Theorem is a Hodge isometry, 
Then $g$ is induced by an auto-equivalence
$\Phi:D(S_1)\rightarrow D(S_2)$
%, or by the composition of such a $\Phi$ 
%with the duality functor 
(see \cite{mirror-symmetry-k3} or Theorem \ref{thm-HLOY} below). 
Furthermore, $\Phi$ is determined by an object $F$ in
the derived category of $S_1\times S_2$ \cite{orlov}.
The class in $K_{alg}[\M_{H_1}(v_1)\times\M_{H_2}(v)]$ representing 
$\gamma(g,v_1)$ is constructed in terms of $F$ in part 
\ref{lemma-item-Chow-theoretic-formula-for-gamma} of Lemma
\ref{lemma-Chow-theoretic-formula-for-gamma}. 
%when $g$ is orientation preserving, and in section
%\ref{sec-the-serre-duality-automorphism-of-a-moduli} in general.

Proof of part \ref{thm-item-invariance-of-chern-classes} 
of Theorem \ref{thm-action}:
The Chern classes $c_k(TM_H(v))$ are invariant under the duality involution
$D_\M$, given in (\ref{eq-duality-M}), as they vanish for odd $k$.
Hence, part \ref{thm-item-invariance-of-chern-classes} holds for all
morphisms $g$ of $\G_{Mukai}$. Part 
\ref{thm-item-invariance-of-chern-classes} holds for morphisms $g$ of 
$\G_{mon}$, as the Chern classes of $TX$ determine global flat sections
of the local system of integral cohomology over the base of any smooth 
deformation of a smooth and compact variety $X$. 
Part \ref{thm-item-invariance-of-chern-classes} follows, since 
$\G_{Mukai}$ and $\G_{mon}$ generate $\G$. 
\EndProof

\begin{rem}
%\label{rem-G-is-generated-by-G-Mukai-and-G-mon}
{\rm
The above proof shows that the groupoid $\G$ is generated by the
sub-groupoids $\G_{mon}$ and $\G_{Mukai}$. 
This fact seems to be related to the conjecture
that equivalences of derived categories of $K3$ surfaces
induce orientation preserving isometries of Mukai lattices. 
We show that an equivalence of derived categories, which violates this 
conjecture, can not induce an isomorphism between any two moduli spaces of
stable sheaves of dimension $\geq 4$ 
(Lemma \ref{lemma-on-the-sign-conjecture}).
}
\end{rem}

%*************************************************
% Outline of the proof of monodromy Theorem 
%*************************************************
\subsubsection{Reduction of the proof of Theorem 
\ref{introduction-thm-Gamma-v-acts-motivicly}}
%{\bf Outline of the proof of Theorem 
%\ref{introduction-thm-Gamma-v-acts-motivicly}:}
We first prove that $\monrep(\Gamma_v)$ 
is a subgroup of $\Mon(\M_H(v))$. 
This follows from the proof of Theorem \ref{thm-action} as follows.
The character $cov$ of the isometry group of $K_{top}(S)$, used in 
equation (\ref{introduction-eq-monodromy-is-gamma-times-cov}),
extends to a functor 
\begin{equation}
\label{eq-functor-cov}
cov \ : \G \ \ \ \rightarrow \ \ \ \Integers/2\Integers,
\end{equation}
where $\Integers/2\Integers$ is the category with one object, whose
automorphism group is $\Integers/2\Integers$ 
(Remark \ref{rem-natural-orientation}). We use the functor $cov$ to
modify the functor $\H$ in (\ref{eq-the-representation-functor-H}) 
and obtain a second functor 
\begin{equation}
\label{eq-functor-H-mon}
\H_{mon} \ : \ \G \ \ \ \rightarrow \ \ \ \A.
\end{equation}
Like $\H$, the functor $\H_{mon}$ sends the object $(S,v,H)$ of $\G$ to 
the algebra $H^*(\M_H(v),\Integers)_{free}$.
The isometry 
$g \in \Hom_{\G}((S_1,v_1,H_1),((S_2,v_2,H_2))$ 
is sent by $\H_{mon}$ to the isomorphism 
\begin{equation}
\label{eq-image-of-g-by-H-mon} 
(D_{\M_{H_2}(v_2)})^{cov(g)} \circ \gamma_{g,v_1} 
\ : \ H^*(\M_{H_1}(v_1),\Integers)_{free} \ \ \ \rightarrow \ \ \ 
H^*(\M_{H_2}(v_2),\Integers)_{free},
\end{equation}
which is denoted $\monrep(g)$.
The representation $\monrep$ given in 
(\ref{introduction-eq-monodromy-is-gamma-times-cov})
is thus the restriction of
$\H_{mon}$ to the automorphism group 
$\Aut_{\G}((S,v,H))$. 
The monodromy operator $\monrep(g)$
corresponds to the class 
\begin{equation}
\label{eq-class-mon-g-v}
\monrep(g,v_1) \ \ := \ \ 
\left\{
\begin{array}{ccc}
\gamma(g,v_1) & \mbox{if} & cov(g)=0,
\\
c_m\left(-\pi_{13_!}\left\{
\pi_{12}^!
\left[
(1\otimes Dg)(e_{v_1})
\right]
\otimes \pi_{23}^!(e_{v_2})
\right\}
\right) & \mbox{if} & cov(g)=1,
\end{array}
\right.
\end{equation}
where we kept the notation used in 
(\ref{eq-intro-class-gamma-g-v}).

The proof of Theorem \ref{thm-action} shows, that each isomorphism 
$\monrep(g)$ in 
(\ref{eq-image-of-g-by-H-mon}) is the composition of 
a finite set of isomorphisms, each either 
induced by an isomorphism of moduli spaces corresponding to a
morphism in $\G_{Mukai}$, 
or a monodromy operator corresponding to a morphism in
$\G_{mon}$. An automorphism of an
irreducible holomorphic symplectic manifold induces 
a monodromy operator (Lemma \ref{lemma-equivariant-family-implies-monodromy}).
It follows that each isomorphism
$\monrep(g)$  in (\ref{eq-image-of-g-by-H-mon}) is a monodromy operator. 

We use results of Verbitsky to prove that $\monrep(\Gamma_v)$ is a normal 
subgroup of finite index (Lemma \ref{lemma-bounding-K}).
\EndProof

\subsubsection{The structure of the rest of the paper}
We summarize the highlights of each of the following sections. 
More detailed summaries appear at the beginning of each section.
In section 
\ref{section-open-questions} we define the notion of a 
{\em local monodromy operator}.
We then relate and discuss two open questions: 
1) a local Torelli type question \ref{question-local-monodromy},
and 2) the question whether the monodromy operator $\monrep(g)$ 
of Theorem \ref{introduction-thm-Gamma-v-acts-motivicly} is local, 
whenever $g$ is a Hodge-isometry. Partial results are
summerized in Theorem \ref{thm-image-of-mu-is-generated-by-local-monodromies}.
Section \ref{sec-automorphisms-of-cohomology-ring} 
begins with background material
(sections \ref{sec-universal-sheaves}, \ref{sec-weight-2-hodge-structure},
and \ref{sec-class-of-diagonal}).
%universal sheaves
%(section \ref{sec-universal-sheaves}), 
%the Mukai lattice and 
%the weight 2 Hodge structure of moduli spaces of sheaves (section
%\ref{sec-weight-2-hodge-structure}), and 
%the class of the diagonal of a moduli space 
%(section \ref{sec-class-of-diagonal}).
We then 
restate Theorem \ref{thm-action}, about symmetries of the collection 
of moduli spaces of sheaves on $K3$ surfaces, without using the 
$K$-theoretic universal classes 
(Theorem \ref{thm-trancendental-reflections}).
We also prove that the isomorphisms 
$\gamma_{g,v}$ compose as stated in part
\ref{thm-item-multiplicative-properties} of Theorem \ref{thm-action}
(Lemma \ref{lemma-conjecture-holds-for-compositions}). 
%we review the relationship between the Mukai lattice of a K3 surface $S$ and 
%the weight 2 Hodge structure of moduli spaces of sheaves on $S$. 
%We then 
In section \ref{sec-monodromy-of-hilbert-schems} 
we relate the representation $\monrep$ of 
Theorem \ref{introduction-thm-Gamma-v-acts-motivicly} 
to a Hodge-theoretic representation of Verbitsky
(sections \ref{sec-monodromy-operators-acting-trivially-on-H-2} 
and \ref{sec-comparison-with-verbitsky}).
We use Verbitsky's results to show that the reflection group
$\W$ of $H^2(\M(v),\Integers)$ maps via the homomorphism
$\mu$ of Corollary \ref{cor-W-consists-of-monodromies} 
onto a normal subgroup of finite index 
in the full monodromy group $\Mon(\M(v))$
(Lemma \ref{lemma-bounding-K}). 
We bound the index $[\mu(\W):\Mon(\M(v))]$ in section
\ref{sec-chern-classes-of-tangent-bundle} 
using an efficient set of generators for the cohomology ring
$H^*(\M(v))$ and the relation of the Chern classes of $\M(v)$ 
to these generators
(Lemma \ref{lemma-chern-classes-project-to-generators}).
%\ref{sec-Mon-2-is-large}. 
%state the results about the monodromy of Hilbert 
%schemes $S^{[n]}$ of a K3 surface $S$ (???). These include the monodromy 
%representation of the stabilizer $\Gamma_v$ 
%(Theorem \ref{introduction-thm-Gamma-v-acts-motivicly}).
%We prove, that the image of $\Gamma_v$ has finite index 
%in the whole monodromy group (Lemma \ref{lemma-bounding-K}).
%In addition, we study constraints on the monodromy group, 
%imposed by the invariance of the Chern classes of $S^{[n]}$ 
%(section \ref{sec-chern-classes-of-tangent-bundle}). 
In section \ref{seq-equivalences-of-derived-categories} 
we lift the topological formula 
(\ref{eq-class-mon-g-v})
%(\ref{eq-intro-class-gamma-g-v}), 
for the monodromy operator $\monrep(g)$, to a Chow theoretic formula,
%(\ref{eq-Chow-theoretic-formula-for-gamma}),
whenever $\monrep(g)$ maps the weight $2$ Hodge structure of 
$\M(v)$ to that of $\M(g(v))$
(Lemma \ref{lemma-Chow-theoretic-formula-for-gamma}). 
This lift is carried out using known results about equivalences of 
derived categories of K3 surfaces. 
In section \ref{sec-conceptual-interpretation-of-the-formula-gamma-g} 
we define the groupoid $\G_{Mukai}$, whose morphisms come from 
Fourier-Mukai functors preserving the stability of all sheaves
in a moduli space. 
In section \ref{sec-the-small-monodromy} we prove the main Theorems
\ref{thm-action} and \ref{introduction-thm-Gamma-v-acts-motivicly} 
in the case of monodromy operators of $S^{[n]}$, arising from deformations
of the surface $S$. 
In section \ref{sec-examples-of-hodge-isometries}
we verify Theorems \ref{thm-action} and 
\ref{introduction-thm-Gamma-v-acts-motivicly} in examples, which do {\em not} 
arise from deformations of the surface $S$. The two key examples are
Proposition 
\ref{reflection-of-Hilbert-schemes-with-respect-to-topological-lb} and 
Theorem \ref{thm+2-reflection-induces-an-isomorphism}, which suffice
for the proof of Theorems \ref{thm-action} and 
\ref{introduction-thm-Gamma-v-acts-motivicly}. 
%we review work of Seidel-Thomas on the reflections of the derived category 
%with respect to spherical objects. 
%We then state Theorems 
%\ref{thm-class-of-correspondence-in-stratified-elementary-trans} 
%and \ref{thm-reflection-sigma-satisfies-main-conj} for the two sequences 
%of isometries $\tau_n$ and $\sigma_n$.
%We recall the orientation character of
%the isometry group of the Mukai lattice. We also remark on 
%the orientation reversing isometry induced by the duality operator.
In section \ref{sec-structure-of-stabilizer} we find a set of generators
for the stabilizer $\Gamma_v$ of the Mukai vector $v$ of $S^{[n]}$. 
This set of generators consist of isometries in $\Gamma_v$, 
for which Theorems \ref{thm-action} and 
\ref{introduction-thm-Gamma-v-acts-motivicly} 
were verified earlier. 
%in sections
%\ref{sec-the-small-monodromy} and 
%\ref{sec-examples-of-hodge-isometries}.

\smallskip
{\em Acknowledgments:} It is a pleasure to acknowledge fruitful conversations 
with Andrei C\u{a}ld\u{a}raru, Igor Dolgachev, Jim Humphreys, 
Daniel Huybrechts, Eduard Looijenga, Vikram Mehta, Shigeru Mukai, 
and Hiraku Nakajima.
Special thanks are due to Kota Yoshioka for
several invaluable conversations and correspondences. 
The author is grateful for the three referees for their detailed reviews and
insightful suggestions, which led to a substantial improvement of
this paper.

In the first submitted version of this paper, Theorem 
\ref{thm-action} was a conjecture and 
Theorem \ref{introduction-thm-Gamma-v-acts-motivicly}
was proven using {\em local} monodromy
operators (Theorem \ref{thm-image-of-mu-is-generated-by-local-monodromies}).
The proof relied on Theorems 
\ref{thm-class-of-correspondence-in-stratified-elementary-trans} and 
\ref{thm-reflection-sigma-satisfies-main-conj} below, which are 
both proven in \cite{markman-part-two}. 
The current proofs of Theorems \ref{thm-action}  and 
\ref{introduction-thm-Gamma-v-acts-motivicly}
are independent of 
Theorems \ref{thm-class-of-correspondence-in-stratified-elementary-trans} and 
\ref{thm-reflection-sigma-satisfies-main-conj} 
(and thus of the results of \cite{markman-part-two}). 
We thank one of the referees for suggesting greater reliance on
Yoshioka's results.

%***************************************************************
% Open questions about local and global monodromy operators
%***************************************************************
\section{Torelli questions and local monodromy operators}
\label{section-open-questions}

Let $X$ be an irreducible holomorphic symplectic manifold 
deformation equivalent to the
Hilbert scheme $S^{[n]}$ of length $n$ subschemes of a $K3$ surface $S$.
The homomorphism $\mu$ in Corollary \ref{cor-W-consists-of-monodromies}
lifts many isometries of $H^2(X,\Integers)$ to monodromy operators in
$\Aut H^*(X,\Integers)_{free}$. When the isometry is of Hodge type, 
another such lift is the local monodromy operator
(Definition \ref{def-local-monodromy-operator}), defined when 
an affirmative answer is known to a strong version of the 
local Torelli Question \ref{question-local-monodromy}. 
We define this second lift, conjecture the equality of the two lifts, 
and describe examples where equality is known to hold (Theorem
\ref{thm-image-of-mu-is-generated-by-local-monodromies}).

\begin{conj}
\label{conj-Mon-injects-into-Mon-2}
The restriction homomorphism
$\Mon(X)\rightarrow \Mon^2(X)$ is injective.
\end{conj}

The conjecture implies that the monodromy group $\Mon(\M_H(v))$ 
is equal to the image $\monrep(\Gamma_v)$ of the homomorphism 
(\ref{introduction-eq-monodromy-is-gamma-times-cov}) (use also 
Corollary \ref{cor-W-consists-of-monodromies} 
and the containment $\Mon^2\subset \W$ proven in 
\cite{markman-integral-generators}). 
The conjecture holds when $\dim(X)=2$ or $4$, because then
$H^2(X,\RationalNumbers)$ generates the cohomology ring of $X$.
An upper-bound for the order of the kernel of the
homomorphism $\Mon(X)\rightarrow \Mon^2(X)$, when $\dim(X)>4$,
is provided in section \ref{sec-chern-classes-of-tangent-bundle}.
Conjecture \ref{conj-Mon-injects-into-Mon-2}
is motivated by the following result of Beauville and Verbitsky:

\begin{prop}
\label{prop-beauville-automorphisms}
%(\cite{beauville-automorphisms} Proposition 10)
Let $X$ be an irreducible holomorphic symplectic manifold 
deformation equivalent to 
the Hilbert scheme $S^{[n]}$, $n\geq 1$, of a $K3$ surface $S$.
A holomorphic automorphism of $X$, 
which acts as the identity on $H^2(X,\Integers)$,  is the identity
automorphism.
\end{prop}

\noindent
{\bf Proof:}
The case $X=S^{[n]}$ was proven by Beauville 
(\cite{beauville-automorphisms} Proposition 10).
The statement follows, for deformations of $S^{[n]}$,
by Corollary 6.9 in \cite{kaledin-verbitsky}.
\EndProof

\medskip
We describe below explicit examples of pairs, of a moduli space $\M_H(v)$
and a monodromy operator $\tilde{g}\in \Aut[H^*(\M_H(v),\Integers)]$, 
inducing an isometry $g$ in $\W\subset \Mon^2(\M_H(v))$. In these examples 
the equality $\tilde{g}=\mu(g)$, implied by 
Conjecture \ref{conj-Mon-injects-into-Mon-2}, is either not known to the 
author (in the very interesting Example 
\ref{example-fiber-square-of-Hilbert-scheme-over-symmetric-product}), 
or is proven with further effort in \cite{markman-part-two}
(Theorems \ref{thm-class-of-correspondence-in-stratified-elementary-trans} 
and \ref{thm-reflection-sigma-satisfies-main-conj} below). 

Let $g$ be an isometry in the subgroup $\W\subset OH^2(X,\Integers)$ given in
(\ref{eq-W}). Assume that  $g$ preserves the Hodge srtucture. 
Then it may be possible to lift 
$g$ to a {\em local} monodromy operator in $\Mon(X)$. 
We formulate this statement as Question \ref{question-local-monodromy} below. 
Recall that the differential of the period map is an isomorphism 
$H^1(X,TX)\rightarrow H^{2,0}(X)^*\otimes H^{1,1}(X)$
\cite{beauville-varieties-with-zero-c-1}. 
Hence, $g$ acts on the infinitesimal deformation space of $X$. 
If $g$ has finite order, 
then we can choose a $g$-invariant simply connected subset 
$B\subset H^1(X,TX)$, open in the complex topology, and a 
local universal family $\pi:\X\rightarrow B$,
with special fiber $X$ over the point $0\in B$. 
If the order of $g$ is infinite, we may not assume $B$ to be $g$-invariant,
but the discussion below generalizes,
replacing $B$ by the union $B'\cup g(B')$ 
of an open subset $B'\subset B$, containing $0$,
such that $g(B')$ is contained in $B$. We  assume, for
simplicity, that $g$ has finite order. 
The local systems $R^2_{\pi_*}\Integers$ and $g^*R^2_{\pi_*}\Integers$
over $B$ are trivial, 
as $B$ is assumed simply connected. The trivializations conjugate
any flat isomorphism
$R^2_{\pi_*}\Integers\cong g^*R^2_{\pi_*}\Integers$ 
to an automorphism of $H^2(X,\Integers)$.
An affirmative answer to the following question would constitute
a version of local Torelli, stronger than the known local Torelli theorem 
\cite{beauville-varieties-with-zero-c-1}:

\begin{question}
\label{question-local-monodromy}
Is there always a $g$-invariant connected open neighborhood $U$ of $0$ 
in $B$, a proper $g$-invariant analytic subset $A\subset U$, 
and a fiber preserving isomorphism 
$f:\restricted{\X}{[U\setminus A]}\rightarrow 
\restricted{g^*\X}{[U\setminus A]}$ 
between the restrictions
of the families $\X$ and $g^*\X$ to $U\setminus A$, inducing an isomorphism
$(R^2_{\pi_*}\Integers\restricted{)}{[U\setminus A]}\cong 
(g^*R^2_{\pi_*}\Integers\restricted{)}{[U\setminus A]}$, which  
is conjugated to the automorphism $g$ of $H^2(X,\Integers)$ 
via the trivializations?
\end{question}

When the answer is affirmative we set $U_0:= U\setminus (A\cup U^g)$,
where $U^g$ is the subset of fixed-points.
We get the quotient family
$(\restricted{\X}{U_0})/f\rightarrow U_0/g$ and 
$g$ generates its monodromy subgroup in $\Mon^2(X)$.
%(we may further assume that $U$ is simply connected). 
The induced isomorphism of local systems
$f_*:(R^*_{\pi_*}\Integers)_{free}\rightarrow 
g^*(R^*_{\pi_*}\Integers)_{free}$
extends over the whole of $U$, due to the triviality of both local systems. 
Denote by $\tilde{g}\in \Aut[H^*(X,\Integers)]$
the value of $f$ at the fiber over the point $0\in U$. 
Then $\tilde{g}$ is a monodromy operator, which restricts to
$g$ on $H^2(X,\Integers)$ 
(Lemma \ref{lemma-equivariant-family-implies-monodromy}). 
Conjecture \ref{conj-Mon-injects-into-Mon-2} 
implies that $\tilde{g}$ is equal to 
the monodromy operator $\mu(g)$ in
Corollary \ref{cor-W-consists-of-monodromies} 
\begin{equation}
\label{eq-g-tilde-equal-mu-g}
\tilde{g} \ \ \ = \ \ \ \mu(g).
\end{equation}

If the isomorphism $f$ satisfying the conditions of
Question \ref{question-local-monodromy} exists, then it is unique 
by Proposition \ref{prop-beauville-automorphisms}. 
The uniquesness justifies the following:

\begin{defi}
\label{def-local-monodromy-operator}
{\em
%Let $\tilde{g}$ be a graded ring automorphism of $H^*(X,\Integers)$,
%whose restriction to $H^2(X,\Integers)$ is a
%Hodge-isometry of finite order. 
%We call $\tilde{g}$ a {\em local monodromy operator}, if 
%it is induced, via the above construction, 
%from an isomorphism $f$ satisfying the conditions of
%Question \ref{question-local-monodromy}.
%
The automorphism $\tilde{g}\in \Aut[H^*(X,\Integers)]$ 
is {\em the local monodromy operator, which restricts to
$g$ on $H^2(X,\Integers)$}. 
Let $\ell$ be an automorphism of $H^*(X,\Integers)$
and denote its restriction to
$H^2(X,\Integers)$ by $\ell_2$. 
We call $\ell$ a
{\em local monodromy operator}, if the local monodromy operator, 
whose restriction to $H^2(X,\Integers)$ is $\ell_2$, exists and is
equal to $\ell$.
}
\end{defi}

Let $\M_H(v)$ be a moduli space as in Corollary 
\ref{cor-W-consists-of-monodromies}. 
An isometry of $H^2(\M_H(v),\Integers)$ is called a 
{\em surface monodromy operator}, 
if it arises from deformations of $\M_H(v)$ as a moduli space of sheaves on 
a $K3$ surface (Definition \ref{def-deformation-equivalent-triples}). 
When $c_1(v)$ is a multiple of an ample class, then 
surface monodromy operators are in a natural bijection with 
signed isometries of $H^2(S,\Integers)$, leaving invariant the class $c_1(v)$
(Theorem \ref{thm-conj-holds-for-hilbert-schemes-and-isometries-of-H-2}). 
If $g$ is a surface monodromy operator, preserving the Hodge structure of
$H^2(\M_H(v),\Integers)$, then an affirmative answer to Question 
\ref{question-local-monodromy} follows from the Torelli Theorem for $K3$ 
surfaces, 
and the equality (\ref{eq-g-tilde-equal-mu-g}) follows from the 
proof of Theorem \ref{introduction-thm-Gamma-v-acts-motivicly}. 

The examples in which equality (\ref{eq-g-tilde-equal-mu-g}) was proven 
suffice to conclude the following Theorem.
Let $n\geq 2$, $\W$ the reflection group (\ref{eq-W}) of 
$H^2(S^{[n]},\Integers)$, and 
and $\mu : \W \rightarrow \Mon(S^{[n]})$ the homomorphism in
Corollary \ref{cor-W-consists-of-monodromies}.

\begin{thm} \cite{markman-part-two}
\label{thm-image-of-mu-is-generated-by-local-monodromies}
The image $\mu(\W)$ in $\Mon(S^{[n]})$ 
is generated by monodromy operators, inducing 
reflections in $+2$ or $-2$ classes of $H^2(S^{[n]},\Integers)$,
%(Definition \ref{def-local-monodromy-operator}), 
each of which is a local monodromy operator for a suitable choice 
of a complex structure on the $K3$ surface $S$.
\end{thm}

Theorem \ref{thm-image-of-mu-is-generated-by-local-monodromies} 
follows from Theorems 
\ref{thm-class-of-correspondence-in-stratified-elementary-trans} and 
\ref{thm-reflection-sigma-satisfies-main-conj} below 
(both proven in \cite{markman-part-two}) via Proposition 
\ref{prop-compare-Gamma-v}
and Theorem \ref{thm-conj-holds-for-hilbert-schemes-and-isometries-of-H-2}
(see Steps 1 and 2 of the proof of Theorem \ref{thm-action}).

O'Grady conjectured 
an affirmative answer to Question \ref{question-local-monodromy}, 
when $g$ is a reflection (\ref{eq-W}) 
by a the Chern class $c_1(L)$ of a line-bundle
$L$ with $(c_1(L),c_1(L))=2$ \cite{ogrady-involutions}. 
He further conjectured, that after a small deformation of the pair $(X,L)$,
the reflection $g$ is induced by a {\em regular involution} of $X$
(the {\em I-conjecture}). 
%This is equivalent to the statement, that 
%the subset $A$ may be chosen not to
%contain the codimention $1$ locus $U^g$ fixed under $g$, because 
%$U^g$ is the locus where the 
%class $c_1(L)$ remains of type $(1,1)$. 
In dimension $2$ the simplest example is 
the regular Galois involution of a double cover of $\PP^2$,
branched along a sextic.

Reflections by negative effective classes can not be induced by a regular 
involution of $X$. 
%(so the subset $A$ in Question \ref{question-local-monodromy}
%must contain the fixed locus $U^g$). 
Following is a simple proof. 
Assume that the line bundle $L$ is effective,    
$(c_1(L),c_1(L))<0$, 
and the isometry $g$ is the reflection with respect to $c_1(L)$. 
Then $g(c_1(L))=-c_1(L)$. If $\kappa$ is a K\"{a}hler class, then
$(\kappa,c_1(L))>0$, because $c_1(L)$ is effective 
(\cite{huybrechts-basic-results} 1.11). Hence, $g(\kappa)$ is not a
K\"{a}hler class and $g$ is not induced by a regular involution. 
A positive answer to Question
\ref{question-local-monodromy} does imply, that the local monodromy 
$\tilde{g}$ is induced by a self-correspondence in
$X\times X$ (see the proof of
\cite{huybrechts-basic-results} Theorem 4.3). 
O'Grady's I-Conjecture is replaced by the question: 

\smallskip
\noindent
{\bf Question:}
{\em Describe the structure of the
self-correspondence for a generic pair $(\X_t,L_t)$,
$t\in U^g$, deformation
equivalent to $(X,L)$.
Use the description to prove the equality (\ref{eq-g-tilde-equal-mu-g}).
}

\medskip
We consider two sequences of reflections by negative Hodge classes, 
which are monodromy operators but {\em not} surface monodromy operators
(Theorem \ref{thm-class-of-correspondence-in-stratified-elementary-trans}
and  Example 
\ref{example-fiber-square-of-Hilbert-scheme-over-symmetric-product}). 

\begin{example}
\label{example-fiber-square-of-Hilbert-scheme-over-symmetric-product}
{\rm
Let $S$ be a $K3$ surface and $n$ an integer $\geq 2$. 
The second cohomology $H^2(S^{[n]},\Integers)$
is isometric to $H^2(S,\Integers)\oplus {\rm span}_\Integers\{\delta\}$,
where $\delta$ is half the class of the big diagonal in $S^{[n]}$
and $(\delta,\delta)=2-2n$. Let $v\in K_{top}(S)$  be the class of the
ideal sheaf of a length $n$ subscheme and $D\in \Gamma_v$ 
the duality involution (\ref{eq-duality-involution}). 
Then the monodromy operator 
$\monrep(D)$, given in (\ref{introduction-eq-monodromy-is-gamma-times-cov}), 
acts on $H^2(S^{[n]},\Integers)$ as 
the reflection $\rho_\delta(x) := x + [(x,\delta)/(n-1)]\delta$. 
%Let us explain why Conjecture \ref{conj-Mon-injects-into-Mon-2} identifies 
%$\monrep(D)$ with a local monodromy operator, as well as with 
%the automorphism induced by a concrete self-correspondence. 
Let $\pi:S^{[n]}\rightarrow S^{(n)}$ be the Hilbert-Chow morphism from the
Hilbert scheme to the symmetric product. The fiber product 
$Z:=[S^{[n]}\times_{\pi} S^{[n]}]\subset [S^{[n]}\times S^{[n]}]$
is reducible but reduced of pure dimension $2n$ 
(part \ref{lemma-item-fiber-product-is-reduced} of Lemma 
\ref{lemma-fiber-product-of-chow-morphism}). 
$Z_*$ is the local monodromy operator 
restricting to $\rho_\delta$ on $H^2(S^{[n]})$ 
(Definition \ref{def-local-monodromy-operator}),
%defined as was the left hand side in equation
%(\ref{eq-g-tilde-equal-mu-g})
by Lemma \ref{lemma-fiber-product-of-chow-morphism}
part \ref{lemma-item-Z-is-the-local-monodromy}.
Conjecture \ref{conj-Mon-injects-into-Mon-2} implies the equality 
$Z_*=\monrep(D)$. 
Equivalently, $Z$ is Poincare dual to the
class $\monrep(D,v)$ given in (\ref{eq-class-mon-g-v}).

The equality $Z_*=\monrep(D)$ 
would follow (by part 
\ref{thm-item-characterization-of-gamma-g} of Theorem \ref{thm-action}), 
independently of Conjecture \ref{conj-Mon-injects-into-Mon-2}, 
if one proves instead that $D_{S^{[n]}}\circ Z_*$ and $D$ satisfy equation  
(\ref{eq-characterization-of-gamma-g}).
Let $e$ be the class in $K_{top}(S^{[n]}\times S)$ 
of the ideal sheaf of the universal subscheme. 
Denote the dual class by $e^\vee$.
Equation (\ref{eq-characterization-of-gamma-g}) 
translates in our case to the equality 
\[
(Z_!\otimes 1)(e) \ \ \ = \ \ \ e^\vee,
\]
where $Z_!\in \End[K_{top}(S^{[n]})]$ is given by
$Z_!(x)=\pi_{2_!}(\pi_1^!(x)\otimes \StructureSheaf{Z})$, and 
$\pi_i$ is the projection from $S^{[n]}\times S^{[n]}$ onto the
$i$-th factor, $i=1,2$ (compare with the equality in part 
\ref{thm-item-Z-maps-universal-class-to-such}
of Theorem \ref{thm-class-of-correspondence-in-stratified-elementary-trans}).
}
\end{example}

\begin{new-lemma}
\label{lemma-fiber-product-of-chow-morphism}
\begin{enumerate}
\item
\label{lemma-item-fiber-product-is-reduced}
The fiber product 
$Z:=[S^{[n]}\times_{\pi} S^{[n]}]$ is reduced, of pure dimension $2n$, and
its irreducible components are in bijection with ordered partitions of $n$. 
\item
\label{lemma-item-Z-is-the-local-monodromy}
$Z_*$ is the local monodromy operator 
restricting to $\rho_\delta$ on $H^2(S^{[n]})$ 
(Definition \ref{def-local-monodromy-operator}).
\end{enumerate}
\end{new-lemma}

\noindent
{\bf Proof:} 
\ref{lemma-item-fiber-product-is-reduced})
Let $\lambda:=(n_1\leq n_2 \leq \dots \leq n_k)$ 
be a partition of $n$ and
$S^{(n)}_\lambda$ the corresponding locally closed subset of $S^{(n)}$.
The dimension of $S^{(n)}_\lambda$ is $2k$. 
The fibers of $\pi$ over points of $S^{(n)}_\lambda$ 
are reduced and irreducible of dimension $n-k$, by
Proposition 2.10 in \cite{haiman}. Hence each partition contributes to 
the fiber product a $2n$-dimensional irreducible component.
Each fiber of the two natural projections $p_i:Z\rightarrow S^{[n]}$ is 
isomorphic to a fiber of $\pi$ and is hence reduced and irreducible. 
$Z$ is reduced, since $p_i$ is a projective morphism with reduced fibers 
from $Z$ to the reduced scheme $S^{[n]}$.
%
%Lemma: Let $X$, $Y$ be noetherian schemes, with $Y$ reduced, and
%$f:X\rightarrow Y$ a projective morphism with reduced fibers.
%Then $X$ is reduced.
%
%{\bf Proof:} The question is local, so we may assume that 
%$f$ is an affine morphism, $X$ is $Spec(\A)$, where $\A$ is a sheaf of
%$\StructureSheaf{Y}$-algebras, and $Y$ is affine. Let $a$ be a nilpotent 
%section of $\A$, $y\in Y$ a closed point, and 
%$m_y\subset \StructureSheaf{Y}$ the ideal sheaf of $y$. Then
%$a$ belongs to $m_y\cdot\A$, because the fiber $f^{-1}(y)$ is reduced.
%We conclude that $a$ belongs to the intersection
%$\cap_{y\in Y}m_y\cdot\A$, which is zero, since $Y$ is reduced

\ref{lemma-item-Z-is-the-local-monodromy})
The homomorphism 
$Z_*\in \End(H^*(S^{[n]},\Integers))$ and the involution 
$\monrep(D)$ both act via the reflection $\rho_\delta$ on 
$H^2(S^{[n]},\Integers)$. 
$\rho_\delta$ acts as an involution of the local deformation space $B$ of 
$S^{[n]}$, 
fixing the divisor $Hilb\subset B$ of deformations as 
Hilbert schemes of K\"{a}hler $K3$ surfaces.
Let $\X\rightarrow B$ be the local universal family 
and $q:B\rightarrow B/\rho_\delta$ the quotient morphism.
An affirmative answer to Question \ref{question-local-monodromy} 
is known in the case of $\rho_\delta$, 
since the Hilbert-Chow morphism is a contraction of type $A_1$ 
as in Namikawa's work \cite{namikawa-deformations}. 
The quotient $B/\rho_\delta$ is the local deformation space of $S^{(n)}$,
the generic fiber of the universal family $\overline{\X}\rightarrow 
B/\rho_\delta$ is smooth, 
and $\X$ is a resolution of the pullback of $\overline{\X}$ to $B$.
It follows that $\X$ and $\rho_\delta^*\X$
restrict to isomorphic families over an open subset $U$ of $B$. 
The open subset $U$ is in fact the complement $B\setminus Hilb$
(possibly after replacing $B$ by a smaller open neighborhood of $0\in B$, 
see the proof of Claim 3 in \cite{namikawa-deformations}). 
%  $S^{(n)}$ is normal and Cohen-Macaulay, by Hochster's Theorem on invariant
% theory. $\pi$ has rational singularities by the Grauert-Riemenschnider 
% vanishing theorem and the triviality of the canonical line-bundle of 
% $S^{[n]}$, and the canonical sheaf of $S^{(n)}$ is trivial by 
% Lemma 5.12 in Kollar-Mori's book ``Birational geometry of algebraic
% varieties'' page 156. It follows that $S^{(n)}$ has Gorenstein singularities
Let $\Z\subset \X\times_B\rho_\delta^*\X$ be the closure of the graph
of the isomorphism. 
The fiber $\Z_0$ of $\Z$ over $0\in B$ is a subscheme of $Z$, since 
$\Z$ is the reduced induced subscheme of the fiber product 
$[\X\times_{q^*\overline{\X}}\rho_\delta^*\X]$ and $Z$ is the fiber of the 
latter over $0$. $Z$ is reduced and hence isomorphic to $\Z_0$.
Consequently $Z_*$ is the local monodromy operator 
restricting to $\rho_\delta$ on $H^2(S^{[n]})$. 
%(Definition \ref{def-local-monodromy-operator}).
\EndProof

\section{Automorphisms of the cohomology ring}
\label{sec-automorphisms-of-cohomology-ring}

We construct in section \ref{sec-universal-sheaves}
a canonical normalization $u_v$ of the Chern character of a universal
(possibly twisted) sheaf over moduli spaces $\M_H(v)$ 
of dimension $\geq 4$ (equation 
(\ref{eq-invariant-normalized-class-of-chern-character-of-universal-sheaf})). 
In section \ref{sec-weight-2-hodge-structure}
we restate the relationship (Theorem \ref{thm-irreducibility}), 
between the Mukai lattice of a K3 surface $S$ and 
the weight 2 Hodge structure of moduli spaces of sheaves on $S$,
using the integral cohomology $H^*(S,\Integers)$, rather than the
$K$-ring $K_{top}(S)$. 
We recall in section \ref{sec-class-of-diagonal},  
that the K\"{u}nneth factors, of the Chern classes of a universal sheaf,
are generators for the cohomology ring of a moduli space of sheaves on
a K3 or abelian surface. In section 
\ref{sec-conjecture} we restate Theorem \ref{thm-action}, 
relating the cohomology rings of pairs of 
moduli spaces of sheaves on K3 surfaces. 
The new statement (Theorem \ref{thm-trancendental-reflections}) 
is in terms of the class $u_v$ rather than the integral class $e_v$. 
%The Theorem leads to an action, of a group of integral isometries of
%the second cohomology of a moduli space, on the full cohomology ring 
%of the moduli space (Corollary \ref{cor-group-homomorphism}). 
%The Theorem has a simple analogue, for sheaves on an elliptic curve,
%which is easily verified in section \ref{sec-elliptic-curve}. 
In section \ref{sec-compositions} we prove that the isomorphisms 
$\gamma_{g,v}$ compose as stated in part
\ref{thm-item-multiplicative-properties} of Theorem \ref{thm-action}.

%******************************************************
% Universal sheaves
%******************************************************
\subsection{Universal sheaves}
\label{sec-universal-sheaves}

Let $S$ be a $K3$ surface, $\LB$ an ample line bundle on $S$, and
$v\in K_{top}(S)$ an effective class (Definition \ref{def-effective-class}). 
Let $\M:=\M_{\LB}(v)$ be the moduli space of $\LB$-stable sheaves with class 
$v$.
We use  the $\LB$-stability of sheaves due to Gieseker, Maruyama, and Simpson
\cite{huybrechts-lehn-book}. 
Recall the construction of semi-universal  families over $\M\times S$
(\cite{mukai-hodge} Theorem A.5).
We start with an open covering $\U:=\{U_i\}$ of $\M$, in the 
complex or \'{e}tale topologies, local universal families 
$\E_i$ over $U_i\times S$, and isomorphisms
\begin{equation}
\label{eq-g-i-j}
g_{ij}:(\E_j\restricted{)}{U_{ij}\times S} \ \ \ \longrightarrow \ \ \ 
(\E_i\restricted{)}{U_{ij}\times S}
\end{equation}
over the intersections $U_{ij}:=U_i\cap U_j$. The co-boundary
$(\delta g)_{ijk}:=g_{ij}g_{jk}g_{ki}$ 
consists of central automorphisms of $(\E_i\restricted{)}{U_{ijk}\times S}$, 
since each $\E_i$ is a family of stable, and hence simple, sheaves.
It follows, that  $(\delta g)_{ijk}$
must be the multiplication by the pull-back $\tilde{\alpha}$ to 
$\M\times S$ of a $2$-cocycle $\alpha\in Z^2(\U,\StructureSheaf{\M}^*)$. 
We get the direct image vector bundles 
$V_i:=p_*(\E_i\otimes \LB^n)$ over $U_i$,  
upon a choice of a sufficiently high power $n$. 
The isomorphisms $g_{ij}$ induce gluing transformations 
$\psi_{ij}$ among the $V_i$'s and $(\delta \psi)_{ijk}$ is
multiplication by $\alpha$. 
Note, that
the $\rho$-th power of $\alpha$ is the co-boundary 
$\delta (\Wedge{\rho}{\psi})$, where $\rho$ is the common rank 
of all the $V_i$. 
The families $\E_i\otimes V_i^*$ glue naturally to a global semi-universal 
family $\F$ of {\em similitude} $\rho$.
The vector bundles $V_i^{\otimes \rho}\otimes (\Wedge{\rho}V_i^*)$
and their subbundles $\Sym^\rho(V_i)\otimes (\Wedge{\rho}V_i^*)$ 
glue to a global vector bundle $W$ and its subbundle $W_+$. 
%Denote by
%$^\rho\sqrt{ch(W^*)}$ the $\rho$-th root in
%$H^*(\M\times S,\RationalNumbers)$, with $1$ as the coefficient of degree $0$.
%Set
%\begin{equation}
%\label{eq-chern-character-of-rational-universal-class}
%ch(\E) \ \ \ := \ \ \ \frac{ch(\F)}{^\rho\sqrt{ch(W^*)}} \ \ 
%\in  \ \ H^*(\M\times S,\RationalNumbers).
%\end{equation}

The projective bundles $\PP{V}_i$ glue naturally to a global 
$\PP^{\rho-1}$-bundle $p:\PP\rightarrow \M$. 
We get a universal quotient sheaf $\widetilde{\E}$ of $p^*\F$. 
The sheaf $\widetilde{\E}$ restricts as $E\otimes \StructureSheaf{}(1)$
to the $\PP^{\rho-1}\times S$ fiber over $E\in \M$. We have also 
a line-bundle $\StructureSheaf{\PP}(\rho)$ over $\PP$, 
which restricts as $\StructureSheaf{}(\rho)$ to each fiber of $p$ and such 
that $p_*\StructureSheaf{}(\rho)$ is $W^*_+$. 
%The following equality holds (see \cite{markman-diagonal} section 3)
There exists a unique class 
\begin{equation}
\label{eq-chern-character-of-rational-universal-class}
ch(\E)
\end{equation}
in $H^*(\M\times S,\RationalNumbers)$ 
satisfying the following equality
\begin{equation}
\label{eq-chern-character-of-universal-sheaf-over-PP}
(p\times id_S)^*ch(\E) \ \ \ = \ \ \ 
ch(\widetilde{\E})\exp[c_1(\StructureSheaf{}(-\rho))/\rho],
\end{equation}
since the right hand side restricts to a trivial class on each fiber
of $\PP\times S$ over a point of $\M\times S$.
The classes $ch(\E)$ and $\widetilde{\E}$ are canonical, 
up to a product by the pullback of a class 
$\exp(\ell)\in H^*(\M,\RationalNumbers)$, 
for $\ell\in H^{1,1}(\M,\RationalNumbers)$.

We introduce next a {\em canonical normalization} of the
Chern character of a universal sheaf $\E$ over $\M_{\LB}(v)\times S$. 
Assume, that $\dim\M_{\LB}(v)>2$. 
Let $\eta_v\in H^2(\M_{\LB}(v),\RationalNumbers)$ be the class
satisfying $(v,v)\cdot \eta_v=\theta_v(v)$, where
$\theta_v$ is given in (\ref{eq-introduction-mukai-homomorphism}). Set
\begin{equation}
\label{eq-invariant-normalized-class-of-chern-character-of-universal-sheaf}
u_v \ \ := \ \ 
ch(\E)\cdot \pi_1^*\exp(\eta_v)\cdot \pi_2^*\sqrt{td_S}.
\end{equation}
Then $u_v$ is a natural class in $H^*(\M_{\LB}(v)\times S,\RationalNumbers)$,
which is independent of the choice of a universal sheaf. 
The appearance of the factor $\sqrt{td_S}$ is explained in section
\ref{sec-weight-2-hodge-structure}.
If a universal sheaf $\E$ does not exist over $\M_{\LB}(v)$, 
replace $ch(\E)$ in 
(\ref{eq-invariant-normalized-class-of-chern-character-of-universal-sheaf})
by the class (\ref{eq-chern-character-of-rational-universal-class}). 
Note, that $u_v$ is the normalization of the Chern character of the 
universal sheaf, which produces the equality $\theta_v(v)=0$,
when $u_v$ is substituted for $ch(\E)$ in the cohomological
translation (\ref{eq-mukai-homomorphism}) of the definition
(\ref{eq-introduction-mukai-homomorphism}) of $\theta_v$ 
(because $(v,v)=-\chi(v^\vee\otimes v)$).

{\bf Twisted sheaves:}
The data $(\E_i,g_{ij})$ above is an example of a twisted sheaf.

\begin{defi}
\label{def-twisted-sheaves}
{\rm
Let $X$ be a scheme or a complex analytic space,
$\U:=\{U_i\}_{i\in I}$ a covering, open in the complex or \'{e}tale topology,
and $\alpha\in Z^2(\U,\StructureSheaf{X}^*)$ a \v{C}ech $2$-cocycle.
An {\em $\alpha$-twisted sheaf} consists of sheaves $\E_i$ of
$\StructureSheaf{U_i}$-modules over $U_i$,
for all $i\in I$, and isomorphisms $g_{ij}$ as in 
(\ref{eq-g-i-j}) satisfying the conditions:

(1) $g_{ii}=id$, \ \ \ 
(2) $g_{ij}=g_{ji}^{-1}$,  \ \ \ 
(3) $g_{ij}g_{jk}g_{ki}=\alpha_{ijk}\cdot id.$ 

\noindent
The $\alpha$-twisted sheaf is {\em coherent}, if the $\E_i$ are.
}
\end{defi}

The class of $\alpha$-twisted sheaves, together with the obvious notion of 
homomorphism, is an abelian category, denoted by ${\frak Mod}(X,\alpha)$,
and referred to as the 
{\em category of $\alpha$-twisted sheaves}. The analogous 
category of $\alpha$-twisted coherent sheaves is denoted by 
${\frak Coh}(X,\alpha)$. 
The categories ${\frak Mod}(X,\alpha)$ and ${\frak Coh}(X,\alpha)$ depend, 
up to equivalence,  only on the class of $\alpha$ in the cohomology group
$H^2(X,\StructureSheaf{X}^*)$, using the analytic or \'{e}tale topology 
\cite{calduraru-thesis}. 
We will always 
assume, that $\alpha$ represents a class in the {\em Brauer group} of $X$,
i.e., in the image of the connecting homomorphism
$H^1(X,{\rm PGL}(n))\rightarrow H^2(X,\StructureSheaf{X}^*)$ of the 
short exact sequence
\[
0 \rightarrow \StructureSheaf{X}^* \rightarrow {\rm GL}(n) \rightarrow 
{\rm PGL}(n)\rightarrow 0,
\]
for some $n\geq 2$. 
%Let $V_i$ be the vector bundle on $U_i$, introduced 
%in the construction of the 
%class (\ref{eq-chern-character-of-rational-universal-class}).
%The projective bundles bundles $\PP{V}_i$ glue naturally to a global 
%$\PP^{\rho-1}$-bundle $\PP$ over the moduli space $\M$. 
The class $[\alpha]$, 
of the cocycle $\alpha$ of the twisted universal sheaf, 
belongs to the Brauer group, since it 
is the image of the class of the projective bundle $\PP$ defined
above over the moduli space $\M$. 
The bounded derived category of complexes of 
$\alpha$-twisted sheaves on $X$ with coherent cohomology is denoted by 
$D^b_{coh}(X,\alpha)$ (see \cite{calduraru-thesis} for its construction). 

The twisted universal sheaf defines an object in
$D^b_{coh}(\M_H(v)\times S,\pi_1^*\alpha)$, which we denote by $\E_v$. 
We associate  next to $\E_v$ a universal class $e_v$ in 
$K_{top}(\M_H(v)\times S)$. 
%We can always cover the moduli space $\M$ by open subsets $U_i$, 
%in the analytic topology, such that a universal sheaf $\E_i$ 
%exists over each $U_i\times S$. The $\E_i$ always glue to a
%twisted universal sheaf
%(Definition \ref{def-twisted-sheaves}). 
%The twisting is captured by a class $\alpha$ in the Brauer group;
%a subgroup of certain torsion elements of the \v{C}ech cohomology
%$H^2_{an}(\M\times S,\StructureSheaf{\M_H(v)\times S}^*)$, 
%calculated using the analytic topology 
%(see section \ref{sec-universal-sheaves}). 
The $\alpha$-twisted universal sheaf $\E_v$ determines a class in 
the $K$-group of $\alpha$-twisted holomorphic vector bundles 
$K_{hol}(\M_H(v)\times S)_{\alpha}$. The latter, 
and its topological analogue $K_{top}(\M_H(v)\times S)_{\alpha}$,
are defined as in the untwisted case. 
$K_{top}(\M_H(v)\times S)_{\alpha}$ depends only on the image 
$\delta(\alpha)$ of $\alpha$
in $H^3(\M_H(v)\times S,\Integers)$, under the connecting homomorphism of the
exponential sequence. We proved that the cohomology 
$H^i(\M_H(v),\Integers)$ vanishes for odd $i$
\cite{markman-integral-generators}.
It follows that $\delta(\alpha)$ is trivial and 
$K_{top}(\M_H(v)\times S)_{\alpha}$ is isomorphic to 
the untwisted group $K_{top}(\M_H(v)\times S)$,
canonically up to tensorization of the latter with a topological line-bundle. 

\begin{defi}
\label{def-e-v}
The universal class $e_v$ is the image of $\E_v$ under the composition
\[
K_{hol}(\M_H(v)\times S)_{\alpha}\rightarrow
K_{top}(\M_H(v)\times S)_{\alpha}\rightarrow
K_{top}(\M_H(v)\times S).
\]
\end{defi}

\begin{rem}
\label{rem-comparizon-between-universal-class-and-twisted-universal-sheaf}
{\rm
When $v$ is a primitive class, 
then the co-cycle $\alpha$ maps to a co-boundary, once
we consider \v{C}ech cohomology with coefficients in the sheaf of complex 
valued smooth functions \cite{markman-integral-generators}. 
%Choose a $1$-co-chain $\psi:=(\psi_{ij})$ satisfying $\delta\psi=\alpha$. 
Furthermore, 
%It follows that 
the class of the
semi-universal sheaf $\F$ in $K_{top}(\M\times S)$
is the product of a class $e_v$ and the pull-back of the class
of a topological vector bundle $B$ on $\M$ 
of rank $\rho$ satisfying 
$W^*\cong B^{\otimes \rho}\otimes (\Wedge{\rho}B^*)$.  
The Chern character 
$ch(e_v)$ and the class $ch(\E)$ in
(\ref{eq-chern-character-of-rational-universal-class})
are related by the equation
\begin{equation}
\label{eq-comparizon-between-universal-class-and-twisted-universal-sheaf}
ch(\E)\exp[c_1(B)/\rho] \ \ \ = \ \ \
ch(e_v).
\end{equation}
The proof uses the equation $[ch(\F)/ch(\E)]^\rho=ch(W^*)$ proven in 
section 3 of \cite{markman-diagonal}.
}
\end{rem}

%******************************************************
% The Mukai pairing on $H^*(S,\Integers)$
%******************************************************
\subsection{The Mukai pairing on $H^*(S,\Integers)$}
\label{sec-weight-2-hodge-structure}
Let $S$ be a K3 surface and $\LB$ an ample line bundle on $S$. 
The Todd class of $S$ is $1+2\omega$, where $\omega$ is the fundamental class
in $H^4(S,\Integers)$. Its square root is $1+\omega$. 
Given a coherent sheaf $F$ on $S$ of rank $r$, we denote by 
\[
v(F)\ := \ ch(F)\sqrt{td_S} = (r,c_1(F),\chi(F)-r) 
\] 
its {\em Mukai vector} in  
\[
H^*(S,\Integers) \ = \ H^0(S,\Integers) \oplus H^2(S,\Integers) 
\oplus H^4(S,\Integers). 
\] 
We used above the standard isomorphisms $H^0(S,\Integers)\cong \Integers$ and
$H^4(S,\Integers)\cong \Integers$ to write the corresponding entries of 
$v(F)$ as integers.
Mukai endowed the cohomology group $H^*(S,\Integers)$ 
with a weight 2 polarized Hodge structure. 
The bilinear form is
\begin{eqnarray*}
\langle (r',c',s'), (r'',c'',s'') \rangle & = & c'c''-r's''-r''s',
\ \ \mbox{or equivalently}, \\
\langle \alpha, \beta\rangle & = & -\int_S \alpha^\vee\wedge\beta, 
\end{eqnarray*}
where $\alpha^\vee$ is the image of $\alpha$ under the duality operator 
\begin{equation}
\label{eq-serre-duality-hodge-isometry}
D \ : \ H^*(S,\Integers) \ \rightarrow \ H^*(S,\Integers)
\end{equation}
acting by $-1$ on the second cohomology (sending the Mukai 
vector $(r,c_1,s)$ to $(r,-c_1,s)$). 
The Hodge filtration is induced by that of $H^2(S,\Integers)$. 
In other words, $H^{2,0}(S)$ is defined to be also the $(2,0)$-subspace of 
the complexified Mukai lattice. The isomorphism
\begin{eqnarray}
\label{eq-isometry-between-K-top-and-Mukai-lattice}
K_{top}(S) & \rightarrow & H^*(S,\Integers),
\\
\nonumber
x & \mapsto & ch(x)\sqrt{td_S},
\end{eqnarray}
pulls back the Mukai pairing to the 
pairing in equation (\ref{eq-mukai-pairing-on-K-top}).
Riemann-Roch's Theorem yields the equation
\begin{equation}
\label{eq-Mukai-pairing-in-terms-of-sheaves}
\langle E, F\rangle := -\chi(E^\vee,F) := -\sum (-1)^i \dim \Ext^i(E,F),
\end{equation}
for the classes of coherent sheaves $E$ and $F$ on $S$. 

Denote by $\M(v):=\M_\LB(v)$ the moduli space of 
$\LB$-stable sheaves with Mukai vector $v$ of non-negative rank $r(v)\geq 0$. 
Mukai constructed the natural homomorphism
\begin{equation}
\label{eq-theta-v-from-v-perp}
\theta_v \ : \ v^\perp \ \ \rightarrow \ \ H^2(\M_\LB(v),\Integers)
\end{equation}
given by
\begin{equation}
\label{eq-mukai-homomorphism}
\theta_v(x) \ := \ \frac{1}{\rho}\left[p_{\M_*}((ch\E)\cdot\sqrt{td_S}\cdot 
\pi_S^*(x^\vee))
\right]_1,
\end{equation}
where $\E$ is a quasi-universal family of similitude $\rho$. 
%In this paper, $\E$ will always be a universal family and $\rho=1$.
Note, that the homomorphism $\theta_v$ extends to the whole Mukai lattice, but
the extension depends on the choice of $\E$. 

The isometry (\ref{eq-isometry-between-K-top-and-Mukai-lattice}) 
conjugates Mukai's homomorphism 
(\ref{eq-introduction-theta-v-from-v-perp})
to the homomorphism (\ref{eq-theta-v-from-v-perp}).
Indeed, the homomorphism $\theta_v$ does not change, if instead of 
the Chern character of a quasi-universal family $\E$ we let $ch(\E)$ 
be the class (\ref{eq-chern-character-of-rational-universal-class}),
considered with similitude $1$ (by the projection formula). The equality 
(\ref{eq-comparizon-between-universal-class-and-twisted-universal-sheaf})
allows us to further replace $ch(\E)$ by $ch(e_v)$ when $v$ is primitive.
The homomorphism $\theta_v$ in equation
(\ref{eq-theta-v-from-v-perp}) does not change, if we replace 
$ch(\E)\sqrt{td_S}$ in equation
(\ref{eq-mukai-homomorphism}) by $u_v$.

The moduli space $\M_{\LB}(v)$, of $\LB$-stable sheaves with Mukai 
vector $v$, depends on the polarization $\LB$.

\begin{defi}
\label{def-v-suitable}
{\rm 
An ample line bundle $\LB$ is said to be $v$-{\em suitable}, 
if every $\LB$-semi-stable sheaf with Mukai vector $v$ is $\LB$-stable. 
}
\end{defi}

{\em Throughout the paper, we will consider only moduli spaces $\M_{\LB}(v)$, 
where $\LB$ is $v$-suitable.}

We translate next the definition of the groupoid 
(\ref{eq-groupoid}), replacing $K_{top}(S)$ by $H^*(S,\Integers)$, 
using the isomorphism (\ref{eq-isometry-between-K-top-and-Mukai-lattice}).
The objects of the groupoid
\begin{equation}
\label{eq-cohomological-groupoid}
\G
\end{equation}
are triples $(S,v,\LB)$ consisting of a $K3$ surface $S$,
a primitive class $v=(r,c,s)$ in $H^*(S,\Integers)$ with
$c\in H^{1,1}(S,\Integers)$, and a $v$-suitable ample line bundle
$\LB$ on $S$. A morphism $g\in \Hom_\G((S_1,v_1,\LB_1)(S_2,v_2,\LB_2))$
is an isometry
$g:H^*(S_1,\Integers)\rightarrow H^*(S_2,\Integers)$ satisfying 
$g(v_1)=v_2$. 

\subsection{The class of the diagonal }
\label{sec-class-of-diagonal}

We will study automorphisms of the cohomology ring of a moduli $\M_{\LB}(v)$
via correspondences in $\M_{\LB}(v)\times \M_{\LB}(v)$. 
We have already dealt with
the case of the identity automorphism in another paper.
Let $\pi_{ij}$ be the projection from
$\M_{\LB}(v)\times S\times \M_{\LB}(v)$ onto the product of the
$i$-th and $j$-th factors.

\begin{thm} \cite{markman-diagonal}
\label{thm-graph-of-diagonal-in-terms-of-universal-sheaves}
Let $\E_v'$, $\E_v''$ be any two universal families of sheaves over the 
$m$-dimensional moduli space $\M_\LB(v)$. Assume, that 
$\LB$ is $v$-suitable.
The class of the diagonal, in the Chow ring of
$\M_\LB(v)\times \M_\LB(v)$, is identified by
\begin{equation}
\label{eq-class-of-diagonal}
c_m\left[- \ 
\pi_{13_!}\left(
\pi_{12}^*(\E'_v)^\vee\stackrel{L}{\otimes}\pi_{23}^*(\E''_v)
\right)
\right],
\end{equation}
where $\pi_{13_!}$ is the K-theoretic push-forward and both the dual
$(\E'_v)^\vee$ and the tensor product are taken in the derived category.
Furthermore, the following vanishing holds
\begin{equation}
\label{eq-c-m-1-vanishes}
c_{m-1}\left[- \ 
\pi_{13_!}\left(
\pi_{12}^*(\E'_v)^\vee\stackrel{L}{\otimes}\pi_{23}^*(\E''_v)
\right)
\right] \ \ \ = \ \ \ 0.
\end{equation}
\end{thm}

An immediate corollary of the Theorem is:
\begin{cor}
\label{cor-kunneth-factors-generate}
The K\"{u}nneth factors of the Chern classes of any universal sheaf $\E$
on $\M_\LB(v)\times S$ generate the cohomology ring 
$H^*(\M_\LB(v),\RationalNumbers)$. 
\end{cor}

\begin{rem}
\label{rem-after-prop-identification-of-diagonal}
{\rm
The geometric identification of the class (\ref{eq-class-of-diagonal}) 
implies its independence of the choice of the universal families. 
This independence follows also formally from the vanishing 
(\ref{eq-c-m-1-vanishes}) and Lemma
\ref{lemma-18-in-integral-generators}.
A universal family may not exist, in general, over the moduli space
$\M_\LB(v)$. Nevertheless, Theorem 
\ref{thm-graph-of-diagonal-in-terms-of-universal-sheaves} 
holds in general, provided we replace the classes
$ch(\E_1)$ and $ch(\E_2)$, in the topological translation
(\ref{eq-gamma-g-E1-E2}) of (\ref{eq-class-of-diagonal}), 
by the universal class 
in equation (\ref{eq-chern-character-of-rational-universal-class})
%in
%$H^*(\M_\LB(v)\times S,\RationalNumbers)$, naturally defined 
%in 
(\cite{markman-diagonal} Section 3).
}
\end{rem}

\begin{new-lemma}
\label{lemma-18-in-integral-generators}
(\cite{markman-integral-generators}, Lemma 18)
Let $X$ be a topological space, $x$ a class of rank $r\geq 0$ in $K_{top}(X)$,
and $L$ a complex line-bundle on $X$. Then 
$c_{r+1}(x\otimes L)=c_{r+1}(x)$ and 
$c_{r+2}(x\otimes L)=c_{r+2}(x)-c_{r+1}(x)c_1(L)$.
\end{new-lemma}

\subsection{The classes $\gamma(g,v)$}
\label{sec-conjecture}

Given a projective variety $M$, we denote by 
\begin{eqnarray*}
\ell \ : \ \oplus_{i}H^{2i}(M,\RationalNumbers) & \longrightarrow &
\oplus_{i}H^{2i}(M,\RationalNumbers)
\\
(r+a_1+a_2+\cdots ) & \mapsto & 1+a_1 + (\frac{1}{2}a_1^2-a_2) + \cdots 
\end{eqnarray*}
the universal polynomial map, which takes the exponential Chern character
of a complex of sheaves to its total Chern class. 

Given two K3 surfaces $S_1$ and $S_2$ and two $m$-dimensional Mukai vectors
$v_i\in H^*(S_i,\Integers)$, we
denote by $\pi_i$ the projection from 
$\M(v_1)\times S_2 \times \M(v_2)$ on the $i$-th factor.
Given classes $\alpha_i \in H^*(\M(v_i)\times S_i,\RationalNumbers)$ 
and a homomorphism 
$g: H^*(S_1,\Integers)\rightarrow H^*(S_2,\Integers)$, 
set
\begin{eqnarray}
\nonumber
x & := & \left[
\ell\left(
\pi_{13_*}\left\{
\pi_{12}^*[(id\otimes g)(\alpha_1)]^\vee
\cdot\pi_{23}^*(\alpha_2)
\right\}
\right)
\right]^{-1},
\\
\label{eq-gamma-delta}
\gamma(g,\alpha_1,\alpha_2) & := & c_m(x),
\\
\label{eq-gamma-prime}
\gamma'(g,\alpha_1,\alpha_2) & := & c_{m-1}(x).
\end{eqnarray}
Then $\gamma(g,\alpha_1,\alpha_2)$ is a class in 
$H^{2m}(\M(v_1)\times \M(v_2),\RationalNumbers)$.

%We denote by $\Delta$ the class of the diagonal in $S_1\times S_1$.
%Given an isometry 
%$g:H^*(S_1,\Integers)\rightarrow H^*(S_2,\Integers)$, 
%we get the class $(1\times g)(\Delta)$ in $H^*(S_1\times S_2,\Integers)$
%inducing $g$.  
%Set
%\begin{equation}
%\label{eq-gamma-g}
%\gamma(g,\alpha,\beta) \ \ := \ \ 
%\gamma((1\times g)(\Delta),\alpha,\beta).
%\end{equation}
If $\E_i$ is a complex of sheaves on $\M(v_i)\times S_i$, we set
\begin{equation}
\label{eq-gamma-g-E1-E2}
\gamma(g,\E_1,\E_2) \ \ := \ \ 
\gamma(g,ch(\E_1)\cdot\sqrt{td_{S_1}},ch(\E_2)\cdot\sqrt{td_{S_2}}).
\end{equation}
When $S_1=S_2$ and $g=id$, Grothendieck-Riemann-Roch yields the equality
$\gamma(id,\E_1,\E_2)  =  c_m\left\{-\pi_{13_!}\left(\pi_{12}^*(\E_1)^\vee
\stackrel{L}{\otimes}
\pi_{23}^*(\E_2)\right)\right\}.
$
%where on the right hand side 
%$\pi_{ij}$ is the projection from $\M(v_1)\times S\times \M(v_2)$.
See section \ref{seq-Fourier-Mukai-transformations} 
for a Chow-theoretic identification of 
$\gamma(g,\E_1,\E_2)$, when $g$ is the Hodge isometry of an
auto-equivalence of the derived category of the surface. 
A conceptual interpretation of formula 
(\ref{eq-gamma-delta}) is provided in section 
\ref{sec-conceptual-interpretation-of-the-formula-gamma-g}. 
%See Lemma \ref{lemma-two-classes-are-equal} part 
%\ref{lemma-item-class-of-reflection-in-terms-of-universal-sheaves} 
%for a Chow-theoretic identification of 
%$\gamma_g(\E_1,\E_2)$ when $g$ is a reflection. When
%$g$ is the Serre's Duality Hodge isometry, the identification is provided in 
%(\ref{eq-chow-class-of-serre-duality-correspondence}). When $g=-id$, 
%the Chow-theoretic identification is given in Section \ref{sec-minus-id}.

When $(v,v)>0$, we set
\begin{equation}
\label{eq-class-gamma-g-v}
\gamma(g,v) \ \ := \ \ \gamma(g,u_v,u_{g(v)}),
\end{equation}
where $u_v$ is given in equation 
(\ref{eq-invariant-normalized-class-of-chern-character-of-universal-sheaf}).
Given universal sheaves $\E_v$ and $\E_{g(v)}$ over the 
moduli spaces, we have the equality 
$\gamma(g,v)=\gamma(g,\E_v,\E_{g(v)})$. 
The latter equality, as well as the more general equality of
the two definitions (\ref{eq-intro-class-gamma-g-v}) and 
(\ref{eq-class-gamma-g-v}) of $\gamma(g,v)$, 
is proven in Lemma \ref{lemma-two-definitions-of-gamma-g-v}
(using Theorem \ref{thm-trancendental-reflections}).
When $(v,v)=0$, set $\gamma(g,v):=\gamma(g,\E_v,\E_{g(v)})$,
when universal sheaves $\E_v$ and $\E_{g(v)}$ exist.
Use instead a universal class $e_v$ or $e_{g(v)}$ (Definition \ref{def-e-v}), 
in the absence of a universal sheaf.
The existence of $e_v$ is clear for two dimensional moduli spaces.
The class $\gamma(g,v)$ is independent of the choice of
universal sheaves or classes in this case as well
(Lemma \ref{lemma-conjecture-holds-for-isotropic-mukai-vectors}). 
%(see Remark \ref{rem-after-prop-identification-of-diagonal}).

%Theorem \ref{thm-graph-of-diagonal-in-terms-of-universal-sheaves} states 
%that $\gamma(id,v)$ is Poincare-dual 
%to the diagonal in $\M_\LB(v)\times \M_\LB(v)$. 

We identify $H^*(\M(v)\times \M(g(v)),\RationalNumbers)$ with 
$H^*(\M(v),\RationalNumbers)^*\otimes H^*(\M(g(v)),\RationalNumbers)$
via the K\"{u}nneth and Poincare-Duality theorems. 
The homomorphism corresponding to the class $\gamma(g,v)$ is denoted by
\[
\gamma_{g,v} \ : \ H^*(\M(v),\RationalNumbers) \ \ \ 
\longrightarrow \ \ \ H^*(\M(g(v)),\RationalNumbers).
\]
We set $\gamma_g:=\gamma_{g,v}$, when the context identifies the 
Mukai vector $v$. 

\bigskip
We can now state parts \ref{thm-item-gamma-g-v-is-an-isomorphism} 
and \ref{thm-item-invariance-of-the-universal-class} of Theorem
\ref{thm-action} without assuming the existence of the universal class
$e_v$:
%a trancendental version of Conjecture 
%\ref{conj-the-total-class-of-the-reflection}. 

\begin{thm}
\label{thm-trancendental-reflections}
Let $(S_i,v_i,\LB_i)$, $i=1,2$, be objects of the groupoid
$\G$ given in (\ref{eq-cohomological-groupoid}) and
$g\in \Hom_\G((S_1,v_1,\LB_1),(S_2,v_2,\LB_2))$ a morphism 
(which need not preserve the Hodge structures). 
Assume that $(v,v)>0$. 
\begin{enumerate}
\item
\label{conj-item-gamma-g-induces-isomorphism-of-cohomology-rings}
The homomorphism 
\[
\gamma_g : H^*(\M_{\LB_1}(v_1),\Integers)_{\rm free} \ \rightarrow \ 
H^*(\M_{\LB_2}(v_2),\Integers)_{\rm free}
\]
is an isomorphism of cohomology rings.
%where $H^*(\M_{\LB_1}(v_1),\Integers)_{\rm free}$ is the torsion free summand
%of the integral cohomology.
%If $g$ is a Hodge isometry, then
%$\gamma_g$ is an isomorphism of Hodge structures. 
\item
\label{conj-item-chern-character-of-univ-family-is-taken-to-same}
The isomorphism 
\[
(\gamma_g\otimes g) \ : \ H^*(\M(v_1)\times S_1,\RationalNumbers) 
\ \LongIsomRightArrow \ 
H^*(\M(v_2)\times S_2,\RationalNumbers)
\]
takes the class 
$u_{v_1}$ 
to the class 
$u_{v_2}$.  
\end{enumerate}
\end{thm}

%See Theorem \ref{thm-Gamma-v-acts-motivicly},  
%Theorem \ref{thm-conj-holds-for-hilbert-schemes-and-isometries-of-H-2}, 
%and Theorem \ref{thm-class-of-correspondence-in-stratified-elementary-trans}
%for verifications of special cases of the conjecture. 
%See Remark
%\ref{rem-auto-equivalences-of-the-derived-category-of-moduli-spaces},
%for a speculative extension of the Theorem, at the level of
%derived categories. 

\begin{rem}
\label{rem-other-base-calabi-yau-varieties}
{\rm
Part \ref{conj-item-gamma-g-induces-isomorphism-of-cohomology-rings} 
of the Theorem makes sense, even when the Mukai vectors $v_i$ are 
isotropic and $\M_{\LB_i}(v_i)$ are K3 surfaces. 
In that case, 
part \ref{conj-item-gamma-g-induces-isomorphism-of-cohomology-rings} 
follows easily from the work of Mukai
(see Lemma \ref{lemma-conjecture-holds-for-isotropic-mukai-vectors}).

Part \ref{conj-item-gamma-g-induces-isomorphism-of-cohomology-rings} 
of the Theorem makes sense, even when we replace the K3 surface 
by an abelian surface or an elliptic curve; 
the two other Calabi-Yau varieties, for which Theorem 
\ref{thm-graph-of-diagonal-in-terms-of-universal-sheaves} is known
\cite{beauville-diagonal,markman-diagonal}. 
The abelian surface case seems very interesting, 
and merits further consideration. 
%We proceed here to comment, that 
The elliptic curve case boils down to an old result of Atiyah.
This is explained in section
\ref{sec-elliptic-curve}.
}
\end{rem}

%The significance of part
%\ref{conj-item-chern-character-of-univ-family-is-taken-to-same}
%of the Theorem, for our purpose, is that it plays the role of a
%``\mbox{\em canonical} \mbox{\em basis}'' for 
%$H^*(S,\RationalNumbers)\otimes H^*(\M(v),\RationalNumbers)$,
%in view of Corollary \ref{cor-kunneth-factors-generate}. 
%It plays a central role in the proof of Theorem 
%\ref{thm-Gamma-v-acts-motivicly}. Further 
%Applications include corollaries \ref{cor-group-homomorphism},
%\ref{cor-generating-subspaces-are-sub-representations}. 
Part \ref{conj-item-chern-character-of-univ-family-is-taken-to-same}
of the Theorem implies the commutativity of the following diagram
\begin{equation}
\label{eq-g-gamma-g-conjugate-theta-to-theta}
\begin{array}{ccc}
v_1^\perp & \LongRightArrowOf{\theta_{v_1}} & 
H^2(\M_{\LB_1}(v_1),\Integers)
\\
g \ \downarrow \hspace{2ex}
& & \hspace{3ex} \downarrow \ \gamma_g
\\
v_2^\perp & \LongRightArrowOf{\theta_{v_2}} &
H^2(\M_{\LB_2}(v_2),\Integers).
\end{array}
\end{equation}
The horizontal arrows
are isomorphisms by Theorem \ref{thm-irreducibility}.

%We will use the following terminology to state special cases of
%the Theorem. 
%Let $\Gamma'$ be a subset of the set of isometries, 
%from the Mukai lattice of 
%$S_1$ to that of $S_2$, and
%$\Sigma$ a subset of the Mukai lattice of $S_1$.
%****************
% Definition
%****************
%\begin{defi}
%\label{def-conjecture-holds}
%We say that 
%{\em Theorem \ref{thm-trancendental-reflections} 
%holds for $(S_1,S_2,\Sigma,\Gamma')$,
%}
%if its statements hold for 
%every $v\in \Sigma$, every $g\in \Gamma'$, and for
%all suitable  polarizations $\LB_i$, $i=1,2$. 
%When $S_1=S_2=S$, we would say that 
%{\em 
%Theorem \ref{thm-trancendental-reflections} 
%holds for $(S,\Sigma,\Gamma')$.
%} 
%Given an ample line bundle $H$ on $S$, we say that
%{\em
%Theorem \ref{thm-trancendental-reflections} 
%holds for $(S,H,\Sigma,\Gamma')$,
%}
%if its statements hold for $S_i=S$, $i=1,2$, $\LB_i=H$, $i=1,2$,
%every $v\in \Sigma$, and every $g\in \Gamma'$. 
%\end{defi}

%*********************************************************************
% Compositions
%*********************************************************************
\subsection{Compositions}
\label{sec-compositions}
We show in Lemma \ref{lemma-conjecture-holds-for-compositions} that 
the validity of Theorem \ref{thm-trancendental-reflections} 
is closed under compositions. In particular, the homomorphisms  
$\gamma_{g,v}$ compose as expected, yielding a representation of the
groupoid $\G$. 
We then discuss the cohomology ring of a fixed 
moduli space as a representation of the automorphism group of 
the corresponding object of $\G$. 

We begin with a characterization of the class $\gamma(g,v)$.
Let $\Delta_i$ be the diagonal in $\M_{\LB_i}(v_i)\times \M_{\LB_i}(v_i)$,
$i=1,2$. Denote by $[\Delta_i]$ the cohomology class Poincare dual to
$\Delta_i$. 

%We generalize in Lemma \ref{lemma-recovering-f} the discussion of 
%Lemma \ref{lemma-the-class-of-f-is-gamma-E1-E2}
%by dropping the assumption that
%$f$ and $g$ arise from a deformation. This result is used in the proof of
%lemma \ref{lemma-conjecture-holds-for-compositions}
%to show, that the classes
%$\gamma_g(\E_v,\E_{g(v)})$ compose as expected.

\begin{new-lemma}
\label{lemma-recovering-f}
Suppose that $f: H^*(\M_{\LB_1}(v_1),\RationalNumbers) \rightarrow 
H^*(\M_{\LB_2}(v_2),\RationalNumbers)$ 
is a {\em ring} isomorphism, 
$g : H^*(S_1,\RationalNumbers) \rightarrow H^*(S_2,\RationalNumbers)$ 
a linear homomorphism, 
and the pair $(f,g)$ satisfies 
\begin{equation}
\label{eq-f-g-takes-ch-to-ch}
(f\otimes g)\left(u_{v_1}\right) \ = \ 
u_{v_2}. 
\end{equation}
Then 
\begin{equation}
\label{eq-f-is-recovered}
(1\otimes f)([\Delta_1]) \ = \ 
\gamma(g,u_{v_1}\exp(\ell_1),u_{v_2}\exp(\ell_2)),
\end{equation}
for any two classes $\ell_i\in H^2(\M_{\LB_i}(v_i),\RationalNumbers)$,
$i=1,2$. 
In particular, given $g$, a ring isomorphism $f$ satisfying
(\ref{eq-f-g-takes-ch-to-ch}) is {\em unique} (if it exists).
Furthermore, the class $\gamma'(g,u_{v_1},u_{v_2})$,
given in (\ref{eq-gamma-prime}), vanishes.
\end{new-lemma}

\noindent
{\bf Proof:}
We have the equalities
\begin{eqnarray*}
[\Delta_2] & = & \gamma(id,u_{v_2}\exp(f(\ell_1)),u_{v_2}\exp(\ell_2)) 
\ \stackrel{(\ref{eq-f-g-takes-ch-to-ch})}{=} \ 
\gamma(id,(f\otimes g)[u_{v_1}\exp(\ell_1)],u_{v_2}\exp(\ell_2))  
\\
& = &
(f\otimes 1)[\gamma(g,u_{v_1}\exp(\ell_1),u_{v_2}\exp(\ell_2))].
\end{eqnarray*}
Theorem \ref{thm-graph-of-diagonal-in-terms-of-universal-sheaves}
yields the first equality
(see also Remark \ref{rem-after-prop-identification-of-diagonal}). 
The last equality is due to the fact, that 
evaluation of $(f\otimes 1)$ commutes with $\gamma$ because $f$
is a {\em ring} isomorphism. 
Evaluating $(f\otimes 1)^{-1}$ on both sides, we conclude the equality
(\ref{eq-f-is-recovered}). 

The vanishing of $\gamma'(g,u_{v_1},u_{v_2})$ follows similarly from 
the corresponding vanishing (\ref{eq-c-m-1-vanishes}) in 
Theorem \ref{thm-graph-of-diagonal-in-terms-of-universal-sheaves} 
and Lemma \ref{lemma-18-in-integral-generators}.
\EndProof

\begin{new-lemma}
\label{lemma-conjecture-holds-for-compositions}
Let $x_i:=(S_i,v_i,\LB_i)$, $i=1,2,3$, be objects of the groupoid
$\G$, given in (\ref{eq-cohomological-groupoid}),
$g\in \Hom_\G(x_1,x_2)$, and $f\in \Hom_\G(x_2,x_3)$.
If $g$ and $f$ satisfy the statements of Theorem 
\ref{thm-trancendental-reflections}, 
then so do $g^{-1}$ and $f\circ g$. 
Moreover, the following equality holds
\begin{equation}
\label{eq-gamma-maps-compositions-to-compositions}
\gamma_{fg,v_1} \ \ = \ \ 
\gamma_{f,v_2}\circ \gamma_{g,v_1}.
\end{equation}
\end{new-lemma}

\medskip
\noindent
{\bf Proof:}
The assumptions imply, that $\gamma_{g,v_1}$ is a
ring isomorphism and $\gamma_{g,v_1}\otimes g$  maps the class $u_{v_1}$ to
$u_{v_2}$. Similarly, the class 
$\gamma_{f,v_2}$ has the analogous properties. 
Hence, the composition 
$\phi:=\gamma_{f,v_2}\circ \gamma_{g,v_1}$ 
is a ring isomorphism, and $\phi\otimes (fg)$ maps $u_{v_1}$ to 
$u_{v_3}$. Lemma \ref{lemma-recovering-f} 
implies the equality 
$\phi=\gamma_{fg,v_1}$. The latter is precisely equality
(\ref{eq-gamma-maps-compositions-to-compositions}). 
The proof for $g^{-1}$ is similar.
\EndProof

Note: The classes $\gamma(g,v)$ and $\gamma(g^{-1},g(v))$
are equal, under the natural identification of 
$\M(v)\times \M(g(v))$ with $\M(g(v))\times \M(v)$. 
A direct proof, without the hypothesis that $g$ satisfies the 
assumptions of 
Theorem \ref{thm-trancendental-reflections}, 
is not hard (see Lemma 4.4 part 3 in the eprint version math.AG/0305042 v2 
of this paper). A second proof of the equality 
$(\gamma_{g,v})^{-1}=\gamma_{g^{-1},g(v)}$ then follows, once one shows that 
$\gamma_{g,v}$ is a ring isomorphism.

Let $\Gamma$ be the isometry group of the Mukai lattice of a K3 surface $S$. 
As a corollary of Theorem \ref{thm-trancendental-reflections} and 
lemma \ref{lemma-conjecture-holds-for-compositions}, 
we get that the stabilizer $\Gamma_v$, of a Mukai vector $v$, 
acts on the cohomology of the corresponding moduli space:

\begin{cor}
\label{cor-group-homomorphism}
Let $(S,v,\LB)$ be an object of the groupoid $\G$.
Assume, that $(v,v)\geq 2$. 
Then the natural map
\begin{eqnarray}
\label{eq-homomorphism-gamma-from-stabilizer-to-ring-auto}
\gamma \ : \ \Gamma_v & \longrightarrow & 
\Aut\left[H^*(\M_\LB(v),\Integers)_{\rm free}\right] 
\\
\nonumber
g & \mapsto & \gamma_{g,v}
\end{eqnarray}
is an injective group homomorphism. 
\end{cor}

\noindent
{\bf Proof:}
The homomorphism $\gamma$ is injective since diagram 
(\ref{eq-g-gamma-g-conjugate-theta-to-theta}) is commutative.
\EndProof

\begin{new-lemma}
\label{lemma-two-definitions-of-gamma-g-v}
Let $g\in \Hom_{\G}((S_1,v_1,\LB_1),(S_2,v_2,\LB_2))$.
The two definitions (\ref{eq-intro-class-gamma-g-v}) and 
(\ref{eq-class-gamma-g-v}) of the class $\gamma(g,v_1)$
agree.
\end{new-lemma}

\noindent
{\bf Proof:}
The two definitions are identical, when $(v_1,v_1)=0$. Assume $(v_1,v_1)>0$.
Theorem \ref{thm-trancendental-reflections} 
proves that the assumptions of Lemma \ref{lemma-recovering-f}
are satisfied for the operator $f:=\gamma_{g,v_1}$ corresponding to 
definition (\ref{eq-class-gamma-g-v}) of the
class $\gamma(g,v_1)$. The class $ch(e_{v_i})\sqrt{td_{S_i}}$ is equal to 
$u_{v_i}\cdot \exp(\ell_i)$, for some class $\ell_i$ in 
$H^2(\M_{\LB_i}(v_i),\RationalNumbers)$, $i=1,2$, by equations 
(\ref{eq-invariant-normalized-class-of-chern-character-of-universal-sheaf}) 
and (\ref{eq-comparizon-between-universal-class-and-twisted-universal-sheaf}). 
Lemma \ref{lemma-recovering-f} implies the equality
$\gamma(g,u_{v_1},u_{v_2})=
\gamma(g,ch(e_{v_1})\sqrt{td_{S_1}},ch(e_{v_2})\sqrt{td_{S_2}})$.
The left hand side appears in the topological definition
(\ref{eq-class-gamma-g-v}) of the class $\gamma(g,v_1)$, while the right hand 
side is the topological translation of the K-theoretic definition 
(\ref{eq-intro-class-gamma-g-v}) via Grothendieck-Riemann-Roch.
\EndProof

%\smallskip
%\noindent
%{\bf Proof:} 
%The property in part 
%\ref{conj-item-chern-character-of-univ-family-is-taken-to-same}
%of Conjecture \ref{thm-trancendental-reflections} 
%certainly behaves well under compositions. The class 
%$\gamma_g(u_v,u_v)$ 
%is characterised by this property (see 
%Lemma \ref{lemma-recovering-f}). 
%\EndProof

\medskip
The set of universal classes, generating the cohomology ring
of $\M_\LB(v)$, is invariant with respect to the representation
(\ref{eq-homomorphism-gamma-from-stabilizer-to-ring-auto}). 
This is another corollary of the Theorem 
\ref{thm-trancendental-reflections}.
Let $B$ be the image of the Mukai lattice in 
$H^*(\M_\LB(v),\RationalNumbers)$, 
via the $\Gamma_v$-equivariant homomorphism
%via the $\Gamma_v$-invariant class $u_v$.
\begin{eqnarray}
\label{eq-homomorphism-induced-by-u-v}
H^*(S,\Integers) & \longrightarrow & H^*(\M_\LB(v),\RationalNumbers)
\\
\nonumber
x & \mapsto & -p_{\M_*}(u_v\cdot \pi_S^*(x^\vee))
\end{eqnarray}
(compare with equation (\ref{eq-mukai-homomorphism}) defining $\theta_v$).
The above homomorphism is the image of $u_v$ under the 
identification of the Mukai lattice $H^*(S,\Integers)$ with its dual.
We use the Mukai pairing for the identification, 
rather than Poincare Duality, because the latter is not
$\Gamma$-invariant. This identification sends
$x\in H^*(S,\Integers)$ to 
$-\int_S x^\vee\cup (\bullet)$ in $H^*(S,\Integers)^*$.
Corollary \ref{cor-kunneth-factors-generate} implies, that 
the cohomology ring $H^*(\M_\LB(v),\RationalNumbers)$ is generated by
the projections $B_i$  of $B$ into 
$H^{2i}(\M_\LB(v),\RationalNumbers)$. Part
\ref{conj-item-chern-character-of-univ-family-is-taken-to-same} of Theorem
\ref{thm-trancendental-reflections} implies:

\begin{cor}
\label{cor-generating-subspaces-are-sub-representations}
Each of the vector subspaces $B_i\subset H^{2i}(\M_\LB(v),\RationalNumbers)$,
generating the cohomology ring of $\M_\LB(v)$,
is $\Gamma_v$-invariant. 
\end{cor}

%\ref{cor-item-generators-for-the-cohomology-ring})
%The class $\gamma(\E_v,\E_v)$, in the cohomology of
%$\M_\LB(v)\times\M_\LB(v)$, represents the diagonal. 
%The K\"{u}nneth factors of the diagonal generate $H^*(\M_\LB(v),\Integers)$. 
%The K\"{u}nneth factors of $\gamma(\E_v,\E_v)$ are in the subring 
%generated by the K\"{u}nneth factors of $ch(\E_v)$. 
%%The  class $ch(\eta_v\otimes\E_v)\cdot\pi_2^*\sqrt{td_S}$ 
%%is $\Gamma_v$-invariant
%%(Part \ref{conj-item-chern-character-of-univ-family-is-taken-to-same}
%%of Conjecture \ref{thm-trancendental-reflections}). 
%%Hence, the image $B_i$ of the Mukai lattice in 
%%$H^{2i}(\M_\LB(v),\RationalNumbers)$ 
%%is a sub-$\Gamma_v$-module. 

Note, that $B\otimes\RationalNumbers$ 
consists of, at most, two irreducible $\Gamma_v$ sub-representations:
the trivial character and 
$v^\perp\otimes\RationalNumbers$, each with multiplicity  $\leq 1$. 

The following lemma will be cited in subsequent sections.
Denote by $\R^{[n]}$ the weighted polynomial ring, 
formally generated by the vector spaces $B_i$, $1\leq i\leq n-1$, 
in Corollary \ref{cor-generating-subspaces-are-sub-representations},
with vectors in $B_i$ having degree $2i$. 
Let $h:\R^{[n]}\rightarrow H^*(S^{[n]},\RationalNumbers)$ be the
natural ring homomorphism and $I^d$ its kernel in degree $d$. 

\begin{new-lemma}
\label{lemma-10-in-markman-diagonal}
(Lemma 10 in \cite{markman-diagonal})
$B_1$ is $23$-dimensional and $B_i$ is $24$-dimensional,
for $2\leq i \leq n/2$. If $n$ is odd, then $B_{(n+1)/2}$ is
either $23$ or $24$ dimensional. 
The homomorphism $h$ is surjective. It is injective
in degree $\leq n$. If $n$ is odd, then the dimension of 
%the degree $(n+1)$ summand 
$I^{n+1}$ is $\dim(B_{(n+1)/2})-23$. 
\end{new-lemma}

%*****************************************************************
% Monodromy of Hilbert schemes:
%*****************************************************************
\section{The monodromy group}
\label{sec-monodromy-of-hilbert-schems}

In section \ref{sec-orientation-character} we recall the orientation 
character of the isometry group of the Mukai lattice. 
This character arises in formula  
(\ref{introduction-eq-monodromy-is-gamma-times-cov}) for the monodromy
representation. 
In sections \ref{sec-monodromy-operators-acting-trivially-on-H-2} and 
\ref{sec-splitting} we show, that the monodromy 
subgroup $\W$, constructed in Corollary 
\ref{cor-W-consists-of-monodromies}, 
maps onto a finite index subgroup of the monodromy group 
$\Mon(S^{[n]}) \subset \Aut[H^*(S^{[n]},\Integers)]$. Furthermore,
the subgroup $K$ of the monodromy group, 
which acts trivially on the second cohomology of $S^{[n]}$, is central, 
finite, and of exponent $2$. 
In section \ref{sec-chern-classes-of-tangent-bundle} 
we improve the upper bound for the order of $K$,
using the monodromy-invariance of the Chern classes $c_{2i}(TS^{[n]})$.
In section \ref{sec-Mon-2-is-large} 
we provide several characterizations of the subgroup $\W$, of the monodromy
group $\Mon^2$, constructed in Corollary 
\ref{cor-W-consists-of-monodromies}.
%we discuss the index of the monodromy
%group, in the isometry group of the weight $2$ cohomology. 
In section \ref{sec-comparison-with-verbitsky}
we compare the monodromy representation of 
Theorem  \ref{introduction-thm-Gamma-v-acts-motivicly} with 
a related representation constructed by Verbitsky
(Lemma \ref{lemma-comparing-two-representations}).  
In section 
\ref{sec-monodromy-invariance-of-the-normalized-chern-character}
we study the equivariance properties, of the Chern 
character of the universal sheaf, 
with respect to the monodromy representation.

%****************************************************************
% The covariance character
%****************************************************************
\subsection{The orientation and covariance characters} 
\label{sec-orientation-character}

\begin{defi}
\label{def-isometry-groups}
Denote by $\Gamma$ the group of isometries of the Mukai lattice 
$H^*(S,\Integers)$ of a fixed K3 surface $S$. Let 
$G$ be the subgroup of (integral) Hodge-isometries of the Mukai lattice. 
Denote by $\Gamma_0$ the group of isometries of the second cohomology 
$H^2(S,\Integers)$ and let $G_0$ be the subgroup of 
integral Hodge-isometries of the weight 2 Hodge structure of $S$. 
\end{defi}

$\Gamma_0$ embeds in $O(3,-19)$ and $\Gamma$ embeds in $O(4,-20)$. 
The Hodge decompositions 
\begin{eqnarray*}
H^2(S,\ComplexNumbers) & = &
[H^{2,0}(S)\oplus H^{0,2}(S)] \ \oplus \ H^{1,1}(S)
\\
H^*(S,\ComplexNumbers) & = &
[H^{2,0}(S)\oplus H^{0,2}(S)] \ \oplus \ 
[H^{1,1}(S)\oplus H^0(S)\oplus H^4(S)]
\end{eqnarray*}
imply that $G_0$ embeds naturally in $Aut(H^{2,0}) \times O(1,-19)$ and 
$G$ embeds naturally in $Aut(H^{2,0}) \times O(2,-20)$.
The group $O(n,-m)$, $n,m> 1$,  has, in addition to the determinant character, 
the {\em orientation character} 
\begin{equation}
\label{eq-top-homology-character}
O(n,-m) \ \ \rightarrow \ \ 
\Aut(H^{n-1}(S^{n-1},\Integers)) \cong \Integers/2\Integers.
%= \{1,-1\} \ \subset \ComplexNumbers^\times.
\end{equation}
The cone 
${\cal C} \ = \ \{v \ : \ (v,v) > 0 \}$ is an 
$(\RealNumbers^{n}\setminus \{0\})$-bundle over $\RealNumbers^{m}$ and
is hence homotopic to the sphere $S^{n-1}$. 
The character (\ref{eq-top-homology-character}) 
is the action of $O(n,-m)$ on the top cohomology 
group of $\cal{C}$. 

In the case of $G_0$, 
the cone ${\cal C}\subset H^{1,1}(S,\RealNumbers)$ has two
components (the case $n=1$) and a Hodge isometry is said to be a
{\em signed isometry}  if it maps each component to itself. 
All automorphisms of $S$, as well as reflections along $-2$ curves on $S$,
are signed isometries. 

In the case of $G$, the cone 
${\cal C}\subset H^{1,1}(S,\RealNumbers)\oplus 
H^{0}(S,\RealNumbers) \oplus H^{4}(S,\RealNumbers)$ 
is homotopic to $S^1$. In this case, we will call 
(\ref{eq-top-homology-character}) the 
{\em covariance} character and denote it by 
\begin{equation}
\label{eq-covariance-character}
cov \ : \ G \ \ \rightarrow \ \ 
\Aut(H^{n-1}({\cal C},\Integers)) \cong \Integers/2\Integers.
%= \{1,-1\} \ \subset \ComplexNumbers^\times.
\end{equation}
We have the natural embedding of $G_0$ into $G$
%\[
%e \ : \ G_0 \ \hookrightarrow G
%\] 
via extension by the identity on $H^0\oplus H^4$. 
This embedding pulls back the character 
$cov$ on $G$ to the same character on $G_0$. 
Similarly, the natural embeddings $G_0\subset \Gamma_0\subset \Gamma$ and
$G_0\subset G \subset \Gamma$ are all compatible with respect to
the orientation character (\ref{eq-top-homology-character}). 

\begin{defi}
\label{def-covariant-subgroups}
Denote by $\Gamma^{cov}$, $\Gamma_0^{cov}$, $G^{cov}$, and $G_0^{cov}$,
the respective kernels of the orientation character. 
\end{defi}

The covariance character sends the action of the Picard 
group, 
%the isometries induced by the Fourier-Mukai functors, 
$\Aut(S)$, 
and reflection in (-2)-vectors, to $0$. 
The isometry $-id$, acting via multiplication by $-1$ on the
whole Mukai lattice, is also covariant. Conjecturally, 
any auto-equivalence, of the derived category of a K3 surface, is
sent by the covariance character to $0$
(see \cite{szendroi,mirror-symmetry-k3}).
Many cases of this conjecture were proven in \cite{huybrechts-stellari}
Proposition 5.5.
On the other hand, the Duality Hodge isometry 
(\ref{eq-serre-duality-hodge-isometry}) 
is contravariant, $cov(D)=1$. 
Reflections in (+2)-vectors are sent to $1$ as well.  
This can be seen as follows. The determinant character is the product 
of the orientation character of the positive cone,
with the orientation character of the negative cone. 
A reflection in a (+2)-vector has determinant $-1$, and it
preserves the orientation of the negative cone. 

\begin{rem}
\label{rem-natural-orientation}
{\rm
The positive cone in $H^{1,1}(S,\RealNumbers)$ has a distinguished 
component; the one containing the K\"{a}hler cone. 
Similarly, the positive cone ${\cal C}$ in $H^{2}(S,\RealNumbers)$ 
has a distinguished {\em orientation}; i.e., a generating class in 
$H^2({\cal C},\Integers)$. Such an orientation is equivalent to 
the determination of an orientation of one
$3$-dimensional subspace of $H^{2}(S,\RealNumbers)$, 
to which the pairing restricts 
as a positive definite one. Let $\sigma$ be a holomorphic 
symplectic structure on $S$ and $\kappa$ a K\"{a}hler form.
Then $\{{\rm Re}(\sigma),{\rm Im}(\sigma),\kappa\}$ is a basis for
such a subspace. The orientation of ${\cal C}$ is independent of the
choices of $\kappa$ and $\sigma$. 
We can further assign a distinguished orientation to the Mukai lattice
of $S$, by adding to the above basis the $+2$ Mukai vector $(1,0,-1)$
as a fourth vector. These four vectors span a 
$4$-dimensional positive definite sublattice of the Mukai lattice. 
This orientation defines an extension of the 
orientation character 
to a functor (\ref{eq-functor-cov})
from the groupoid $\G$ to $\Integers/2\Integers$.
Note that the Mukai vector $(1,0,-1)$ is effective 
(Definition \ref{def-effective-class}).
}
\end{rem}

\subsection{Monodromy operators acting trivially on $H^2$}
\label{sec-monodromy-operators-acting-trivially-on-H-2}

Let $X$ be an irreducible holomorphic symplectic manifold. 
The theory of Lefschetz-modules, developed by
Verbitsky and Looijenga-Lunts, introduces an action of the group
$Spin H^2(X,\RealNumbers)$ on the cohomology of $X$. 
The action is monodromy-equivariant and preserves the ring structure. 
We will use this theory to study the subgroup $K\subset Mon(X)$
of monodromy operators, which act trivially on $H^2(X,\Integers)$
(Corollary \ref{cor-monodromy-invariant-decomposition}). 

Let $A_{k}\subset H^*(X,\Integers)_{\rm free}$, $k\geq 0$,   
be the graded subalgebra generated by 
$\oplus_{i=0}^k H^{i}(X,\Integers)_{\rm free}$. 
Set $(A_{k})^j:=A_{k}\cap H^j(X,\Integers)_{\rm free}$.
Define $A_{k}(\RationalNumbers)$ and $[A_{k}(\RationalNumbers)]^j$
using rational coefficients. 
The theory of Lefschetz-modules will enable us to construct a natural 
subspace $C_k(\RationalNumbers)$ of $H^k(X,\RationalNumbers)$,
leading to a {\em monodromy-invariant} decomposition
\begin{equation}
\label{eq-monodromy-invariant-decomposition-of-H-k}
H^k(X,\RationalNumbers) \ \ \ = \ \ \ 
[A_{k-2}(\RationalNumbers)]^k \oplus C_k(\RationalNumbers). 
\end{equation}

Set $d:=\dim_\ComplexNumbers(X)$. Let $h\in \End[H^*(X,\RealNumbers)]$ 
be the linear endomorphism acting via multiplication by $(i-d)$ on 
$H^i(X,\RealNumbers)$. 
Given a class $a\in H^2(X,\RealNumbers)$, denote by 
$e_a\in \End[H^*(X,\RealNumbers)]$ 
the operator obtained by cup product with $a$.
The class $a$ is called of {\em Lefschetz type}, if there exists 
$f_a\in  \End[H^*(X,\RealNumbers)]$ satisfying the $\LieAlg{sl}_2$ 
commutation relations
\[
[e_a,f_a]=h, \ \ \ 
[h,e_a]=2e_a, \ \ \ 
[h,f_a]=-2f_a.
\]
Such $f_a$ is unique, if it exists. 
The triple $(e_a,h,f_a)$ is called a 
{\em Lefschetz triple}.
The set of classes $a\in H^2(X,\RealNumbers)$ of 
Lefschetz type is a Zariski dense open subset. 

Let $\LieAlg{g}(X)$ be the graded Lie subalgebra of 
$\End[H^*(X,\RealNumbers)]$
generated by 
$\{e_a,f_a\}$  for all
$a$ in $H^2(X,\RealNumbers)$ of Lefschetz type. 
Denote by $\LieAlg{g}_k(X)$ the subspace of grade $k$.
The primitive subspace $Prim^k(X)\subset H^k(X,\RealNumbers)$ 
is the set of classes killed by $\LieAlg{g}_{<0}(X)$ and
$Prim(X):=\oplus_{k}Prim^k(X)$. 
Set $b_2:=\dim H^2(X,\RealNumbers)$.
The following Theorem 
was proven by Verbitsky \cite{verbitsky-announcement}, 
and refined in the above language 
by Looijenga and Lunts \cite{looijenga-lunts}, who proved also a more 
general version for all compact K\"{a}hler manifolds.

\begin{thm}
\label{thm-so-action-determines-hodge-structure}
%(\cite{verbitsky-announcement} and \cite{looijenga-lunts})
Let $X$ be an irreducible holomorphic symplectic manifold. 
\begin{enumerate}
\item
\label{thm-tem-LL-Proposition-4.5}
(\cite{looijenga-lunts} Proposition 4.5)
$\LieAlg{g}(X)\cong \LieAlg{so}(4,b_2-2,\RealNumbers)$. 
$\LieAlg{g}(X)$ is defined over $\RationalNumbers$. 
The degree-zero summand
$\LieAlg{g}_0(X)$ is isomorphic to 
$\LieAlg{so}(H^2(X,\RealNumbers))\oplus \RealNumbers h$. 
The homomorphism 
\[
e \ : \ H^2(X,\RealNumbers) \ \ \ \rightarrow \ \ \ 
\End[H^*(X,\RealNumbers)],
\]
sending $a$ to $e_a$, is injective with image $\LieAlg{g}_2(X)$. 
$\LieAlg{g}_k(X)=0$, for $k\not\in\{-2,0,2\}$.
In particular, $\LieAlg{g}_2(X)$ and $\LieAlg{g}_{-2}(X)$ are abelian 
subalgebras.
\item
(\cite{verbitsky-mirror-symmetry})
The action of the semisimple part of $\LieAlg{g}_0(X)$, which is 
isomorphic to $\LieAlg{so}(H^2(X,\RealNumbers))$, integrates to an action of 
$Spin(H^2(X,\RealNumbers))$ via ring isomorphisms. The action on the even 
cohomology factors through a representation 
\begin{equation}
\label{eq-verbitskys-representation}
\rho \ : \ SO(H^2(X,\RealNumbers)) \ \ \longrightarrow \ \ 
\Aut[H^{even}(X,\RealNumbers)].
\end{equation}
\item
\label{thm-item-invariance-of-Poincare-pairing}
(\cite{looijenga-lunts} Proposition 1.6)
$\LieAlg{g}(X)$ preserves, infinitesimally, the Poincare pairing on
$H^*(X,\RealNumbers)$.
\item
\label{thm-item-prim-generates-as-A-2-module}
(\cite{looijenga-lunts} Corollary 2.3)
$H^*(X,\RealNumbers)$ is generated, as an $A_2$-module, by $Prim(X)$.
\item
\label{thm-item-irreducible-submodules-are-generated-by-their-prim}
(\cite{looijenga-lunts} Corollary 1.13)
Let $W$ be an irreducible $\LieAlg{g}_0(X)$-submodule of
$Prim^k(X)$. Then the $A_2$-submodule generated by $W$ is an irreducible 
$\LieAlg{g}(X)$-submodule. 
Conversely, all irreducible $\LieAlg{g}(X)$-submodules 
are of this type.
\item
\label{thm-item-so-action-determines-hodge-structure}
(\cite{verbitsky-mirror-symmetry})
Let $I$ be the complex structure of $X$. Denote by $ad_I$ the 
semisimple endomorphism
of $H^*(X,\ComplexNumbers)$, with $H^{p,q}(X)$ an eigenspace with 
eigenvalue $\sqrt{-1}(p-q)$. Then $ad_I$ is an element of 
$\LieAlg{g}_0(X)\otimes_\RealNumbers\ComplexNumbers$.
\end{enumerate}
\end{thm}
 
%The integrality properties of the $SO(H^2(X,\Integers))$ action, 
%on the cohomology of weight $i>2$, seem not to have been studied. 
We could not determine directly from Theorem
\ref{thm-so-action-determines-hodge-structure}, how large 
is the intersection of the image of 
$\rho$ with the monodromy group. One seems to need 
some Torelli type result (see section \ref{sec-monodromy-for-K3}). 
In that respect, Theorems 
\ref{introduction-thm-Gamma-v-acts-motivicly} and 
\ref{thm-so-action-determines-hodge-structure}
seem to complement each other nicely, in case $X=S^{[n]}$. 
Roughly, Theorem 
\ref{introduction-thm-Gamma-v-acts-motivicly} implies, that
the image of $\Mon\rightarrow \Mon^2$ is large 
(section \ref{sec-Mon-2-is-large}), while Theorem 
\ref{thm-so-action-determines-hodge-structure}
implies, that the kernel $K$ is small (Lemma \ref{lemma-bounding-K}).
The action (\ref{introduction-eq-monodromy-is-gamma-times-cov}) 
is compared to
Verbitsky's (\ref{eq-verbitskys-representation})
in Lemma \ref{lemma-comparing-two-representations}.

Let $A'_i$ be the $\LieAlg{g}(X)$-submodule of
$H^*(X,\RealNumbers)$ generated by $Prim(X)\cap A_i$. 
If $i\geq 2$, then $A'_i$ is the maximal $\LieAlg{g}(X)$-submodule of
$H^*(X,\RealNumbers)$, which is contained in $A_i$,
by parts \ref{thm-item-prim-generates-as-A-2-module} 
and \ref{thm-item-irreducible-submodules-are-generated-by-their-prim} 
of the Theorem. 
Clearly, we have the equalities $A'_2=A_2$ and 
$
(A'_i)^k =  H^k(X,\RealNumbers), 
\ \mbox{for} \ k\leq i.
% \ \mbox{and} \ k\geq 2d-i.
$
\begin{new-lemma}
\label{lemma-A'-i-equal-A-i}
$(A'_{i-2})^i \ \ \ = \ \ \ (A_{i-2})^i$, for $i\geq 4$.
\end{new-lemma}

\noindent
{\bf Proof:}
%We prove only the less obvious inclusion $(A'_{i-2})^i\supset (A_{i-2})^i$.
Let $x\in H^i(X,\RealNumbers)$. Then $x$ can be written in the form
${\displaystyle x=
\sum_{j=0}^{\lfloor\frac{i}{2}\rfloor} a_jx_j,
}$
where $a_j\in (A_2)^{2j}$, 
$x_j\in Prim^{i-2j}(X)$, 
by part \ref{thm-item-prim-generates-as-A-2-module} 
of the Theorem. Note that $Prim^2(X)=0$.
For $j\geq 1$, $x_j$ belongs to $(A'_{i-2})^{i-2j}$,
by the equality $(A'_{i-2})^{i-2j}=H^{i-2j}(X,\RealNumbers)$.
Hence, $a_jx_j$ belongs to $(A'_{i-2})^i$, for $j\geq 1$.
Thus, $x$ belongs to $(A_{i-2})^i$ (respectively $(A'_{i-2})^i$), 
if and only if 
$x_0$ belongs to $Prim^i(X)\cap (A_{i-2})^i$ 
(respectively $Prim^i(X)\cap (A'_{i-2})^i$). 
But $Prim^i(X)\cap (A_{i-2})^i=Prim^i(X)\cap (A'_{i-2})^i$,
by definition of $A'_{i-2}$. Hence, 
$x$ belongs to $(A_{i-2})^i$, 
if and only if it belongs to $(A'_{i-2})^i$.
\EndProof

\smallskip
The Poincare pairing restricts, as a {\em non-degenerate pairing}, 
to each irreducible $\LieAlg{g}(X)$-submodule,
by part \ref{thm-item-invariance-of-Poincare-pairing} of the Theorem
and the fact that $\LieAlg{g}(X)$ is semi-simple.
In particular, it restricts to a non-degenerate pairing on $A'_k$. 
Set $C_2:=H^2(X,\Integers)$, $C_3:=H^3(X,\Integers)_{free}$ and let 
\begin{equation}
\label{eq-C-k}
C_k\subset H^k(X,\Integers)_{free}, \ \ \  k\geq 4,
\end{equation}
be the weight $k$ summand of the orthogonal complement $(A'_{k-2})^\perp$.
Set $C_k(\RationalNumbers):=C_k\otimes_\Integers\RationalNumbers$. 
Let $K\subset Mon(X)$ be the subgroup of monodromy operators, 
which act trivially on $H^2(X,\Integers)$.

\begin{cor}
\label{cor-monodromy-invariant-decomposition}
\begin{enumerate}
\item
\label{cor-item-decomposition}
$H^k(X,\RationalNumbers)$ admits the monodromy-invariant 
and $\LieAlg{g}_0(X)$-invariant decomposition
(\ref{eq-monodromy-invariant-decomposition-of-H-k}), for $k\geq 2$.
\item
\label{cor-item-K-has-exponent-leq-2}
Assume, that each irreducible 
$\LieAlg{g}_0(X)$-submodule of $C_k(\RationalNumbers)$ appears 
in $C_k(\RationalNumbers)$ with multiplicity at most $1$, for all $k\geq 2$. 
Then $K$ is a finite subgroup of exponent $\leq 2$. 
\item
\label{cor-item-K-is-central}
Assume further, that none of the $\LieAlg{g}_0(X)$-modules
$C_k(\RationalNumbers)$, $k\geq 2$, contains two 
irreducible $\LieAlg{g}_0(X)$-submodules with isomorphism classes, 
which are conjugate under the 
automorphism group of $\LieAlg{g}_0(X)$.
Then $K$ is contained in the center of $Mon(X)$.
\end{enumerate}
\end{cor}

\noindent
{\bf Proof:}
\ref{cor-item-decomposition}) $C_k(\RationalNumbers)$ and 
$[A'_{k-2}(\RationalNumbers)]^k$ are complementary, by definition, so 
(\ref{eq-monodromy-invariant-decomposition-of-H-k}) is a direct sum 
decomposition, by Lemma \ref{lemma-A'-i-equal-A-i}.
The $\LieAlg{g}_0(X)$-invariance is clear. 
$\LieAlg{g}(X)$ is a $Mon(X)$-invariant Lie subalgebra of 
$\End[H^*(X,\RealNumbers)]$, by its definition.
Thus the decomposition (\ref{eq-monodromy-invariant-decomposition-of-H-k}) 
is $Mon(X)$-invariant as well.

\ref{cor-item-K-has-exponent-leq-2})
Let $f$ be an element of $K$. Then $f$ commutes with $\LieAlg{g}_0(X)$,
by definition of $\LieAlg{g}_0(X)$. 
Thus, $f$ acts on each irreducible $\LieAlg{g}_0(X)$-submodule of
$C_k(\RationalNumbers)$ via multiplication by $\pm 1$, by the
multiplicity assumption. The subspaces $C_k(\RationalNumbers)$, $k\geq 2$,
generate $H^*(X,\RationalNumbers)$. Hence, $K$ is finite and 
has exponent $\leq 2$. 

\ref{cor-item-K-is-central})
$Mon(X)$ acts on $\LieAlg{g}_0(X)$ via Lie algebra automorphisms. 
Each irreducible $\LieAlg{g}_0(X)$-submodule of
$C_k(\RationalNumbers)$ is $Mon(X)$-invariant, by the
assumed absence of $\Aut\LieAlg{g}_0(X)$-conjugate 
$\LieAlg{g}_0(X)$-submodules. Thus $K$ is in the center of $Mon(X)$.
\EndProof

\smallskip
The assumptions in parts 
\ref{cor-item-K-has-exponent-leq-2} and \ref{cor-item-K-is-central}
of Corollary \ref{cor-monodromy-invariant-decomposition} hold 
for the Hilbert schemes $S^{[n]}$, of a $K3$ surface $S$
(Lemma \ref{lemma-Mon-invariant-irreducible-generating-subspaces} below).
These assumptions fail to hold
for generalized Kummers of dimension $2n\geq 4$. 
Indeed, $K$ has exponent $>2$ and is not abelian
in that case. This can be seen as follows.
The subgroup, of automorphisms of an abelian surface, 
generated by multiplication by $-1$ and by translations by points 
of order $n+1$, acts trivially on the second cohomology of the generalized 
Kummer (\cite{beauville-automorphisms} proposition 9), but it acts 
faithfully on the higher cohomology \cite{salamon}. 
Hence it injects into $K$.
%****************************************************************
% A splitting of the monodromy group
%****************************************************************
\subsection{A splitting of the monodromy group}
\label{sec-splitting}
In view of the local Torelli Theorem for 
irreducible holomorphic symplectic manifolds 
\cite{beauville-varieties-with-zero-c-1}, it is natural to compare
the monodromy of the total cohomology ring of the Hilbert scheme $S^{[n]}$
and the monodromy for the weight $2$ cohomology. 
Our model is the case $n=1$; in which the monodromy group for K3 
surfaces is an index 2 subgroup of
the isometry group of the lattice $H^2(S,\Integers)$. 
It consists of isometries, which preserve the orientation of the 
positive cone.
This identification of the monodromy is a consequence of the Torelli Theorem
for K3 surfaces and the surjectivity of the period map
(Corollary \ref{cor-Mon-S}). 

Let $\Mon$ be the group of automorphisms of the ring
$H^*(S^{[n]},\Integers)_{\rm free}$, which is generated by monodromy
(Definition \ref{def-monodromy}). 
Denote by $\Mon^2$ the image of $\Mon$ in
the isometry group of $H^2(S^{[n]},\Integers)$ and let
$K$ be the kernel
\begin{equation}
\label{eq-exact-seq-of-monodromies}
0 \ \rightarrow \ K \ \rightarrow \ \Mon \ \rightarrow \ 
\Mon^2 \ \rightarrow \ 0.
\end{equation}
%If a version of the Torelli theorem holds for birational classes of 
%compact hyperk\"{a}hler varieties, which are 
%deformation equivalent to $S^{[n]}$, 
%we would expect that 1) $\Mon^2$ is a large subgroup of the
%isometry group of the lattice $H^2(S^{[n]},\Integers)$
%and 2) the restriction homomorphism 
%$\Mon \rightarrow \Mon^2$ is an isomorphism 
%(Conjecture \ref{conj-Mon-injects-into-Mon-2}). 
%Theorem \ref{introduction-thm-Gamma-v-acts-motivicly} 
%verifies the first prediction.
%We elaborate further on the size of $\Mon^2$ in section 
%\ref{sec-Mon-2-is-large}. 
%Theorem \ref{introduction-thm-Gamma-v-acts-motivicly} verifies 
%also a weaker version of the second prediction, namely, that
The exact sequence (\ref{eq-exact-seq-of-monodromies})
{\em splits naturally} in case  $n-1$ is a prime power,
by Theorem \ref{introduction-thm-Gamma-v-acts-motivicly}. 
More generally, the pullback of (\ref{eq-exact-seq-of-monodromies}),
via the inclusion\footnote{The equality $\W=\Mon^2$ is proven in 
\cite{markman-part-two}} 
of $\W$ in $\Mon^2$, splits (see Corollary
\ref{cor-W-consists-of-monodromies} for the inclusion $\W\subset \Mon^2$).
%(see Lemma \ref{lemma-the-index-of-Gamma-v}). 
The splitting is natural, since the image of $\Gamma_v$ in $\Mon$ 
is a normal subgroup, which is isomorphic to $\W$
(Lemmas \ref{lemma-bounding-K} and \ref{lemma-EO-is-maximal}). 
Conjecture \ref{conj-Mon-injects-into-Mon-2} suggests that $K$ is trivial.
%
%Generalized Kummer varieties admit automorphisms, 
%which act trivially on their second cohomology 
%\cite{beauville-automorphisms}.
%These automorphisms act non-trivially on the higher cohomology
%\cite{salamon}.
%Consequently, the sequence analogous to (\ref{eq-exact-seq-of-monodromies}) 
%above has non-trivial kernel. 
%Beauville proved, that 
%any automorphism of $S^{[n]}$, which acts
%trivially on $H^2(S^{[n]})$, is the identity (Proposition
%\ref{prop-beauville-automorphisms}). 
The rest of this section is dedicated to the study of $K$. 
We start with a rough estimate of $K$ in 
the following lemma.

\begin{new-lemma}
\label{lemma-bounding-K}
$K$ is contained in the center of $\Mon$. Furthermore, 
$K$ is isomorphic to a subgroup of $(\Integers/2\Integers)^{2n-4}$. 
Consequently, the image of the monodromy representation 
(\ref{introduction-eq-monodromy-is-gamma-times-cov}) 
of $\Gamma_v$ is a normal subgroup of finite index in $\Mon$.
\end{new-lemma}

Further constraints on $K$ 
are discussed in section \ref{sec-chern-classes-of-tangent-bundle}.
Lemma \ref{lemma-bounding-K} is a corollary of the following lemma.
%Let $A_{k}\subset H^*(S^{[n]},\Integers)_{\rm free}$, $k\geq 1$,be the 
%subalgebra generated by 
%$\oplus_{i=1}^k H^{2i}(S^{[n]},\Integers)_{\rm free}$.
%The Poincare pairing restricts 
%to a non-degenerate pairing on $A_{k}$, 
%being an $A_{1}$-submodule of $H^*(S^{[n]},\Integers)_{\rm free}$. 
%The non-degeneracy follows from the definiteness of the polarization 
%of the Hodge structure, induced by a K\"{a}hler form. 
%Set $C_1:=H^2(S^{[n]},\Integers)$ and let 
%$C_{k+1}\subset H^{2k+2}(S^{[n]},\Integers)_{\rm free}$, 
%$k\geq 1$,  be the weight $2k+2$ summand of the orthogonal
%complement $A_{k}^\perp$. 
Let $C_k\subset H^k(S^{[n]},\Integers)_{free}$ be the 
$Mon$-submodule defined in (\ref{eq-C-k}).
Then $C_2$, $C_4$, \dots, $C_{2n-2}$ generate the cohomology ring 
with rational 
coefficients (Lemma \ref{lemma-10-in-markman-diagonal}).
%(part 3 of Lemma 10 in \cite{markman-diagonal}). 
%Clearly, $C_i$ is $\Mon$-invariant, and 
Thus $\Mon$ is determined by its action on the $C_i$'s. 
%Set $C_i(\RationalNumbers):=C_i\otimes_\Integers\RationalNumbers$.
Let $\monrep^{2i}:\Gamma_v\rightarrow GL(C_{2i}(\RationalNumbers))$
be the composition of $\monrep$, given in 
(\ref{introduction-eq-monodromy-is-gamma-times-cov}),
with the restriction homomorphism from 
$\Mon$ to $GL(C_{2i}(\RationalNumbers))$.
Set $\theta_{v,1}:=\theta_v$, 
given in (\ref{eq-theta-v-from-v-perp}), and let 
$\theta_{v,i}:H^*(S,\RationalNumbers)\rightarrow C_{2i}(\RationalNumbers)$, 
$i\geq 2$, be the composition of the homomorphism
(\ref{eq-homomorphism-induced-by-u-v}) with the projection 
$H^*(S^{[n]},\RationalNumbers)\rightarrow C_{2i}(\RationalNumbers)$.

\begin{new-lemma}
\label{lemma-Mon-invariant-irreducible-generating-subspaces}
$C_{2i}(\RationalNumbers)$ admits a $\Mon$-invariant decomposition
\begin{equation}
\label{eq-C-k-decomposes-into-at-most-two-irreducibles}
C_{2i}(\RationalNumbers) \ \ = \ \ 
C'_{2i}(\RationalNumbers) \ \oplus \ C''_{2i}(\RationalNumbers).
\end{equation}
$C'_{2i}(\RationalNumbers)$ either vanishes, or it is a
one-dimensional character $\chi'_i$ of $\Mon$. The summand 
$C''_{2i}(\RationalNumbers)$ either vanishes, or it is an
irreducible $\Mon$ representation, which is isomorphic to
$H^2(S^{[n]},\RationalNumbers)\otimes \chi''_i$, for a 
character $\chi''_i:\Mon\rightarrow \{\pm 1\}$. 
The pullback of $\chi'_i$ to $\Gamma_v$, via
(\ref{introduction-eq-monodromy-is-gamma-times-cov}), 
is the orientation character $cov$,
if $C'_{2i}(\RationalNumbers)$ does not vanish and $i$ is odd. 
The pullback of $\chi'_i$ to $\Gamma_v$ is trivial, when $i$ is even. 
The pullback of $\chi''_i$ to $\Gamma_v$ is the orientation character $cov$,
when $C''_{2i}(\RationalNumbers)$ does not vanish and $i$ is even. 
The pullback of $\chi''_i$  to $\Gamma_v$ is trivial, when $i$ is odd.  
The summand $C'_{2i}(\RationalNumbers)$ does not vanish,
for $2\leq i\leq n/2$, and $C''_{2i}(\RationalNumbers)$ does not vanish, 
for $0<i\leq (n+1)/2$. 
Consequently, $\theta_{v,i}$ is an isomorphism, for $1\leq i\leq n/2$, and
\[
\monrep^{2i}(g) \ \ \ = \ \ \ (-1)^{i\cdot cov(g)}
[\theta_{v,i}\circ g\circ \theta_{v,i}^{-1}].
\]
\end{new-lemma}

{\bf Proof of Lemma \ref{lemma-bounding-K}:}
The kernel $K$ of $\Mon\rightarrow \Mon^2$
acts on $C''_{2i}(\RationalNumbers)$ diagonally, via the 
character $\chi''_i$. 
We conclude, that the action of $K$ is determined by the $2n-4$
characters $\{\chi'_i,\chi''_i\}$, $2\leq i\leq {n-1}$. 
The image of the monodromy representation 
(\ref{introduction-eq-monodromy-is-gamma-times-cov}), in the
isometry group of $H^2(S^{[n]},\Integers)$, 
is a normal subgroup of finite index in $\Mon^2$ 
(Lemma \ref{lemma-the-index-of-Gamma-v}).
Hence, the image of the monodromy representation 
(\ref{introduction-eq-monodromy-is-gamma-times-cov}) is  a normal subgroup of 
finite index in $\Mon$.
\EndProof

{\bf Proof of lemma 
\ref{lemma-Mon-invariant-irreducible-generating-subspaces}:}
The decomposition 
(\ref{eq-C-k-decomposes-into-at-most-two-irreducibles})
of $C_{2i}(\RationalNumbers)$ is constructed in three steps:
%1) We determine the possible Hodge numbers of $C_{2i}(\RationalNumbers)$
%using its 
%We define it first as the decomposition 
%into irreducibles with respect to the monodromy representation 
%(\ref{introduction-eq-monodromy-is-gamma-times-cov}) of $\Gamma_v$. 
%2) We establish the $\LieAlg{g}_0(S^{[n]})$-invariance of the decomposition.
%3) We conclude the $Mon$-invariance of the decomposition.
Step 1: Let us first determine the possible Hodge numbers of 
$C_{2i}(\RationalNumbers)$. The homomorphism
$\theta_{v,i}:H^*(S,\Integers)\rightarrow C_{2i}(\RationalNumbers)$ 
is surjective, 
%Each subspace $B_i$ in corollary 
%\ref{cor-generating-subspaces-are-sub-representations}
%projects onto $C_{2i}(\RationalNumbers)$,
by Corollary \ref{cor-kunneth-factors-generate}, and  
%This projection 
is $\Gamma_v^{cov}$-equivariant with respect to the
action (\ref{introduction-eq-monodromy-is-gamma-times-cov}), by Theorem 
\ref{introduction-thm-Gamma-v-acts-motivicly}. 
Hence, $C_{2i}(\RationalNumbers)$  contains at most two irreducible
$\Gamma_v^{cov}$-modules.
The rank of $\theta_{v,i}$
%the projection $B_i\rightarrow C_{2i}(\RationalNumbers)$
can be $24$, $23$, $1$, or $0$. 
Moreover, $\theta_{v,i}$ is 
%the projection induces 
a rational Hodge structure homomorphism.
Consequently, the possible Hodge numbers of $C_{2i}$ are as follows.
$[h^{i-1,i+1},h^{i,i},h^{i+1,i-1}]$ 
is $(1,22,1)$, $(1,21,1)$, $(0,1,0)$, or $(0,0,0)$,
and all the other $h^{k,2i-k}$ vanish. 

Step 2: The decomposition 
(\ref{eq-C-k-decomposes-into-at-most-two-irreducibles})
of $C_{2i}(\RationalNumbers)$ is defined to be 
the decomposition into irreducible
$\LieAlg{g}_0(S^{[n]})$-submodules. 
If $h^{i-1,i+1}$ does not vanish, $C_{2i}(\RealNumbers)$ 
can not be a trivial $\LieAlg{so}(3,20,\RealNumbers)$-representation, because
the Hodge structure is determined by the 
$\LieAlg{so}(3,20,\RealNumbers)$-module structure
(Theorem \ref{thm-so-action-determines-hodge-structure}
part \ref{thm-item-so-action-determines-hodge-structure}).
$\LieAlg{so}(3,20,\RealNumbers)$ does not have an irreducible representation
of dimension $24$. 
The minimal dimension of a non-trivial irreducible representation
of $\LieAlg{so}(3,20,\RealNumbers)$ is the $23$-dimensional
standard representation. 
Summarizing, 
%Theorem \ref{thm-so-action-determines-hodge-structure}  implies, that 
$C_{2i}(\RealNumbers)$ is a direct sum 
$
C'_{2i}(\RealNumbers) \ \oplus \ C''_{2i}(\RealNumbers)
$
with respect to Verbitsky's $SO(H^2(S^{[n]},\RealNumbers))$
action. The decomposition is defined over $\RationalNumbers$, because 
$\LieAlg{so}(3,20,\RealNumbers)$ is defined over $\RationalNumbers$,
by Theorem \ref{thm-so-action-determines-hodge-structure} part 
\ref{thm-tem-LL-Proposition-4.5}.
The summand $C'_{2i}$ is either zero or
the trivial representation and $C''_{2i}$ is either zero or 
the standard representation. The non-vanishing of 
$C'_{2i}$ and $C''_{2i}$, in the stated ranges of $i$, follow 
from Lemma \ref{lemma-10-in-markman-diagonal}.
%Lemma 10 in \cite{markman-diagonal}. 

Step 3: When non-zero, 
$C'_{2i}(\RealNumbers)$ and $C''_{2i}(\RealNumbers)$ are 
irreducible subrepresentations of $\LieAlg{g}_0(S^{[n]})$ 
of different dimensions.
Hence, each is $Mon$-invariant, and $K$ is finite, central,
and of exponent $\leq 2$, by Corollary
\ref{cor-monodromy-invariant-decomposition}.

It remains to define the character $\chi''_i$. 
Consider the representation 
$H^2(S^{[n]},\RationalNumbers)\otimes C''_{2i}(\RationalNumbers)$
of $\Mon$. Restrict it to a representation
of $(\Gamma_v)^{cov}$ via 
(\ref{introduction-eq-monodromy-is-gamma-times-cov}). 
This restriction is isomorphic to the tensor square of
$H^2(S^{[n]},\RationalNumbers)$. It decomposes as a direct sum of 
three pairwise non-isomorphic and irreducible representations; 
the trivial character $V$ of $(\Gamma_v)^{cov}$,
$\Wedge{2}H^2(S^{[n]},\RationalNumbers)$, and another irreducible
representation (\cite{looijenga-lunts} Proposition 2.14). 
The image of $(\Gamma_v)^{cov}$ in
$\Mon$ is a normal subgroup, by Lemma \ref{lemma-the-index-of-Gamma-v} 
and the centrality of $K$. 
Hence, each of the three subrepresentations is $\Mon$-invariant.
$V$ determines the character $\chi''_i$.
\EndProof

\subsection{The Chern classes of $TS^{[n]}$}
\label{sec-chern-classes-of-tangent-bundle}

%One may be able to bound by $2$ the order $\Abs{K}$ of the kernel $K$ 
%of $\Mon\rightarrow \Mon^2$, 
%using further knowledge of the ring structure of 
%$H^*(S^{[n]},\RationalNumbers)$ (see \cite{lehn-sorger} 
%for a determination of the ring structure). 
%One may even be able to prove the injectivity of $\Mon\rightarrow \Mon^2$,
%using further knowledge of the integral structure of
%$H^*(S^{[n]},\Integers)$. 
%We continue the study of the kernel $K$ of $\Mon\rightarrow \Mon^2$.
%We formulate here conditions, that imply the bound $\Abs{K}\leq 2$,
%and prove some of them. 

We continue the study of the kernel $K$ of $\Mon\rightarrow \Mon^2$
and improve the upper bound for its order $\Abs{K}$, given in
Lemma \ref{lemma-bounding-K}. 
We conjectured $K$ to be trivial (Conjecture
\ref{conj-Mon-injects-into-Mon-2}).
$K$ is embedded in $(\Integers/2\Integers)^{2n-4}$ via the 
restriction of the characters 
$\chi'_i$ and $\chi''_i$, $2\leq i \leq n-1$ (Lemma \ref{lemma-bounding-K}). 
We prove in this section 
that the character $\chi'_i$ is trivial, for even $i$
(Lemma \ref{lemma-chern-classes-project-to-generators}). 
The restriction of $\chi''_2$ to $K$ is trivial for $n\geq 4$, by
part 1 of Theorem 34 of
\cite{markman-integral-generators}.
The equality $\chi''_{2i}=\chi''_2$ is proven 
in part 3 of Theorem 34 of
\cite{markman-integral-generators} for
integer values of $i$ in the range
$1\leq i\leq n/4$ mentioned in Lemma 
\ref{lemma-Mon-invariant-irreducible-generating-subspaces}.

In case $n=2$, the kernel is trivial
(Lemma \ref{lemma-bounding-K}). In case $n=3$, 
$\Abs{K}$ is bounded by $2$, by
the vanishing of $C'_{2i}$, for all $i$,
and of $C''_{2i}$, for $i>2$ (Example 14 in \cite{markman-diagonal}).
Question 1 in \cite{markman-diagonal} asks, if $C'_{2i}$ and 
$C''_{2i}$ vanish,
for $i$ outside the corresponding ranges in Lemma 
\ref{lemma-Mon-invariant-irreducible-generating-subspaces}. 
This suggested range of vanishing seems to be a close approximation of the 
truth \cite{lehn-sorger-generators}. 
When $n=4$, the triviality of $K$ follows from an affirmative 
answer to Question 1 in \cite{markman-diagonal} and
Lemma \ref{lemma-chern-classes-project-to-generators}. 

Let $v=(1,0,1-n)$ be the Mukai vector of $S^{[n]}$, 
$\pi_\M$ and $\pi_S$ the two projections from
$\M(v)\times S$, and $\E_v$ a universal sheaf. 
Consider the decomposition 
$H^{4i}(\M(v),\RationalNumbers)=(A_{4i-2})^{4i}(\RationalNumbers)\oplus
C'_{4i}(\RationalNumbers)\oplus C''_{4i}(\RationalNumbers)$, $i\geq 1$,
introduced in Lemma 
\ref{lemma-Mon-invariant-irreducible-generating-subspaces}. 
Let $E$ be a sheaf with Mukai vector $v$. 
%Let $C'_{4i}$ be the summand, in the decomposition of
%$H^{4i}(S^{[n]},\RationalNumbers)$, 
%introduced in Lemma 
%\ref{lemma-Mon-invariant-irreducible-generating-subspaces}.

\begin{new-lemma}
\label{lemma-chern-classes-project-to-generators}
The projections to $C'_{4i}(\RationalNumbers)$, of $c_{2i}(TS^{[n]})$ 
and $2c_{2i}[-\pi_{\M_!}(\E_v\otimes p_S^!(E^\vee))]$, are equal and 
span $C'_{4i}(\RationalNumbers)$. Consequently, we get: 
\begin{enumerate}
\item 
If $4\leq 2i\leq n/2$, then $c_{2i}(TS^{[n]})$ does not belong to the 
subalgebra of $H^*(S^{[n]},\RationalNumbers)$ generated by 
classes in $H^k(S^{[n]},\RationalNumbers)$, $k<4i$. 
\item 
The character $\chi'_{2i}$ is trivial, for all $i$.
\end{enumerate}
\end{new-lemma}

\noindent
{\bf Proof:}
Denote the $\Gamma_v$-equivariant homomorphism 
(\ref{eq-homomorphism-induced-by-u-v})
by $u_v$ as well. 
The projection of $ch_{2i}[-\pi_{\M_!}(\E_v\otimes p_S^!(E^\vee))]$
to $C'_{4i}$ is equal to the projection of $u_v(v)$ 
and hence spans $C'_{4i}$. 
The coset $c_{2i}(x)+(A_{4i-2})^{4i}$,
$x\in K_{top}(\M(v))$, is a multiple of $ch_{2i}(x)+(A_{4i-2})^{4i}$ 
by a universal constant. Hence, 
it suffices to prove the equality of the projections  
to $C'_{4i}(\RationalNumbers)$ of $ch_{2i}(TS^{[n]})$ and $2u_v(v)$.

%Each of the statements concerns a one dimensional 
%subrepresentation $\bar{\ell}$ in $H^{4j}$ or its quotient.
%We need to show, that $\bar{\ell}$ is the trivial character of $\Mon$.
%We show, that the projection of the Chern class $c_{2j}(TS^{[n]})$
%to $\bar{\ell}$ is equal (modulo classes in $A_{2i-1}$ in part
%\ref{lemma-item-chi-prime-2i} and classes in $H^2\cup A_{2i}$ in part
%\ref{lemma-item-chi-double-prime-2i-plus-1}) to the contribution 
%of the K\"{u}nneth factors of the Chern classes of the universal sheaf.
%Since these  K\"{u}nneth factors generate the ring (Corollary
%\ref{cor-kunneth-factors-generate}), 
%then $\bar{\ell}$ is generated by the projection of 
%$c_{2j}(TS^{[n]})$. It follows, that 
%$\bar{\ell}$ is the trivial character of $\Mon$. 
%
We express next the Chern classes of $TS^{[n]}$ in terms of the 
K\"{u}nneth factors of the Chern classes of the universal sheaf. 
The tangent bundle is isomorphic to the relative extension sheaf
$\RelExt^1_{\pi_\M}(\E_v,\E_v)$, while the sheaves
$\RelExt^i_{\pi_\M}(\E_v,\E_v)$, $i=0,2$, are the trivial line bundles.
We get the equality of Chern polynomials 
\[
c(T\M(v)) \ \ = \ \ 
c\left(-\pi_{\M_!}\left[\E_v^\vee\stackrel{L}{\otimes}\E_v\right]\right).
\]
The class, on the right hand side, is independent of the choice of a
universal sheaf. Grothendieck-Riemann-Roch yields the equality
\[
ch(T\M(v)-2\StructureSheaf{\M(v)}) \ \ = \ \ 
-\pi_{\M_*}(u_v^\vee\cup u_v). 
\]

Let $\Lambda:=H^*(S,\Integers)$ be the Mukai lattice and 
$\{v_1, v_2, \dots, v_{24}\}$ an orthonormal basis of 
$\Lambda\otimes\RationalNumbers$.
%, such that $v_1$ is in ${\rm span}\{v\}$. 
%The dual basis, with respect to the cup product pairing, is $\{-D(v_i)\}$.
We have the equalities
$u_v=\sum_i v_i\otimes u_v(v_i)$ and 
\[
%\begin{equation}
%\label{eq-decomposition-of-chern-character-of-tangent-bundle}
\pi_{\M_*}(u_v^\vee\cup u_v) \ = \ 
-\sum_i\sum_j(v_i,v_j)\cdot u_v(v_i)\cup u_v(v_j)^\vee \ = \ 
-\sum_iu_v(v_i)\cup u_v(v_i)^\vee.
%\end{equation}
\]
More invariantly, 
the class $-\pi_{\M_*}(u_v^\vee\cup u_v)$ is the image of the  class 
$q=\sum_iv_i\otimes v_i$, 
inverse to the Mukai pairing, 
in the tensor square of the Mukai lattice,
under the composition
\[
%H^*(S)\otimes H^*(S) 
%\LongRightArrowOf{D\otimes D}
\Lambda\otimes \Lambda
\LongRightArrowOf{u_v\otimes u_v} 
H^*(\M(v))\otimes H^*(\M(v)) 
\LongRightArrowOf{1\otimes D_{\M(v)}}
H^*(\M(v))\otimes H^*(\M(v)) 
\RightArrowOf{\cup} 
H^*(\M(v)).
\]

Let $u_{2k}:\Lambda\rightarrow H^{2k}(\M(v),\RationalNumbers)$ 
be the composition of $u_v$ with the projection onto the graded summand.
Then  $ch_{2i}(T\M(v))$ is the image of $q$ under
\begin{equation}
\label{eq-sum-of-u-2j}
\sum_{j=0}^{2i}
u_{2j}\cup D_{\M(v)}\circ u_{4i-2j} \ : \ 
\Lambda\otimes \Lambda \ \ \ \longrightarrow\ \ \ 
H^{4i}(\M(v),\RationalNumbers). 
\end{equation}
Let $\bar{u}_{2k}:\Lambda\rightarrow 
H^{2k}(\M(v),\RationalNumbers)/
[A_{2k-1}(\RationalNumbers)^{2k}+C_{2k}''(\RationalNumbers)]
\cong C'_{2k}(\RationalNumbers)$ be the composition of
$u_{2k}$ with the projection to $C'_{2k}(\RationalNumbers)$, $k\geq 1$.
Set $\bar{u}_{0}:=u_0$. 
The projection of $ch_{2i}(T\M(v))$ to
$C'_{4i}(\RationalNumbers)$ is the image of $q$ under the composition
\[
\Lambda\otimes \Lambda \longrightarrow 
\left([\Lambda/v^\perp]\otimes[\Lambda/v^\perp]\right)
\oplus \left([\Lambda/v^\perp]\otimes[\Lambda/v^\perp]\right)
\LongRightArrowOf{(\bar{u}_{4i}\cup u_0)+(u_0\cup \bar{u}_{4i})}
C'_{4i}.
\]
We ignored above the intermediate summands of 
(\ref{eq-sum-of-u-2j}), since their images are in $A_{4i-2}$.

Assume, as we may, that $n\geq 2$, so that $(v,v)\neq 0$. The images 
of $q$ and $\frac{v\otimes v}{(v,v)}$ in 
$[\Lambda/v^\perp]\otimes[\Lambda/v^\perp]$ are equal and 
$u_0(v)=(v,v)\cdot 1\in H^0(\M(v),\Integers)$. 
Thus, the image of $q$ is $2\bar{u}_{4i}(v)$. 
We conclude, that the image of $ch_{2i}(T\M(v))$ 
in $C'_{4i}(\RationalNumbers)$ is equal to $2\bar{u}_{4i}(v)$. 
%Modulo $A_{4i-2}$, $\bar{u}_{4i}(v)$ is 
%equal to the image in $C'_{4i}(\RationalNumbers)$ of
%$ch_{2i}[-\pi_{\M_!}(\E_v\otimes p_S^!(v^\vee))]$. 
\EndProof

\subsection{
The monodromy subgroup of $O(H^2(S^{[n]},\Integers))$
}
\label{sec-Mon-2-is-large}

%\begin{question}
%Is the image $\Mon^2$, of the monodromy group,
%equal to the kernel $O_+(H^2(S^{[n]},\Integers))$ of
%the orientation character (\ref{eq-top-homology-character}) of
%$O(H^2(S^{[n]},\Integers))$?
%\end{question}
%
%An affirmative answer to the above question, seems to be a necessary 
%condition for a weight $2$ birational Torelli theorem to hold. 
%The index of $\Mon^2$ in $O_+(H^2(S^{[n]},\Integers))$ is an
%approximation for the ``degree of the period map''.
%We address the question next. 

We describe the image of $\Gamma_v$ in $O(H^2(S^{[n]},\Integers))$, 
under the monodromy representation of Theorem
\ref{introduction-thm-Gamma-v-acts-motivicly}.
Consider the homomorphism 
\begin{equation}
\label{eq-homomorphism-from-Gamma-v-to-EO}
f \ : \ O(v^\perp) \ \ \longrightarrow \ \ O_+(v^\perp),
\end{equation}
into the kernel of the orientation character 
(\ref{eq-top-homology-character})
of $O(v^\perp)$, sending $g$ to $(-1)^{cov(g)}\cdot g$. 
%Set $EO:=f(\Gamma_v)$. Then $EO$ 

\begin{new-lemma}
\label{lemma-EO-is-maximal}
The following subgroups of $O_+(v^\perp)$ are all equal. 
\begin{enumerate}
\item
\label{lemma-item-image-under-monodromy-rep}
The image via the monodromy representation 
$\monrep$ of $\Gamma_v$ in $O(H^2(S^{[n]},\Integers))$,
under the identification $v^\perp\cong H^2(S^{[n]},\Integers)$ of
Theorem \ref{thm-irreducibility}. 
\item
\label{lemma-item-W}
The subgroup $\W(v^\perp)$ given in (\ref{eq-W}).
\item
\label{lemma-item-extendible-isometries}
The subgroup of orientation preserving isometries, 
which can be extended to isometries of the whole Mukai lattice.
\item
\label{lemma-item-image-under-f}
The image $f(\Gamma_v)$, of the stabilizer $\Gamma_v$ of $v$ in the
Mukai lattice, under the homomorphism 
(\ref{eq-homomorphism-from-Gamma-v-to-EO}).
\item
\label{lemma-item-inverse-image-of-minus-1}
The inverse image of $\{1,-1\}$, under the natural homomorphism 
$O_+(v^\perp) \rightarrow \Aut[(v^\perp)^*/v^\perp]$.
\end{enumerate}
\end{new-lemma} 
 
\noindent 
{\bf Proof:} 
\ref{lemma-item-image-under-monodromy-rep} $=$
\ref{lemma-item-image-under-f}: Clear from the equality
$\theta_v\circ f(g)\circ \theta_v^{-1}=\monrep(g)$, $g\in\Gamma_v$, 
where $\theta_v$ is the 
isomorphism (\ref{eq-mukai-homomorphism}) used in Theorem 
\ref{thm-irreducibility}.

\ref{lemma-item-extendible-isometries} $=$
\ref{lemma-item-image-under-f}:
The inclusion \ref{lemma-item-image-under-f} $\subset$
\ref{lemma-item-extendible-isometries} is clear. 
We prove the reverse inclusion.
Let $\tilde{g}$ be an isometry of the whole Mukai lattice, which
leaves $v^\perp$ invariant and restricts to an orientation
preserving isometry $g$ of $v^\perp$. Then $\tilde{g}(v)=v$ or $-v$. 
If $\tilde{g}(v)=v$, then $\tilde{g}$ is an orientation preserving 
isometry in $\Gamma_v$ and $f(\tilde{g})=g$. 
If $\tilde{g}(v)=-v$, then 
$\tilde{g}$ reverses the orientation of the positive cone of the Mukai 
lattice, since $v$ belongs to that cone. 
Consequently, $-\tilde{g}$ belongs to $\Gamma_v$ and is orientation reversing.
Hence, $f(-\tilde{g})=g$. 

\ref{lemma-item-W} $=$ \ref{lemma-item-image-under-f}: 
The inclusion $\W(v^\perp)\subset f(\Gamma_v)$ is clear. 
%We prove the reverse inclusion. 
%Wall proved, that the isometry group of $H^2(S,\Integers)$ is
%generated by reflections, with respect to 
%$+2$ and $-2$ vectors (\cite{wall} Theorem 4.8). 
%It follows from Wall's Theorem and Proposition \ref{prop-compare-Gamma-v},
%that $\Gamma_v$ is generated by reflections, with respect to 
%$+2$ and $-2$ vectors in $v^\perp$. Now 
The reverse inclusion follows from Corollary 
\ref{cor-stabilizer-is-generated-by-reflections}
and the fact, that $f$ leaves invariant each reflection 
with respect to a $-2$ vector, and multiplies by $-1$ each 
reflection with respect to a $+2$ vector.

\ref{lemma-item-image-under-f} $=$
\ref{lemma-item-inverse-image-of-minus-1}:
Follows immediately from  Lemma \ref{lemma-the-index-of-Gamma-v}. 
\EndProof

\medskip
When $n-1$ is a prime power, $\W(v^\perp)$ is the whole of
$O_+(v^\perp)$, by Lemma
\ref{lemma-the-index-of-Gamma-v}. 
%The subgroup $EO$ is contained in the 
%kernel of the orientation character (\ref{eq-top-homology-character})
%of $O(v^\perp)$. 
%The notation $EO$ stands for extendible isometries; $EO$ 
%satisfies the following maximality property.
%
%Extendibility to the whole Mukai lattice is an integral obstruction, 
%for elements of $O(H^2(S^{[n]},\Integers))$, to belong to $EO$. 
%$EO$ may be a proper subgroup of the image $\Mon^2$ of the monodromy group. 
%On the other hand, similar integral obstructions may exist, 
%for elements of $O(H^2(S^{[n]},\Integers))$, to belong to $\Mon^2$.
%These obstructions are imposed by the integral subrepresentations
%$C_i$, introduced in the proof of
%Lemma \ref{lemma-bounding-K}. We proved, that $C_i$ is isogenous to the
%Mukai lattice, if $i$ is in the range $2 \leq i \leq n/2$ 
%(Lemma \ref{lemma-10-in-markman-diagonal}).

%*****************************************************************
% Comparison with Verbitsky's representation
%*****************************************************************
\subsection{Comparison with Verbitsky's representation}
\label{sec-comparison-with-verbitsky}
We study first the Zariski closures of the two actions of
the group $\Gamma_v$ in $GL(H^*(S^{[n]},\ComplexNumbers))$ (Lemma 
\ref{lemma-extension-of-gamma-representation-to-all-isometries-of-H-2}). 
We then compare the monodromy action to Verbitsky's related representation
(Lemma \ref{lemma-comparing-two-representations}).

The group $O(n,m;\RealNumbers)$, of a non-definite and non-degenerate
bilinear form, has four connected components in the classical topology.
Its Zariski closure $O(n+m,\ComplexNumbers)$, in the standard representation, 
has two connected components. 
%The Zariski closure of the image $\gamma(\Gamma_v)$ has two connected 
%components, while the Zariski closure of the monodromy group has, at least, 
%{\em four} connected components, if $n\geq 3$. 

%\begin{new-lemma} 
%\label{lemma-connected-components-of-zariski-closure}
%The Zariski closure of $\monrep(\Gamma_v)$ 
%in $GL(H^*(S^{[n]},\ComplexNumbers))$ is isomorphic to
%$O(H^2(S^{[n]},\ComplexNumbers))\times \Integers/2\Integers$,
%if $n\geq 3$. It has four connected components.
%\end{new-lemma}
%
%The lemma is proven in section 
%\ref{sec-monodromy-invariance-of-the-normalized-chern-character}. 

\begin{new-lemma}
\label{lemma-extension-of-gamma-representation-to-all-isometries-of-H-2}
\begin{enumerate}
\item
\label{lemma-item-Zariski-closure-of-image-of-gamma}
The Zariski closure of the image $\gamma(\Gamma_v)$ in
$GL(H^*(S^{[n]},\ComplexNumbers))$ is isomorphic to the group
$O(H^2(S^{[n]},\ComplexNumbers))$ of isometries of 
$H^2(S^{[n]},\Integers)$ with respect to Beauville's pairing.
\item
\label{lemma-item-Zariski-closure-of-image-of-mon}
The Zariski closure of $\monrep(\Gamma_v)$ 
in $GL(H^*(S^{[n]},\ComplexNumbers))$ is isomorphic to
$O(H^2(S^{[n]},\ComplexNumbers))\times \Integers/2\Integers$,
if $n\geq 3$. It has four connected components.
\item
\label{lemma-item-extension-of-gamma-and-mon}
Both the $\gamma$ and the $\monrep$ actions, of $\Gamma_v$ 
on $H^*(S^{[n]},\Integers)_{\rm free}$, extend naturally to an action 
of the group 
$O(H^2(S^{[n]},\Integers))$. 
$O(H^2(S^{[n]},\Integers))$ acts on $H^*(S^{[n]},\RationalNumbers)$ 
via ring automorphisms. The extended  
$\gamma$-action leaves invariant the class $u_v$, given in
(\ref{eq-invariant-normalized-class-of-chern-character-of-universal-sheaf}). 
\end{enumerate}
\end{new-lemma}

\noindent
{\bf Proof:} 
\ref{lemma-item-extension-of-gamma-and-mon})
The character $cov$ extends to the orientation character 
(\ref{eq-top-homology-character}) of $O(H^2(S^{[n]},\Integers))$. 
Hence, an extension of $\gamma$ yields an extension of 
$\monrep$. We may restrict attention to the representation $\gamma$. 
The group $\Gamma_v$ is a normal subgroup, of finite index, 
in the group $O(H^2(S^{[n]},\Integers))$ 
(see Lemma \ref{lemma-the-index-of-Gamma-v}). 
Let ${\cal R}$ be the weighted polynomial ring, generated by 
the subspaces of generators $B_i$, introduced in 
Lemma \ref{lemma-10-in-markman-diagonal}.
%introduced in Corollary 
%\ref{cor-generating-subspaces-are-sub-representations}. 
The invariance with respect to the representation $\gamma$, 
of the subspaces $B_i$, 
gives rise to a $\Gamma_v$-action on the ring ${\cal R}$. 
Moreover, the surjective homomorphism 
$q:{\cal R}\rightarrow H^*(S^{[n]},\RationalNumbers)$ 
is $\Gamma_v$-equivariant. 
Hence, the relations ideal in ${\cal R}$ is $\gamma(\Gamma_v)$-invariant. 
The $B_i$ are $O(H^2(S^{[n]},\Integers))$
representations, via the isometry $v^\perp \cong H^2(S^{[n]},\Integers)$. 
The $\Gamma_v$ action on ${\cal R}$ extends to an algebraic action of
$O(H^2(S^{[n]},\ComplexNumbers))$ on the complexification of ${\cal R}$. 
$\Gamma_v$ is Zariski dense in $O(H^2(S^{[n]},\ComplexNumbers))$. 
Hence, the relations ideal in ${\cal R}$ 
is $O(H^2(S^{[n]},\Integers))$-invariant. 
The extension of the $\gamma$-action to an $O(H^2(S^{[n]},\Integers))$-action
follows. 

The invariance, of the class 
(\ref{eq-invariant-normalized-class-of-chern-character-of-universal-sheaf}),
with respect to the Zariski dense subgroup $\Gamma_v$, implies
the invariance with respect to 
$O(H^2(S^{[n]},\ComplexNumbers))$ and thus with respect to
$O(H^2(S^{[n]},\Integers))$ as well.

\ref{lemma-item-Zariski-closure-of-image-of-gamma}) 
The subspaces $B_i$, generating ${\cal R}$, are defined in terms of 
the $\Gamma_v$-invariant class $u_v$. Hence
the Zariski closure, of the image $\gamma(\Gamma_v)$ of $\Gamma_v$ in 
$GL({\cal R}\otimes_\RationalNumbers\ComplexNumbers)$, 
is isomorphic to $O(H^2(S^{[n]},\ComplexNumbers))$.
The latter Zariski closure maps {\em injectively} into 
$GL(H^*(S^{[n]},\ComplexNumbers))$, via the homomorphism $q$.
The Zariski closure in 
$GL({\cal R}\otimes_\RationalNumbers\ComplexNumbers)$ must map to
the Zariski closure in $GL(H^*(S^{[n]},\ComplexNumbers))$, since
the $\gamma$-action of $\Gamma_v$ on $H^*(S^{[n]},\ComplexNumbers)$ 
factors through $q$. 

\ref{lemma-item-Zariski-closure-of-image-of-mon})
When $n\geq 3$, the $\Gamma_v$-representation $v^\perp\otimes cov$ 
appears as a subrepresentation of
$H^4(S^{[n]},\RationalNumbers)$, by Lemma
\ref{lemma-Mon-invariant-irreducible-generating-subspaces}. 
Hence, the direct sum
$v^\perp\oplus [v^\perp\otimes cov]$ appears as a subrepresentation of
$H^2(S^{[n]},\RationalNumbers)\oplus H^4(S^{[n]},\RationalNumbers)$.
The Zariski closure, of the image of $\Gamma_v$ in the 
complexification of 
$v^\perp\oplus [v^\perp\otimes cov]$, is isomorphic to 
$O(H^2(S^{[n]},\ComplexNumbers))\times \Integers/2\Integers$. 
%When $n\geq 6$, the additional connected components are detected also via 
%the appearance of the orientation character 
%(\ref{eq-top-homology-character}) 
%as a subrepresentation of 
%$H^6(S^{[n]},\RationalNumbers)$ 
%(see Lemma
%\ref{lemma-normalized-chern-character-is-not-monodromy-invariant}).
\EndProof

\begin{rem}
\label{rem-other-extensions}
{\rm
Though natural, the extensions in part
\ref{lemma-item-extension-of-gamma-and-mon} 
of Lemma 
\ref{lemma-extension-of-gamma-representation-to-all-isometries-of-H-2}
are not unique, if $(n-1)$ is not a prime power. 
The quotient $O(H^2(S^{[n]},\Integers))/\Gamma_v$ is a finite 
abelian group of exponent $2$ (section \ref{sec-structure-of-stabilizer}), 
and other extensions of the action are
parametrized by characters of the quotient. 
%It follows, when $n-1$ is a prime power, 
%that the group of isometries of $H^2(S^{[n]},\Integers)$ 
%lifts to a subgroup of the automorphisms of the cohomology ring of
%$S^{[n]}$ preserving Beauville's pairing
%(see Lemma \ref{lemma-the-index-of-Gamma-v}). 
}
\end{rem}

We may regard the monodromy representation $H^*(S^{[n]},\Integers)$ 
as a representation of the subgroup $\W$, given in (\ref{eq-W}), 
since the representation $\monrep$ of $\Gamma_v$ factors through 
the injective homomorphism
$
\mu  :  \W  \longrightarrow \Mon
$
given in (\ref{eq-homomorphism-mu-from-W-to-Mon}).
Assume $n\geq 3$. Then the isomorphism 
$f:\Gamma_v\rightarrow \W$ conjugates the orientation character, of
$\Gamma_v$, to the character 
\begin{equation}
\label{eq-character-phi}
\phi\ : \ \W  \ \ \ \rightarrow \ \ \ 
\W/[\W\cap\Gamma_v] \ \ \ \cong \ \ \ \Integers/2\Integers.
\end{equation} 
The unique extension of $g\in \W$ to $\Gamma$ sends
$v$ to $\pm v$, where the sign is determined by $\phi(g)$. 
Note, that the image of $\W$ under 
the homomorphism (\ref{eq-projection-to-finite-isometry-group})
has order $2$ and thus, the character $\phi$ is the restriction of 
(\ref{eq-projection-to-finite-isometry-group}). 
Hence, the character $\phi$ is precisely a character of the type mentioned 
in remark \ref{rem-other-extensions}. 
Rewriting the formula 
(\ref{introduction-eq-monodromy-is-gamma-times-cov}) for the 
monodromy representation in terms of $\W$, we use the character $\phi$,
instead of $cov$.
%see that the 
%monodromy representation detects, 
%whether an element of $\W$ extends to an integral
%isometry of the Mukai lattice. 

Let $S\W\subset \W$ be the index $2$ subgroup 
$\W\cap SO(H^2(S^{[n]},\Integers))$. 
Let us compare the two representations 
of $S\W$ in $\Aut[H^*(S^{[n]},\RealNumbers)]$; the monodromy representation
$\mu$ and the restriction of Verbitsky's representation $\rho$,
given in (\ref{eq-verbitskys-representation}). 
Define $\eta : S\W \rightarrow \Aut[H^*(S^{[n]},\RealNumbers)]$ by
$\eta(f):=\rho(f)^{-1}\mu(f)$. 

\begin{new-lemma}
\label{lemma-comparing-two-representations}
Verbitsky's representation $\rho$ agrees with the monodromy 
representation $\mu$ on the intersection of $S\W$ with the kernel of $\phi$. 
Furthermore, 
the map $\eta$ is a homomorphism from $S\W$
into the center of the subgroup of $\Aut[H^*(S^{[n]},\RealNumbers)]$,
generated by $\Mon$ and the image of $\rho$. 
The image of $\eta$ is isomorphic to $\Integers/2\Integers$, if $n\geq 3$.
%a subgroup of $(\Integers/2\Integers)^{2n-4}$. 
%The homomorphism $\eta$ is non-trivial, if $n\geq 3$. 
\end{new-lemma}

\noindent
{\bf Proof:}
Given a monodromy operator $g\in \Mon$, denote by 
$\bar{g}$ its restriction to $H^2(S^{[n]})$. 
The $\Mon$-equivariance of $\rho$ 
means the equality $g\rho(f)g^{-1}=\rho(\bar{g}f\bar{g}^{-1})$, 
for $g\in \Mon$ and $f\in SO(H^2(S^{[n]},\RealNumbers))$. 
Extend $\rho$ to a representation of 
$O(H^2(S^{[n]},\RealNumbers))$, by letting $-id$ act as the Duality 
operator $D_{S^{[n]}}$ on $H^*(S^{[n]})$ 
(see (\ref{introduction-eq-monodromy-is-gamma-times-cov})). 
The extended $\rho$ remains $\Mon$-equivariant.
Then $\rho(\bar{g})^{-1}g$ commutes with $\rho(f)$, for all 
$g\in\Mon$ and $f\in SO(H^2(S^{[n]},\RealNumbers))$.
The map $\eta(g):=\rho(\bar{g})^{-1}g$ is thus
a representation of $\Mon$. 
%The proof of lemma 
%\ref{lemma-Mon-invariant-irreducible-generating-subspaces} implies, that 

$C'_{2i}(\RealNumbers)$ and $C''_{2i}(\RealNumbers)$ are 
$\Mon$-invariant and $\rho$-invariant, by Lemma 
\ref{lemma-Mon-invariant-irreducible-generating-subspaces}. 
Thus, they are also $\eta$-invariant. 
Furthermore, for each $g\in \Mon$, the element 
$\eta(g)$ must act on $C'_{2i}(\RealNumbers)$ and 
$C''_{2i}(\RealNumbers)$ diagonally via a character of $\Mon$
(since $\eta(g)$ is an intertwining operator of a zero or 
irreducible subrepresentation of $\rho$).
%which we denote by
%$\eta''_i$.
%We conclude, that the group $\eta(\Mon)$ is contained in the
%center of the subgroup of 
%$\Aut(H^*(S^{[n]},\Integers))$, generated by 
%$\Mon$ and the image of $\rho$. 
%The restriction of $\eta$ to the kernel $K$ is independent of 
%Verbitsky's representation $\rho$ and is equal to the monodromy
%representation of $K$. 
%Hence, $K$ belongs to the center of $\Mon$. 
%Each of $C'_{2i}(\RealNumbers)$ and $C''_{2i}(\RealNumbers)$
%comes from a rational subrepresentation of $C_{2i}$,
%because the $\LieAlg{so}(3,20,\RealNumbers)$-representation is defined over 
%$\RationalNumbers$. 
%Thus, $\chi'_i$ and the restrictions to $K$, of the 
%character $\eta''_i$, have values in $\{1,-1\}$. 
%Furthermore, each
%$\eta(f)$, $f\in S\W$, 
%acts diagonally on the subrepresentations $C'_{2i}(\RealNumbers)$ 
%and $C''_{2i}(\RealNumbers)$ via a character of $S\W$. 

The character group ${\rm Char}(S\W)$ is generated by 
the restriction of $\phi$,
by Corollary \ref{cor-stabilizer-is-generated-by-reflections}.
${\rm Char}(S\W)$ is isomorphic to $\Integers/2\Integers$, if $n\geq 3$,
and ${\rm Char}(S\W)$ is trivial, if $n=2$. 
% We have the short exact sequence of abelian groups:
% \[ 0 \rightarrow S\W/[\W,\W] \rightarrow \W/[\W,\W] \rightarrow 
%      \Integers/2\Integers\rightarrow 0 \]
% Take $Hom(\bullet,\ComplexNumbers^\times)$ and use 
% the vanishing of $Ext^1(\Integers/2\Integers,\ComplexNumbers^\times)$.
%%Every character of $S\W$ is of finite order, because the subgroup of $S\W$,
%%generated by reflections, is of finite index. 
%%Thus, a real valued character has values in $\{1,-1\}$. 
Hence, the restriction of $\eta$ to the subgroup 
$\mu(\ker(\phi)\cap S\W)$ of $Mon$ is trivial. 

If $n\geq 3$, the non-triviality of $\eta$ follows, 
by comparing the Zariski closures of the images $\mu(S\W)$ and $\rho(S\W)$. 
%The orientation character (\ref{eq-top-homology-character}) of 
%$O(v^\perp)$ restricts to a non-trivial
%character of $S\W$. Hence, 
The Zarisky closure of $\mu(S\W)$ in
$GL(H^*(S^{[n]},\ComplexNumbers))$ has two connected components
(Lemma 
\ref{lemma-extension-of-gamma-representation-to-all-isometries-of-H-2}). 
On the other hand, 
the Zariski closure of the image of $\rho$ has one connected
component, because $\rho$ can be extended to a representation 
of $SO(H^2(S^{[n]},\ComplexNumbers))$. 
\EndProof

%****************************************************************
% Monodromy-invariance of the normalized Chern character
%****************************************************************
\subsection{Monodromy-equivariance of the Chern character}
\label{sec-monodromy-invariance-of-the-normalized-chern-character}

The normalized Chern character  $u_v$, given in 
(\ref{eq-invariant-normalized-class-of-chern-character-of-universal-sheaf}),
is invariant with respect to the diagonal 
action of the stabilizer $\Gamma_v$, provided we use on the moduli
factor the representation $\gamma$, given in 
(\ref{eq-homomorphism-gamma-from-stabilizer-to-ring-auto}) 
(Theorem \ref{thm-trancendental-reflections}). 
We show in Lemma 
\ref{lemma-normalized-chern-character-is-not-monodromy-invariant}, that
the class $u_v$ is {\em not} invariant, in general, 
if we replace $\gamma$ by the monodromy representation
(\ref{introduction-eq-monodromy-is-gamma-times-cov}). 
Let us first identify the diagonal action of an appropriate subgroup of 
$O(v^\perp)$. Theorem \ref{thm-irreducibility} yields the isomorphism 
$v^\perp\cong H^2(\M(v),\Integers)$. 
Consequently, we can regard the class $u_v$ as an element of
\begin{equation}
\label{eq-vhs-containing-chern-character}
[H^2(\M(v),\RationalNumbers)\oplus \RationalNumbers\cdot v] \ \otimes \
H^*(\M(v),\RationalNumbers). 
\end{equation}
%The isometry group $O(v^\perp)$ embeds in $O(3,-20)$. 
%We have the orientation character
%(\ref{eq-top-homology-character}).
%The isometry $-id_{v^\perp}$ is orientation reversing. 
%Let $\langle \Gamma_v, -id_{v^\perp}\rangle$ be the subgroup of 
%$O(v^\perp)$ generated by $\Gamma_v$ and $-id_{v^\perp}$. 
%We will consider two index $2$ subgroups of this group, $\Gamma_v$ and $\W$.
%The subgroup $\W$ is the kernel of the restriction 
%$
%cov \ : \ \langle \Gamma_v, -id_{v^\perp}\rangle \ \ \rightarrow \ 
%\Integers/2\Integers,
%$
%of the orientation character (\ref{eq-top-homology-character}) of 
%$O(v^\perp)$. The subgroup 
%$\W$ is precisely the image of the monodromy representation 
%$\monrep$, given in 
%(\ref{introduction-eq-monodromy-is-gamma-times-cov}). 
%Additional characterizations of $\W$ are given in Lemma 
%\ref{lemma-EO-is-maximal}. 
Let $\W:=\W(v^\perp)$ be the subgroup of $O(v^\perp)$ 
described in Lemma \ref{lemma-EO-is-maximal}. 
Consider the diagonal action of $\W$
on (\ref{eq-vhs-containing-chern-character}), regarding 
$\RationalNumbers\cdot v$ as the trivial character, and identifying 
$H^2(\M(v),\RationalNumbers)$ with $v^\perp\otimes_\Integers\RationalNumbers$ 
via the isomorphism in Theorem \ref{thm-irreducibility}.

The case of $S^{[2]}$ is exceptional. 
Invariance holds for the monodromy action 
as well. 

\begin{new-lemma}
\label{lemma-chern-character-is-invariant-in-case-n-equal-2}
When $v=(1,0,-1)$ is the Mukai vector of $S^{[2]}$, 
the normalized Chern character $u_v$ 
%(\ref{eq-invariant-normalized-class-of-chern-character-of-universal-sheaf})
is invariant under the diagonal action of $\W$.
\end{new-lemma}

%Let us first recall Verbitski's result about the subring
%$A$ generated by $H^2(\M(v),\Integers)$:
%
%\begin{prop}
%\label{prop-verbitzki}
%\cite{verbitski,bogomolov}
%Let $X$ be a K\"{a}hler symplectic manifold of dimension
%$2n$, and let $A$ be the subalgebra of $H^*(X,\RationalNumbers)$ 
%spanned by $H^2(X,\RationalNumbers)$. Then 
%$H^*(X,\RationalNumbers)=A\oplus A^\perp$, and $A$ is  the quotient of 
%$\Sym^*H^2(X,\RationalNumbers)$ by the ideal generated by the elements 
%$x^{n+1}$ for all isotropic $x\in H^2(X,\RationalNumbers)$.
%\end{prop}

\noindent
{\bf Proof:} 
%of Lemma \ref{lemma-chern-character-is-invariant-in-case-n-equal-2}:}
Theorems \ref{thm-trancendental-reflections} and 
\ref{introduction-thm-Gamma-v-acts-motivicly} imply, that $u_v$ is invariant
with respect to the monodromy representation of 
$(\Gamma_v)^{cov}$. It suffices to show, that
$(\Gamma_v)^{cov}$ surjects onto $\W$. This follows from the fact,
that the kernel of the homomorphism $f:\Gamma_v\rightarrow \W$,
given in (\ref{eq-homomorphism-from-Gamma-v-to-EO}), 
contains (and is generated by)
the orientation reversing isometry $-\sigma_v$.
The Mukai vector $v=(1,0,-1)$ of $S^{[2]}$ 
is a $+2$ vector. $\Gamma_v$ surjects onto $O(v^\perp)$. 
The reflection $-\sigma_v$ is in $\Gamma_v$
and it restricts as $-id_{v^\perp}$. Hence, $f(-\sigma_v)=1$. 
%
%We know, that 
%$\gamma_{-\sigma_v}$ is a ring automorphism, 
%which is equal to the composition 
%$\gamma_D\circ\gamma_{\tau_{v_0}}$, where $v_0=(1,0,1)$. 
%Hence, $\gamma_{-\sigma_v}$ acts by $-1$ on $H^2(S^{[2]},\Integers)$. 
%The ring $H^*(S^{[2]},\Integers)$ is generated by $H^2(S^{[2]},\Integers)$
%(Proposition \ref{prop-verbitzki}).
%Hence, $\gamma_{-\sigma_v}$ is equal to the Duality automorphism
%$D_{S^{[2]}}$. The $\Gamma_v$-invariance of $u_v$ 
%implies, in particular, its $-\sigma_v\otimes D_{S^{[2]}}$ invariance.
\EndProof

\medskip
We show next, that the class $u_v$ 
is not invariant, when $v=(1,0,1-n)$ is the Mukai vector of 
$S^{[n]}$ and $n\geq 3$. In that case, the homomorphism 
(\ref{introduction-eq-monodromy-is-gamma-times-cov}) 
induces an isomorphism between 
$\Gamma_v$ and $\W$. It conjugates the orientation character of $\Gamma_v$ 
to the character $\phi$ of $\W$ given in (\ref{eq-character-phi}).
%\[
%\phi \ : \ \W \ \ \longrightarrow \ \ \W/[\W\cap \Gamma_v] \ \cong \ 
%\{\pm 1\}.
%\]
The kernel $\W\cap \Gamma_v$ of $\phi$ is precisely $\Gamma_v^{cov}$. 
The character $\phi$ 
is trivial if $(v,v)=2$, but non-trivial if $(v,v)\geq 4$.

Set $u_v:=u'_v+u''_v$ to be the decomposition of $u_v$ with respect 
to the direct sum decomposition of the left factor of 
(\ref{eq-vhs-containing-chern-character}), where $u'_v$ is in
${\rm span}\{v\}\otimes H^*(S^{[n]},\RationalNumbers)$. 
Let $u'_v:=u'_v({\rm odd})+u'_v({\rm even})$ be the decomposition of
$u'_v$ according to the decomposition of
the right hand factor $H^*(S^{[n]},\RationalNumbers)$ 
into weights congruent to $2$ modulo $4$ (the odd case)
and weights congruent to $0$ modulo $4$ (the even case). 
Let $u''_v:=u''_v({\rm odd})+u''_v({\rm even})$ be the analogous 
decomposition of $u''_v$.
Set $u_v^+ := u'_v({\rm even})+u''_v({\rm odd})$ and 
$u_v^- := u'_v({\rm odd})+u''_v({\rm even})$. 
Theorem \ref{thm-trancendental-reflections} implies, that the class $u_v$ 
is invariant with respect to the diagonal action of $\W\cap \Gamma_v$. 
Furthermore, it implies the equivariance of 
$u_v^+$ and $u_v^-$ in the following Lemma.

%Consider the diagonal action of 
%$\langle \Gamma_v, -id_{v^\perp}\rangle$ 
%on (\ref{eq-vhs-containing-chern-character}), regarding 
%$\RationalNumbers\cdot v$ as the trivial character. 
%We let the element $-id_{v^\perp}$ act on the moduli factor via
%the natural $\monrep$-extension of Lemma 
%\ref{lemma-extension-of-gamma-representation-to-all-isometries-of-H-2}. 
%The automorphism $\gamma(-id_{v^\perp})$ acts as multiplication by $-1$
%on the image 
%$B_{2i}'$ of $v^\perp$ in $H^{2i}(\M(v),\RationalNumbers)$.
%On the image $B_{2i}''$ of ${\rm span}\{v\}$, the automorphism 
%$\gamma(-id_{v^\perp})$ acts as the identity. 
%Now, $-id_{v^\perp}$ is orientation reversing, 
%when extended to a rational isometry of the Mukai lattice,
%via the identity on $v$.
%Thus, $\monrep(-id_{v^\perp})$ acts on $B_{2i}'$ via
%$(-1)^{1+i}$ and on $B_{2i}''$ via $(-1)^i$.
%
%The homomorphism 
%(\ref{eq-projection-to-finite-isometry-group}) restricts to 
%another character $\phi$
%\begin{equation}
%\label{eq-extention-character}
%\phi \ : \
%\langle \Gamma_v, -id_{v^\perp}\rangle \ \rightarrow \ 
%\langle \Gamma_v, -id_{v^\perp}\rangle/\Gamma_v
%\hookrightarrow \ \{\pm 1\}.
%\end{equation}

\begin{new-lemma}
\label{lemma-normalized-chern-character-is-not-monodromy-invariant}
The class $u_v$, given in 
(\ref{eq-invariant-normalized-class-of-chern-character-of-universal-sheaf}),
decomposes as a sum of the $\W$ invariant class 
$u_v^+$
and the class 
$u_v^-$,
which is acted upon by the character $\phi$ of $\W$. 
The class 
$u_v^+$ 
does not vanish for $n\geq 1$. 
If $n\geq 3$, then the summand 
$u_v^-$
does not vanish as well.
\end{new-lemma}

\noindent
{\bf Proof:}
Only the non-vanishing of $u_v^-$ remains to be proven.
Let $u_v(2i)$ be the summand of degree $2i$,
with respect to the weight decomposition of the right hand factor of
(\ref{eq-vhs-containing-chern-character}). 
Lemma \ref{lemma-10-in-markman-diagonal} 
%Lemma 10 in \cite{markman-diagonal} 
implies, that $u''_v(2i)$ does not vanish for even weights in the range 
$0< 2i\leq n+1$. 
Similarly, $u'_v(2i)$ does not vanish for even weights  in the range 
$4\leq 2i\leq n$. 
It follows, that $u''_v({\rm even})$ does not vanish, for $n\geq 3$. 
Furthermore, $u'_v({\rm odd})$ does not vanish, for $n\geq 6$. 
%old proof
%for $n\geq 3$, that
%the projection to 
%$H^2(S^{[n]},\RationalNumbers)\otimes H^4(S^{[n]},\RationalNumbers)$, 
%of the $\Gamma_v^{cov}$-invariant class $u_v$
%must induce an injective homomorphism from $H^2(S^{[n]})^*$ into 
%$H^4(S^{[n]})$. Hence, $u''_v({\rm even})$ does not vanish. 
%Furthermore, if $n\geq 6$, then the class $u_v$ 
%embeds ${\rm span}\{v\}^*$ in $H^6(S^{[n]},\RationalNumbers)$ as
%the character $\phi$ of $\W$
%(Lemma 10 in \cite{markman-diagonal}). Hence, $u'_v({\rm odd})$ 
%does not vanish.
%
%old old proof
%G\"{o}ttsche's formula
%yields that the Betti number $b_2$ of $S^{[n]}$ is $23$, while 
%\[
%b_4(S^{[n]}) \ \ = \ \ \left\{
%\begin{array}{ccc}
%276 & \mbox{if} & n=2,
%\\
%299 & \mbox{if} & n=3,
%\\
%300 & \mbox{if} & n\geq 4.
%\end{array}
%\right.
%\]
%The multiplication homomorphism 
%$\Sym^2H^2(S^{[n]})\rightarrow H^4(S^{[n]})$ is injective
%(Proposition \ref{prop-verbitzki}). 
%When $n=2$, $b_4$ is precisely the dimension of the second symmetric product
%of $H^2(S^{[2]})$. When $n=3$ (resp. $n\geq 4$), 
%the group $H^4(S^{[n]},\RationalNumbers)$
%contains $\Sym^2H^2(S^{[n]})$ and a complementary $23$-dimensional 
%(resp. $24$-dimensional) subspace. 
%Corollary \ref{cor-kunneth-factors-generate} implies, that the latter is a 
%copy of the irreducible representation $v^\perp$ of $\Gamma_v^{cov}$,
%when $n=3$, and the sum $\RationalNumbers\cdot v \oplus v^\perp$,
%when $n\geq 4$. Consequently, 
\EndProof

\begin{example}
\label{example-two-components-of-automorphisms-of-grassmannian}
{\rm
The near-invariance, of the normalized Chern character $u_v$ of the universal
sheaf, has the following simple analogue. 
Let $V$ be a $2n$-dimensional 
vector space, $n\geq 2$. The automorphism group $\Aut(G(n,V))$
of $G(n,V)$ has two connected components \cite{chow}. 
The identity component is $PGL(V)$. 
The character $\Aut(G(n,V))\rightarrow \Aut(G(n,V))/PGL(V)$ 
plays a role analogous to that of the character $cov$ 
in the monodromy representation 
(\ref{introduction-eq-monodromy-is-gamma-times-cov}). 
The Chern character $ch(\tau)$, of the tautological subbundle, 
is $PGL(V)$-invariant but not $\Aut(G(n,V))$-invariant. 
Choose any isomorphism $f:V\rightarrow V^*$
and denote by $f$ also the induced isomorphism $G(n,V)\rightarrow G(n,V^*)$,
sending a subspace $W$ to $f(W)$. Let $\iota: G(n,V^*)\rightarrow G(n,V)$
be the natural isomorphism. The composition $\iota\circ f$ is an automorphism
of $G(n,V)$, given by $W\mapsto f(W)^\perp$. Clearly, $\iota^*\tau$ 
is isomorphic to the dual $q^*_{G(n,V^*)}$ of the universal quotient bundle 
$q_{G(n,V^*)}$. Hence, the pullback 
$(\iota\circ f)^*c_2(\tau)=f^*(c_2(q^*_{G(n,V^*)}))=c_2(q^*_{G(n,V)})
=c_2(q_{G(n,V)})=c_1(\tau)^2-c_2(\tau)$
is different from $c_2(\tau)$. 

The proof of Theorem \ref{thm-reflection-sigma-satisfies-main-conj} 
relates the above example to the monodromy representation 
(\ref{introduction-eq-monodromy-is-gamma-times-cov}). A certain 
grassmannian bundle embeds in $S^{[n]}$ and plays an important role in the 
construction of the monodromy operator in 
Theorem \ref{thm-reflection-sigma-satisfies-main-conj}.
}
\end{example}

\section{Equivalences of Derived Categories}
\label{seq-equivalences-of-derived-categories}

In section \ref{seq-Fourier-Mukai-transformations} 
we recall the language of Fourier-Mukai functors. 
In section \ref{sec-Chow-theoretic-formula-for-gamma} 
we use this language
to provide a conceptual and motivic formula 
for the monodromy operator $\monrep(g)$, given in 
(\ref{eq-image-of-g-by-H-mon}), associated to 
a Hodge-isometry $g$ of Mukai lattices
(Lemma \ref{lemma-Chow-theoretic-formula-for-gamma}). 
The topological formula (\ref{eq-gamma-g-E1-E2}) for the class 
$\gamma(g,\E_v,\E_{g(v)})$ was 
motivated by the ideal case discussed in section
\ref{sec-conceptual-interpretation-of-the-formula-gamma-g}. 
This is the case, where an auto-equivalence induces an isomorphism of 
moduli spaces. 
%Section \ref{sec-the-induced-representation} 
%discusses a representation, of the group of auto-equivalences of the 
%derived category of a K3 surface $S$. This representation is a candidate, 
%for the action of the group of auto-equivalences 
%on the cohomology of all moduli spaces of sheaves on $S$ 
%(with a primitive Mukai vector). 

%*********************************************************************
% Fourier-Mukai transformations
%*********************************************************************
\subsection{Fourier-Mukai transformations}
\label{seq-Fourier-Mukai-transformations} 

Given a projective variety $X$, 
denote by $D(X):= D^b_{Coh}(X)$ the bounded derived category of coherent
sheaves on $X$. 
Given projective varieties $X_1$ and $X_2$, together with an object $F$ 
in $D(X_1\times X_2)$, we get the functor
\[
\Phi^F_{X_1\rightarrow X_2} \ : \ D(X_1) \ \ \longrightarrow \ \ D(X_2),
\] 
defined by
\begin{equation}
\label{eq-fourier-mukai-functor}
\Phi^F_{X_1\rightarrow X_2}(x) \ := \ R_{\pi_{2,*}}
\left(
F\stackrel{L}{\otimes}\pi_1^*(x)\right),
\end{equation}
where $\pi_i:X_1\times X_2 \rightarrow X_i$ is the projection, for $i=1,2$. 
If the functor is an equivalence, it is called a 
{\em Fourier-Mukai transformation} and 
the object $F$ is called its {\em kernel}. 
A theorem of Orlov implies, that every equivalence, 
of the derived categories of two projective 
varieties, is a Fourier-Mukai transformation. Moreover, the kernel 
is unique up to isomorphism in the derived category of their product 
\cite{orlov}. 

Let $S_1$ and $S_2$ be K3 surfaces. 
An equivalence $\Phi: D(S_1)\rightarrow D(S_2)$ 
determines a class
\[
\Z_\Phi \ \ := \ \ \pi_1^*\sqrt{td_{S_1}}\cdot  ch(F) \cdot 
\pi_2^*\sqrt{td_{S_2}}
\]
in $H^*(S_1\times S_2,\Integers)$, 
where $F$ is the isomorphism class of the kernel associated to 
$\Phi$ by Orlov's Theorem. 
%Orlov showed, that $S_2$ is isomorphic to a moduli space of stable sheaves on 
%the K3 surface $S_1$, and 
%the object $F$ in $S_1\times S_2$ can be chosen to be a universal sheaf.
The integrality of the class $\Z_\Phi$ 
was proven in \cite{mirror-symmetry-k3}, 
using ideas of Mukai \cite{mukai-hodge}. 

For a K3 surface $S$, 
let $\Aut D(S)$ be the group of auto-equivalences of its derived category. 
Denote by $G$ the group of Hodge isometries of the Mukai lattice of $S$.
Results of Mukai and Orlov yield the following Theorem.

\begin{thm}
\label{thm-HLOY}
(\cite{mirror-symmetry-k3} Theorem 1.6)
\begin{enumerate}
\item
The class $\Z_\Phi$ induces a Hodge isometry of Mukai lattices, denoted by 
$ch(\Phi) : H^*(S_1)\rightarrow H^*(S_2)$. 
\item
The resulting map 
\begin{equation}
\label{eq-ch}
ch \ : \ \Aut D(S) \ \ \longrightarrow \ \ G
\end{equation}
is a group homomorphism. Its image has index at most $2$. 
Moreover, $G$ is generated by the image of $ch$ and
the isometry $D$ of the Mukai lattice given in 
(\ref{eq-serre-duality-hodge-isometry}).
%\begin{eqnarray}
%\label{eq-serre-duality-hodge-isometry}
%D : H^*(S) & \longrightarrow & H^*(S) 
%\\
%\nonumber
%(r,c_1,s) & \mapsto & (r,-c_1,s).
%\end{eqnarray}
\end{enumerate}
\end{thm}

Let us prove part
\ref{conj-item-gamma-g-induces-isomorphism-of-cohomology-rings} of Theorem  
\ref{thm-trancendental-reflections} for two dimensional moduli spaces.
%under the assumption,
%that universal families $\E_v$ and $\E_{g(v)}$ exist.

\begin{new-lemma}
\label{lemma-conjecture-holds-for-isotropic-mukai-vectors}
Let $v$ be a Mukai vector of a K3 surface $S_1$ and 
$g:H^*(S_1,\Integers)\rightarrow H^*(S_2,\Integers)$ an isometry of 
Mukai lattices as in Theorem \ref{thm-trancendental-reflections}.
Assume, in addition, that $v$ is isotropic. 
Then part \ref{conj-item-gamma-g-induces-isomorphism-of-cohomology-rings} 
of Theorem  
\ref{thm-trancendental-reflections} holds.
\end{new-lemma}

\noindent
{\bf Proof:}
The proof below assumes that 
universal families $\E_v$ and $\E_{g(v)}$ exist.
The general proof is similar, using universal classes $e_v$
and $e_{g(v)}$ (Definition \ref{def-e-v}). 
The moduli spaces $\M(v)$ and $\M(g(v))$
are K3 surfaces, by the work of Mukai. 
We note first a simplification, in the formula
(\ref{eq-gamma-delta}) for $\gamma(g,\E_v,\E_{g(v)})$, 
when the dimension $m$ of the moduli spaces is $2$. 
Given a class $f$, in the K-group of 
$\M(v)\times\M(g(v))$, with $c_1(f)=0$, then the following equation
holds in $H^4(\M(v)\times\M(g(v)),\Integers)$.
\begin{equation}
\label{eq-c-2-equals-ch-2}
c_2(-f) \ \ = \ \ ch_2(f).
\end{equation}
This equality extends to an equality between the
topological formula (\ref{eq-gamma-delta}) and the one obtained from it by 
replacing 
$c_2(\{\ell(\bullet)\}^{-1})$ with $(\bullet)_2$, provided 
$(\bullet)_1=0$. 

We claim, that the operator 
$\gamma_{g,v}:H^*(\M(v),\RationalNumbers)\rightarrow
H^*(\M(g(v)),\RationalNumbers)$, associated to the 
class $\gamma(g,\E_v,\E_{g(v)})$, 
is the grade $0$ summand $h_0$, in the weight decomposition of 
the homomorphism 
\[
h \ := \
ch\left(\Phi^{\E_{g(v)}}\right)\circ (D\circ g \circ D) \circ 
ch\left(\Phi^{\E_v}\right)^{-1}.
\]
Let $\Delta$ be the diagonal in $S_1\times S_1$. 
The above claim follows from the definition (\ref{eq-gamma-g-E1-E2})
of $\gamma(g,\E_v,\E_{g(v)})$, the vanishing of $h_{-2}$ proven below,
equation (\ref{eq-c-2-equals-ch-2}), the identity 
\[
ch(\Phi^{\E_v}_{S_1\rightarrow\M(v)}) \ \ = \ \ 
\left[ch(\Phi^{\E_v^\vee}_{\M(v)\rightarrow S_1})\right]^{-1}
\]
(Theorem 4.9 in \cite{mukai-hodge}), 
and the identity
\[
[(1\otimes g)(ch(\E_v))]^\vee\ = \ 
(D_\M\otimes D)(1\otimes g)(ch(\E_v)) \ = \ 
%(1\otimes DgD)(D_\M\otimes D)(ch(\E_v)) \ = \ 
(1\otimes DgD)(ch(\E_v^\vee)).
\]
%\[
%[(1\otimes g)(\Delta)]^\vee \ = \
%(D\otimes D)(1\otimes g)(\Delta) \ = \ 
%(D\otimes Dg)(\Delta) \ = \
%(1\otimes DgD)(\Delta).
%\]
The class $\gamma(g,\E_v,\E_{g(v)})$ is integral, since $h$ is a composition 
of integral isometries.
The homomorphism $ch\left(\Phi^{\E_v}\right):H^*(S,\Integers)\rightarrow 
H^*(\M(v),\Integers)$ 
maps the Mukai vector $D(v)$ to the fundamental class $w_v$ 
in $H^4(\M(v),\Integers)$,
by Lemma 4.11 of \cite{mukai-hodge}. 
Consequently, $h$
is an integral isometry of Mukai lattices, which 
maps $w_v$ to the fundamental class $w_{g(v)}$. 
Considering orthogonal complements, we see that $h$ maps the direct sum 
$H^2(\M(v),\Integers)\oplus H^4(\M(v),\Integers)$
to $H^2(\M(g(v)),\Integers)\oplus H^4(\M(g(v)),\Integers)$.
In other words, the negative summands, in 
the weight decomposition of $h$, vanish and $h=h_0+h_2+h_4$. 
We have the following equality
\[
(h_0(1),w_{g(v)}) \ = \ (h(1),w_{g(v)}) \ = \ 
(h(1),h(w_v)) \ = \ (1,\omega_v) \ = \ 1. 
\]
Hence, $h_0(1)=1$. Moreover, $(h_0(x),h_0(y))=(h(x),h(y))=(x,y)$,
for all $x$ and $y$ in $H^2(\M(v))$. 
It follows, that $\gamma_{g,v}$ is a ring isomorphism.
%
%The isometries $ch\left(\Phi^{\E_v}\right)$
%and $ch\left(\Phi^{\E_{g(v)}}\right)$ are Hodge-isometries \cite{mukai-hodge}.
%Hence, so are $h$ and $\gamma_{g,v}$, 
%provided $g$ is a Hodge isometry.
\EndProof

%***********************************************************
% A Chow-theoretic formula for monodromy operators
%***********************************************************
\subsection{A Chow-theoretic formula for monodromy operators}
\label{sec-Chow-theoretic-formula-for-gamma}
Given an object $F$ of $D(S_1\times S_2)$ and a projective variety $M$, we 
construct the object 
\[
F_M \ := \ (\pi_{24}^* F)\stackrel{L}{\otimes} 
(\pi_{13}^*\StructureSheaf{\Delta}) 
\]
in $D(M\times S_1\times M \times S_2)$, 
where $\Delta$ is the diagonal of $M\times M$ and $\pi_{ij}$ is the
projection from $M \times S_1\times M \times S_2$ onto the product
of the $i$-th and $j$-th factors.
Given a Fourier-Mukai transformation $\Phi_{S_1\rightarrow S_2}^F$,
we get a natural functor
\[
\Phi^{F_M} \ : \ 
D(M\times S_1) \ \ \longrightarrow \ \ D(M\times S_2), 
\]
using formula (\ref{eq-fourier-mukai-functor}) with the kernel $F_M$. 
The functor $\Phi^{F_M}$ is an equivalence of derived categories as well
\cite{orlov-abelian-varieties}. 

\begin{defi}
\label{def-relative-fourier-mukai-transform}
{\rm
Let $\M(v)$ be a moduli space of stable sheaves on $S_1$ and 
$\E_v$ a universal sheaf over $\M(v)\times S_1$. 
We will refer to the class $\Phi^{F_{\M(v)}}(\E_v)$ as the
{\em relative Fourier-Mukai-transform}
of the universal sheaf. 
}
\end{defi}

The construction of the Fourier-Mukai functor $\Phi^{F_M}$
generalizes to define an equivalence of derived categories
of twisted sheaves
$\Phi^{F_M}: D^b_{coh}(M\times S_1,\pi_M^*\alpha) \rightarrow 
D^b_{coh}(M\times S_2,\pi_M^*\alpha)$, for any class $\alpha$ in the 
Brauer group of $M$ (see section \ref{sec-universal-sheaves}). 

\medskip
Let $(S_i,v_i,\LB_i)$ be two objects of the groupoid $\G$. 
Let $p_i:\PP_{v_i}\rightarrow \M_{\LB_i}(v_i)$ be the projective bundle 
defined in section \ref{sec-universal-sheaves} and $\widetilde{\E}_{v_i}$
the universal sheaf over $\PP_{v_i}$
(see equation (\ref{eq-chern-character-of-universal-sheaf-over-PP})). 
When a universal sheaf $\E_{v_i}$ exists over
$\M_{\LB_i}(v_i)$, we may assume that $p_i$ is the identity and 
$\widetilde{\E}_{v_i}=\E_{v_i}$.

%The Grothendieck-Riemann-Roch theorem yields the following equality:
\begin{new-lemma}
\label{lemma-Chow-theoretic-formula-for-gamma}
\begin{enumerate}
\item
\label{lemma-item-the-chern-character-of-the-relative-FM-transform}
The following equality holds in $H^*(\PP_{v_1}\times S_2)$.
\[
%\label{eq-the-chern-character-of-the-relative-FM-transform}
ch\left(
\Phi^{F_{\PP_{v_1}}}(\widetilde{\E}_{v_1})
\right)\sqrt{td_{S_2}} \ \ = \ \ 
\left(id_{\PP_{v_1}}\otimes ch(\Phi^F)\right)
(ch(\widetilde{\E}_{v_1})\sqrt{td_{S_1}}). 
\]
\item
\label{lemma-item-Chow-theoretic-formula-for-gamma}
%Set $w:=ch(\Phi^F)(v)$. Assume, that $\LB_1$ is $v$-suitable and
%$\LB_2$ is $w$-suitable. 
Assume that $v_2=ch(\Phi^F)(v_1)$. 
The following equality holds in 
the cohomology ring of $\PP_{v_1}\times \PP_{v_2}$. 
\begin{equation}
\label{eq-Chow-theoretic-formula-for-gamma}
(p_1\times p_2)^*\gamma(ch(\Phi^F),v_1) \ \ = \ \ 
c_m\left[- \ 
\pi_{13_!}
\left(
\pi_{12}^*(\Phi^{F_{\PP_{v_1}}}[\widetilde{\E}_{v_1}])^\vee
\stackrel{L}{\otimes}
\pi_{23}^*(\widetilde{\E}_{v_2})
\right)
\right],
\end{equation}
where $m$ is the dimension of $\M(v_i)$ and 
$\pi_{ij}$ is the projection
from $\PP_{v_1} \times S_2 \times \PP_{v_2}$.
\item
\label{lemma-item-Chow-theoretic-formula-for-mon}
If $v_2=[ch(\Phi^F)(v_1)]^\vee$ and 
the isometry $ch(\Phi^F)$ is orientation-preserving, 
then
\[
(p_1\times p_2)^*\monrep(D\circ ch(\Phi^F),v_1) \ \ = \ \ 
c_m\left[- \ 
\pi_{13_!}
\left(
\pi_{12}^*(\Phi^{F_{\PP_{v_1}}}[\widetilde{\E}_{v_1}])
\stackrel{L}{\otimes}
\pi_{23}^*(\widetilde{\E}_{v_2})
\right)
\right],
\]
where $\monrep(g,v_1)$ is given in
(\ref{eq-class-mon-g-v}).
%\item
%\label{lemma-item-classes-are-adjoint}
%Set $g:=ch(\Phi^F)$. 
%The homomorphisms
%\begin{eqnarray*}
%\gamma_{g,v_1} & : &
%H^*(\M(v_1)) \ \ \longrightarrow \ \ H^*(\M(v_2)),
%\\
%\gamma_{g^{-1},v_2} & : &
%H^*(\M(v_2)) \ \ \longrightarrow \ \ H^*(\M(v_1))
%\end{eqnarray*}
%are adjoint, with respect to the Poincare-Duality pairings on
%$H^*(\M(v_i))$, $i=1,2$. 
\end{enumerate}
\end{new-lemma}

\begin{rem}
\label{rem-auto-equivalences-of-the-derived-category-of-moduli-spaces}
{\rm
It is tempting to speculate, that the object
\[
R\pi_{13_!}
\left(
\pi_{12}^*(\Phi^{F_{\M(v_1)}}[\E_{v_1}])^\vee\stackrel{L}{\otimes}
\pi_{23}^*(\E_{v_2})
\right), 
\]
on the right hand side of equation 
(\ref{eq-Chow-theoretic-formula-for-gamma}),
may play a role in the study of equivalences of 
the derived category of moduli spaces of sheaves. 
One might even relax the assumption, that $v_2=ch(\Phi^F)(v_1)$. 
We take in (\ref{eq-Chow-theoretic-formula-for-gamma})
the middle dimensional Chern class of the object, as we are studying
graded automorphisms of the cohomology ring. 
See lemma \ref{lemma-conjecture-holds-for-isotropic-mukai-vectors}, 
for two dimensional moduli spaces, and 
remark \ref{rem-other-base-calabi-yau-varieties}, 
for the elliptic curve case. 
It should be noted, however, that some adjustments to the above object 
would be necessary, when $m>2$. When $\Phi$ is the identity, $v=v_i$,
and $\E_v=\E_{v_i}$, $i=1,2$, 
the above object fits in a distinguished triangle, 
involving also the relative extension sheaves 
$\RelExt^i_{\pi_{13}}(\pi_{12}^*\E_v,\pi_{23}^*\E_v)$, $i=1,2$. 
The sheaf 
$\RelExt^2_{\pi_{13}}(\pi_{12}^*\E_v,\pi_{23}^*\E_v)$ 
is the structure sheaf of the diagonal in $\M(v)\times \M(v)$, while
the sheaf $\RelExt^1_{\pi_{13}}(\pi_{12}^*\E_v,\pi_{23}^*\E_v)$ 
is either $0$, when $\M(v)$ is two dimensional, or 
a torsion free sheaf of rank $m-2$, when $m>2$ 
\cite{mukai-hodge,markman-diagonal}. One would need to
eliminate the contribution of 
$\RelExt^1_{\pi_{13}}(\pi_{12}^*\E_v,\pi_{23}^*\E_v)$. 
Theorem \ref{thm-graph-of-diagonal-in-terms-of-universal-sheaves}
is equivalent to the statement, that 
$c_m(\RelExt^1_{\pi_{13}}(\pi_{12}^*\E_v,\pi_{23}^*\E_v))$
is $1-(m-1)!$ times the class of the diagonal \cite{markman-diagonal}.
}
\end{rem}

\noindent
{\bf Proof of lemma \ref{lemma-Chow-theoretic-formula-for-gamma}:}
\ref{lemma-item-the-chern-character-of-the-relative-FM-transform})
Let $\iota:\Delta\hookrightarrow \PP_{v_1}\times \PP_{v_1}$, 
be the closed immersion of the diagonal and
$[\Delta]$ the cohomology class Poincare dual to the diagonal.
Let $\pi_{ij}$ be the projection from 
$\PP_{v_1}\times S_1\times \PP_{v_2} \times S_2$ onto the product of the 
$i$-th and $j$-th factors.
\begin{eqnarray*}
ch\left(
\Phi^{F_{\PP_{v_1}}}(\widetilde{\E}_{v_1})
\right)\cdot \sqrt{td_{S_2}}
& = &
ch\left(
\pi_{34_!}
\left\{
\pi_{12}^*\widetilde{\E}_{v_1}\otimes \pi_{24}^* F \otimes 
\pi_{13}^*\StructureSheaf{\Delta}
\right\}
\right)\cdot \sqrt{td_{S_2}}
\\
& = &
\pi_{34_*}\left(
ch\left\{
\pi_{12}^*\widetilde{\E}_{v_1}\otimes \pi_{24}^* F \otimes 
\pi_{13}^*\StructureSheaf{\Delta}
\right\}
\pi_2^*td_{S_1}\cdot
\pi_1^*td_{\PP_{v_1}}
\right) \cdot \sqrt{td_{S_2}}
\\
& = &
\pi_{34_*}\left(
\pi_{12}^*[ch(\widetilde{\E}_{v_1})\sqrt{td_{S_1}}]\cdot 
\pi_2^*\sqrt{td_{S_1}}\cdot \pi_{24}^*ch(F) \cdot 
\pi_{13}^*[\Delta]
\right)\cdot \sqrt{td_{S_2}}
\\
& = &
\left(ch(\Phi^F)\otimes id_{\PP_{v_1}}\right)
(ch(\widetilde{\E}_{v_1})\sqrt{td_{S_1}}). 
\end{eqnarray*}
The first equality is the definition of $\Phi^{F_{\PP_{v_1}}}$. The second 
follows from Grothendieck-Riemann-Roch. 
The third follows via Grothendieck-Riemann-Roch again, 
applied to compute the Chern character of the pushforward of 
$\StructureSheaf{\Delta}$ via the closed immersion $\iota$,
\[
ch(\iota_*\StructureSheaf{\Delta}) \ = \ 
\iota_*(td(N_{\Delta/[\PP_{v_1}\times \PP_{v_1}]})^{-1}) \ = \ 
[\Delta]\cup \pi_1^*td_{\PP_{v_1}}^{-1}
\]
(\cite{fulton} Example 15.2.15).
The fourth equality follows from the projection formula and the
definition of $ch(\Phi^F)$. 

\ref{lemma-item-Chow-theoretic-formula-for-gamma}) 
Set $g:=ch(\Phi^F)$. 
Assume first that universal sheaves $\E_{v_i}$, $i=1,2$, exist.
Then $\gamma(g,v_1)=\gamma(g,\E_{v_1},\E_{v_2})$,
by Lemma \ref{lemma-two-definitions-of-gamma-g-v}, and 
part \ref{lemma-item-Chow-theoretic-formula-for-gamma} 
follows from part 
\ref{lemma-item-the-chern-character-of-the-relative-FM-transform} via
Grothendieck-Riemann-Roch.

We sketch the proof of the general case. 
Let $\gamma'(g,\widetilde{\E}_{v_1},\widetilde{\E}_{v_1})$
be the class obtained by taking instead the Chern class $c_{m-1}$
on the right hand side of (\ref{eq-Chow-theoretic-formula-for-gamma}).
Equations
(\ref{eq-chern-character-of-universal-sheaf-over-PP}) and 
(\ref{eq-invariant-normalized-class-of-chern-character-of-universal-sheaf}) 
%relating $u_{v_i}$ and $ch(\widetilde{\E}_{v_i})$,
imply that $ch(\widetilde{\E}_{v_i})\sqrt{td_{S_i}}=
(p_i\times id_{S_i})^*(u_{v_i})\cdot 
\exp(\ell_i)$, for some class $\ell_i\in H^2(\PP_{v_i},\RationalNumbers)$.
Set $\ell:=\pi_2^*(\ell_2)-\pi_1^*(\ell_1)
\in H^2(\PP_{v_1}\times \PP_{v_2},\RationalNumbers)$.
The difference between the two sides of equation 
(\ref{eq-Chow-theoretic-formula-for-gamma}) is the cup product 
$\ell\cup \gamma'(g,\widetilde{\E}_{v_1},\widetilde{\E}_{v_1})$,
by part \ref{lemma-item-the-chern-character-of-the-relative-FM-transform},
Grothendieck-Riemann-Roch, and
Lemma \ref{lemma-18-in-integral-generators}. 
Now $\gamma'(g,\widetilde{\E}_{v_1},\widetilde{\E}_{v_1})$
is equal to the pull-back of the
class $\gamma'(g,u_{v_1},u_{v_2})$, defined in 
(\ref{eq-gamma-prime}), by Grothendieck-Riemann-Roch and
Lemma \ref{lemma-18-in-integral-generators}.
The class $\gamma'(g,u_{v_1},u_{v_2})$
vanishes by Lemma \ref{lemma-recovering-f}.

The proof of part \ref{lemma-item-Chow-theoretic-formula-for-mon}
is similar.
\EndProof
\subsection{The ideal case of isomorphic moduli spaces}
\label{sec-conceptual-interpretation-of-the-formula-gamma-g}

The language of Fourier-Mukai transforms provides a 
conceptual interpretation of formula (\ref{eq-gamma-delta})
for the class $\gamma(g,v)$. When $g$ is a Hodge isometry, 
Theorem \ref{thm-trancendental-reflections} says, roughly, the
following. 
{\em 
The homomorphism 
$\gamma_{g,v}: H^*(\M(v))\rightarrow H^*(\M(g(v)))$ 
behaves, as if it is pushforward by an isomorphism 
$\M(v)\rightarrow \M(g(v))$ of moduli spaces, 
induced by a Fourier-Mukai transformation $\Phi$, with $ch(\Phi)=g$. 
} 
%the relative Fourier-Mukai transform
%of a universal sheaf (Definition \ref{def-relative-fourier-mukai-transform})
%behaves, 
%as if it is the pullback of a universal sheaf on another moduli space,
%by an isomorphism of moduli spaces
We explain the above statement in this section.

We define first the sub-groupoid
\begin{equation}
\label{eq-G-Mukai}
\G_{Mukai}
\end{equation}
of the groupoid $\G$ given in (\ref{eq-cohomological-groupoid}). 
The objects of $\G_{Mukai}$ are the same as those of  $\G$.
Let  $x_i:=(S_i,v_i,H_i)$, $i=1,2$, be objects of $\G$. 
The morphisms in $\Hom_{\G_{Mukai}}(x_1,x_2)$ come in two flavors and 
depend on a Fourier-Mukai transformation 
$\Phi_{S_1\rightarrow S_2}^F$ 
with kernel $F$, satisfying the two conditions below.

\noindent
(1) The isometry $g:=ch(\Phi^F)$ is orientation-preserving 
(Definition \ref{def-covariant-subgroups}), 
with respect to the orientation of the
Mukai lattices in Remark \ref{rem-natural-orientation},
where $ch$ is the homomorphism (\ref{eq-ch}). 

\noindent
(2) One of the following holds for some integer $i$:
\begin{enumerate}
\item[(2a)]
$g(v_1)=v_2$ and 
$\Phi^F(E_1)$ is equivalent in $D(S_2)$ to the shift
$E_2[i]$ of some $H_2$-stable sheaf $E_2$, for every 
$H_1$-stable sheaf $E_1$ with Mukai vector $v_1$.
%the relative Fourier-Mukai transform
%$\Phi^{F_{\M(v_1)}}(\E_{v_1})$ is represented by a
%$\pi_{\M(v_1)}^*\alpha$-twisted sheaf $\E'$ on $\M_{H_1}(v_1)\times S_2$,
%which is a family of $H_2$-stable sheaves on $S_2$. 
\item[(2b)]
$g(v_1)=(v_2)^\vee$ and the dual 
$[\Phi^F(E_1)]^\vee$ is equivalent in $D(S_2)$ to the shift
$E_2[i]$ of some $H_2$-stable sheaf $E_2$, for every 
$H_1$-stable sheaf $E_1$ with Mukai vector $v_1$.
%$[\Phi^{F_{\M(v_1)}}(\E_{v_1})]^\vee$, of the 
%relative Fourier-Mukai transform,
%is represented by a
%$\pi_{\M(v_1)}^*\alpha$-twisted sheaf $\E'$ on $\M_{H_1}(v_1)\times S_2$,
%which is a family of $H_2$-stable sheaves on $S_2$. 
\end{enumerate}

\noindent
If conditions (1) and (2a) hold, then $g$ belongs to
$\Hom_{\G_{Mukai}}(x_1,x_2)$. If conditions (1) and (2b) hold, 
then $D\circ g$ belongs to $\Hom_{\G_{Mukai}}(x_1,x_2)$,
where $D$ is the duality involution
(\ref{eq-serre-duality-hodge-isometry}).

\begin{new-lemma}
\label{lemma-G-Mukai-is-in-G-OK}
Let $x_i:=(S_i,v_i,H_i)$, $i=1,2$, be objects of $\G$. 
If $f$ is a morphism in $\Hom_{\G_{Mukai}}(x_1,x_2)$, then
Theorem \ref{thm-trancendental-reflections} holds for $f$. 
More specifically,
$\gamma_{f,v_1}$ is an isomorphism of cohomology rings and
$(\gamma_{f,v_1}\otimes f)(u_{v_1}) = (u_{v_2})$.
\end{new-lemma}

\noindent
{\bf Proof:} 
Let $\alpha$ be a class in the Brauer group of $\M_{H_1}(v_1)$ and 
$\E_{v_1}$ a $\pi_{\M(v_1)}^*\alpha$-twisted sheaf over
$\M_{H_1}(v_1)\times S_1$. 
When condition (2a) holds, then 
the relative Fourier-Mukai transform
$\Phi^{F_{\M(v_1)}}(\E_{v_1})$ is represented by a
$\pi_{\M(v_1)}^*\alpha$-twisted sheaf $\E'$ on $\M_{H_1}(v_1)\times S_2$,
which is a family of $H_2$-stable sheaves on $S_2$. This follows 
from a flatness result for Fourier-Mukai functors (\cite{mukai-applications} 
Theorem 1.6). 
When condition (2b) holds, then 
$[\Phi^{F_{\M(v_1)}}(\E_{v_1})]^\vee$ 
is represented by such a twisted sheaf $\E'$.
The classifying morphism
\begin{equation}
\label{eq-kappa}
\kappa \ : \ \M_{H_1}(v_1) \ \ \longrightarrow \ \ \M_{H_2}(v_2),
\end{equation}
associated to the twisted sheaf $\E'$, is an isomorphism.
This is seen as follows. 
$\kappa$ pulls back the holomorphic symplectic $2$-from of 
$\M_{H_2}(v_2)$ to that of $\M_{H_1}(v_1)$, by the construction of 
this $2$-from in \cite{mukai-symplectic-structure}. Thus 
$\kappa$ is \'{e}tale. It is surjective, since both moduli spaces
are compact and irreducible (Theorem \ref{thm-irreducibility}). 
It is injective, since $\M_{H_2}(v_2)$ is simply-connected 
(Theorem \ref{thm-irreducibility}).

Assume first condition (2a) with $g=f$. 
The pushforward 
$\E_{v_2}:=(\kappa\times id_{S_2})_*(\E')$ is a
$\pi_{\M(v_2)}^*(k_*(\alpha))$-twisted universal sheaf
over $\M(v_2)\times S_2$. 
Let $\PP_{v_1}$ be the projective bundle over $\M(v_1)$, as in 
Lemma \ref{lemma-Chow-theoretic-formula-for-gamma}, and 
$\tilde{\kappa}:\PP_{v_1}\rightarrow k_*\PP_{v_1}$ 
the natural isomorphism. Then $(\kappa\times id_{S_2})_*\E'$ lifts
to a universal sheaf $\widetilde{\E}_{v_2}$ over $k_*\PP_{v_1}$.
Lemma \ref{lemma-Chow-theoretic-formula-for-gamma} part
\ref{lemma-item-the-chern-character-of-the-relative-FM-transform}
yields the equality 
$(id\otimes g)(ch(\widetilde{\E}_{v_1})\sqrt{td_{S_1}}) = 
(\tilde{\kappa}\otimes id_{S_2})^*(ch(\widetilde{\E}_{v_2}))\sqrt{td_{S_2}}$
in $H^*(\PP_{v_1}\times S_2)$. 
The equality 
$(id_{\M(v_1)}\otimes g)(u_{v_1})=(\kappa\times id_{S_2})^*(u_{v_2})$
follows from the latter, equations 
(\ref{eq-chern-character-of-universal-sheaf-over-PP}) and 
(\ref{eq-invariant-normalized-class-of-chern-character-of-universal-sheaf}), 
and the injectivity of 
$p_1^*:H^*(\M(v_1))\rightarrow H^*(\PP_{v_1})$. 
Hence, $\kappa_*$ is a ring isomorphism satisfying
$(\kappa_*\otimes g)(u_{v_1})=u_{v_2}$. 
The equality $\kappa_*=\gamma_{g.v_1}$ follows from Lemma 
\ref{lemma-recovering-f}. 

Assume next condition (2b) with $g=D\circ f$. 
The equality $(D_{\M(v_1)}\otimes D\circ g)(u_{v_1})=
(\kappa\times id_{S_2})^*(u_{v_2})$ follows from 
Lemma \ref{lemma-Chow-theoretic-formula-for-gamma}.
Hence, $\kappa_*\circ D_{\M(v_1)}$ is a ring isomorphism
satisfying the equality 
$(\kappa_*\circ D_{\M(v_1)}\otimes D\circ g)(u_{v_1})=u_{v_2}$.
Lemma 
\ref{lemma-recovering-f} implies the equality 
$\gamma_{D\circ g,v_1}=\kappa_*\circ D_{\M(v_1)}$. 
\EndProof

%\medskip
%We will later show, that condition (1) above follows from (2) 
%(Lemma \ref{lemma-on-the-sign-conjecture}). 
%In the meanwhile, we will need the following definition.
%
%\begin{defi}
%\label{def-equivalence-functor-induces-an-isomorphism}
%{\rm
%Let $(S_i,v_i,H_i)$, $i=1,2$, be objects of the groupoid $\G$
%and $\Phi:D^b_{Coh}(S_1)\rightarrow D^b_{Coh}(S_2)$ an equivalence
%of derived categories. We say that
%{\em $\Phi$ induces an isomorphism from
%$\M_{H_1}(v_1)$ to $\M_{H_2}(v_2)$,
%}
%if conditions (2a) and (3) hold.
%We say that 
%{\em 
%the composition of $\Phi$ with the functor of taking dual induces
%an isomorphism from
%$\M_{H_1}(v_1)$ to $\M_{H_2}(v_2)$,
%}
%if conditions (2b) and (3) hold.
%}
%\end{defi}

\medskip
Condition (1) above follows from (2), when $(v_1,v_1)>0$. We prove
the implication 
(2a)$\Rightarrow$(1) using Theorems \ref{thm-action} and 
\ref{introduction-thm-Gamma-v-acts-motivicly}.
The proof of the implication (2b)$\Rightarrow$(1) is similar. 
Let 
$\Phi:D(S_1)\rightarrow D(S_2)$ be an equivalence,
$v_1\in K_{top}(S_1)$ a primitive and effective class satisfying 
$(v_1,v_1)>0$,
and $H_1$ a $v_1$-suitable polarization. Assume that the isometry
$\phi:K_{top}(S_1)\rightarrow K_{top}(S_2)$ is orientation 
{\em reversing} (Remark \ref{rem-natural-orientation}). 
Set $v_2:=\phi(v_1)$.
%Assume, further, that the class $v_2:=\phi(v_1)$ (???)
%is effective (???) (otherwise, replace $\Phi$ by its shift $[1]\circ \Phi$). 

\begin{new-lemma}
\label{lemma-on-the-sign-conjecture}
For every $v_2$-suitable polarization $H_2$ and for every integer $i$, 
there exists an $H_1$-stable sheaf $E_1$ with Mukai 
vector $v_1$, such that $\Phi(E_1)[i]$ is not equivalent in
$D(S_2)$ to any $H_2$-stable sheaf $E_2$ with Mukai vector 
$(-1)^i(v_2)$.
%Condition (2) above fails for all 
%$v_2$-suitable polarizations $H_2$. Consequently, 
%the equivalence $\Phi$ does not induce an isomorphism from
%$\M_{H_1}(v)$ to $\M_{H_2}(\phi(v))$,  
%%(Definition \ref{def-equivalence-functor-induces-an-isomorphism}), 
%for any $v_2$-suitable polarization $H_2$.
\end{new-lemma}

\noindent
{\bf Proof:} 
Assume otherwise that a pair $(H_2,i)$ exists, 
which satisfies condition (2a).
%such that 
%$\Phi(E_1)[i]$ is equivalent in
%$D(S_2)$ to some $H_2$-stable sheaf $E_2$ with Mukai vector 
%$(-1)^i(v_2)$, for every $E_1\in \M_{H_1}(v_1)$.
%We may assume that $i=0$, after replacing $\Phi$ by $[i]\circ\Phi$. 
%Then condition (2a) is satisfied. 
%by some twisted sheaf $\E'$, 
%by a flatness result for Fourier-Mukai functors (\cite{mukai-applications} 
%Theorem 1.6). 
Let 
$\kappa:\M_{H_1}(v_1)\rightarrow \M_{H_2}(v_2)$ be the isomorphism 
(\ref{eq-kappa}). 
Then 
%$\kappa_*\otimes \bar{\phi}$ satisfies 
%equation (\ref{eq-characterization-of-gamma-g}). Thus 
$\kappa=\gamma_{\phi,v_1}$ by the proof of Lemma
\ref{lemma-G-Mukai-is-in-G-OK}.
%(part \ref{thm-item-characterization-of-gamma-g} of 
%Theorem \ref{thm-action}).
On the other hand, the isomorphism 
$\monrep(\phi):=D_{\M_{H_2}(v_2)}\circ \gamma_{\phi,v_1}$ given in
(\ref{eq-image-of-g-by-H-mon}) was shown to be a monodromy operator
in the proof of Theorem \ref{introduction-thm-Gamma-v-acts-motivicly}. 
It follows that the automorphism $D_{\M_{H_2}(v_2)}$
is a monodromy operator of $\M_{H_2}(v_2)$. The lattice 
$H^2(\M_{H_2}(v_2),\Integers)$ is isometric to $v_2^\perp$,
by the assumption that $(v,v)\geq 2$
(Theorem \ref{thm-irreducibility}). $D_{\M_{H_2}(v_2)}$
acts on $H^2(\M_{H_2}(v_2),\Integers)$
as multiplication by $-1$, which is orientation reversing.
This contradicts the fact, that the distinguished orientation 
is monodromy invariant.
\EndProof

\section{The surface-monodromy representation}
\label{sec-the-small-monodromy}

Let $S_1$ and $S_2$ be two $K3$ surfaces. 
A signed isometry from $H^2(S_1,\Integers)$ to $H^2(S_2,\Integers)$ 
is an isometry, which maps 
the distinguished orientation of the positive cone of $H^2(S_1,\RealNumbers)$
to that of $H^2(S_2,\RealNumbers)$ 
(Remark \ref{rem-natural-orientation}).
The set of signed isometries of the second cohomologies
of the two surfaces naturally embeds 
in ${\rm Isom}[H^*(S_1,\Integers),H^*(S_2,\Integers)]$. 
An isometry is extended to $H^i(S_1,\Integers)$, for $i=0$ and $4$, 
by sending the Poincare duals of a point and of $S_1$ to 
the corresponding classes in $H^i(S_2,\Integers)$. 
%Let 
%\begin{equation}
%\label{eq-G-n}
%\G_{[n]}
%\end{equation}
%be the following sub-groupoid of the groupoid
%$\G$ given in (\ref{eq-cohomological-groupoid}). 
%The objects of $\G_{[n]}$ consist of triples
%$(S,v,H)$, where $v=(1,0,1-n)$ is the Mukai vector of the ideal sheaf of 
%a length $n$ subscheme. 
%A morphism in $\Hom_{\G_{[n]}}(x_1,x_2)$,
%between two objects $x_i=(S_i,(1,0,1-n),H_i)$ of $\G_{[n]}$, is an isometry 
%in ${\rm Isom}[H^*(S_1,\Integers),H^*(S_2,\Integers)]$, arising from a 
%signed isometry of the second cohomologies. 
%
The main result of this section is the following special case 
of Theorem \ref{thm-trancendental-reflections}. 

\begin{thm}
\label{thm-conj-holds-for-hilbert-schemes-and-isometries-of-H-2}
Let $x_i:=(S_i,v_i,H_i)$, $i=1,2$, be two objects of the groupoid $\G$.
Assume that either each $c_1(v_i)$ is a non-zero multiple of an ample class,
or that both $v_i$ are the Mukai vectors $(1,0,1-n)$
of $S_i^{[n]}$, $n\geq 1$.
Let $g\in \Hom_\G(x_1,x_2)$ be a morphism arising from a signed isometry
of the second cohomologies of the two surfaces. 
In particular, both $v_i$ have the same rank and Euler characteristic.
Then $\gamma_{g,v_1}$ is a monodromy operator.
If $(v_i,v_i)>0$, then $(\gamma_{g,v_1}\otimes g)(u_{v_1})=u_{v_2}$.
\end{thm}

The Hilbert scheme case of the Theorem 
is used in the reduction of the proof of Theorem \ref{thm-action}
in section \ref{sec-reduction-of-proof-of-thm-action}.
The proof of Proposition 
\ref{reflection-of-Hilbert-schemes-with-respect-to-topological-lb} 
uses Theorem 
\ref{thm-conj-holds-for-hilbert-schemes-and-isometries-of-H-2}
for Mukai vectors $v_i$ with rank $0$.
%Theorem \ref{thm-conj-holds-for-hilbert-schemes-and-isometries-of-H-2}
%follows from  Proposition \ref{prop-torelli-and-monodromy} 
%part \ref{prop-item-hypothesis-hold} and Lemma
%\ref{lemma-the-class-of-f-is-gamma-E1-E2}.
%The automorphism group of an object of $\G_{[n]}$ is $\Gamma_0^{cov}$ 
%(Definition \ref{def-covariant-subgroups}). 
We will prove that the morphism $g$ in Theorem 
\ref{thm-conj-holds-for-hilbert-schemes-and-isometries-of-H-2}
is a morphism of the sub-groupoid
$\G_{mon}$, used already in 
section \ref{sec-reduction-of-proof-of-thm-action} 
%in the reduction of the proof of Theorem \ref{thm-action}, 
and defined below (Definition 
\ref{def-deformation-equivalent-triples}).
We refer to the operator $\gamma_{g,v_1}$ in the Theorem as a
{\em surface-monodromy\/} operator because it 
arises from deformations of the surface $S_1$ to $S_2$. 
The irreducible holomorphic symplectic manifold 
$S^{[n]}$ admits deformations, which do not arise from deformations 
of the $K3$ surface $S$. These deformations give rise to a 
larger monodromy group 
(Theorem \ref{introduction-thm-Gamma-v-acts-motivicly}). 

The Hilbert schemes case of 
Theorem \ref{thm-conj-holds-for-hilbert-schemes-and-isometries-of-H-2}
is proven in section \ref{sec-monodromy-for-K3}. 
The case $c_1(v_i)\neq 0$ of Theorem 
\ref{thm-conj-holds-for-hilbert-schemes-and-isometries-of-H-2}
is proven in section \ref{sec-monodromy-for-polarized-K3}.
The main ingredients are the Strong Torelli Theorem and 
the surjectivity of the period map for K3 surfaces. 
%An analogue of Theorem 
%\ref{thm-conj-holds-for-hilbert-schemes-and-isometries-of-H-2},
%for more general moduli spaces of sheaves, is proven in
%Proposition \ref{prop-torelli-and-monodromy}
%part \ref{prop-item-full-monodromy}.

%*********************************************************************
% Deformation equivalence
%*********************************************************************
\subsection{Deformation equivalence}
\label{sec-deformation-equivalence}

We prove next that the homomorphism
$\gamma_{g,v_1}:H^*(\M_{H_1}(v_1),\RationalNumbers)\rightarrow 
H^*(\M_{H_1}(v_1),\RationalNumbers)$ is a monodromy operator, whenever
the moduli space $\M_{H_1}(v_1)$ on one $K3$ surface 
is related to a moduli space 
$\M_{H_2}(v_2)$ on another surface 
via a deformation of the surfaces (Lemma 
\ref{lemma-the-class-of-f-is-gamma-E1-E2}).

%Theorem \ref{thm-graph-of-diagonal-in-terms-of-universal-sheaves}
%can be viewed as a factorization of the class of the diagonal. 
%It enables us to deform the class of the diagonal, 
%by holding one factor fixed and deforming the other factor. 
%An immediate consequence is a formula 
%lifting a monodromy representation
%from $H^*(S,\Integers)$ to $H^*(\M_\LB(v),\Integers)$
%(Lemma \ref{lemma-the-class-of-f-is-gamma-E1-E2}). 

Let $T$ be a connected algebraic space, 
$p:\S\rightarrow T$ a smooth family of $K3$ surfaces, 
and $\LB$ a line bundle on $\S$.
Set $v:=(r,d\LB,s)$ to be the flat section of the local system
$R^*p_*\Integers$, with $\gcd(r,d,s)=1$. 
Let $H$ be a section of $R^2p_*\Integers$, consisting of ample classes,
such that $H_t$ is $v_t$-suitable, for every $t\in T$.
We do {\em not} assume $H$ to be continuous.
Let $Spl_{\S/T}\rightarrow T$ be the moduli space of simple sheaves $E$ on 
$\S_t$, $t\in T$, with Mukai vector equal $v_t$ 
(see \cite{altman-kleiman} Theorem 7.4 for 
its existence). 
Assume, that there exists an algebraic subspace 
${\frak M}(v)\subset Spl_{\S/T}$ and a smooth and proper morphism 
${\frak M}(v)\rightarrow T$, such that for every point $t\in T$, 
${\frak M}(v)_t$ is isomorphic to the moduli space 
$\M_{H_t}(v_t)$ of $H_t$-stable sheaves on $\S_t$ with Mukai vector $v_t$. 

We define next the subgroupoid $\G_{mon}$ of the groupoid $\G$
given in (\ref{eq-cohomological-groupoid}). The object of $\G_{mon}$
are the same. The morphisms arise from deformations, which we now describe. 
Let $S_i$ be a projective $K3$ surface with an ample line bundle $H_i$,
$i=1,2$. Let $v_i:=(r,dL_i,s)$ be an algebraic Mukai vector
in $H^*(S_i,\Integers)$. Assume, that $H_i$ is $v_i$-suitable. 

\begin{defi}
\label{def-deformation-equivalent-triples}
{\rm
\begin{enumerate}
\item
The objects $(S_1,v_1,H_1)$ and $(S_2,v_2,H_2)$ of $\G$ 
are said to be {\em deformation-equivalent}, if there exists data 
$(T,\S,\LB,H,{\frak M}(v))$ as above, and two points
$t_1,t_2\in T$, such that $(S_i,v_i,H_i)=(S_{t_i},(r,d\LB_{t_i},s),H_{t_i})$,
$i=1,2$.
If $v_i=(1,0,1-n)$, $i=1,2$, $n\geq 1$, 
then we allow $H$ to be trivial on $T\setminus\{t_1,t_2\}$, 
$p:\S\rightarrow T$ to be a smooth and proper analytic family of
K\"{a}hler $K3$ surfaces over a connected analytic space $T$,
and we let ${\frak M}(v)\rightarrow T$ be the relative Douady space 
of ideal sheaves of length $n$ zero-dimensional subschemes
of the fibers of $p$.
\item 
A morphism $g\in \Hom_{\G}((S_1,v_1,H_1),(S_2,v_2,H_2))$ is a morphism 
in $\G_{mon}$, if the two objects are deformation-equivalent by data  
$(T,\S,\LB,H,{\frak M}(v))$ as above, and $g$ 
is a monodromy operator, corresponding to a homotopy class of paths in 
the base $T$ from $t_1$ to $t_2$.
\end{enumerate}
}
\end{defi}

$\G_{mon}$ is easily seen to be a sub-groupoid of $\G$. 

\begin{prop}
\label{prop-yoshioka-5.1}
(\cite{yoshioka-abelian-surface} Proposition 5.1)
Assume that the base $T$, of the smooth family $p:\S\rightarrow T$,
is a smooth connected quasi-projective curve, the line bundle
$\LB$ is relatively ample, and $\Pic(\S_t)$ is generated by $\LB_t$, 
for some $t\in T$. Then $\LB_t$ is $v_t$-suitable for all
but finitely many $t\in T$. Choose a section $H$ of $R^2_{p_*}\Integers$  
satisfying $H_t=\LB_t$, for all but finitely many $t\in T$, and 
such that $H_t$ is $v_t$-suitable, for every $t\in T$.
Then there exists an algebraic subspace 
${\frak M}(v)\subset Spl_{\S/T}$ and a smooth and proper morphism 
${\frak M}(v)\rightarrow T$, such that for every point $t\in T$, 
${\frak M}(v)_t$ is isomorphic to the moduli space 
$\M_{H_t}(v_t)$ of $H_t$-stable sheaves on $\S_t$ with Mukai vector $v_t$. 
%Note: the existence of ${\frak M}(v)$ as above is proven in 
%\cite{yoshioka-abelian-surface} Proposition 5.1, under the 
%additional assumptions that
%$T$ is a smooth curve, $\LB$ is relatively ample, 
%$\Pic(\S_t)$ is generated by $\LB_t$, for some $t\in T$,
%provided we choose $H_t=\LB_t$ for all but finitely many $t\in T$
%(it is further proven that $\LB_t$ is $v_t$-suitable for all
%but finitely many $t\in T$).
\end{prop}

If $H_i$, $i=1,2$, are two $v$-suitable polarizations and either
$c_1(v)$ is ample or $-c_1(v)$ is ample, 
then the objects $(S,v,H_1)$ and $(S,v,H_2)$ are deformation-equivalent,
by  Proposition \ref{prop-yoshioka-5.1},
and the identity isometry is a morphism in 
$\Hom_{\G_{mon}}((S,v,H_1),(S,v,H_2))$. More generally, we have:

\begin{cor}
\label{cor-same-r-s-k-imply-deformation-equiv}
Let $(S_i,v_i,H_i)$, $i=1,2$, be two objects of $\G$ with Mukai vectors 
$v_i=(r,kA_i,s)$, with the same $r$, $s$, $k$. Assume that 
$c_1(A_1)^2=c_1(A_2)^2$ 
and that both $A_i$ are ample and indivisible.
Then $(S_i,v_i,H_i)$, $i=1,2$, are deformation equivalent. 
\end{cor}

\noindent
{\bf Proof:} The case where $H_i=A_i$ follows from the proof of
\cite{ogrady-hodge-str} Proposition 2.3. Each $(S_i,v_i,H_i)$
is deformation equivalent to an object $(S'_i,v'_i,H'_i)$,
where $v'_i=(r,kA'_i,s)$, and $H'_i=A'_i$, by 
Proposition \ref{prop-yoshioka-5.1}.
\EndProof

\begin{new-lemma}
\label{lemma-the-class-of-f-is-gamma-E1-E2}
Let $g\in \Hom_{\G_{mon}}((S_1,v_1,H_1),(S_2,v_2,H_2))$ be a morphism. 
Then 
\[
\gamma_{g,v_1}  \ : \ H^*(\M_{H_1}(v_1),\RationalNumbers)  \rightarrow 
H^*(\M_{H_2}(v_2),\RationalNumbers)
\]
is a monodromy operator. Furthermore, 
$(\gamma_{g,v_1}\otimes g)(u_{v_1})= u_{v_2}$.
\end{new-lemma}

\noindent
{\bf Proof:} We prove the case $(v_i,v_i)\geq 2$. The case $(v_i,v_i)\leq 0$
is similar. Let 
$\pi: {\frak M}(v)\times_T\S\rightarrow T$ be the fiber product. 
There exists a $2$-cocycle $\alpha$ of the sheaf 
$\StructureSheaf{{\frak M}(v)}^*$, defining a class in the Brauer group of
${\frak M}(v)$, and an $\alpha$-twisted universal sheaf $\E$ 
over ${\frak M}(v)\times_T\S$ (by the proof of Theorem A.5 in 
\cite{mukai-hodge}). The normalization of the
Chern character of a twisted universal sheaf, given in equations
(\ref{eq-chern-character-of-rational-universal-class}) and 
(\ref{eq-invariant-normalized-class-of-chern-character-of-universal-sheaf}),
produces a flat section $u_v$ of $R^*_{\pi_*}\RationalNumbers$, with value
$u_{v_t}$ over $t\in T$ given by 
(\ref{eq-invariant-normalized-class-of-chern-character-of-universal-sheaf}).

Let $\eta$ be a homotopy class of paths in $T$ from $t_1$ to $t_2$. 
%Let $\E_i$ be the restriction of $\E$ to 
%$S_{t_i}\times {\frak M}(v_i)_{t_i}$.
We get isomorphisms (which need not preserve Hodge structures) 
\begin{eqnarray}
\label{eq-g-and-f}
g \ : \ H^*(S_1,\RationalNumbers) & \rightarrow &  H^*(S_2,\RationalNumbers)
\\
\nonumber
f \ : \ H^*(\M_{H_1}(v_1),\RationalNumbers) & \rightarrow &  
H^*(\M_{H_2}(v_2),\RationalNumbers)
\end{eqnarray}
satisfying $g(v_1)=v_2$  and 
\begin{equation}
\label{eq-monodromy-takes-chern-character-to-chern-character}
(f\otimes g)(u_{v_1}) \ = \ u_{v_2}. 
\end{equation}
The equality 
$f = \gamma_{g,v_1}$
follows  from Lemma \ref{lemma-recovering-f}.
\EndProof

\subsection{Monodromy for K\"{a}hler K3 surfaces and Torelli}
\label{sec-monodromy-for-K3}

We prove the case $v_i=(1,0,1-n)$ of Theorem
\ref{thm-conj-holds-for-hilbert-schemes-and-isometries-of-H-2} 
in this section.

%\begin{prop}
%\label{prop-torelli-and-monodromy}
%\begin{enumerate}
%\item
%\label{prop-item-hypothesis-hold}
%The sub-groupoid $\G_{[n]}$ of $\G$, given in (\ref{eq-G-n}),
%is a sub-groupoid of $\G_{mon}$ as well.
%\item
%\label{prop-item-stabilizer-of-ample-lb-consists-of-monodromies}
%Let $(S,v,H)$ be an object of $\G$, $v=(r,kc_1(A),s)$, $k$ an integer,
%$r\geq 0$, and
%$A$ an ample line-bundle on $S$ with $c_1(A)$ an indivisible class. 
%Then the stabilizer of $c_1(A)$ 
%in $\Gamma_0^{cov}$ is contained in 
%$\Aut_{\G_{mon}}(S,v,H)$.
%\item
%\label{prop-item-Mon-S}
%The group $\Gamma_0^{cov}$
%is equal to the monodromy group of $S$.
%\end{enumerate}
%\end{prop}
%
%Proposition \ref{prop-torelli-and-monodromy} is proven below using 
%the relationship between the Torelli theorem for (K\"{a}hler) 
%K3 surfaces, the surjectivity of the period map, 
%and the monodromy of the cohomology of the surface. 

Let $\Gamma$ be a subgroup of the automorphism group of a lattice $L$, 
$\pi:\S\rightarrow M$ a family of smooth compact complex varieties over
a connected base $M$, and
$\Phi : R^i_{\pi_*}\Integers \IsomRightArrow L$ a trivialization of the
weight $i$ local system of integral cohomology groups.
Denote by $\phi_m$ the value of $\Phi$ at $m\in M$. 
Let $Ad_{\phi_m}:\Gamma\rightarrow \Aut H^i(S_m,\Integers)$
be the homomorphism sending $\gamma\in \Gamma$ to $\phi_m^{-1}\gamma\phi_m$. 

\begin{new-lemma}
\label{lemma-equivariant-family-implies-monodromy}
Assume that the following conditions hold:
\begin{enumerate}
\item
\label{assumption-automorphisms-act-faithfully}
The automorphism group $Aut(S_m)$ acts faithfully on 
$H^i(S_m,\Integers)$ for every $m\in M$.
\item
\label{assumption-Gamma-orbits-represent-isomorphic_K3}
$\Gamma$ acts on the base $M$ in such a way, 
that $(S_m,\gamma\circ\phi_m)$ and
$(S_{\gamma(m)},\phi_{\gamma(m)})$ are isomorphic, for every 
$\gamma\in \Gamma$. In other words, there exists a  unique
isomorphism $f:S_m\rightarrow S_{\gamma(m)}$ making the following diagram 
commutative 
\begin{equation}
\label{eq-commutative-diagram-of-trivializations}
\begin{array}{ccc}
H^i(S_m,\Integers) & \LongRightArrowOf{\phi_m} &
L
\\
\uparrow \ f^* & & \downarrow \gamma
\\
H^i(S_{\gamma(m)},\Integers) & \LongRightArrowOf{\phi_{\gamma(m)}} & L.
\end{array}
\end{equation}
\item
\label{assumption-non-trivial-action}
The action of $\Gamma$ on $M$ is non-trivial. I.e., 
there exists a pair $\gamma\in \Gamma$ and $m\in M$, 
such that $\gamma(m)\neq m$.
\end{enumerate}
Then the subgroup $Ad_{\phi_m}(\Gamma)$ of 
$\Aut H^i(S_m,\Integers)$ consists of 
monodromy operators. 
%and automorphisms induced by automorphisms of $S_m$. 
%In particular, if $m_0\in M$ is a point whose 
%stabilizer subgroup in $\Gamma$ 
%is trivial, then $\Gamma$ acts on $H^i(S_{m_0},\Integers)$ via monodromy. 
\end{new-lemma}

\medskip
\noindent
{\bf Proof:} 
The existence of the trivialization $\Phi$ implies, that $\pi_1(M)$
acts trivially on $R^i_{\pi_*}\Integers$. 
Let $m\in M$, $\gamma\in \Gamma$, and 
choose $\eta:[0,1]\rightarrow M$ a path from $m$ to $\gamma(m)$. 
As we deform $S_{\eta(t)}$ from $S_m$ to $S_{\gamma(m)}$,
the class of the diagonal in $S_m\times S_m$ deforms flatly, 
through classes in the cohomology of $S_m\times S_{\eta(t)}$,
to a class $F$ inducing a homomorphism
$F: H^i(S_{\gamma(m)}) \rightarrow H^i(S_m)$ satisfying
$
\phi_{\gamma(m)} = \phi_m\circ F. 
$
Let $f:S_m \rightarrow S_{\gamma(m)}$ be the isomorphism 
(\ref{eq-commutative-diagram-of-trivializations}) determined by
the equality $\gamma\circ\phi_m\circ f^* =  \phi_{\gamma(m)}$.
We get the equality
\[
(\phi_m^{-1}\circ\gamma\circ\phi_m) \ = \ F \circ (f^{-1})^*.
\]
If $\gamma(m)=m$, then $F$ is the class of the diagonal, 
and the automorphism $\phi_m^{-1}\circ\gamma\circ\phi_m$ of 
$H^i(S_m,\Integers)$ is induced by the automorphism $f$. 
If $\gamma(m)\neq m$, we can use the isomorphism $f$ to glue
the pulled back family $\eta^{-1}\S$ to a family $(\eta^{-1}\S)/f$ 
over a circle. 
The automorphism $\phi_m^{-1}\circ\gamma\circ\phi_m$ of 
$H^i(S_m,\Integers)$ is then the monodromy operator of the family
$(\eta^{-1}\S)/f$. 
%Hence, the automorphism $\phi_m^{-1}\circ\gamma\circ\phi_m$ of 
%$H^i(S_m,\Integers)$ is induced by monodromy and the automorphisms
%(\ref{eq-commutative-diagram-of-trivializations}) 
%of $S$, which are determined by elements $\gamma$ in the stabilizer of $m$.
Fix a point $m'\in M$ with a non-trivial $\Gamma$-orbit 
(assumption \ref{assumption-non-trivial-action}). 
Then $\Gamma$ is generated by the elements $\gamma$ satisfying 
$\gamma(m')\neq m'$. Hence, $Ad_{\phi_{m'}}(\Gamma)$ consists of monodromy
operators. The same holds for every $m\in M$, since $M$ is connected.
\EndProof

\medskip
Denote by $L$ the $K3$ lattice
\begin{equation}
\label{eq-K3-lattice}
-E_8\oplus -E_8 \oplus H \oplus H \oplus H, 
\end{equation}
where $H$ is the even rank 2 hyperbolic lattice. 
Set $L_\ComplexNumbers:=L\otimes_\Integers\ComplexNumbers$.
Let $Q\subset \PP(L\otimes_\Integers\ComplexNumbers)$ be the quadric defined
by the pairing. Let $\Omega\subset Q$ be the analytic open subset
$\{\omega \ : \ (\omega,\bar{\omega})>0\}$. 
%(the moduli of marked K3's). 
A {\em marked triple}
$(S,\kappa,\phi)$ consists of a $K3$ surface $S$, 
a K\"{a}hler form $\kappa$, and a marking by an isometry 
$\phi:H^2(S,\Integers)\rightarrow L$, mapping the distinguished orientation 
of the positive cone of $H^2(S,\RealNumbers)$ to a fixed orientation of
that of $L\otimes_\Integers\RealNumbers$ 
(see Remark \ref{rem-natural-orientation}). 
The moduli of marked triples is 
\begin{eqnarray*}
K\Omega^0 & := & \{([w],\kappa) \ \mid \ [w]\in \Omega, \ 
\kappa\in L\otimes_\Integers\RealNumbers\cap[\omega]^\perp, \ 
(\kappa,\kappa)>0, \\
& & \ \ \ \  
(\kappa,d)\not=0 \ \forall d\in L\cap [\omega]^\perp
\ \mbox{with} \ (d,d)=-2\}
\end{eqnarray*}
(\cite{BHPV} chapter VIII Theorems 12.3 and 14.1). 
$K\Omega^0$ is a connected $60$-dimensional real manifold
(\cite{BHPV} chapter VIII Corollary 9.2).
Let $\Gamma_0^{cov}$ be the group of orientation preserving isometries
of $L$ (Definition \ref{def-covariant-subgroups}). 
$K\Omega^0$ is a $\Gamma_0^{cov}$ invariant subset of 
$\PP(L_\ComplexNumbers)\times L_\ComplexNumbers$. 
There exists a universal family $\pi:\S\rightarrow K\Omega^0$ 
of $K3$ surfaces and a trivialization 
$\Phi:R^2_{\pi_*}\Integers\rightarrow L$, satisfying:
(i) Each $\phi_{([w],\kappa)}$ is an isometry,
(ii) $\phi_{([w],\kappa)}[H^{2,0}(S_{([w],\kappa)})]=[w]$, and 
(iii) $\phi_{([w],\kappa)}^{-1}(\kappa)$ is a K\"{a}hler class in
$H^{1,1}(S_{([w],\kappa)})$. 

Note: we avoid taking the quotient of $K\Omega^0$ by $\Gamma_0^{cov}$
because some pairs $(S,\kappa)$ do have non-trivial automorphisms. 
The automorphism group is equal to the subgroup of $\Gamma_0^{cov}$ 
stabilizing $([w],\kappa)$. In particular, the
universal family of K3's does not descend to the quotient.

\medskip
\noindent
{\bf Proof of the case $v_i=(1,0,1-n)$ of Theorem
\ref{thm-conj-holds-for-hilbert-schemes-and-isometries-of-H-2}:}
%{\bf Proof of Proposition \ref{prop-torelli-and-monodromy}:} 
%\ref{prop-item-hypothesis-hold}) 
The hypothesis of Lemma \ref{lemma-equivariant-family-implies-monodromy}
hold for K3 surfaces with $i=2$, $\Gamma=\Gamma_0^{cov}$, 
and $M$ is the fine moduli space $K\Omega^0$ of marked triples.
Assumptions 
\ref{assumption-automorphisms-act-faithfully} and 
\ref{assumption-Gamma-orbits-represent-isomorphic_K3} 
of Lemma \ref{lemma-equivariant-family-implies-monodromy}
follow from the Strong 
Torelli Theorem. Assumption
\ref{assumption-non-trivial-action} is clear. 
Lemma \ref{lemma-equivariant-family-implies-monodromy}
implies, that the isometry $g$ in Theorem
\ref{thm-conj-holds-for-hilbert-schemes-and-isometries-of-H-2} 
is a monodromy operator 
for some family $p:\S\rightarrow T$ of $K3$ surfaces. 
(The total space $\S$ is assumed to be complex analytic and $p$ is smooth
and proper).
There is a relative Douady family $\S^{[n]}\rightarrow T$,
parametrizing ideal sheaves of length $n$ subschemes in the
fibers of the  family. 
Theorem \ref{thm-conj-holds-for-hilbert-schemes-and-isometries-of-H-2} now 
follows from Lemma \ref{lemma-the-class-of-f-is-gamma-E1-E2}
in the Hilbert schemes case. 
\EndProof 

The special case $v_1=v_2=(1,0,0)$ of 
Theorem \ref{thm-conj-holds-for-hilbert-schemes-and-isometries-of-H-2} is
the well known:

\begin{cor}
\label{cor-Mon-S}
The group $\Gamma_0^{cov}$ is equal to the monodromy group of  a 
$K3$ surface $S$.
\end{cor}

\noindent
{\bf Proof:}
The above proof shows that the monodromy group $\Mon(S)$ contains 
$\Gamma_0^{cov}$. 
$\Mon(S)$ is contained in the index $2$ subgroup $\Gamma_0^{cov}$ of
$\Gamma_0$,
since the distinguished orientation of the positive cone is 
monodromy invariant (Remark \ref{rem-natural-orientation}).
\EndProof

%\ref{prop-item-Mon-S})
%The monodromy group is contained in the index $2$ subgroup $\Gamma_0^{cov}$
%(see \cite{donaldson} Proposition 6.2). 
%Part \ref{prop-item-Mon-S} follows from part
%\ref{prop-item-hypothesis-hold} via Lemma 
%\ref{lemma-equivariant-family-implies-monodromy}. 

%**************************************************************
% Monodromy for polarized K3 surfaces and Torelli
%**************************************************************
\subsection{Monodromy for polarized K3 surfaces and Torelli}
\label{sec-monodromy-for-polarized-K3}
We prove in this section the case of Theorem 
\ref{thm-conj-holds-for-hilbert-schemes-and-isometries-of-H-2}, where
$c_1(v_i)$ is a non-zero multiple of an ample class.
The proof actually works when $c_1(v_i)$ is zero, under 
the additional assumption that the isometry $g$ 
in Theorem 
\ref{thm-conj-holds-for-hilbert-schemes-and-isometries-of-H-2}
maps some ample class of $S_1$ to an ample class of $S_2$.
We already know, that the objects $x_1$ and $x_2$ 
are related by some morphism $g$ in $\G_{mon}$, by Corollary
\ref{cor-same-r-s-k-imply-deformation-equiv}. Furthermore,  
$\gamma_{g,v_1}$ has the properties 
stated in Theorem 
\ref{thm-conj-holds-for-hilbert-schemes-and-isometries-of-H-2},
by Lemma \ref{lemma-the-class-of-f-is-gamma-E1-E2}.
%\ref{prop-item-stabilizer-of-ample-lb-consists-of-monodromies})
Hence, it suffices to prove the
statement for automorphisms $g$ of 
one object $(S,v,H)$ with $v=(r,kA,s)$ and $A$ ample, 
for each value of $r$, $s$, $k$, and for $c_1(A)^2=2d$, $d\geq 1$.
We need only consider automorphisms arising from 
orientation-preserving isometries of $H^2(S,\Integers)$ stabilizing $A$, 
and prove that they belong to $\Aut_{\G_{mon}}(S,v,H)$.

We keep the notation of section
\ref{sec-monodromy-for-K3}.
Fix a primitive element $h\in L$ satisfying $(h,h)=2d$, $d\geq 1$.
Let $\Gamma_h$ be the subgroup of $\Gamma_0$ stabilizing $h$ and
set $\Gamma_h^{cov}:=\Gamma_h\cap \Gamma_0^{cov}$.
Set $\Omega'_h:=\{\omega\in\Omega \ : \ (\omega,h)=0\}$. 
$\Omega'_h$ has two connected components, which are interchanged by 
complex conjugation. The two components are distinguished by the orientation 
of the positive cone of $L_\RealNumbers$ determined by the basis 
$\{{\rm Re}(\omega),{\rm Im}(\omega),h\}$ of the positive $3$-dimensional
subspace (see Remark \ref{rem-natural-orientation}).
Choose  one of the connected components of 
$\Omega'_h$ and denote it by $\Omega_h$. 
Then $\Omega_h$ is $\Gamma_h^{cov}$-invariant. 
Set
\[
K\Omega_h^0 \ \ := \ \ 
\left\{
\omega\in \Omega_h \ : \ (\omega,h) \in K\Omega^0
\right\}.
\]
$K\Omega_h^0$ is a dense open subset of $\Omega_h$.  
$K\Omega_h^0$ admits a natural embedding as a closed subset of
$K\Omega^0$. 

$\Gamma_h^{cov}$ acts on $\Omega_h$ properly-discontinuously
and hence the quotient $\Omega_h/\Gamma_h^{cov}$ has a canonical structure 
of a normal analytic space. $\Omega_h/\Gamma_h^{cov}$ is moreover a (normal) 
quasi-projective variety \cite{baily-borel}.
%The quotient $Y'_{2d}:=K\Omega_h^0/\Gamma_h^{cov}$ is a Zariski dense open 
%subset of $\Omega_h/\Gamma_h^{cov}$.
The complement $Y_{2d}$ in $\Omega_h/\Gamma_h^{cov}$, 
of the branch locus of the quotient map, is a Zariski dense smooth open 
subset of $\Omega_h/\Gamma_h^{cov}$ (\cite{peters} Proposition 2.2.2).

Let $\widetilde{Y}_{2d}$ be the complement in $\Omega_h$
of the ramification locus of the quotient map. 
Then $\widetilde{Y}_{2d}$ is an open and dense subset of $K\Omega_h^0$. 
The universal family of marked $K3$ surfaces over $K\Omega^0$ restricts to one
over $\widetilde{Y}_{2d}$, which we denote by
$\tilde{\pi}:\widetilde{\S}\rightarrow \widetilde{Y}_{2d}$.
The natural section $h$ of $R^2_{\tilde{\pi}_*}\Integers$ corresponds to an 
ample line-bundle $L_t$ on the fiber $\S_t$ 
(\cite{BHPV} chapter VIII Corollary 3.9). 
Furthermore, the identity is the only automorphism of $\S_t$ 
leaving the class of $L_t$ invariant. 
Hence, the family $\tilde{\pi}$ descends to $Y_{2d}$. 
The above construction is summarized in the following Proposition.

\begin{prop} 
\label{prop-algebraic-family-of-K3-of-degree-2d}
(A special case of \cite{peters} Theorem 3.4.1 and Corollary 
3.4.2)
\begin{enumerate}
\item
There exists a smooth and proper family of projective $K3$ surfaces 
$\pi:\S\rightarrow Y_{2d}$, 
over the smooth quasi-projective variety $Y_{2d}$,
whose pullback to $\widetilde{Y}_{2d}$ is isomorphic to the universal family 
$\tilde{\pi}:\widetilde{\S}\rightarrow \widetilde{Y}_{2d}$.
%and a (???) holomorphic line-bundle $\LB$ over $\S$,
%restricting to each fiber of $\pi$ as 
%an ample line bundle $\LB_t$ of degree $2d$. 
%Each class $c_1(\LB_t)$ is indivisible. 
\item
\label{prop-item-full-monodromy}
Choose a point $(\omega,h)$ in the subset $\widetilde{Y}_{2d}$ 
of $K\Omega^0_h$. Let $t$ be the corresponding $\Gamma_h^{cov}$-orbit in 
$Y_{2d}$. We get an associated marking 
$\phi_{(\omega,h)}:H^2(\S_t,\Integers)\rightarrow L$.
The monodromy group of the local system $R^2\pi_*\Integers$
is conjugated, via the marking $\phi_{(\omega,h)}$,
to the subgroup $\Gamma_h^{cov}$ of $O(L)$.
\end{enumerate}
\end{prop}

The construction of a smooth and proper relative moduli space of sheaves
need not be possible over the whole family $\pi:\S\rightarrow Y_{2d}$. 
We replace the above family by one, whose base is a smooth curve,
for which Yoshioka's construction (Proposition
\ref{prop-yoshioka-5.1}) applies.
There exists a smooth connected quasi-projective
curve $C\subset Y_{2d}$, for which the homomorphism
$\pi_1(C)\rightarrow \pi_1(Y_{2d})$ is surjective. This follows 
from a version of Lefschetz' Hyperplane Section Theorem
(\cite{deligne-LHST}, see also \cite{fulton-lazarsfeld} Theorem 1.1). 
Denote by $p:\restricted{\S}{C}\rightarrow C$ the restriction
of the family $\pi:\S\rightarrow Y_{2d}$.
We may and do further assume, that $C$ passes through a point $y\in Y_{2d}$, 
corresponding to a $K3$ surface $\S_y$ with a cyclic Picard group. 
The tautological section of $R^2_{\tilde{\pi}_*}\Integers$ over
$\widetilde{Y}_{2d}$ comes from a section $h$ of 
$R^2_{\pi_*}\Integers$, which restricts to a section 
$\restricted{h}{C}$ of $R^2_{p_*}\Integers$.

We would like to apply Proposition \ref{prop-yoshioka-5.1} 
with the family $p:\restricted{\S}{C}\rightarrow C$, the
section $v:=(r,k\restricted{h}{C},s)$ of the local system of Mukai lattices,
%the line bundle $\LB$ of Claim \ref{claim-algebraic-LB}, 
and a $v$-suitable
choice $H$ of a section of $R^2p_*\Integers$. 
As in Proposition \ref{prop-yoshioka-5.1}, the section $H$ is assumed to be 
equal to $\restricted{h}{C}$ away from finitely many points of $C$.
We are missing however a relatively ample line bundle $\LB$ on 
$\restricted{\S}{C}$ inducing the section $\restricted{h}{C}$. 
We can find instead an open \'{e}tale covering $\{U_i\}$ of $C$, and 
such a line bundle $\LB_i$ over each $\restricted{\S}{U_i}$.
This follows from the algebraic construction of the moduli space of 
polarized $K3$ surfaces (\cite{viehweg} Theorem 1.13)
and the existence of a universal family over the Hilbert scheme of
polarized $K3$ surfaces (\cite{viehweg} section 1.7). 
Proposition \ref{prop-yoshioka-5.1} produces the moduli spaces 
${\frak M}(v)_i\rightarrow U_i$, which glue to the desired algebraic space 
${\frak M}(v)$, which is smooth and proper over $C$. 

Let $v_t:=(r,k\cdot h(t),s)$
be the value of $v$ at a point $t$ of $C$, 
%$H_t$ a $v_t$-suitable ample line-bundle over $\S_t$, 
and $\eta$ an element of $\pi_1(C,t)$ corresponding to a monodromy
operator $g$ of $H^*(\S_t,\Integers)$. Then 
$\eta$ induces 
a surface-monodromy element $g$ in $\Aut_{\G_{mon}}(\S_t,v_t,H_t)$,
by Definition \ref{def-deformation-equivalent-triples}. 
Every element of $\Gamma_h^{cov}$ is obtained this way,
by the surjectivity of the homomorphism 
$\pi_1(C,t)\rightarrow \pi_1(Y_{2d},t)$
and by part \ref{prop-item-full-monodromy} of Proposition 
\ref{prop-algebraic-family-of-K3-of-degree-2d}.
This completes the proof of Theorem 
\ref{thm-conj-holds-for-hilbert-schemes-and-isometries-of-H-2}.
\EndProof

%*********************************************************************
% Examples of Hodge isometries
%*********************************************************************
\section{Examples}
\label{sec-examples-of-hodge-isometries}

We verify in this section special cases of Theorem
\ref{introduction-thm-Gamma-v-acts-motivicly}, 
inducing monodromy operators of the Hilbert scheme $S^{[n]}$, 
which do not come from deformations of the K3 surface $S$. 
In section \ref{sec-elliptic-k3} we consider reflections of the
Mukai lattice of $S$ with respect to the class of
a {\em topological} line-bundle on $S$. The automorphisms
$\gamma_{g}$, associated to such reflections $g$, 
are shown to be monodromy involutions of the cohomology
of Hilbert schemes $S^{[n]}$ (Proposition 
\ref{reflection-of-Hilbert-schemes-with-respect-to-topological-lb}). 
In section \ref{sec-reflections-by-spherical-objects}
we review work of Seidel-Thomas, on the reflections of the derived category 
with respect to spherical objects. 
We then state an analogue of Proposition 
\ref{reflection-of-Hilbert-schemes-with-respect-to-topological-lb},
when the reflection $g$ is with respect to the class of an 
{\em algebraic} line bundle on $S$
(Theorem 
\ref{thm-class-of-correspondence-in-stratified-elementary-trans}).
%corresponding to reflections $\gamma_{\tau}$ 
%of Hilbert schemes with respect to the class of an 
%{\em algebraic} line bundles on $S$. 
The automorphisms
$\gamma_{g}$, associated to such reflections $g$, 
are shown to be {\em local} monodromy operators of $S^{[n]}$.
The isometry $-1$ of the Mukai lattice is treated in section
\ref{sec-minus-id}. 
%In section \ref{sec-the-serre-duality-automorphism-of-a-moduli}
%we remark on the orientation reversing isometry induced by the 
%duality operator. 
In section
\ref{sec-reflections-with-respect-to-plus-2-vectors}
we verify cases of Theorem
\ref{introduction-thm-Gamma-v-acts-motivicly}
for certain orientation reversing isometries (Theorems
\ref{thm+2-reflection-induces-an-isomorphism} and 
\ref{thm-reflection-sigma-satisfies-main-conj}).

%****************************************************************
% Elliptic K3 surfaces
%****************************************************************
\subsection{Elliptic K3 surfaces}
\label{sec-elliptic-k3}
Let $S$ be a $K3$ surface admitting an elliptic fibration 
$\pi:S\rightarrow\PP^1$ with integral fibers and with a section.
Denote by $\sigma$ and $f$ the classes of the section and the fiber in
$H^2(S,\Integers)$. Let $v:=(1,0,1-n)$ be the Mukai vector of the
ideal sheaf of a length $n$ subscheme, $n\geq 2$. 
Fix a class $\beta\in H^2(S,\Integers)$, which is orthogonal to both
$\sigma$ and $f$, and satisfies $(\beta,\beta)=2n-4$.
We do {\em not} assume that $\beta$ is of Hodge-type $(1,1)$.
Then $v_0:=(1,\beta-f,n-1)$ is 
the Mukai vector in $v^\perp$ of a topological line-bundle 
and $(v_0,v_0)=-2$. 
Let $g$ be the reflection of $H^*(S,\Integers)$ with respect to $v_0$.

\begin{prop}
\label{reflection-of-Hilbert-schemes-with-respect-to-topological-lb}
The endomorphism $\gamma_{g,v}\in\End[H^*(S^{[n]},\Integers)_{free}]$
is a monodromy operator and $(\gamma_{g,v}\otimes g)(u_v)=u_v$,
where $u_v$ is the class 
(\ref{eq-invariant-normalized-class-of-chern-character-of-universal-sheaf}).
\end{prop}

\noindent
{\bf Proof:} 
We consider below an auto-equivalence $\Phi$ of $D(S)$, inducing 
a Hodge isometry $\phi$ of $H^*(S,\Integers)$ as well as an isomorphism 
$\M(v)\cong\M_H(\phi(v))$, such that $\phi g\phi^{-1}$
is a {\em surface-monodromy operator} of $S$ and
an element of $\Aut_{\G}(S,\phi(v),H)$. 
Theorem 
\ref{thm-conj-holds-for-hilbert-schemes-and-isometries-of-H-2}
proves the desired properties for $\gamma_{\phi g\phi^{-1},\phi(v)}$,
from which 
Proposition 
\ref{reflection-of-Hilbert-schemes-with-respect-to-topological-lb}
then easily follows.

We recall first the set-up of Theorem 3.15 in \cite{yoshioka-abelian-surface}.
Set $H:=\sigma+kf$, $k>>0$. Then $H$ is a $(0,f,0)$-suitable polarization, 
and the moduli space $\M_H(0,f,0)$, which is a compactified relative 
Jacobian, is isomorphic to $S$. We choose a
universal sheaf on $S\times_{\PP^1}\M_H(0,f,0)$ and
regard it as a sheaf $\P$ on $S\times S$.
Let $p_i:S\times S\rightarrow S$, $i=1,2$, be the two projections. 
The Fourier-Mukai functor $\Phi^\P:D(S)\rightarrow D(S)$
with kernel $\P$, defined in (\ref{eq-fourier-mukai-functor}), 
was studied in \cite{bridgeland-elliptic-surfaces}.
Set $\Phi:=[1]\circ \Phi^\P$, where $[1]$ is the shift auto-equivalence. 
Then $\Phi(\StructureSheaf{S})$ is represented by a line-bundle on 
the zero section.
Replacing $\P$ by $\P\otimes p_2^*\StructureSheaf{\PP^1}(i)$, 
for some $i\in \Integers$, we may assume that 
$\chi(\Phi(\StructureSheaf{S}))=1$. Set $e_1:=(1,0,0)$ and $e_2:=(0,0,1)$. 

\begin{new-lemma}
\label{lemma-matrix-of-phi}
The Hodge-isometry $\phi$, associated to $\Phi$, acts by $-1$ on the 
sublattice of $H^*(S,\Integers)$ orthogonal to 
$\Lambda:={\rm span}_\Integers\{e_1,\sigma,f,e_2\}$ and the matrix of
the restriction of $\phi$ to $\Lambda$, in the given basis, is
\begin{scriptsize}
$
\left(
\begin{array}{cccc}
0 & -1 & 0 & 0
\\
1 & 0 & 0 & 0
\\
1 & -1 & 0 & -1
\\ 
1 & -1 & 1 & 0
\end{array}
\right).
$
\end{scriptsize}
\end{new-lemma}
{\bf Proof:}
When the Picard lattice of $S$ is spanned by $\sigma$ and $f$, then
$\Lambda^\perp$ is the transcendental lattice in $H^2(S,\Integers)$.
In this case $\Lambda^\perp$ is $\phi$ invariant.
For a generic such $S$, the only Hodge isometry of $\Lambda^\perp$
is $1$ or $-1$. Hence, $\phi$ acts on $\Lambda^\perp$ by $1$ or $-1$. 
The latter conclusion follows, 
without the assumption on the Picard number, by a standard deformation 
argument. Columns $1$, $3$, and $4$ of the above matrix were calculated in 
equation (3.18) in \cite{yoshioka-abelian-surface}. 
It is also shown there, that $\phi$ has a $-1$ eigenvector in $\Lambda^\perp$
(the class $\sigma-\tau+d(\tau)f$ in the notation of equation (3.18) in
\cite{yoshioka-abelian-surface}).
The second column of the matrix is determined by the other three, since
$\phi$ is an isometry. 
\EndProof

$\Phi$ induces an isomorphism 
%(Definition \ref{def-equivalence-functor-induces-an-isomorphism})
\begin{equation}
\label{eq-isomorphism-of-moduli-spaces-due-to-yoshioka}
\M_H(1,0,1-n) \ \ \ \longrightarrow \ \ \
\M_H(0,\sigma + nf,1),
\end{equation}
%for every positive integer $\ell$, 
by Theorem 3.15 in \cite{yoshioka-abelian-surface}.
Set $w:=(0,\sigma + nf,1)$. 
Consequently, $\phi$ is a morphism from 
$(S,v,H)$ to $(S,w,H)$ in the groupoid $\G_{Mukai}$ 
given in (\ref{eq-G-Mukai}). 
The homomorphism $\gamma_{\phi,v}$ is induced by the isomorphism
(\ref{eq-isomorphism-of-moduli-spaces-due-to-yoshioka}) of the
moduli spaces, by Lemma \ref{lemma-G-Mukai-is-in-G-OK}.
Lemma \ref{lemma-matrix-of-phi} yields the equality
$\phi(1,\beta-f,n-1)=(0,\alpha,0)$, with  $\alpha:=\sigma+(2-n)f-\beta$.
Note, that $(\alpha,\alpha)=-2$ and $(0,\alpha,0)$ is orthogonal to $w$.

We claim that $\gamma_{\rho,w}$ is a monodromy operator and 
$(\gamma_{\rho,w}\otimes \rho)(u_w)=u_w$. 
When  $n\geq 3$, then $c_1(w)=\sigma+nf$ is ample and 
the claim follows 
from Theorem \ref{thm-conj-holds-for-hilbert-schemes-and-isometries-of-H-2}.
%part \ref{prop-item-stabilizer-of-ample-lb-consists-of-monodromies} 
%of Proposition \ref{prop-torelli-and-monodromy}. 
When $n=2$,  $c_1(w)$ is not ample. Nevertheless, $c_1(w)^2=2$ 
and $c_1(w)$ is effective. These conditions imply that 
the object $(S,w,H)$ is deformation equivalent to an object
$(S',w',H')$, with $c_1(w')^2=2$, where $c_1(w')$ is ample and 
$\Pic(S')$ is cyclic generated by $c_1(w')$
(this is proven by Yoshioka in the unpublished note proving Theorem
\ref{thm-irreducibility} in the case $\rank(v)=0$ and $c_1(v)$ not ample).
Let $\G_{mon}$ be the groupoid 
given in Definition \ref{def-deformation-equivalent-triples}.
Choose any morphism $\eta$ of $\Hom_{\G_{mon}}((S,w,H),(S',w',H'))$.
The conjugate $\eta\circ\rho\circ\eta^{-1}$ of $\rho$
is an element of $\Aut_{\G_{mon}}(S',w',H')$, 
by Theorem \ref{thm-conj-holds-for-hilbert-schemes-and-isometries-of-H-2}.
We conclude that $\rho$ is an element of $\Aut_{\G_{mon}}(S,w,H)$.
Lemma \ref{lemma-the-class-of-f-is-gamma-E1-E2} implies that 
$\gamma_{\rho,w}$ is a monodromy operator in 
$H^*(\M_H(w),\Integers)_{free}$ as claimed. 
 
The equality $g=\phi^{-1}\circ \rho\circ \phi$ implies 
the equality
$\gamma_{g,v}=\gamma_{\phi,v}^{-1}\circ \gamma_{\rho,w}\circ \gamma_{\phi,v}$,
by Lemma \ref{lemma-conjecture-holds-for-compositions}.
Hence, $\gamma_{g,v}$ is a monodromy operator
and $(\gamma_{g,v}\otimes g)(u_v)=u_v$. 
This completes the proof of Proposition 
\ref{reflection-of-Hilbert-schemes-with-respect-to-topological-lb}.
\EndProof

%****************************************************************
% Reflections by spherical objects
%****************************************************************
\subsection{Reflections by spherical objects}
\label{sec-reflections-by-spherical-objects}

Fourier-Mukai functors, which induce $-2$ reflections of the Mukai lattice,
often arise from spherical 
objects on the K3 surface $S$.

\begin{defi}
{\rm
An object $x$ of $D(S)$ is {\em spherical}, if
it satisfies the following conditions
\begin{eqnarray*}
& \Hom^0(x,x) & \ \cong \ \ComplexNumbers \ \cong \ \Hom^2(x,x),
\\
& \Hom^i(x,x) & \ = \ 0, \ \mbox{for} \ i\notin \{0,2\}.
\end{eqnarray*}
}
\end{defi}

The Mukai vector $v_0$, of a spherical object, 
is a $-2$ vector $(v_0,v_0)=-2$. 
Examples of spherical objects include sheaves, which are both 
simple and rigid. 
In particular, any line bundle on $S$ is a spherical object. 
In addition, if $\iota:\PP^1\hookrightarrow S$ is an embedding, 
then the pushforward $\iota_*\StructureSheaf{\PP^1}(n)$,
of any line bundle on $\PP^1$, is a spherical object. 

Given a complex $E_\bullet$ of locally free sheaves, 
the cone $C(E_\bullet)$ of $E_\bullet$ is the object in
$D(S\times S)$ represented by the complex
\begin{equation}
\label{eq-cone-complex}
\pi_1^*(E_\bullet^\vee)\stackrel{L}{\otimes} \pi_2^*(E_\bullet) \ \ 
\LongRightArrowOf{ev} \ \ \Delta_*\StructureSheaf{S},
\end{equation}
with $\Delta_*(\StructureSheaf{S})$ in degree $0$. 
Above, $\Delta:S\hookrightarrow S\times S$ is the diagonal embedding
and $ev$ is the natural evaluation homomorphism. 
If we take, for example, $E_\bullet$ to be the trivial line bundle,
then $C(\StructureSheaf{S})$ is the complex
\begin{equation}
\label{eq-cone-complex-of-structure-sheaf}
\StructureSheaf{S\times S} \ \ \longrightarrow \ \ \Delta_*\StructureSheaf{S}.
\end{equation}

\begin{thm}
\label{thm-reflection-by-a-sherical-object}
\cite{seidel-thomas}
If the complex $E_\bullet$ represents a spherical object,
then the Fourier-Mukai functor $\Phi_{S\rightarrow S}^{C(E_\bullet)}$
is an auto-equivalence of $D(S)$. 
\end{thm}

We will refer to the functor $\Phi_{S\rightarrow S}^{C(E_\bullet)}$
as the {\em reflection functor with respect to $E_\bullet$}. 
On the level of K-theory, it sends the class of a sheaf $F$, to 
\[
F-\chi(E_\bullet,F)\cdot E_\bullet,
\]
where $\chi(E_\bullet,F) := \sum (-1)^i \dim\Hom^i(E_\bullet,F)$. 
On the level of the Mukai lattice, it induces the reflection
$\tau_{v_0}$ with respect to the Mukai vector $v_0$ of $E_\bullet$
\[
\tau_{v_0}(v) \ \ := \ \ v + (v_0,v)\cdot v_0. 
\]
Yoshioka made extensive use of reflection functors in his study of
the geometry of moduli spaces of sheaves on K3 and abelian surfaces
\cite{yoshioka-examples-of-reflections,yoshioka-irreducibility,
yoshioka-abelian-surface}.

%\begin{rem}
%{\rm 
%We can interpret the construction of the cone complex
%(\ref{eq-cone-complex}) in the spirit of our interpretation of 
%Theorem \ref{thm-action} 
%(see section \ref{sec-conceptual-interpretation-of-the-formula-gamma-g}).
%Let us start with the case, where the spherical object is the 
%trivial line bundle. 
%Observe, that the complex 
%(\ref{eq-cone-complex-of-structure-sheaf})
%is quasi-isomorphic to $I_\Delta[1]$, the shift of the ideal sheaf of the 
%diagonal. 
%Moreover, the reflection $\tau_{v_0}(0,0,1)=(-1,0,0)$, of the Mukai vector 
%$(0,0,1)$ of the sky-scraper sheaf, 
%is $-1$ times the Mukai vector $(1,0,0)$ of the ideal sheaf of a point. 
%The sheaf $I_\Delta$ is the universal sheaf over $S\times \M(1,0,0)$, 
%and 
%$I_\Delta[1]$ is a universal sheaf over $S\times \M(-1,0,0)$,
%according to our convention (Definition \ref{def-M-minus-v}). 
%
%More generally, let $w:=\tau_{v_0}(0,0,1)$. 
%We can view the cone complex (\ref{eq-cone-complex}) 
%as a substitute for the pullback of the universal sheaf 
%$\E_w$ over $S\times \M(w)$, via 
%$id\times \kappa:S\times S\rightarrow S\times \M(w)$, where 
%$\kappa:S\rightarrow \M(w)$ is the isomorphism lifting the Hodge isometry 
%\[
%H^2(S,\Integers) \ \cong \ (0,0,1)^\perp/\mbox{span}\{(0,0,1)\} 
%\ \ \LongRightArrowOf{\tau_{v_0}} \ \ 
%w^\perp/\mbox{span}\{w\} \ \cong \ H^2(\M(w),\Integers). 
%\]
%In general, the isomorphism $\kappa$ need not exist, but the
%Hodge isometry does determine a composition of an isomorphism,
%with reflections by $-2$ curves (The Strong Torelli Theorem).  
%}
%\end{rem}

\begin{example}
\label{example-relative-fourir-mukai-transform-wrt-E-0}
{\rm
Assume, that the spherical object $E_\bullet$ is a simple sheaf $E_0$.
Choose a moduli space $\M(v)$ of sheaves on $S$ with 
a universal sheaf  $\E_v$ over $\M(v)\times S$. 
Then the relative Fourier-Mukai transform of $\E_v$ 
(Definition \ref{def-relative-fourier-mukai-transform}), 
with respect to the reflection functor of $E_0$, determines the class  
in K-theory, which is represented also by 
\[
\tau(\E_v) \ \ := \ \ 
\E_v \ - \ 
\pi_2^*E_0\stackrel{L}{\otimes}  
\pi_1^!\pi_{1_!}(\pi_2^*E_0^\vee\stackrel{L}{\otimes}\E_v), 
\] 
where 
%$q_i$ is the projection from $\M(v)\times S\times S$ to the
%$i$-th factor and 
$\pi_i$ is the projection
from $\M(v)\times S$. 
When $E_0$ is the trivial line bundle, we have
\[
\tau(\E_v) \ \ := \ \ 
\E_v \ - \ 
\pi_1^!R_{\pi_{1,!}}\E_v.
\]
}
\end{example}

Assume, that $v=(r,\LB,-r)$, $r\geq 1$, 
$\LB$ satisfies condition
\ref{condition-minimality} with respect to the ample line-bundle $H$, 
and $\mbox{gcd}(r,\deg{\LB})=1$
(so that there is a universal family $\E_v$ over $\M(v)\times S$). 

\begin{condition}
\label{condition-minimality}
\begin{enumerate}
\item 
\label{cond-part-minimality}
$\LB$ is an effective Cartier divisor with {\em minimal}
degree in the sense that the subgroup $H\cdot \Pic(S)$ of $H^4(S,\Integers)$
is generated by $c_1(H)\cdot c_1(\LB)$. 
In particular, all curves in the linear system
$\linsys{\LB}$ are reduced and irreducible. 
\item 
The base locus of 
$\linsys{\LB}$ is either empty or zero-dimensional. 
\item \label{cond-part-generic-curve-is-smooth}
The generic curve in $\linsys{\LB}$ is smooth. 
\item
$H^1(S,\LB)=0$.
\end{enumerate}
\end{condition}

The condition implies, that the polarization $H$ is $v$-suitable
(Definition \ref{def-v-suitable}), by \cite{markman-reflections} Lemma 3.5. 
Let $v_0=(1,0,1)$ be the Mukai vector of $\StructureSheaf{S}$.
%Note, that $v=\tau(v)$ and the Brill-Noether locus $\M(v)^1$ is a divisor. 
Denote by $\M(v)^t$, $t\geq 0$, the subvariety of $\M(v)$ consisting of 
sheaves $F$ with $h^1(F)\geq t$. Then $\M(v)\setminus\M(v)^1$ 
is a dense Zariski
open subset and $\M(v)^1$ is a divisor, whose
class is $\theta(-v_0)$, where $\theta$ is given in 
(\ref{eq-theta-v-from-v-perp}), by  \cite{markman-part-two} Lemma 4.11.
Serre's Duality yields the isomorphism
$H^1(S,F)\cong \Ext^1(\StructureSheaf{S},F)^*$. Given $F\in \M_H(v)^t$,
we get a tautological extension 
\[
0\rightarrow H^1(S,F)\otimes \StructureSheaf{S}\rightarrow E(F) \rightarrow F
\rightarrow 0,
\]
and $E(F)$ is $H$-stable with $h^1(E(F))=0$ (\cite{markman-reflections}
Lemma 3.7). Set $w_t:=(r+t,\LB,-r+t)$. We get a regular morphism
\[
\M_H(v)^t\setminus \M_H(v)^{t+1} \ \ \ \longrightarrow \ \ \ 
\M_H(w_t)\setminus \M_H(w_t)^1,
\]
which is a Grassmannian fibration with $G(t,2t)$ fibers
\cite{markman-reflections}. Note, that $\dim \M(w_t)=\dim\M(v)-2t^2$ and
thus $\M(v)^t$ is non-empty, if and only if $t^2\leq (v,v)/2+1$.
We conclude, that the subset
\[
\Z \ \ := \ \
\{(F_1,F_2) \ : \ h^1(F_1)=h^1(F_2) \ \mbox{and} \
E(F_1)\cong E(F_2)\}
\]
is of pure dimension equal to $\dim(\M(v))$. $\Z$ is in fact a Zariski-closed
subset, which we endow with the reduced subscheme structure, and each 
of its irreducible components is the closure of the fiber square of
$\M_H(v)^t\setminus \M_H(v)^{t+1}$ over $\M_H(w_t)$, $t\geq 0$
\cite{markman-part-two}.

The Hilbert scheme $S^{[n]}$ is a special case of 
$\M_H(v)$, when $r=1$, $c_1(\LB)^2=2n-4$, $n\geq 1$,
and the following choices are made to satisfy 
Condition \ref{condition-minimality}:
(a) $n\geq 3$ and $H=\LB$ is a generator of
$\Pic(S)$. (b) $n=2$, $\varphi:S\rightarrow \PP^1$ is an elliptic 
$K3$ with reduced and irreducible fibers, 
$\LB=\varphi^*\StructureSheaf{\PP^1}(1)$, and 
$H=\StructureSheaf{S}(C)\otimes \varphi^*\StructureSheaf{\PP^1}(k)$,
where the curve $C$ is either a section of $\varphi$ or
a multi-section of minimal degree, and $k>>0$.
(c) $n=1$,  $\LB=\StructureSheaf{S}(\Sigma)$, where $\Sigma\subset S$ 
is a smooth rational curve, and $H$ restricts to $\Sigma$ 
as a generator of the image of $\Pic(S)$ in $\Pic(\Sigma)$. 
The Brill-Noether loci $(S^{[n]})^t:=\M(v)^t$ parametrize ideal sheaves 
$I_A$ of length $n$ subschemes $A\subset S$, such that 
$h^1(I_A\otimes\LB)\geq t$. These loci 
admit a geometric description in terms of the morphism
$\varphi:S\rightarrow \linsys{\LB}^*\cong\PP^{n-1}$,
which is an embedding when $n\geq 4$. When $n\geq 2$,
the stratum $(S^{[n]})^t\setminus (S^{[n]})^{t+1}$ 
consists of ideal sheaves $I_A$, such that
$\varphi(A)$ spans a $\PP^{n-1-t}$. When $n=1$, then $(S^{[1]})^1=\Sigma$
and
$\Z\subset S\times S$ is the union $\Delta\cup [\Sigma\times \Sigma]$, 
where $\Delta$ is the diagonal.

The following Theorem is proven in 
%part II of the paper
\cite{markman-part-two}. 

%***********************************
% Theorem: Stratified -2 reflection:
%***********************************
\begin{thm}
\label{thm-class-of-correspondence-in-stratified-elementary-trans}
\begin{enumerate}
\item
The homomorphism 
$\Z_*:H^*(\M_H(v),\Integers)\rightarrow H^*(\M_H(v),\Integers)$,
induced by the self correspondence $\Z$, is the local monodromy operator
inducing the reflection of $H^2(\M_H(v),\Integers)$ with respect to the
$-2$ class of $\M_H(v)^1$ (Definition \ref{def-local-monodromy-operator}).
%, providing an affirmative answer to Question
%\ref{question-local-monodromy} in this case.
\item
\label{thm-item-Z-maps-universal-class-to-such}
Set 
$\Z_!:=\delta_{2,!}\delta_1^!:K_{alg}(\M_H(v))\rightarrow K_{alg}(\M_H(v))$,
where  $\delta_i:\Z\rightarrow \M_H(v)$, $i=1,2$, are the two projections. 
Then 
\[
\Z_!(\E_v-\pi_1^!\pi_{1,!}\E_v) \equiv \E_v(\M_H(v)^1\times S).
\] 
\item
The isomorphism $\Z_*$ and the reflection $\tau_{v_0}$
of $K_{top}(S)$ with respect to $\StructureSheaf{S}$ satisfy equation
(\ref{eq-characterization-of-gamma-g}) with
$\ell=\StructureSheaf{\M_H(v)}(\M_H(v)^1)$.
Consequently, $\Z$ is Poincare-Dual to the cohomology class 
$\gamma(\tau_{v_0},v)$, 
the homomorphism $\gamma_{\tau_{v_0},v}$ is the monodromy operator $\Z_*$,
and $(\gamma_{\tau_{v_0},v}\otimes \tau_{v_0})(u_v)=u_v$.
%the statements of 
%Theorem \ref{thm-trancendental-reflections} hold 
%for the automorphism $\tau_{v_0}$ of the object
%$(S,v,H)$ of the groupoid $\G$ given in (\ref{eq-cohomological-groupoid}). 
\end{enumerate}
\end{thm}

\subsection{The Hodge isometry $-id$}
\label{sec-minus-id}

The isometry $-id$ fixes only the zero vector. 
Moreover, if a non-zero Mukai vector $v$ represents a sheaf,
then $-v$ does not. Nevertheless, it is convenient to 
consider $-v$ as representing shifted complexes $F[1]$, where $F$ is a
coherent sheaf with Mukai vector $v$. 

\begin{defi}
\label{def-M-minus-v}
The moduli spaces $\M_\LB(v)$ and $\M_\LB(-v)$ are defined to be the same
projective variety. A universal ``sheaf'' $\E_{-v}$ over $\M_\LB(-v)\times S$
is the K-theoretic class of the complex
\[
\E_{-v} \ := \ \E_v[d], \ \ d \ \mbox{odd},
\]
where $\E_v$ is a universal sheaf over $\M_\LB(v)\times S$. 
We set $u_{-v}:=-u_v$.
\end{defi}

%Definition \ref{def-M-minus-v} enables us to lift the Hodge 
%isometry $-id$ to a {\em ring} isomorphism 
Let 
\[
\gamma_{-id} : H^*(\M(v),\Integers) \rightarrow H^*(\M(-v),\Integers)
\]
be the ring isomorphism, induced by the natural isomorphism 
$\M(v)\cong \M(-v)$. 
%Then the isomorphism 
%\[
%(\gamma_{-id} \otimes -id)\ : \ H^*(\M(v)\times S,\Integers) \ \rightarrow \ 
%H^*(\M(-v)\times S,\Integers)
%\]
%is $-id$. 
The identity 
$ch(\E_v[1])=-ch(\E_v)$ translates to the equality
\[
(\gamma_{-id} \otimes -id)_*(ch(\E_v)) \ \ = \ \ ch(\E_{-v}).
\]
Moreover, Lemma 
\ref{lemma-Chow-theoretic-formula-for-gamma} holds, 
with $\Phi$ being the shift auto-equivalence. 
\subsection{Stratified reflections with respect to $+2$ vectors}
\label{sec-reflections-with-respect-to-plus-2-vectors}

Let $u_0=(1,0,-1)$ be the $+2$ Mukai vector of the ideal sheaf of
two points. 
Denote by $\sigma_{u_0}$ the reflection 
of the Mukai lattice with respect to $u_0$
\begin{equation}
\label{eq-sigma}
\sigma_{u_0}(w) \ = \ w - (w,u_0)u_0.
\end{equation}
Let $v_0$ be the $-2$ vector $(1,0,1)$ of the trivial line-bundle. 
The corresponding reflections, $\sigma_{u_0}$ and $\tau_{v_0}$,
commute and satisfy the relation 
\[
-(\sigma_{u_0}\circ \tau_{v_0}) \ = \ D. 
\]
In particular, $\sigma_{u_0} = -(D \circ \tau_{v_0})$. 

\begin{thm}
\label{thm+2-reflection-induces-an-isomorphism}
Let $(S,(r,\LB,s),H)$ be an object of the groupoid $\G$, 
where the polarization $H$ and the line-bundle $\LB$
satisfy Condition \ref{condition-minimality}.
Assume that $r\geq 0$, $s\geq 0$, and
\begin{equation}
\label{eq-inequality-implying-empty-BN-locus}
c_1(\LB)^2 \ \ \ < \ \ \ 2[r+s+rs].
\end{equation}
Then the composition $\Phi:=([1]\circ \Psi)^\vee$ 
of the reflection $\Psi:D(S)\rightarrow D(S)$ with respect to the
spherical object $\StructureSheaf{S}$ 
(constructed in Theorem \ref{thm-reflection-by-a-sherical-object}),
the shift auto-equivalence, and the functor of taking dual, 
induces a regular isomorphism
\[
f_{(r,\LB,s)} \ : \ \M_H(r,\LB,s) \ \ \ \rightarrow \ \ \ \M_H(s,\LB,r).
\]
On the level of cohomology rings 
$f_{(r,\LB,s)_*}=D_{\M}\circ \gamma_{\sigma_{u_0},(r,\LB,s)}$. 
Furthermore, $f_{(r,\LB,s)}\circ f_{(s,\LB,r)}=id$.
\end{thm}

\noindent
{\bf Proof:}
The isomorphism $f_{(r,\LB,s)}$ is constructed in
\cite{markman-reflections} Theorem 3.21 (use also the inequality
(\ref{eq-inequality-implying-empty-BN-locus}) above and
\cite{markman-reflections} Corollary 3.16 to conclude the emptyness of the 
Brill-Noether loci).
Related results were proven independently in 
\cite{yoshioka-examples-of-reflections}.
\EndProof

If $r=s$ then $f_{(r,\LB,r)}$ is a regular involution of the moduli space
$\M_H(r,\LB,r)$.
When the inequality (\ref{eq-inequality-implying-empty-BN-locus}) 
does not hold, the functor induces a non-regular birational isomorphism. 
The following theorem is proven in 
%part II of the paper
\cite{markman-part-two}.

\begin{thm}
\label{thm-reflection-sigma-satisfies-main-conj}
Set $v:=(1,\LB,1)$. 
The statements of 
Theorem \ref{thm-trancendental-reflections} 
holds for the automorphism $\sigma_{u_0}$ of the object 
$(S,v,H)$ of the groupoid
(\ref{eq-cohomological-groupoid}),
where the polarization $H$ and
the line-bundle $\LB$ satisfy condition 
\ref{condition-minimality}. 
Moreover, the class 
of the composition $D_\M\circ\gamma_{\sigma_{u_0},v}$ 
is the local monodromy operator 
restricting as 
$-\theta_v\circ \sigma_{u_0}\circ \theta_v^{-1}$ on
$H^2(\M_H(v),\Integers)$
(Definition \ref{def-local-monodromy-operator}).
\end{thm}

Note, that $\M(1,\LB,1)$ is isomorphic to $S^{[n]}$, where
$n:=\frac{1}{2}c_1(\LB)^2$. 
Associated to the reflection $\sigma_{u_0}$ is
a Lagrangian correspondence $\Z \subset \M(v)\times \M(v)$ 
(see \cite{markman-reflections} section 3.2). 
The class $\monrep(\sigma_{u_0},v)$, given in (\ref{eq-class-mon-g-v}),
%$(1\otimes D_\M)\gamma(\sigma_{u_0},v)$
is Poincare-dual to the class of $\Z$. This fact is proven in the course 
of the proof of Theorem 
\ref{thm-reflection-sigma-satisfies-main-conj}. 
$\Z$ is the graph of a regular involution, for $n=1$ or $2$
(Theorem \ref{thm+2-reflection-induces-an-isomorphism}). 
When $c_1(\LB)^2=2$ and $n=1$, then $\Z$ is the graph of the Galois
involution of the double cover 
$\M(1,\LB,1)\cong S\rightarrow \linsys{\LB}^*\cong\PP^2$. 
$\Z$ is reducible for $n\geq 3$. 
One of its irreducible components is the graph of a birational involution.

%*********************************************************************
% The structure of the stabilizer $\Gamma_v$
%*********************************************************************
\section{Generators for the stabilizer $\Gamma_v$}
\label{sec-structure-of-stabilizer}

Let us study the structure of the stabilizer $\Gamma_v$ of the
Mukai vector $v=(1,0,-m)$ of the Hilbert scheme $S^{[m+1]}$, $m\geq 1$. 
Lemma \ref{lemma-the-index-of-Gamma-v} relates $\Gamma_v$ and 
$O(v^\perp)$. 
Proposition \ref{prop-compare-Gamma-v}
states, that $\Gamma_v$ is generated by the subgroup
$\Gamma_0$, in Definition \ref{def-isometry-groups}, and 
reflections in Mukai vectors of line bundles. 
The character group of $\Gamma_v$ is calculated in Corollary 
\ref{cor-stabilizer-is-generated-by-reflections}.

The signature $(\ell_+,\ell_-)$ of $H^2(S,\Integers)$ is 
$(3,19)$, however we will use in this section only the 
inequalities $\ell_+\geq 3$ and $\ell_-\geq 3$, so that the results
hold in case $S$ is an abelian surface (though $v$ is then the Mukai vector of
$S^{[m]}$).

\begin{new-lemma}
\label{lemma-primitive-embeddings-of-rank-2-lattices}
Let $\Lambda$ be an even unimodular lattice 
of signature $(\ell_+,\ell_-)$, satisfying $\ell_+\geq 3$ and $\ell_-\geq 3$. 
1) Let $M=
\left[
\begin{array}{cc}
2a & b \\
b & 2d
\end{array}
\right]$ 
be a symmetric matrix with $a,b, d\in \Integers$, and  
$\lambda_1\in \Lambda$  a primitive element with 
$(\lambda_1,\lambda_1)=2a$.
Then there exists a primitive rank $2$ sublattice $\Sigma\subset \Lambda$, 
containing $\lambda_1$, and an element $\lambda_2\in \Sigma$, such that 
$\{\lambda_1,\lambda_2\}$ is a basis for $\Sigma$, and 
$M$ is the matrix of the bilinear form of 
$\Sigma$ is this basis. 

2) Assume that $\rank(M)=2$. Let $\{\lambda_1',\lambda_2'\}$ 
be a basis for another primitive sublattice of $\Lambda$,
having the same matrix $M$. 
Then there exists an isometry $g$ of $\Lambda$ satisfying 
$g(\lambda_i)=\lambda_i'$, $i=1,2$.

3) The kernel $O(\Lambda)^{cov}$ of the orientation character
(\ref{eq-top-homology-character}) acts transitively on 
the set of primitive integral classes $A\in \Lambda$, 
with fixed ``squared-length'' $(A,A)$. 
\end{new-lemma}

\noindent
{\bf Proof:}
1) Let $\Sigma'$ be the rank $2$ lattice with basis $\{e_1,e_2\}$,
whose bilinear form is given by the matrix $M$.
Then there exists a primitive embedding 
$\eta: \Sigma'\hookrightarrow \Lambda$, by Theorem 1.14.4 of
\cite{nikulin}, when $\rank(M)=2$, and by Proposition
1.17.1 of \cite{nikulin}, when $\rank(M)<2$. 
Furthermore, there exists an isometry $g$ of $\Lambda$, satisfying 
$g(\eta(e_1))=\lambda_1$, again by Theorem 1.14.4 of
\cite{nikulin}. Take $\Sigma=(g\circ \eta)(\Sigma')$ and
$\lambda_2=g(\eta(e_2))$. 

\noindent
2) Follows from Theorem 1.14.4 of \cite{nikulin}. 

\noindent
3) We already explained the transitivity of the $O(\Lambda)$-action. 
The transitivity of $O(\Lambda)^{cov}$ would follow from that of $O(\Lambda)$,
once we prove that the stabilizer of $A$ in
$O(\Lambda)$ contains an orientation reversing isometry.
It suffices to prove that $A^\perp$ contains a class $\lambda$ with 
$(\lambda,\lambda)=2$,
because the reflection with respect to $\lambda$ would be 
orientation reversing (see section \ref{sec-orientation-character}).
The existence of such $\lambda$ follows from part 1.
\EndProof

A $-2$ vector $v_0$, whose reflection $\tau_{v_0}$ is in the 
stabilizer of $v$, must be in $v^\perp$. Consequently, $v_0$ 
has the form
\[
v_0 \ = \ (r,\LB,rm), 
\]
where $\LB$ has self-intersection $2r^2m-2$. 
%and, if effective, its genus is $r^2m$. 
The lattice $v^\perp$ is the direct sum 
$H^2(S,\Integers)\oplus {\rm span}_\Integers\{(1,0,m)\}$. 
It is even, but {\em not} unimodular; the primitive vector 
$(1,0,m)$ sends $v^\perp$, via the Mukai pairing, to 
$2m\cdot \Integers$.
Observe, that the stabilizer $\Gamma_v$ embeds in 
the isometry group of $v^\perp$, but the latter  is larger in general. 

The quotient $(v^\perp)^*/v^\perp$ is isomorphic to $\Integers/2m$.
Embed $(v^\perp)^*$ in $v^\perp\otimes_\Integers\RationalNumbers$.
We get the well-defined quadratic form on $(v^\perp)^*/v^\perp$
\begin{eqnarray*}
q \ : \ (v^\perp)^*/v^\perp & \rightarrow & \RationalNumbers/2\Integers, 
\ \ \mbox{determined by}
\\
\frac{(1,0,m)}{2m} & \mapsto & -1/2m, \ \ \mbox{and satisfying}
\\
q\left(\frac{(r,0,rm)}{2m}\right) & = & -r^2/2m.
\end{eqnarray*}
The isometry group $O[(v^\perp)^*/v^\perp]$ is isomorphic to 
$(\Integers/2\Integers)^\rho$, where the Euler number $\rho$ of $m$ is 
the number of distinct primes
in the prime factorization $m=p_1^{e_1}p_2^{e_2}\cdots p_\rho^{e_\rho}$ 
(and $e_i$ are positive integers) (see \cite{oguiso}).
%the subgroup $U[2]$ of square roots of $1$ in the multiplicative group $U$ 
%of units in $\Integers/(2m)$. 
We have the natural homomorphism 
\begin{equation}
\label{eq-projection-to-finite-isometry-group}
O(v^\perp) \ \rightarrow \ O[(v^\perp)^*/v^\perp].
\end{equation}

\begin{new-lemma}
\label{lemma-O-v-perp-surjects-on-finite-orthogonal-group}
%\begin{enumerate}
%\item
%\label{lemma-item-surjectivity-of-O-v-perp}
The  homomorphism (\ref{eq-projection-to-finite-isometry-group})
is surjective.
%\item 
%\label{lemma-item-image-of-1-0-n-1}
%((???) not needed)
%If $\varphi\in O(v^\perp)$, 
%then $\varphi(1,0,n-1)=(r,c,r(n-1))$ where the class $c$ is divisible by
%$2n-2$ and $c\cdot c= (r^2-1)(2n-2)$. 
%%If $c\neq 0$, and we write $c=k\cdot c'$ where 
%%$c'$ is an indivisible class in $H^2(S,\Integers)$, 
%then $2n-2$ divides $k$. 
%Consequently, $r^2 \equiv 1 \ (\mbox{mod} \ 2n-2)$. 
%\end{enumerate}
\end{new-lemma}

\noindent
{\bf Proof:} Follows from  \cite{nikulin} Theorem 1.14.2.

\begin{new-lemma}
\label{lemma-the-index-of-Gamma-v}
$\Gamma_v $ is the kernel of the natural surjective homomorphism 
(\ref{eq-projection-to-finite-isometry-group}). 
In particular, $\Gamma_v=O(v^\perp)$ if and only if $m=1$.
If $m$ is a prime power, then
$\Gamma_v$ is a subgroup of $O(v^\perp)$ of index 2 and
$O(v^\perp)= \Gamma_v\cup (-id_{v^\perp})\cdot \Gamma_v$.
\end{new-lemma}

\noindent
{\bf Proof:} We prove only that $\Gamma_v$ is the kernel $N$ of
(\ref{eq-projection-to-finite-isometry-group}). The rest follows 
from the description of $O[(v^\perp)^*/v^\perp]$.

There is a natural isomorphism between the residue groups
${\rm span}_\Integers\{v\}^*/{\rm span}_\Integers\{v\}$ 
 and $(v^\perp)^*/v^\perp$, each being isomorphic to 
$H^*(S,\Integers)/[\Integers v+(v^\perp)]$.
An element $\varphi\in O(v^\perp)$ extends to an element of $\Gamma$, 
if and only if the image $\bar{\varphi}$ of $\varphi$ in
$O[(v^\perp)^*/v^\perp]$ is 
contained in the image of $O[{\rm span}_\Integers\{v\}]$; i.e., 
if and only if $\bar{\varphi}$ is the identity, or multiplication by $-1$
(Proposition 1.5.1 in \cite{nikulin}).
The image  $\bar{\varphi}$  of $\varphi$ is determined by the image
$\varphi(1,0,m)$ in $(v^\perp)^*/v^\perp$.
Any isometry $\varphi\in O(v^\perp)$ satisfies 
$\varphi(1,0,m) \ = \ (r,c,rm)$
and $c\cdot c= 2m(r^2-1)$. 
We see, that an
isometry $\varphi\in O(v^\perp)$ extends to an element of $\Gamma$, 
if and only if $r \equiv \pm 1 \ (\mbox{mod} \ 2m)$. 
Any element of $\Gamma$, which leaves
$v^\perp$ invariant, must send $v$ to $\pm v$. 
Hence, $N$ is contained in the subgroup $E$ generated by $\Gamma_v$ and 
$-id_{v^\perp}$. 
If $m=1$, then $N=O(v^\perp)$, 
$v$ is a $+2$ vector, and the reflection $\sigma_v$ sends $v$ 
to $-v$ and acts as the identity on $v^\perp$. 
Hence, if $\varphi\in\Gamma$ extends an isometry in 
$O(v^\perp)$, then either 
$\varphi$ or $\sigma_v\circ\varphi$ is an extension in $\Gamma_v$. 
Thus $\Gamma_v=O(v^\perp)$.

If $m>1$, then both $N$ and $\Gamma_v$ have index $2$ in $E$, so
it remains to show that $\Gamma_v$ is contained in $N$, i.e., 
that the additional condition $\varphi(v)=v$
implies, that $r \equiv 1 \ (\mbox{mod} \ 2m)$. 
This is 
%clear if $m=1$ and for $m\geq 2$ it is 
precisely equation (\ref{eq-r-is-congruent-to-1-mod-2n-2}), 
proven below. 
\EndProof

%Assume $r\equiv 1$ (mod $2n-2$). 
%Extend $\varphi$ to an isometry in $\Gamma$. 
%We claim that $\varphi(v)=v$ or $n=2$. 
%Suppose $\varphi(v)\not=v$. Then $\varphi(v)=-v$. 
%We calculate
%\[
%\varphi(2,0,0)  = \varphi(v+(1,0,n-1)) =
%-v+(r,kc',r(n-1)) =
%(r-1,kc',(r+1)(n-1)). 
%\]
%Thus, $\varphi(1,0,0)$ is divisible by $n-1$. 
%We conclude that $n=2$. 
%
%Assume next that $r\equiv-1$ (mod $2n-2$).
%Extend $\varphi$ to an isometry in $\Gamma$. 
%We claim that $\varphi(v)=-v$ or $n=2$. 
%Suppose $\varphi(v)\not=-v$. Then $\varphi(v)=v$. 
%We calculate
%\[
%\varphi(2,0,0)  = \varphi(v+(1,0,n-1)) =
%v+(r,kc',r(n-1)) =
%(r+1,kc',(r-1)(n-1)). 
%\]
%Hence, $2$ is divisible by $2n-2$ and $n=2$. 
%\EndProof

\medskip
Our next goal is to find a manageable set of generators for $\Gamma_v$
(Proposition \ref{prop-compare-Gamma-v}). 

\begin{new-lemma}
\label{lemma-stabilizer-is-generated-by-reflections-and-G-0}
The stabilizer $\Gamma_v$ is generated by $\Gamma_0$ 
(Definition \ref{def-isometry-groups}) and 
reflections in $-2$ vectors in $v^\perp$. 
\end{new-lemma}

\noindent
{\bf Proof:}
Let $g$ be an isometry in $\Gamma_v$ and 
set $w:=(1,0,m)$. 
If $g(w)=w$, then $g$ belongs to $\Gamma_0$.
If $g(w)= -w$, then $g(2,0,0)=g(v+w)=v-w=(0,0,-2m)$,
so $m=1$, $w$ is a $-2$ vector, and the reflection by $w$ multiplies
$g$ into $\Gamma_0$.
We may thus assume that 
\[
g(w) \ = \ (r,k\LB,rm), \ \ r, k\in \Integers \ \mbox{and} \ \LB 
\ \mbox{a primitive class}.
\]
%Observe first that $r\neq 0$. If $r=0$, then $\LB$ is a primitive 
%class in $H^2(S,\Integers)$. There exist classes $A\in H^2(S,\Integers)$
%such that $a\cdot \LB=1$, because the lattice $H^2(S,\Integers)$ 
%is unimodular.
%But the pairing of $w$ with any vector in $v^\perp$ is divisible by
%$2n-2$. 
We claim that the following equalities hold:
\begin{eqnarray}
\nonumber
k & = & 2mc  \ \mbox{for some integer} \ c,
\\
\label{eq-r-is-congruent-to-1-mod-2n-2}
r-1 & = & 2m\rho \ \mbox{for some integer} \ \rho,
\\
\label{eq-c-rho}
c^2(\LB\cdot\LB) & = & 2\rho + 2m\rho^2,
\end{eqnarray}
and consequently, 
$g(w)  =  (1,0,m) + 2m(\rho,c\LB,m\rho)$.

The equalities $g(1,0,-m)=(1,0,-m)$ and 
$g(1,0,m)=(r,k\LB,rm)$ imply:
\[
g(2,0,0) \ = \ (r+1,k\LB,(r-1)m).
\]
Thus, $r$ is odd, $k$ is even, and
\[
g(0,0,m) \ = \ \left(\frac{r-1}{2},\frac{k}{2}\LB,\frac{m(r+1)}{2}
\right).
\]
In particular, $2m$ divides $r-1$ and $k$. 
The equation (\ref{eq-c-rho}) is equivalent to the statement, that
$(g(w),g(w))=-2m$. 
%It implies that $c^2$ divides $\rho(1+(n-1)\rho)$.
As $w$ is primitive, $c$ and $r=1+2m\rho$ are relatively prime. 

%(Not needed ???$\rightarrow$)
%Next, we show that $c$ and $(2n-2)\rho$ are relatively prime as well. 
%Let 
%\begin{eqnarray*}
%a \ : \  H^2(S,\Integers) & \rightarrow & \Integers \ \ \mbox{and}
%\\
%A  \ : \  H^2(S,\Integers) & \rightarrow & H^2(S,\Integers)
%\end{eqnarray*}
%be the linear homomorphisms, such that $g(0,x,0)=(a(x),A(x),a(x)(n-1))$.
%Since $(0,x,0)$ belongs to $w^\perp$, then 
%\begin{eqnarray*}
%0 & = & \left(g(w),(a(x),A(x),a(x)(n-1))\right) \ = \ 
%\\
%& = & (2n-2)\cdot c\cdot (\LB\cdot A(x)) -a(x)(2n-2)-a(x)(2n-2)^2\rho.
%\end{eqnarray*}
%Thus, $a(x)+a(x)(2n-2)\rho \ = \ c\cdot (\LB\cdot A(x))$. 
%Set $h:=\mbox{gcd}(c,(2n-2)\rho)$. We get
%\begin{eqnarray*}
%a(x) & \equiv & 0 \ (\mbox{mod} \ h) \ \mbox{and}
%\\
%\mbox{rank}(g(w)) & \equiv & 1 \ (\mbox{mod} \ h).
%\end{eqnarray*} 
%But $g$ maps $v^\perp$ {\em onto} itself. Hence, $h=1$ and 
%$c$, $(2n-2)\rho$ are relatively prime.
%($\leftarrow$???  Not needed )

Let $\Gamma_v'$ be the subgroup of $\Gamma_v$ generated by 
$\Gamma_0$ and reflections in $-2$ vectors in $v^\perp$. We need to show, 
that the coset $\Gamma_v'g$ is $\Gamma_v'$. 
Given a $-2$ vector $u=(a,A,am)\in v^\perp$, 
%$A\cdot A=2a^2m-2$, 
we get
\begin{equation}
\label{eq-reflection-of-g-w}
\tau_u(g(w)) \ = \ 
g(w) + 2m\left[
c(A\cdot \LB)-ar
\right]\cdot (a,A,am).
\end{equation}
Since $\mbox{gcd}\{c,r\}=1$, there are integers $a$, $b$ such that
$bc-ar=1$. Choose $u$ with such rank $a$, satisfying  
$A\cdot \LB=b$ and with $c\LB+A$ primitive 
(possible by Lemma \ref{lemma-primitive-embeddings-of-rank-2-lattices}). Then 
the first Chern class of $\tau_u(g(w))$ is $2m(c\LB+A)$ with $c\LB+A$ 
primitive.
Thus, any coset of $\Gamma_v'$ contains a representative $g$ with $c=1$.
We may assume, that $g$ is such a representative.
%We may further assume that $\rho<0$. 

If $\rho=-1$, set $v_0 := (-1,\LB,-m)$. Equation 
(\ref{eq-c-rho}) implies, that $v_0$ is a $-2$ vector.
A straightforward calculation shows, that the composition
$\tau_{v_0}\circ g$ fixes both $v$ and $w$. 
Hence, $\tau_{v_0}\circ g$ comes from an isometry of
$\{v,w\}^\perp$, namely of $H^2(S,\Integers)$. 
Consequently, $\Gamma_v'g=\Gamma_v'$. 

When $\rho\neq-1$,  we choose a $-2$ vector $u$ with 
$a=1$ and $A\in H^2(S,\Integers)$ satisfying 
\begin{eqnarray*}
A\cdot \LB & = & r-\rho-1 \ = \ 2m\rho-\rho, \ \mbox{and}
\\
A\cdot A & = & 2m-2
\end{eqnarray*}
($A$ exists, by Lemma \ref{lemma-primitive-embeddings-of-rank-2-lattices}, 
because $\LB$ is primitive).
Equation (\ref{eq-reflection-of-g-w}) implies that
\[
(\tau_u\circ g)(w)=(1,0,m)+2m(-1,\LB-(\rho+1)A,-m).
\]
The invariant $\rho$ of $\tau_u\circ g$ is $-1$. 
Consequently, $\Gamma_v'\tau_u g= \Gamma_v'$. 
But $\Gamma_v'g=\Gamma_v'\tau_u g$. 
%The case $\rho=-1$ implies that $\tau_u\circ g$ in the
%subgroup generated by $\Gamma_0$ and reflections in $-2$ vectors in 
%$v^\perp$.
%Hence, the same applies to $g$.
This completes the proof of Lemma 
\ref{lemma-stabilizer-is-generated-by-reflections-and-G-0}.
\EndProof

%\medskip
%Let $\Gamma_v^-$ be the normal subgroup of $\Gamma_v$ generated by 
%$\Gamma_0$ and reflections in $-2$ vectors from  $A_-$. 

\begin{new-lemma}
\label{lemma-two-normal-subgroups}
The stabilizing subgroup $\Gamma_v$ is generated by 
$\Gamma_0$ and reflections in
$-2$ vectors $(a,A,am)$ in $v^\perp$, with $A$ a primitive class
in $H^2(S,\Integers)$. 
%The subgroup $\Gamma_v^-$ is, in fact, the whole of $\Gamma_v$
\end{new-lemma}

\noindent
{\bf Proof:} If $g\in \Gamma_v$ and $g(w)\neq -w$, 
simply apply the proof of Lemma
\ref{lemma-stabilizer-is-generated-by-reflections-and-G-0}.
The point is that the reflection
$\tau_{u}$ in (\ref{eq-reflection-of-g-w}) may be chosen with 
$u=(a,A,am)$ and $A$ primitive.
If $g(w)= -w$, then $m=1$ (see the proof of Lemma
\ref{lemma-stabilizer-is-generated-by-reflections-and-G-0}).
Take $v_0=(1,\LB,1)$ with $\LB$ primitive and isotropic. 
Then $v_0$ is a $-2$ vector in $v^\perp$, $\tau_{v_0}(w)=-w$, 
and $\tau_{v_0}g$ belongs to $\Gamma_0$.
\EndProof

\begin{prop} 
\label{prop-compare-Gamma-v}
$\Gamma_v$ is generated by $\Gamma_0$ and reflections in $-2$ vectors 
of the form $(1,-\LB,m)$, with $\LB$ 
a primitive class in $H^{2}(S,\Integers)$ of self-intersection $2m-2$.
\end{prop}

Note that $(1,-\LB,m)$ is 
the Mukai vector of the inverse of a line bundle $\LB$, if
 $c_1(\LB)^2=2m-2$. We will need Lemmas
\ref{lemma-relation-between-reflections-and-concatenation-of-rational-curves} 
and \ref{lemma-existence-of-pairs-of-line-bundles} 
for the proof of Proposition \ref{prop-compare-Gamma-v}.

\begin{new-lemma}
\label
{lemma-relation-between-reflections-and-concatenation-of-rational-curves}
Let $\LB_1$  and $\LB_2$ be two classes in $H^2(S,\Integers)$ satisfying
the following conditions:
\begin{eqnarray}
\label{eq-conditions-for-two-line-bundles}
\LB_1\cdot \LB_1 & = & 2a^2m-2,
\\
\nonumber
\LB_2\cdot \LB_2 & = & 2b^2m-2,   \ \  \mbox{and}
\\
\nonumber
\LB_1\cdot \LB_2 & = & 1 + 2abm. 
\end{eqnarray}
Then the Mukai vectors $v_1:=(a,-\LB_1,am)$, 
$v_2:=(b,-\LB_2,bm)$, and $v_0:=v_1+v_2$
are all $-2$ vectors in $v^\perp$.
The subgroup of $\Gamma_v$, generated by 
$\tau_{v_0}$, $\tau_{v_1}$, and $\tau_{v_2}$, is isomorphic to
$\Sym_3$ and is generated by any two out of the three reflections. 
\end{new-lemma}

\noindent
{\bf Proof:}
The equality $(v_1,v_2)=1$ yields 
\begin{eqnarray*}
\tau_{v_1}(v_2) & = & v_2 + (v_2,v_1)\cdot v_1 = v_0 \ \ \mbox{and}
\\
\tau_{v_2}(v_1) & = & v_1  + (v_2,v_1)\cdot v_2 = v_0.
\end{eqnarray*}
The class $\LB_0 := \LB_1 + \LB_2$ 
has self-intersection $2(a+b)^{2}m-2$ and the vector $v_0:=v_1+v_2$ is a
$-2$ vector, whose reflection $\tau_{v_0}$ satisfies the relations: 
\begin{eqnarray}
\label{eq-relation-between-reflections-and-concatenation-of-rational-curves}
\tau_{v_0} & = & \tau_{v_1}\circ\tau_{v_2}\circ\tau_{v_1} \ \ \mbox{and} 
\\
\nonumber
\tau_{v_0} & = & \tau_{v_2}\circ\tau_{v_1}\circ\tau_{v_2}.
\end{eqnarray}
Consequently, the subgroup of $\Gamma_v$, generated by 
$\tau_{v_0}$, $\tau_{v_1}$, and $\tau_{v_2}$, is isomorphic to
$\Sym_3$.
%and is generated by any two out of the three reflections. 
Observe, that $(v_1,v_2)=(-v_0,v_1)=(-v_0,v_2)=1$. 
The ``triangle'' with vertices $-v_0$, $v_1$, and $v_2$ is an extended Dynkin
diagram of type $A_2$. 
\EndProof

%Note that the conditions (\ref{eq-conditions-for-two-line-bundles})
%are equivalent to
%\begin{eqnarray*}
%\LB_1\cdot \LB_1 & = & 2a^2(n-1)-2,
%\\
%(\LB_1-\LB_2)\cdot(\LB_1-\LB_2) & = & 2(a-b)^2(n-1)-6, \ \ \mbox{and}
%\\
%(\LB_1-\LB_2)\cdot\LB_2 & = & 2(ab-b^2)(n-1)+3.
%\end{eqnarray*}
%In the special case where $a=b$,  we get that $\LB_1-\LB_2$ has degree $-6$
%and $(\LB_1-\LB_2)\cdot \LB_2=3$. 
%Two ways to achieve this set-up are:
%\begin{enumerate}
%\item
%Choose $\LB_2$ as an effective curve of genus $b^2(n-1)$ intersecting three
%disjoint rational curves $\Sigma_i$ each transversally at one point. 
%We then set $\LB_1:=\LB_2+\sum_{i=1}^3\Sigma_i$. 
%\item
%The difference $\LB_1-\LB_2$ could be also equivalent to
%the difference $E_1-E_2$ of two elliptic curves such that
%$E_1\cdot E_2=3$. 
%Choose $\LB_2$ as an effective curve of genus $b^2(n-1)$ whose intersection
%with $E_1-E_2$ is $3$. In particular, when $n=2$ and $a=b=1$, we can choose
%$\LB_i$ to be $E_i$. 
%\end{enumerate}
%
%
%\begin{example}
%{\rm
%Let $v_0=(r,-\LB,r(n-1))$, and assume that 
%$\Sigma$ is a smooth rational curve on $S$ satisfying 
%$\LB\cdot \Sigma=1$.  
%The line bundle $\LB_2 := \LB+\Sigma$ is also  primitive 
%(since $\LB_2\cdot \Sigma=-1$). 
%Set $v_1 := (0,\Sigma,0)$. 
%The Mukai vector $v_2:=v_0-v_1 = (r,-\LB_2,r(n-1))$ 
%induces a reflection $\tau_{v_2}$. The relations
%(\ref{eq-relation-between-reflections-and-concatenation-of-rational-curves}) 
%hold among the reflections $\tau_{v_i}$. 
%}
%\end{example}

\begin{new-lemma}
\label{lemma-existence-of-pairs-of-line-bundles}
\begin{enumerate}
\item
\label{lemma-item-decomposition-of-2-vector}
Given an even rank $r\geq 2$ and a 
class $\LB_0$ of self-intersection $2r^{2}m-2$ (not necessarily primitive), 
we can decompose $\LB_0$ as a sum $\LB_0=\LB_1+\LB_2$ 
of two classes satisfying Condition
(\ref{eq-conditions-for-two-line-bundles}) with $a=b=\frac{r}{2}$.
We can further choose $\LB_1$ and $\LB_2$ primitive. 
\item
\label{lemma-item-increasing-by-one-the-rank-of-a-primitive-2-vector}
Fix a rank $a\geq 1$ and let
$\LB_1$ be a primitive class of self-intersection $2a^{2}m-2$. 
Then there exists a primitive class $\LB_2$ of self-intersection $2m-2$,
such that 
$v_1:=(a,-\LB_1,am)$ and  $v_2:=(1,-\LB_2,m)$ satisfy condition
(\ref{eq-conditions-for-two-line-bundles}). 
%\item
%\label{lemma-item-changing-odd-rank-to-even-rank-2-vector}
%Let $a$ be an odd rank and $\LB_1$ a class  of self-intersection 
%$2a^{2}m-2$. 
%If $\LB_1$ is not divisible by $2$, then 
%there exists an odd rank $b$ and a primitive class $\LB_2$ of the degree 
%$2b^2(n-1)-2$, such that the vectors $v_1:=(a,\LB_1,a(n-1))$ and
%$v_2:=(b,\LB_2,b(n-1))$ satisfy: 
%$(v_1,v_2)$ is odd (and hence $\tau_{v_2}(v_1)=v_1+(v_1,v_2)v_2$ 
%has even rank). 
\end{enumerate}
\end{new-lemma}

\noindent
{\bf Proof:}
Part \ref{lemma-item-decomposition-of-2-vector}) 
Say $\LB_0=i\LB'_0$ with $\LB'_0$ primitive.
Set $d:=(\LB'_0\cdot\LB'_0)/2=r^2m-1/i^2$.
Let $\Sigma\subset H^2(S,\Integers)$ be a rank $2$ primitive sublattice 
containing $\LB'_0$ as well as a class $\LB_1$, such that 
$\{\LB'_0,\LB_1\}$ is a basis for $\Sigma$, and whose bilinear form
has matrix
$\left[
\begin{array}{cc}
2d & id\\
id & (r^2m-4)/2
\end{array}
\right]$
(Lemma \ref{lemma-primitive-embeddings-of-rank-2-lattices}).
Set $\LB_2:=\LB_0-\LB_1$. Then $\LB_1$, $\LB_2$ are primitive and they 
satisfy condition (\ref{eq-conditions-for-two-line-bundles}) with $a=b=r/2$.

Part 
\ref{lemma-item-increasing-by-one-the-rank-of-a-primitive-2-vector}
follows immediately from Lemma 
\ref{lemma-primitive-embeddings-of-rank-2-lattices}
\EndProof

\bigskip
\noindent
{\bf Proof of Proposition \ref{prop-compare-Gamma-v}:}
Let $W_v$ be the subgroup of $\Gamma_v$ generated by $\Gamma_0$ 
and reflections in $-2$ vectors of the form $(1,-\LB,m)$ 
with $\LB$ in $H^{2}(S,\Integers)$ a primitive class of 
self-intersection $2m-2$. Lemma 
\ref{lemma-two-normal-subgroups}
reduces the proof of the equality $W_v=\Gamma_v$,
to the proof, that $W_v$ contains the reflection $\tau_{v_0}$,
for every $-2$ vector $v_0=(r,-\LB_0,rm)$, with $\LB_0$ primitive.
Since $\tau_{v_0}=\tau_{-v_0}$,
we may assume that the rank $r$ is non-negative. 
We prove Proposition
\ref{prop-compare-Gamma-v} by induction on the rank
of a $-2$ vector $v_0=(r,-\LB_0,rm)$ in $v^\perp$. 

%We prove first that the reflection $\tau_{v_0}$ is in $W_v$, provided that 
%$\LB_0$ is primitive.
%\begin{equation}
%\label{eq-r-is-even-or-LB-is-primitive}
%r \ \mbox{is even, or} \  \LB_0 \ \mbox{is primitive.} 
%\end{equation}

If $r=1$ (and $\LB_0$ is primitive), or if $r=0$, then  $\tau_{v_0}$ 
is in $W_v$  (by definition of $W_v$). Assume $r\geq 2$ and 
that any $-2$ vector in $v^\perp$, 
%of type 
%(\ref{eq-r-is-even-or-LB-is-primitive}) 
with primitive first Chern class and of rank
$1\leq a \leq r-1$, is in $W_v$. If $r$ is even, 
Part \ref{lemma-item-decomposition-of-2-vector} of Lemma
\ref{lemma-existence-of-pairs-of-line-bundles}
implies, that $\tau_{v_0}$ 
is in $W_v$. If $r$ is odd and $\LB_0$ is primitive, then 
part \ref{lemma-item-increasing-by-one-the-rank-of-a-primitive-2-vector}
of Lemma \ref{lemma-existence-of-pairs-of-line-bundles}
replaces $v_0$ by $w_0$, of rank $r+1$, such that 
$\tau_{v_0}$ is in $W_v$ if and only if $\tau_{w_0}$ is in $W_v$. 
Part \ref{lemma-item-decomposition-of-2-vector} implies that $\tau_{w_0}$ 
is in $W_v$. 
%This completes the proof of Proposition \ref{prop-compare-Gamma-v}.
\EndProof

%Next, we use part 
%\ref{lemma-item-changing-odd-rank-to-even-rank-2-vector} of Lemma
%\ref{lemma-existence-of-pairs-of-line-bundles}
%to replace the primitivity assumption on $\LB$ by
%the weaker assumption that $\LB$ is not divisible by $2$. 
%Let $v_1=(a,-\LB_1,a(n-1))$, with $\LB_1$ of degree $2a^2(n-1)-2$,
%but $\LB_1$ is not primitive. Choose a primitive $\LB_2$ 
%as in part \ref{lemma-item-changing-odd-rank-to-even-rank-2-vector}. 
%Then the vector $v_0:=\tau_{v_2}(v_1)$ has even rank. Thus $\tau_{v_0}$ 
%is in $W_v$. By construction, $\tau_{v_2}$ is in $W_v$ as well. 
%The relation $\tau_{v_1}=\tau_{v_2}\circ\tau_{v_0}\circ\tau_{v_2}$ 
%implies that $\tau_{v_1}$ is also in $W_v$. 
%
%We conclude that $W_v$ is equal to $\Gamma_v^-$.
%Lemma \ref{lemma-two-normal-subgroups} 
%implies the equality $W_v=\Gamma_v$. 
%This completes the proof of Proposition \ref{prop-compare-Gamma-v}.
%\EndProof

%\medskip
%Reflections in $-2$ vectors play a major role in this paper. Hence, 
%we conclude this section with the following lemma,
%though we will not use it. 

\begin{new-lemma}
\label{lemma-two-orbits-of-minus-2-vectors}
The set of $-2$ vectors in $v^\perp$ is the union of two 
$O(v^\perp)$-orbits (which are also  $\Gamma_v$-orbits):
\begin{eqnarray*}
A_+ & := & \{v_0=(r,\LB,rm) \ : \  \ \LB\cdot\LB=2r^{2}m-2 \ 
\mbox{and} \ \LB \ \mbox{is divisible by} \ 2\} \ \ \ \mbox{and} 
\\
A_- & := & \{v_0=(r,\LB,rm) \ : \  \ \LB\cdot\LB=2r^{2}m-2 \ 
\mbox{and} \ \LB \ \mbox{is not divisible by} \ 2\}.
%\mbox{the complement of} \ A_+ \ \mbox{in the set of} \ -2 \
%\mbox{vectors in} \ v^\perp\} 
\end{eqnarray*}
The orbit $A_+$ is non-empty if and only if $m\equiv 1$ (mod $4$). 
Moreover, every  vector $v_0$ in $A_+$ has {\em odd} rank. 
\end{new-lemma}

\noindent
{\bf Proof:}
$A_+$ is contained in the kernel of
$
v^\perp \rightarrow  (v^\perp)^*/2(v^\perp)^*.
$
Equivalently, if $v_0$ is in $A_+$, then 
the reflection $\tau_{v_0}$ is in the kernel of
$
\Gamma_v \ \rightarrow \ \Aut(v^\perp/2v^\perp).
$
Clearly, $A_-$ contains elements which do not have this property. 
%Hence, once we show that $A_-$ is a single $O(v^\perp)$ orbit,
%we will conclude that $A_+$ and $A_-$ are disjoint. 

Let $v_0:=(r,b\LB_0,rm)$ be a $-2$ vector in $v^\perp$, where
$b$ is an integer and $\LB_0$ a primitive class. The equality
\[
-1 \ = \ \frac{1}{2}(v_0,v_0) \ = \ b^2\frac{(\LB_0,\LB_0)}{2}-r^{2}m
\]
implies that $b$ and $rm$ are relatively prime.
%\[
%\mbox{gcd}(b,rm) \ = \ 1.
%\]

If $v_0\in A_-$, then the lattice $\mbox{span}_\Integers\{v,v_0\}$
is a primitive sub-lattice of the Mukai lattice. Indeed,
given $x,y \in \RationalNumbers$, the vector 
$xv+y(v_0-rv)=(x,yb\LB_0,2myr-xm)$ is integral, if and only if
$x$, $yb$ and $2myr$ are integral. Since $v_0\in A_-$, 
$b$ is odd and is thus relatively prime to $2mr$. 
Thus, $y$ must be an integer as well. 

Given two $-2$ vectors in $A_-$, we get two primitive
sublattices $\Delta_0:=\mbox{span}\{v,v_0\}$ and 
$\Delta_1:=\mbox{span}\{v,v_1\}$
of the Mukai lattice. The isomorphism $\Delta_0\rightarrow\Delta_1$,
sending the ordered basis $\{v,v_0\}$ to $\{v,v_1\}$, extends to 
a global isometry $g$ of the Mukai lattice
(Lemma \ref{lemma-primitive-embeddings-of-rank-2-lattices}). 
Hence, $v_0$ and $v_1$ belong to the same
$\Gamma_v$ orbit in $v^\perp$.

Let $v_0:=(r,b\LB_0,rm)$ belong to $A_+$. Then $b$ is even,
$r$ is odd, and  $\mbox{span}_\Integers\{v,v_0\}$
is an index $2$ sublattice of 
$\Delta_0:=\mbox{span}_\Integers\{v, (0,\frac{b}{2}\LB_0,rm)\}$.
It is easy to check that  $\Delta_0$ is a primitive sublattice 
of the Mukai lattice and 
\[
\left\{v,\frac{v_0-v}{2}\right\} \ \ = \ \ 
\left\{(1,0,-m), \ \ 
\left(\frac{r-1}{2}, \frac{b}{2}\LB_0, \frac{m(r+1)}{2}\right)\right\}
\]
is another basis of $\Delta_0$. The intersection form of 
$\Delta_0$ with respect to the latter basis is
$\left(
\begin{array}{cc}
2m & -m
\\
-m & \frac{2m-2}{4}
\end{array}
\right)
$. 
We conclude that $m\equiv 1$ (mod $4$). 
Conversely, if $m\equiv 1$ (mod $4$), then $2m-2\equiv 0$ (mod $8$) 
and there exists a class $\LB_0$ of self-intersection $\frac{2m-2}{4}$. 
The vector $v_0=(1,2\LB_0,m)$ is then a $-2$ vector in $A_+$. 

Given two vectors $v_0$ and $v_1$ in $A_+$, we get the two primitive 
sublattices
$\Delta_0:=\mbox{span}\{v,\frac{1}{2}(v_0-v)\}$ and
$\Delta_1:=\mbox{span}\{v,\frac{1}{2}(v_1-v)\}$. The isomorphism
$\Delta_0\rightarrow \Delta_1$, sending the basis
$\{v,\frac{1}{2}(v_0-v)\}$ to $\{v,\frac{1}{2}(v_1-v)\}$,
extends to an isometry of the Mukai lattice (Lemma 
\ref{lemma-primitive-embeddings-of-rank-2-lattices}). 
Hence, $v_0$ and $v_1$ belong to the same $\Gamma_v$ orbit in $v^\perp$.
%Let $v_0:=(r,\LB,rm)$ belong to $A_+$. 
%If $v_1=(a,\LB_1,am)$ is a $-2$ vector in $v^\perp$ then 
%$\tau_{v_1}(v_0) = v_0 + k\cdot v_1$, where 
%$k=(v_0,v_1)=\LB\cdot\LB_1-2ar(n-1)$ is even. Thus, $\tau_{v_1}(v_0)$ 
%belongs to $A_+$ and $A_+$ is invariant. 
%As $\LB$ is divisible by $2$, its length square $\LB\cdot\LB$ is divisible 
%by $8$.
%If $r$ is even, then $(v_0,v_0)$ is divisible by $8$ and $v_0$ can not
%be a $-2$ vector. In fact, the congruence 
%$2r^2(n-1) \equiv 2$ (mod $8$) is equivalent to 
%\[
%r^2(n-1) \equiv 1 \ \ (\mbox{mod} \ 4).
%\]
%Thus, $n-1$ is a square of an odd number (mod $4$). Hence,
%$n-1\equiv 1$ (mod $4$). 
%Conversely, if $n\equiv 2$ (mod $4$), then $2n-4\equiv 0$ (mod $8$) and we 
%can choose $v_0=(1,2\LB',(n-1))$ where $\LB'$ has degree $\frac{2n-4}{4}$.
\EndProof

\begin{cor}
\label{cor-stabilizer-is-generated-by-reflections}
\begin{enumerate}
\item
\label{cor-item-stabilizer-is-generated-by-reflections}
$\Gamma_v$ is generated by reflections in $+2$ and $-2$ 
vectors in $v^\perp$.
\item
\label{cor-item-characters}
The character group ${\rm Char }(\Gamma_v)$, of homomorphisms
from $\Gamma_v$ to $\ComplexNumbers^\times$, is isomorphic to 
$\Integers/2\times \Integers/2$.
\end{enumerate}
\end{cor}

\noindent
{\bf Proof:}
\ref{cor-item-stabilizer-is-generated-by-reflections})
Wall proved, that $\Gamma_0$ is generated by reflections in $+2$ and $-2$ 
vectors in $H^2(S,\Integers)$ (\cite{wall} Theorem 4.8).
The statement follows from Wall's result and Lemma 
\ref{lemma-stabilizer-is-generated-by-reflections-and-G-0}. 

\ref{cor-item-characters}) 
We construct first an isomorphism 
${\rm Char}(\Gamma_0)\cong \Integers/2\times \Integers/2$.
There is precisely one $\Gamma_0$-orbit of 
$-2$ vectors and one $\Gamma_0$-orbit of $+2$ vectors 
in $H^2(S,\Integers)$ 
(Lemma \ref{lemma-primitive-embeddings-of-rank-2-lattices}).
Reflections in $-2$ vectors, in a fixed orbit, are all conjugate. 
Hence, every character is determined by its values on these two orbits.
Let $v_0$ be a $-2$ vector,  $u_0$ a $+2$ vector, and
$\tau_{v_0}$ and $\sigma_{u_0}$ the reflections. Then
$\det(\tau_{v_0})=-1$ and $\det(\sigma_{u_0})=-1$, while
$cov(\tau_{v_0})=0$ and $cov(\sigma_{u_0})=1$. 
Hence, $\{\det,cov\}$ is a basis for ${\rm Char}(\Gamma_0)$
as a $\Integers/2$-module.

We prove next the equality 
${\rm Char}(\Gamma_v)={\rm Char}(\Gamma_0)$. 
Proposition \ref{prop-compare-Gamma-v} and Lemma
\ref{lemma-two-orbits-of-minus-2-vectors}
imply, that $\Gamma_v$ is generated by $\Gamma_0$ and reflections in $-2$ 
vectors $v_0$, whose $\Gamma_v$-orbit contains $-2$ vectors in 
$H^2(S,\Integers)$. 
Reflections in $-2$ vectors, in a fixed orbit of
$\Gamma_v$, are all conjugate. Hence, the restriction
${\rm Char}(\Gamma_v)\rightarrow{\rm Char}(\Gamma_0)$ is injective. 
The restriction is surjective, since both $cov$ and $\det$ extend to
$\Gamma_v$.
\EndProof


\begin{thebibliography}{B-N-R}

%\bibitem[ACGH]{a-c-g-h} Arbarello E., Cornalba M., Griffiths P., Harris J.: 
%{\em Geometry of Algebraic curves Volume I.\/} Springer 1984. 

\bibitem[At]{atiyah-book} Atiyah, M. F.:
{\em $K$-theory.\/} 
Lecture notes by D. W. Anderson. 
W. A. Benjamin, Inc., New York-Amsterdam 1967. 

\bibitem[AK]{altman-kleiman} Altman, A., Kleiman, S.:
{\em Compactifying the Picard scheme.\/} 
Adv. in Math. 35, 50-112 (1980).

\bibitem[BB]{baily-borel} Baily, W. L., Borel, A.: 
{\em Compactification of arithmetic quotients of bounded symmetric 
domains.\/}  Ann. of Math. (2)  84  1966 442-528.

\bibitem[BFM]{bfm}
Baum, P., Fulton, W., MacPherson, R.:
{\em Riemann-Roch and topological $K$-Theory
for singular varieties.\/}
Acta Math. 143 (1979), no. 3-4, 155-192.

%\bibitem[BPV]{CCS} Barth, W., Peters, C., and Van de Ven, A.:
%{\em Compact Complex Surfaces.\/} Springer-Verlag 1984

\bibitem[BHPV]{BHPV} Barth, W., Hulek, K., Peters, C., and Van de Ven, A.:
{\em Compact Complex Surfaces.\/} Second edition, 
Springer-Verlag, 2004. 

\bibitem[B1]{beauville-varieties-with-zero-c-1}
Beauville, A. {\em Varietes K\"ahleriennes dont la premiere classe de Chern 
est nulle.}  J. Diff. Geom. 18, p. 755-782 (1983).

\bibitem[B2]{beauville-diagonal} Beauville, A. {\em 
Sur la cohomologie de certains espaces de modules de fibr\'{e}s vectoriels.}
Geometry and analysis (Bombay, 1992),
37--40, Tata Inst. Fund. Res., Bombay, 1995. 

\bibitem[B3]{beauville-automorphisms} Beauville, A. {\em Some remarks on 
k\"{a}hler manifolds with $c_1=0$.\/
}
in Classification of algebraic and analytic manifolds (Katata, 1982), 1-26, 
Progr. Math., 39, Birkh\"{a}user Boston, 1983. 

%\bibitem[BD]{beauville-donagi}{Beauville, A., Donagi, R.:}
%{\em La vari\'et\'e des droites d'une hypersurface cubique de 
%dimension 4. \/}
%C. R. Acad. Sci. Paris Ser. I t. 301, 703-706 (1985)

%\bibitem[Bog]{bogomolov}
%Bogomolov, F.: {\em On the cohomology ring of a simple 
%hyper-k\"{a}hler manifold.\/} Geom. Funct. Anal. 6, 612-618 (1996)

\bibitem[BO]{bondal-orlov}
Bondal, A., Orlov, D.: 
{\em 
Reconstruction of a variety from the derived category and groups of 
autoequivalences.\/}
Compositio Math. 125 (2001), no. 3, 327--344.

\bibitem[Br]{bridgeland-elliptic-surfaces} Bridgeland, T.:
{\em Fourier-Mukai transforms for elliptic surfaces.\/} 
J. Reine Angew. Math.  498  (1998), 115--133.

\bibitem[Ca]{calduraru-thesis} C\u{a}ld\u{a}raru, A.:
{Derived categories of twisted sheaves on
Calabi-Yau manifolds.\/} Thesis, Cornell Univ., May 2000.

\bibitem[C]{chow} Chow, W. L.:
{\em On the geometry of algebraic homogeneous spaces.\/}
Ann. of Math. (2) 50, 32-67 (1949) 

\bibitem[CG]{chriss-ginzburg} Chriss, N., Ginzburg, V.:
{\em Representation theory and complex geometry.\/} 
Birkh\"{a}user Boston, Boston, MA, 1997. 

\bibitem[De]{deligne-LHST}
Deligne, P.: {\em Le groupe fondamental du compl\'{e}ment d'une courbe plane 
n'ayant que des points doubles ordinaires est ab\'{e}lien 
(d'apr\`{e}s W. Fulton).\/} 
Bourbaki Seminar, Vol. 1979/80,  pp. 1--10, 
Lecture Notes in Math., 842, Springer, 1981.

%\bibitem[Do]{donaldson} 
%Donaldson, S. K.:
%{\em Polynomial invariants for smooth four-manifolds.\/}
%Topology, Vol. 29 No. 3, 257-315, 1990.

\bibitem[Fu]{fulton} Fulton, W.: {\em Intersection Theory.\/}
Springer-Verlag 1984

\bibitem[FL]{fulton-lazarsfeld} Fulton, W., Lazarsfeld, R.:
{\em  Connectivity and its applications in algebraic geometry\/.}
Algebraic geometry (Chicago, Ill., 1980), pp. 26-92,
Lecture Notes in Math., 862,
Springer, 1981. 

\bibitem[Hai]{haiman} Haiman, M.:
{\em $t,q$-Catalan numbers and the Hilbert scheme.\/} 
Discrete Mathematics 193 (1998) 201-224.

\bibitem[Har]{hartshorne} Hartshorne, R.:
{\em Algebraic geometry.\/} 
Graduate Texts in Mathematics, No. 52.
Springer-Verlag, New York-Heidelberg, 1977. 

\bibitem[HLOY]{mirror-symmetry-k3}
Hosono, S., Lian, B., Oguiso, K., Yau S-T.:
{\em Autoequivelences of derived category of a 
K3 surface and monodromy transformations.\/} 
math.AG/0201047

%\bibitem[Hu1]{huybrechts-birational-implies-deformation-equivalent}
%Huybrechts, D.: 
%{\em Birational symplectic manifolds and their deformations.\/}
%J. Differential Geom. 45 (1997), no. 3, 488--513.

\bibitem[Hu]{huybrechts-basic-results}
Huybrechts, D.: 
{\em Compact Hyperk\"{a}hler Manifolds: Basic results.\/}
Invent. Math. 135 (1999), no. 1, 63-113 and
Erratum in Invent. Math. 152, 209-212 (2003).

%\bibitem[Hu]{huybrechts-norway}
%Huybrechts, D.: {\em Compact hyperk\"{a}hler manifolds.\/}
% Calabi-Yau manifolds and related geometries (Nordfjordeid, 2001),  
%161-225, Universitext, Springer, Berlin, 2003. 

\bibitem[HL]{huybrechts-lehn-book}
Huybrechts, D, Lehn, M.: 
{\em The geometry of moduli spaces of sheaves.\/} 
Aspects of Mathematics, E31. Friedr. Vieweg \& Sohn, Braunschweig, 1997.

\bibitem[HS]{huybrechts-stellari}
Huybrechts, D.; Stellari, P.: 
{\em Equivalences of twisted $K3$ surfaces.\/}
arXiv:math.AG/0409030 v3

\bibitem[K]{karoubi} Karoubi, M.:
{\em $K$-theory. An introduction.\/}
Springer-Verlag,  1978.

%\bibitem[K]{kovacs} Kov\'acs, S.: {\em The cone of curves of a K3 surface.}
%Math. Ann. 300, 681-691 (1994)

\bibitem[KV]
{kaledin-verbitsky} Kaledin, D., Verbitsky, M.:
{\em Partial resolutions of Hilbert type, Dynkin diagrams and generalized 
Kummer varieties.\/} Preprint, math.AG/9812078


\bibitem[LS1]{lehn-sorger} Lehn, M., Sorger, C.: {\em 
The cup product of the Hilbert scheme for $K3$ surfaces.\/}
%math.AG/0012166 
Invent. Math.  152  (2003),  no. 2, 305-329.

\bibitem[LS2]{lehn-sorger-generators} Lehn, M., Sorger, C.: 
Private communication of work in progress.

%\bibitem[L]{le-potier} Le Potier, J.:
%{\em Syst\'{e}mes coh\'{e}rents et structures de niveau.\/}
%Ast\'{e}risque 214, 1993

\bibitem[LL]{looijenga-lunts}  Looijenga, E., Lunts, V.:
{\em A Lie algebra attached to a projective variety.\/} Invent. Math. 129
(1997), no. 2, 361--412. 

%\bibitem[Loj]{lojasiewicz} Lojasiewicz, S.:
%{\em Biholomorphismes des vari\'{e}t\'{e}s grassmanniennes.\/}
%Geometry seminars, 1982--1983, 93--113, Univ. Stud. Bologna, Bologna, 1984. 

%\bibitem[Mag]{magyar} Magyar, P.:
%{\em Affine Schubert Varieties and Circular Complexes.\/}
%Preprint, math.AG/0210151


%\bibitem[Mat]{matsumura} Matsumura, H.: 
%{\em Commutative Algebra.\/} W. A. Benjamin Co., New York (1970)



\bibitem[Ma1]{markman-reflections} Markman, E.:
{\em Brill-Noether duality for moduli spaces of sheaves on K3 surfaces}
J. of Alg. Geom. {\bf 10} (2001), no. 4, 623-694

\bibitem[Ma2]{markman-diagonal} Markman, E.: 
{\em Generators of the cohomology ring of moduli spaces of sheaves on K3 
and Abelian surfaces.\/}  
Journal fur die reine und angewandte Mathematik 544 (2002), 61-82

\bibitem[Ma3]{markman-part-two} Markman, E.: 
{\em On the monodromy of moduli spaces of sheaves on K3 surfaces II.\/} 
Preprint, math.AG/0305043 v4.

\bibitem[Ma4]{markman-integral-generators} Markman, E.: 
{\em Integral generators for the cohomology ring of moduli spaces of 
sheaves over Poisson surfaces.\/} 
Preprint, math.AG/0406016 v1. 

%\bibitem[Mi]{milne} Milne, J.: {\em \'{E}tale Cohomology.\/} 
%Princeton Univ. Press (1980)


%\bibitem[MT]{mehta-trivedi}
%Mehta, V. B. and Trivedi, V.:
%{\em The variety of circular complexes and
%F-splitting.\/}
%Invent. math. 137, 449-460 (1999)


%\bibitem[Mo]{morrison} Morrison, D.:
%{\em On K3 surfaces with large Picard number}
%Invent. Math. 75, 105-121 (1984)

\bibitem[Mu1]{mukai-symplectic-structure}  Mukai, S.:
{\em Symplectic structure of the moduli space of
sheaves on an abelian or K3 surface}, 
Invent. math. 77,101-116 (1984)

\bibitem[Mu2]{mukai-hodge} Mukai, S.: 
{\em On the moduli space of bundles on K3 surfaces I},
Vector bundles on algebraic varieties,
Proc. Bombay Conference, 1984, Tata Institute of Fundamental Research Studies,
no. 11, Oxford University Press, 1987, pp. 341-413.

\bibitem[Mu3]{mukai-fourier} Mukai, S.: 
{\em
Duality between $D(X)$ and $D(\hat X)$ with its application to 
Picard sheaves. 
\/}
Nagoya Math. J. 81 (1981), 153--175.

\bibitem[Mu4]{mukai-applications} Mukai, S.: 
{\em Fourier functor and its application to the moduli of bundles 
on an Abelian variety.\/}
Adv. Studies in Pure Math. 10, 515-550 (1987).

%\bibitem[Mu5]{mukai-survey} Mukai, S.: 
%{\em Moduli of vector bundles on K3 surfaces, and symplectic manifolds.}
%Suguku Expositions, Vol. 1, No 2, (1988)

%\bibitem[Mu6]{mukai-curves-K3-fanos}   S. Mukai, 
%{\em Curves, K3 surfaces and Fano 3-folds of genus $\leq 10$.\/}
%In Algebraic Geometry and Commutative Algebra in honor of
%M. Nagata, pp. 357-377(1987)

\bibitem[Na1]{nakajima-reflections}
Nakajima, H.: {\em
Reflection functors for quiver varieties and Weyl group actions.\/} 
 Math. Ann.  327  (2003),  no. 4, 671--721.

\bibitem[Na2]{nakajima-representations}
Nakajima, H.: {\em 
Convolution on homology groups of moduli spaces of sheaves on K3 surfaces.
\/} 
 Vector bundles and representation theory (Columbia, MO, 2002),  
75-87, Contemp. Math., 322, Amer. Math. Soc., Providence, RI, 2003. 

\bibitem[Na3]{nakajima-book}
Nakajima, H.: {\em
Lectures on Hilbert schemes of points on surfaces.\/} 
University Lecture Series, 18. Amer. Math. Soc., 
Providence, RI, 1999.

%\bibitem[Nam]{namikawa} Namikawa, Y.:
%{\em Mukai flops and derived categories.\/} Preprint, 
%Math.AG/0203287

\bibitem[Nam]{namikawa-deformations} Namikawa, Y.:
{\em Deformation theory of singular symplectic $n$-folds.\/}  
Math. Ann.  319  (2001),  no. 3, 597--623.

\bibitem[Ni]{nikulin} Nikulin, V. V.: 
{\em Integral symmetric bilinear forms and some of their applications.\/}
Math. USSR Izvestija, Vol. 14 (1980), No. 1

\bibitem[OG1]{ogrady-hodge-str} O'Grady, K.: 
{\em The weight-two Hodge structure of moduli spaces of sheaves on a K3 
surface.\/} J. Algebraic Geom. 6 (1997), no. 4, 599-644.

\bibitem[OG2]{ogrady-involutions} O'Grady, K.: 
{\em Involutions and linear systems on holomorphic symplectic manifolds.\/}
Preprint, arXiv.org math.AG/0403519.

\bibitem[O]{oguiso} Oguiso, K.:
{\em K3 surfaces via almost-primes. \/}
Math. Res. Lett. 9 (2002), no. 1, 47--63. 


\bibitem[Or1]{orlov} Orlov, D. O.: 
{\em Equivalences of derived categories and K3 surfaces.\/}
Algebraic geometry, 7. J. Math. Sci. (New York) 84 (1997), no. 5, 1361--1381. 

\bibitem[Or2]{orlov-abelian-varieties} Orlov, D. O.: 
{\em Derived categories of coherent sheaves on abelian varieties 
and equivalences between them.\/
} 
%alg-geom/9712017
(Russian)  Izv. Ross. Akad. Nauk Ser. Mat.  66  (2002),  no. 3, 131-158;  
translation in  Izv. Math.  66  (2002),  no. 3, 569-594

\bibitem[P]{peters} Peters, C.: 
{\em Monodromy and Picard-Fuchs equations for families of 
$K3$-surfaces and elliptic curves.\/}  Ann. Sci. \'{E}cole Norm. Sup. (4)  19  (1986),  
no. 4, 583-607.

\bibitem[Sa]{salamon} Salamon, S. M.:
{\em On the cohomology of K\"{a}hler and hyper-K\"{a}hler manifolds.\/} 
Topology  35  (1996),  no. 1, 137--155. 

%\bibitem[St]{strickland} Strickland, E.: 
%{\em On the conormal bundle of the determinantal variety.\/}
%J. of Algebra 75, 523-537 (1982).

\bibitem[Sz]{szendroi} Szendr\"{o}i, B.:
{\em Diffeomorphisms and families of Fourier-Mukai transforms in mirror 
symmetry.\/}
Applications of algebraic geometry to coding theory, 
physics and computation (Eilat, 2001), 317--337, 
NATO Sci. Ser. II Math. Phys. Chem., 36, 
Kluwer Acad. Publ., Dordrecht, 2001.

\bibitem[ST]{seidel-thomas} 
Seidel, P., Thomas, R.:  {\em 
Braid group actions on derived categories of coherent sheaves.\/} 
Duke Math. J. 108 (2001), no. 1, 37--108. 

%\bibitem[Ty]{tyurin-cycles-curves-surfaces} {Tyurin, A. N.:}
%{\em Cycles, curves and vector bundles on an
%algebraic surface.\/}
%Duke Math. J. Vol. 54 No. 1 (1987) pages 1-26

\bibitem[Ve1]{verbitsky-announcement} Verbitsky, M.: 
{\em 
Cohomology of compact hyper-K\"{a}hler manifolds and its applications.\/},
Geom. Funct. Anal. 6 (1996), no. 4, 601--611. 


\bibitem[Ve2]{verbitsky-mirror-symmetry} Verbitsky, M.: 
{\em Mirror symmetry for hyper-K\"{a}hler manifolds.\/}
Mirror symmetry, III (Montreal, PQ, 1995), 115--156, 
AMS/IP Stud. Adv. Math., 10, Amer. Math. Soc., Providence, RI, 1999. 

\bibitem[Vi]{viehweg} Viehweg, E.:
{\em Quasi-projective Moduli for Polarized Manifolds.\/}
Springer-Verlag (1995).

\bibitem[W]{wall} Wall, C. T. C.: 
{\em On the orthogonal groups of unimodular quadratic forms. II.\/} 
J. Reine Angew. Math. 213  (1963/1964), 122-136. 

\bibitem[Y1]{yoshioka-examples-of-reflections} Yoshioka, K.: 
{\em Some examples of Mukai reflections on K3 surfaces.\/} 
J. Reine Angew. Math. 515 (1999), 97--123.  

\bibitem[Y2]{yoshioka-irreducibility} Yoshioka, K.:
{\em Irreducibility of moduli spaces of vector bundles on K3 surfaces.\/}
math.AG/9907001

\bibitem[Y3]{yoshioka-abelian-surface} Yoshioka, K.:
{\em 
Moduli spaces of stable sheaves on abelian surfaces. \/
}
Math. Ann. 321 (2001), no. 4, 817--884.

\bibitem[Y4]{yoshioka-note-on-fourier-mukai}
 Yoshioka, K.:
{\em A Note on Fourier-Mukai transform.\/}
Eprint arXiv:math.AG/0112267 v3.


\end{thebibliography}
\end{document}